Министерство образования и науки Российской Федерации
Федеральное государственное автономное
образовательное учреждение высшего образования
«Московский физико-технический институт
(государственный университет)»

А. В. Гасников

# СОВРЕМЕННЫЕ ЧИСЛЕННЫЕ МЕТОДЫ ОПТИМИЗАЦИИ.

## МЕТОД УНИВЕРСАЛЬНОГО ГРАДИЕНТНОГО СПУСКА

Учебное пособие

Москва
МФТИ
2018

Р е ц е н з е н т ы :

Директор Вычислительного центра им. А. А. Дородницына ФИЦ ИУ РАН
академик РАН *Ю. Г. Евтушенко*

Директор Института вычислительной математики РАН, заведующий кафедрой
вычислительных технологий и моделирования ВМиК МГУ
академик РАН *Е. Е. Тыртышников*

Рассматривается классический градиентный спуск. Однако изложение ведется на продвинутом уровне. Пособие отличается довольно полным обзором современного состояния методов типа градиентного спуска.

В данном пособии делается акцент не на изложение методов, а на способы получения из старых методов новых с помощью небольшого числа общих приемов.

Учебное пособие используется для преподавания численных оптимизация студентам 3 курса школы ПМИ МФТИ и 3 курса ФКН ВШЭ.

# Оглавление

# Обозначения

$\mathbb{R}^n$ – $n$-мерное вещественное (векторное) пространство.

$\mathbb{R}^n_+ = \left\{x \ge 0:\ x \in \mathbb{R}^n\right\}$ – неотрицательный ортант $\mathbb{R}^n$.

$\mathrm{const}$ – числовая константа, значение которой зависит от контекста.

$\dim x$ – размерность вектора $x$, в частности, $\dim x = n$, если $x \in \mathbb{R}^n$.

$[x]_i$, $x_i$ – $i$-я компонента вектора $x$.

$\|x\|_p = \left(\sum_{i=1}^n |x_i|^p\right)^{1/p}$ – $p$-норма вектора $x = \{x_i\}_{i=1}^n \in \mathbb{R}^n$, $p \ge 1$.

$\langle x, y\rangle = \sum_{i=1}^n x_i y_i$ – скалярное произведение векторов $x, y \in \mathbb{R}^n$. Отметим, что при определении графа используется похожее обозначение $G = \langle V, E\rangle$, имеющее другой смысл.

$\|y\|_* = \sup_{\|x\|\le 1}\langle x, y\rangle$ – сопряженная норма к норме $\|\ \|$. В частности, для $p$-нормы сопряженной будет $q$-норма, где $1/p + 1/q = 1$.

$\lceil a \rceil = \max\{1, a\}$.

$A \Rightarrow B$ – из утверждения (формулы) $A$ следует утверждение (формула) $B$.

$A \Leftrightarrow B$ – утверждение (формула) $A$ эквивалентно (равносильно) утверждению (формуле) $B$, т. е. $A \Rightarrow B$ и $B \Rightarrow A$.

$x \ll y$ – число $y$ много больше числа $x$.

$x \simeq y$, $x \approx y$ – число $x$ приближенно равно числу $y$.

$x \sim y$ – значение выражения $x$ пропорционально значению выражения $y$, например, $V = HR \Rightarrow V \sim R$ или $T = 2\pi\sqrt{l/g} \Rightarrow T \sim \sqrt{l/g}$.

$x := y$ – $x$ присваивается $y$ (пришло из программирования), например, $x := x + 1$.

$\mathrm{Lin}\left\{z^1, .., z^m\right\}$ – линейное пространство (подпространство $\mathbb{R}^n$), «натянутое» на векторы $z^1, .., z^m \in \mathbb{R}^n$, т. е. любой элемент такого пространства можно представить в виде

$$\alpha_1 z^1 + \ldots + \alpha_m z^m,\ \alpha_1, \ldots, \alpha_m \in \mathbb{R}.$$

$A^T$ – матрица, транспонированная к матрице $A = \left\| A_{ij} \right\|_{i,j=1}^{n}$, т. е. $A^T = \left\| A_{ji} \right\|_{i,j=1}^{n}$.

$\operatorname{tr}(A)$ – след квадратной матрицы $A$, т.е. сумма всех ее диагональных элементов.

$\operatorname{rank} A$ – ранг матрицы $A$, т.е. максимальное число линейно независимых столбцов (или строк – не важно).

$I$ – единичная матрица, т. е. $I = \left\| I_{ij} \right\|_{i,j=1}^{n}$, где $I_{ij} = 0$, если $i \neq j$, $I_{ij} = 1$, иначе.

$A^{-1}$ – матрица, обратная к квадратной матрице $A$, т. е. $A^{-1}A = AA^{-1} = I$.

$(\operatorname{Ker} A)^{\perp}$ – ортогональное дополнение подпространства, натянутого на собственные векторы матрицы $A$, отвечающие нулевому собственному значению.

$\sqrt{A}$ – квадратный корень из симметричной неотрицательно определенной матрицы $A$.

◊ Для каждой неотрицательно определенной симметричной матрицы существует такой ортонормированный базис, в котором действие этой матрицы можно понимать как соответствующие растяжение/сжатие/проектирование (задается собственными значениями $\lambda_i$) вдоль ортов. Тогда действие матрицы $\sqrt{A}$ можно понимать как растяжение/сжатие/проектирование (задается собственными значениями $\sqrt{\lambda_i}$) вдоль тех же самых ортов. ◊

$A^j$ – $j$-й столбец матрицы $A$; $A_i$ – $i$-я строка матрицы $A$.

$A \succ 0$ – симметричная матрица $A$ ($A = A^T$) неотрицательно определена, т. е.

$$\forall \ x \in \mathbb{R}^n \to \langle x, Ax \rangle \geq 0.$$

$A \succ B$ – означает, что $A - B \succ 0$.

$1_n = \underbrace{(1,\ldots,1)}_{n}^T$ – вектор из единиц.

$e_i = (\underbrace{0,\ldots,0,1}_{i},0,\ldots,0)^T$ – $i$-орт.

$S_n(1) = \left\{ x \in \mathbb{R}_+^n : \sum_{i=1}^{n} x_i = 1 \right\}$ – единичный симплекс пространства $\mathbb{R}^n$.

$B_{R,Q}(y) = \{x \in Q:\ \|x - y\|_2 \le R\}$ – пересечение евклидова шар радиуса $R$ с центром в точке $y$ и множества $Q$.

$\tilde{x} = \arg\min_{x \in P} F(x)$ – означает, что $F(\tilde{x}) < F(x)$ для всех $x \in P \setminus \tilde{x}$.

$\tilde{x} \in \operatorname{Arg}\min_{x \in P} F(x)$ – означает, что $F(\tilde{x}) \le F(x)$ для всех $x \in P$.

$\pi_Q(x) = \arg\min_{y \in Q} \|x - y\|_2^2$ – евклидова проекция точки (вектора) $x$ на замкнутое выпуклое множество $Q$.

$\log a$ – логарифм положительного числа $a$ по основанию, зависящему от контекста (в частности, $\ln = \log_e$). Если основание логарифма не указано, значит, основание зависит от контекста и это хочется подчеркнуть (см. приложение).

$E_\xi\left[F(x,\xi)\right]$ – математическое ожидание по случайной величине (вектору) $\xi$ от измеримой по $\xi$ (вектор) функции $F(x,\xi)$. Здесь $x$ следует понимать как параметр.

$E_\xi\left[f(\xi,\eta)\middle|\eta\right]$ – условное математическое по случайной величине (вектору) $\xi$ при «замороженной» случайной величине $\eta$ от измеримой по $\xi$ и $\eta$ функции $f(\xi,\eta)$. Условное математическое ожидание является случайной величиной, зависящей от $\eta$.

$\lambda_{\max}(A) = \max\{\lambda:\ \exists\ x \ne 0: Ax = \lambda x\}$,

$\lambda_{\min}(A) = \min\{\lambda:\ \exists\ x \ne 0: Ax = \lambda x\}$.

$\sigma_{\max}(A) = \lambda_{\max}(A^T A) = \lambda_{\max}(AA^T) = \max\{\lambda:\ \exists\ x \ne 0: AA^T x = \lambda x\}$.

$\tilde{\sigma}_{\min}(A) = \min\{\lambda > 0:\ \exists\ x \ne 0: AA^T x = \lambda x\}$.

$\nabla f(x) = \left(\partial f(x)/\partial x_1, \dots, \partial f(x)/\partial x_n\right)^T$ – градиент гладкой функции $f(x)$, так же

$$\nabla_x f(x,y) = \left(\partial f(x,y)/\partial x_1, \dots, \partial f(x,y)/\partial x_n\right)^T.$$

$\partial f(x)$ – субдифференциал выпуклой функции $f(x)$. По определению $g \in \partial f(x)$ ($g$ – субградиент $f$ в точке $x$) тогда и только тогда, когда для всех $y$ имеет место неравенство $f(y) \ge f(x) + \langle g, y - x\rangle$.

$\nabla h(y) = \left\| \partial h_i(y) / \partial y_j \right\|_{i,j=1}^{n,m}$ – матрица Якоби гладкого отображения $h : \mathbb{R}^m \to \mathbb{R}^n$.

$\nabla^2 f(x) = \left\{ \partial \nabla f(x) / \partial x_j \right\}_{j=1}^{n} = \left\| \partial^2 f(x) / \partial x_i \partial x_j \right\|_{i,j=1}^{n}$ – матрица Гессе дважды дифференцируемой функции $f(x)$. Аналогично можно определить

$$\nabla^{r+1} f(x) = \left\{ \partial \nabla^r f(x) / \partial x_j \right\}_{j=1}^{n}.$$

$\nabla^{r+1} f(x)[u] = \sum_{j=1}^{n} \partial \nabla^r f(x) / \partial x_j \cdot u_j$ – тензор ранга $r$ ($u \in \mathbb{R}^n$).

$A(\text{параметры}) = \mathrm{O}(B(\text{параметры}))$ – означает, что существует такая абсолютная числовая константа $C$, не зависящая ни от каких параметров, что

$$A(\text{параметры}) \le C \cdot B(\text{параметры}).$$

$A(\text{параметры}) = \tilde{\mathrm{O}}(B(\text{параметры}))$ – означает, что существует такой множитель $\tilde{C}$, зависящий от параметров не сильнее, чем логарифмическим образом, что

$$A(\text{параметры}) \le \tilde{C} \cdot B(\text{параметры}).$$

$A \overset{\text{def}}{=} B$ – означает $A = B$, и это равенство определяет либо $A$, либо $B$.

# Предисловие

Данное пособие написано по материалам лекций, прочитанных автором в Летней школе «Современная математика» в Ратмино (г. Дубна) в июле 2017 года.

Идея курса состояла в том, чтобы, с одной стороны, рассказать основные приемы, с помощью которых порождается многообразие современных численных методов выпуклой оптимизации первого порядка (рестарты, регуляризация, переход к двойственной задаче, адаптивная настройка на гладкость задачи, минибатчинг, каталист, слайдинг и т. д.). С другой стороны, хотелось провести все рассуждения на строгом математическом языке (с полным обоснованием). Поэтому для наглядности было решено ограничиться изучением только градиентного спуска и его окрестностей.

По данному пособию было прочитано несколько курсов лекций, которые снимались на видео: в осеннем семестре 2018/2019 в КМЦ АГУ [609], в весеннем семестре 2018/2019, 2019/2020 в школе ПМИ МФТИ [605, 610].

Опыт использования пособия при чтении лекций студентам школы ПМИ МФТИ показывает, что в пособие желательно еще добавить подробно разобранные примеры решения практических задач оптимизации. Эту проблему планируется в перспективе решить за счет издания другого пособия. Однако отметим, что даже в текущем варианте пособия можно найти примеры вполне реальных (практических) задач, с которыми мы сталкивались в разное время. Например, в конце § 5 и в примере из приложения рассматривается задача композитной оптимизации, возникшая при решении задачи восстановления матрицы корреспонденций в большой компьютерной сети по замерам потоков на линках (ребрах). Эта обратная задача (обратные задачи являются естественным источником постановок задач в обычной оптимизации [355]) была поставлена компанией Хуавей в 2015 году.

Работа по подготовке пособия в § 1 была выполнена при поддержке Министерства науки и высшего образования Российской Федерации (госзадание) №075-00337-20-03, номер проекта 0714-2020-0005, в § 2 была поддержана грантом РФФИ 18-29-03071 мк, в § 3, 5 – грантом РФФИ 18-31-20005 мол_а_вед, в § 4 – грантом РФФИ 19-31-51001 Научное наставничество, в Приложении – грантом РНФ 18-71-10108.

# Введение

Это пособие написано прежде всего для студентов-математиков, начинающих изучать численные методы оптимизации и желающих впоследствии серьезно погрузиться в данную область.

Пожалуй, основным численным методом современной оптимизации является *метод градиентного спуска*. Метод прекрасно изложен в замечательной книге Б.Т. Поляка [82], вышедшей в 1983 году. В некотором смысле этот метод порождает[1] большинство остальных численных методов оптимизации. Метод градиентного спуска активно используется в вычислительной математике не только для непосредственного решения задач оптимизации (минимизации), но и для задач, которые могут быть переписаны на языке оптимизации [16, 49, 82, 91, 281, 355] (решение нелинейных уравнений, поиск равновесий, обратные задачи и т. д.). Метод градиентного спуска можно использовать для задач оптимизации в бесконечномерных пространствах [487], например, для численного решения задач оптимального управления [14, 48, 49, 82, 83, 281]. Но особенно большой интерес к градиентным методам в последние годы связан с тем, что градиентные спуски и их стохастические/рандомизированные варианты лежат в основе почти всех современных алгоритмов обучения, разрабатываемых в *анализе данных* [39, 78, 156, 171, 180, 196, 320, 343, 386, 413, 470, 508, 522, 537, 579].

◊ Все это также хорошо можно проследить по трем основным конференциям по анализу данных: COLT, ICML, NIPS, которые за последние 10–15 лет частично превратились в конференции, посвященные использованию градиентных методов в решении задач *машинного обучения*. ◊

Не удивительно в этой связи, что подавляющее большинство современных курсов по численным методам оптимизации построено вокруг градиентных методов [13, 46, 72, 154, 158, 176, 180, 319, 320, 386, 413, 452, 472]. Данное пособие, подготовленное по материалам курса, прочитанного в ЛШСМ 2017, также построено по такому принципу. Однако принципиальное методическое отличие предложенного курса от остальных заключается в том, что в данном курсе предпринята попытка на примере только градиентного спуска продемонстрировать основной арсенал приемов, с помощью которых разрабатываются новые численные методы и теоретически исследуется их скорость сходимости. Такое построение курса было обусловлено желанием в первую очередь донести основную идею того или иного приема, не отягощая изложение техническими дета-лями. Градиентный спуск был выбран по нескольким причинам: во-первых, пожалуй, он самый простой, во-вторых, он лежит в основе большинства других методов, и если хорошо разобраться с тем или иным приемом на примере градиентного спуска, то это можно использовать при перенесении на более сложный метод, лучше подходящий для решения конкретной задачи.

[1] Собственно, данное пособие имеет одной из своих целей пояснить смысл этого предложения и слова «порождает» в данном контексте.

Курс начинается со стандартного изложения в § 1 того, что такое градиентный спуск. А именно, исходно сложная минимизируемая (целевая) функция заменяется в окрестности рассматриваемой точки, касающимся её графика в этой точке параболоидом вращения, который по построению должен также мажорировать исходную функцию. Далее исходная задача минимизации заменяется задачей минимизации построенного параболоида. Последняя задача решается явно (осуществляется шаг градиентного спуска). Найденное решение задачи принимается за новую точку (положение метода) и процесс повторяется. В зависимости от того, какими свойствами обладала исходная функция (свойства гладкости, выпуклости), устанавливаются оценки на скорость сходимости описанной процедуры.

Начиная с § 2 изложение заметно усложняется, обрастая деталями. В § 2 рассматриваются задачи выпуклой оптимизации на множествах простой структуры (например, к таким множествам можно отнести неотрицательный ортант) в условиях небольших шумов неслучайной природы (см., например, [82, гл. 4]). Описанная выше процедура переносится на этот случай. Наличие шума играет ключевую роль в достижении одной из главных целей курса – построении *универсального градиентного спуска*. Этот метод сам настраивается на гладкость задачи и не требует параметров на входе.

В § 3 предлагается *концепция модели функции*, заключающаяся в том, что вместо параболоида вращения, аппроксимирующего (касающегося надграфика и мажорирующего) исходную выпуклую функцию в окрестности данной точки, можно использовать какие-то другие функции. Таким образом, например, можно дополнительно переносить «тяжесть» исходной постановки задачи на вспомогательные подзадачи, надеясь, что это ускорит сходимость метода. Понятно, что такое ускорение будет достигнуто за счет того, что каждая итерация станет дороже. Чтобы правильно по задаче выбрать модель функции, нужно иметь оценки того, насколько скорость сходимости внешней процедуры зависит от вида вспомогательных задач, точности их решения, и понимать, как сложность вспомогательных задач зависит от точности их решения. Все это прорабатывается в данном параграфе при достаточно общих условиях.

В § 4 демонстрируется *прямодвойственная* природа обсуждаемых методов для выпуклых задач. Свойство прямодвойственности метода поз-воляет почти бесплатно получать решение задачи, двойственной к данной. Как правило, для большинства оптимизационных задач, приходящих из практики (экономика [17, 455, 465], транспорт [18, гл. 1, 3], проектирование механических конструкций [467] и даже анализ данных [18, гл. 5], [500]), двойственная задача несет в себе дополнительную полезную информацию об изучаемом объекте (явлении), которую также хотелось бы получить в результате оптимизации. Другая не менее важная причина популярности прямодвойственных методов заключается в том, что, имея пару прямая–двойственная задача, можно выбирать, которую из них решать (какая проще). В частности, двойственные задачи, являются задачами выпуклой оптимизации на множествах простой структуры. Если при решении выбранной задачи (прямой или двойственной) использовать прямодвойственный метод, то решив её с некоторой точностью, гарантированно решим с такой же точностью и сопряженную (двойственную) к ней задачу.

◊ Напомним, что при весьма общих условиях [176, гл. 5] двойственной задачей для двойственной к исходной выпуклой задаче будет исходная задача (теорема Фенхеля–Моро [62, п. 1.4, 2.2]). ◊

В § 5 строится *прямодвойственный универсальный градиентный спуск* для задачи выпуклой оптимизации на множестве простой структуры. Концепция универсального метода обобщает известное и популярное на практике правило выбора шага дроблением/удвоением [42, п. 6.3.2], см. также правила Армихо, Вулфа, Голдстейна [13, гл. 5], [53, п. 3.1.2], [56, п. 9.4], [72, п. 1.2.3], [82, гл. 3], [472, гл. 3], выбора шага градиентного спуска. Эта концепция подготавливалась около 30 лет (см., например, [68, 75]), и лишь весной 2013 года была оформлена Ю.Е. Нестеровым сначала в виде препринта, а потом в виде статьи [462]. Статья вызвала большой интерес, и сейчас активно цитируется в оптимизационном сообществе. Отличие универсального подхода от *адаптивного* (к последнему можно отнести методы с выбором шага по отмеченным выше правилам, типа Армихо) заключается в том, что настройка происходит не только на константу гладкости, но и на степень гладкости по шкале: негладкая $\to$ гёльдерова $\to$ гладкая функция. Универсальные прямодвойственные методы, сейчас активно используются при поиске равновесий в больших транспортных сетях [7, 18, 28]. Большая популярность самонастраивающихся оптимизационных процедур в анализе данных, особенно в глубоком обучении[2] [39, 78, 129, 292] (в том числе использование нейросети

[2] Несмотря на огромную популярность этого направления и огромные усилия, затраченные на объяснение успешного практического опыта использования глубоких нейронных сетей в различных приложениях, важно подчеркнуть, что на данный момент, насколько нам известно, ученые по-прежнему достаточно далеки

для выбора величины шага в обучении другой нейросети), определенно указывает на то, что за адаптивными (самонастраивающимися), а по нашей терминологии «универсальными», методами будущее! Все это, безусловно, также сильно сказалось на отборе материала и сделанных в пособии акцентах.

◊ Опыт использования терминов *прямодвойственный* и *универсальный* (следуя [455, 462]) показывает, что оптимизационное сообщество в России принимает эти термины не однозначно. В частности, часто можно было слышать следующие замечания. «Представляется более естественным говорить про просто *двойственный метод* – см., например, метод Эрроу–Гурвица [82, п. 3 § 2, гл. 8], который имеет еще более ярко выраженную прямодвойственную структуру, чем рассматриваемые в пособии, однако относится к классу *двойственных методов*. Словосочетание *универсальный метод* несколько вводит в заблуждение масштабами универсальности. Ведь в данном контексте речь идет только об универсальном по гладкости методе, т. е. методе, который на вход не требует никакой информации о свойствах гладкости задачи (в том числе и константах, характеризующих гладкость). Однако, например, для сильно выпуклых задач такие методы требуют знания константы сильной выпуклости, и никакой самонастройки на эту константу по ходу работы (как в случае с константами, отвечающими за гладкость) уже не происходит.» В целом, несмотря на эти замечания, было решено сохранить термины в неизменном виде, поскольку в англоязычной литературе они уже достаточно прочно успели закрепиться и их исправление может осложнить последующее изучение читателями современной литературы по данной тематике, которая в основном вся на английском языке. ◊

В приложении приводится краткий обзор современного состояния дел в активно развивающейся в последние годы области численных методов выпуклой оптимизации. Материал излагается в контексте результатов, приведенных в основном тексте пособия. Приложение написано, в первую очередь, для читателей, желающих продолжить изучение курса численных методов оптимизации. Надеемся, что приложение поможет сориентироваться читателям и укажет на некоторые новые направления и возможности.

Важную роль в тексте пособия играют замечания и упражнения, которые рекомендуется, как минимум, просматривать, а лучше прорешивать. В частности, таким образом (через замечания и упражнения) вводятся два основных приема (сохраняющих оптимальность методов в смысле

---

от возможности научно все это объяснить, в том числе с точки зрения оптимизации. Более того, здесь имеются и вполне определенные отрицательные результаты [533].

числа обращений к оракулу), позволяющих переходить от выпуклых задач к сильно выпуклым и обратно. Соответственно, *метод регуляризации* и *метод рестартов*. Имея метод, настроенный на сильно выпуклые задачи с помощью регуляризации функционала, можно привести любую задачу к сильно выпуклой и использовать имеющийся метод. Обратно, имея метод, настроенный на выпуклые задачи, можно использовать данный метод для решения сильно выпуклых задач, *рестартуя* (перезапуская) его каждый раз, когда расстояние до решения сокращается в два раза. В упражнениях также обсуждается *ускоренный* градиентный (*быстрый, моментный*) спуск и *теория нижних оракульных оценок* сложности задач выпуклой оптимизации, построенная в конце 70-х годов XX века А.С. Немировским и Д.Б. Юдиным [70]. В замечании 3.3 описывается общий способ (*каталист*) ускорения неускоренных методов любого порядка.

Чтобы приблизить изложение к «живым» лекциям и местами немного «разбавить» достаточно насыщенный формулами материал, в пособии имеется также несколько исторических замечаний и замечаний «второго плана», выделенных следующим образом:

◊ … ◊.

Изложение построено таким образом, что по ходу изучения материала должна появляться интуиция о возможности практически произвольным образом и в любом количестве сочетать различные описанные приемы (конструкции, надстройки) друг с другом, получая, таким образом, все более и более сложные методы, лучше подходящие под решаемую задачу. В этой связи, наверное, можно сказать, что в пособии описаны «структурные блоки», из которых строятся современные градиентные методы. Замечательно, что эти же структурные блоки используются и для ускоренных методов и их стохастических и рандомизированных вариантов, см. приложение, а также [16, 18, 23, 25, 30, 71, 72, 93, 118, 122, 158, 225, 257, 264, 337, 339, 358, 384, 386, 445, 462, 468].

Приведем здесь для удобства основные структурные блоки (приемы) для методов первого порядка (градиентных методов), с указанием частей пособия, в которых они описаны. Эти блоки переносятся и на методы другого порядка, однако, детали вынуждены здесь опустить. Ограничимся, также для простоты, только четырьмя бинарными признаками, характеризующими решаемую задачу оптимизации и используемый метод: 1) задача гладкая / негладкая; 2) задача сильно выпуклая / выпуклая (вырожденная задача выпуклой оптимизации); 3) при решении задачи доступен градиент функционала / стохастический градиент; 4) для решения задачи используется ускоренный метод / неускоренный. Далее (см.

также табл. 2 в приложении и комментарии к ней) будут описаны приемы, которые в совокупности позволяют по (оптимальному) алгоритму, отвечающий конкретному набору этих четырех признаков строить (оптимальный) алгоритм, отвечающий любому из пятнадцати оставшихся наборов этих признаков. Впрочем, необходимости строить по ускоренным методам неускоренные на практике не возникает, поэтому соответствующее описание далее опущено.

1) Негладкая задача может рассматриваться как гладкая за счет искусственного введения неточности в параболическую модель аппроксимации оптимизируемой функции и адаптивной стратегии выбора кривизны параболической модели, см. § 5.

2) Любую выпуклую задачу можно сделать сильно выпуклой с помощью *регуляризации* (см. замечание 4.1), а любой алгоритм, настроенный на решение выпуклой задачи можно использовать для решения сильно выпуклой задачи за счет *рестартов*, см. упражнение 2.3 и конец параграфа 5, см. также приложение.

3) Стохастического оракула, выдающего градиент, можно свести к неточному (с малым шумом), но уже детерминированному оракулу с помощью *минибатчинга*, см. начало приложения. Идея приема: возвращение вместо стохастического градиента оптимизируемой функции в рассматриваемой точке среднее арифметическое независимых реализаций стохастических градиентов в этой же точке.

4) На основе конструкции *каталист* (см. замечание 3.3), в основе которой лежит проксимальный ускоренный градиентный метод, можно ускорять произвольные неускоренные методы, предназначенные для решения задач гладкой сильно выпуклой оптимизации. При этом получаются ускоренные методы, сходящиеся согласно нижним оценкам с точностью до логарифмических множителей. То есть в отличие от конструкций, описанных выше, в данной конструкции согласно теоретическим оценкам все же приходится «заплатить» логарифмический множитель за «общность».

Отметим также, что все описанные выше конструкции могут быть рассмотрены в общности § 4, 5, т.е. с более общей моделью и в прямо-двойственном контексте.

Для более комфортного изучения материала пособия рекомендуется предварительно познакомиться с основами выпуклого анализа, например, в объеме одной из следующих книг [62, 176, 503] и основами (вычислительной) линейной алгебры [91, 597].

Список литературы к пособию включает более 600 источников (при том, что мы далеко не всегда ссылались на первоисточники, в ряде случаев предпочитая более современные статьи и обзоры), поэтому вряд ли можно рассчитывать, что даже хорошо мотивированный читатель сможет ознакомиться с большей его частью. В этой связи для удобства выделим из этого списка учебники, изучение которых вместе с данным пособием можно рекомендовать в первую очередь:

I. *Boyd S., Vandenberghe L.* Convex optimization. – Cambridge University Press, 2004.
II. *Nocedal J., Wright S.* Numerical optimization. – Springer, 2006.
III. *Поляк Б.Т.* Введение в оптимизацию. – М.: URSS, 2014. – 392 с.
IV. *Bubeck S.* Convex optimization: algorithms and complexity // Foundations and Trends in Machine Learning. – 2015. – V. 8, N 3–4. – P. 231–357.

Стэнфордский учебник [I] является наглядным и одновременно строгим введением в выпуклую оптимизацию (теорию двойственности, *принцип множителей Лагранжа*, как следствие теоремы об отделимости гиперплоскостью граничной точки выпуклого множества от этого множества [62, п. 2.1], теоремы о дифференцировании функции максимума и т.п.), основы которой активно используются в настоящем пособии. Учебники [II, III] представляют собой достаточно подробное и хорошо проработанное описание основ численных методов оптимизации (выпуклой и не выпуклой). Во многом на базе именно этих двух учебников происходит обучение студентов основам численных методов оптимизации в большинстве продвинутых учебных заведениях по всему миру. Собранные в этих учебниках материалы отражают развитие данной области в основном в 60-80-е годы XX века. Более современные тенденции, связанные с развитием методов внутренней точки, ускорением методов и различными рандомизациями градиентных методов отражены в Принстонском учебнике [IV]. Этот учебник можно рекомендовать в качестве основного источника для последующего изучения.

Во вторую очередь (для более строгого изучения предмета) можно рекомендовать учебники

a) *Ben-Tal A., Nemirovski A.* Lectures on modern convex optimization analysis, algorithms, and engineering applications. – Philadelphia: SIAM, 2019.
b) *Nesterov Yu.* Lectures on convex optimization. – Springer, 2018.
c) *Lan G.* First-order and stochastic optimization methods for Machine Learning. – Springer, 2020.

◊ С. Бойд [173] является сейчас одним из самых цитируемых и активно публикующихся ученых в области численных методов оптимизации. С. Бойд имеет инженерное образование, и большое внимание в своих исследованиях уделяет практической составляющей, изящно сочетая её с фундаментальной. Оптимизационное сообщество практически едино во мнении, что работы С. Бойда (речь прежде всего о его книгах и документациях к разработанным под его руководством пакетам типа CVX [592]) являются хорошим образцом ясности изложения. Курс [I] является, пожалуй, самым известным (востребованным) в последнее десятилетие курсом по выпуклой оптимизации. ◊

Отметим, что настоящее пособие довольно сильно отличается и по отбору материалы и по форме изложения от подавляющего большинства известных нам учебников по оптимизации, в том числе и от выделенных четырех. Достаточно сказать, что в пособие не были включены ставшие уже классическими разделы про задачи линейного и полуопределенного программирования. Приведем здесь ссылки на то, как эти материалы в 2018/2019 учебном году преподавал А.С. Немировский студентам и аспирантам университета Джорджия в Атланте [437, 440]. Отметим также некоторые недавние достижения в этих областях [398, 400]. С другой стороны, почти половина из материалов пособия, по-видимому, впервые излагается (осмысляется) в учебном контексте.

В пособии имеется большое число ссылок на современную иностранную литературу. После распада Советского Союза «оптимизационный крен» сильно сместился на Запад. Однако считаем важным подчеркнуть определяющую роль российских ученых и научных школ [490] в создании того фундамента, на котором сейчас стоит молодая (чуть больше 60 лет), но бурно развивающаяся область знаний: «Численные методы оптимизации». На Западе даже есть такая вполне серьёзная шутка: «Если ты придумал новый численный метод оптимизации, не торопись радоваться, наверняка его уже знал какой-нибудь русский еще в 60-е годы прошлого века, и опубликовал, конечно, на русском языке». В частности, многое из того, что включено в данное пособие, было придумано нашими соотечественниками.

В 2004–2005 гг. автор, будучи студентом факультета управления и прикладной математики (ФУПМ) МФТИ, на базовой кафедре в ВЦ РАН слушал курс профессора В.Г. Жадана [46] по дополнительным главам численных методов оптимизации, оказавший заметное влияние на последующий интерес к этой области. В целом, стоит отметить большое влияние школы акад. Н.Н. Моисеева на формирование как базового, так и дополнительного цикла оптимизационных дисциплин на ФУПМ [10, 11, 46, 48, 49, 64, 65]. Современный учебный план студентов ФУПМ состоит из сочетания отмеченного опыта школы Н.Н. Моисеева и опыта коллег с ВМиК МГУ [13, 14, 53, 56, 88]. В данном пособии предпринята попытка посмотреть на этот учебный план, формировавшийся в течение полувека, сквозь призму современных достижений в области численных методов выпуклой оптимизации [72, 158, 180] и новых приложений [18, 39, 579]. Отметим также практикумы [595, 603] к упомянутому циклу лекций для студентов ФУПМ.

Автор также постарался учесть и обыграть в пособии некоторые наработки, которыми любезно с ним делились на всевозможных конференциях и семинарах представители различных научных школ: В.П. Булатова (Иркутск), В.Ф. Демьянова (Санкт-Петербург), И.И. Еремина (Екатеринбург), Л.В. Канторовича (Санкт-Петербург, Новосибирск, Москва), М.М. Лаврентьева (Новосибирск), А.А. Милютина (Москва), В.А. Скокова (Москва), А.Н. Тихонова (Москва), Я.З. Цыпкина (Москва), Н.З. Шора (Киев). Особенно, школ Ю.Г. Евтушенко (ВЦ РАН), Б.Т. Поляка (ИПУ РАН) и В.М. Тихомирова (мехмат МГУ). Вот уже более 10 лет автор имеет возможность обсуждать различные, связанные с оптимизацией, вопросы с Е.А. Нурминским, В.Ю. Протасовым, С.П. Тарасовым, С.В. Чукановым и А.А. Шананиным, А.Б. Юдицким.

Серьезное влияние на автора оказало регулярное общение с 2011 года с Б.Т. Поляком, А.С. Немировским и, особенно, с Ю.Е. Нестеровым. В большой части данный курс (пособие) был построен на расшифровке этих бесед. Автор очень благодарен трем оракулам за это.

Хотелось бы отметить важную роль, которую оказала совместная научная работа, выполняемая с А.Ю. Горновым, П.Е. Двуреченским, Ф.С. Стонякиным на данный текст.

Автор также выражает благодарность своему коллеге по кафедре Математических основ управления МФТИ доценту А.Г. Бирюкову за внимательное прочтение данной рукописи и предложенные исправления, а также Ф. Баху, Е.А. Воронцовой, К.В. Воронцову, А.И. Голикову, Н.В. Дойкову, Ю.В. Дорну, С.Э. Парсегову, А.О. Родоманову, Ф.Н. Рыбакову, Н. Сребро, А. Тейлору, Ц. Урибе, Р. Хильдебранду, А.В. Чернову, О. Шамиру за ряд ценных замечаний. На ряд неточностей автору было указано учениками: Артёмом Агафоновым, Кириллом Бобыревым, Эдуардом Горбуновым, Сергеем Гуминовым, Анастасией Ивановой, Дмитрием Камзоловым, Виктором Мишиным, Петром Остроуховым, Дмитрием Пасечнюком, Александром Рогозиным, Антоном Рябцевым, Даниилом Селихано-вичем, Игорем Соколовым, Александром Титовым, Александром Тюриным, Ильнурой Усмановой и Салихом Хабибуллиным.

Особую благодарность за постоянную поддержку хотелось бы выразить своей жене Даше Двинских.

Все возможные ошибки лежат всецело на авторе. В случае обнаружения неточностей просьба присылать информацию на адрес электронной почты <gasnikov.av@mipt.ru>.

На обложке изображен спуск горнолыжников на горе Монблан (февраль 2016 г.). Примеры «горных аналогий» в оптимизации см. в [96].

# § 1. Градиентный спуск

Рассмотрим задачу

$$f(x) \to \min_{x \in \mathbb{R}^n}. \tag{1.1}$$

Далее в этом параграфе приведены классические способы получения/понимания одного из основных инструментов современной вычислительной математики – *метода градиентного спуска*, восходящего к работам О. Коши [195], Л.В. Канторовича [55], Б.Т. Поляка [83]. Более современное изложение, в котором прорабатываются различные тонкие вопросы, начнется со следующего параграфа.

◊ Стоит особо отметить большой (наверное, можно даже сказать, решающий) вклад, который внес Б.Т. Поляк в 60-е годы XX века в развитие градиентных методов. Многие из современных методов и подходов, активно использующихся для решения задач оптимизации больших размеров, восходят к работам Бориса Теодоровича: усреднение Поляка [39, п. 8.7.3]; субградиентный метод Поляка [460]; *метод тяжелого шарика* (импульсный метод), породивший впоследствии целую линейку ускоренных градиентных методов, в частности, очень популярный в последние годы (быстрый, ускоренный, моментный) *градиентный метод Нестерова* (см. указание к упражнению 1.3). Собственно, знакомство с градиентными методами далее в пособии (особенно в § 1, § 2) осуществляется во многом под влиянием отмеченного цикла работ Б.Т. Поляка [82]. ◊

Рассмотрим систему обыкновенных дифференциальных уравнений [195]

$$\frac{dx}{dt} = -\nabla f(x). \tag{1.2}$$

Покажем, что значения функции $W(x) = f(x)$ убывают на траекториях динамической системы (1.2), т. е. $W(x)$ будет *функцией Ляпунова* системы (1.2). Действительно,

$$\frac{dW(x(t))}{dt} = \left\langle \nabla f(x(t)), \frac{dx(t)}{dt} \right\rangle = \left\langle \nabla f(x(t)), -\nabla f(x(t)) \right\rangle = -\left\| \nabla f(x(t)) \right\|_2^2 \le 0,$$

$$\frac{dW(x)}{dt} = 0 \iff \nabla f(x) = 0.$$

Отсюда можно сделать вывод, что любая траектория такой системы должна сходиться к *стационарной точке*[3] функции $f(x)$, вообще говоря, зависящей от точки старта (на рис. 1 рассмотрен случай выпуклой функции). Аналогичного свойства можно ожидать и от дискретизованной по схеме Эйлера версии динамики (1.2)

$$x^{k+1} = x^k - h\nabla f\left(x^k\right), \tag{1.3}$$

в случае достаточно малого шага $h$ [82, гл. 2], см. также [5, 13, 48, 70, 145, 232, 310, 431, 517, 544, 574]. Метод (1.3) обычно называют *методом градиентного спуска* или просто *градиентный спуск* [82], а приведенный здесь способ получения оценки скорости сходимости метода относят ко *второму методу Ляпунова* [82, § 2, гл. 2].

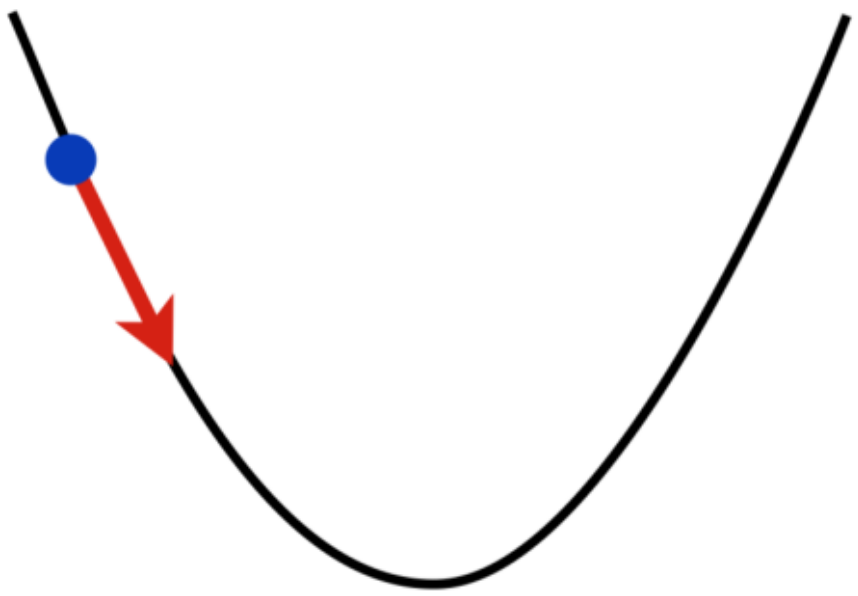

Рис. 1

Чтобы количественно оценить скорость сходимости и получить условие на выбор шага сделаем следующее предположение о *липшицевости градиента* в 2-норме [82, гл. 1]: для любых $x$ и $y$ имеет место неравенство

$$\left\|\nabla f(y) - \nabla f(x)\right\|_2 \le L\left\|y - x\right\|_2. \tag{1.4}$$

Из этого неравенства имеем[4] $\lambda_{\max}\left(\nabla^2 f(x)\right) \le L$, т. е. все собственные значения *матрицы Гессе* $\nabla^2 f(x) = \left\|\partial^2 f(x)/\partial x_i \partial x_j\right\|_{i,j=1}^n$ не больше $L$. По формуле Тейлора с остаточным членом в форме Лагранжа для любых $x$ и $y$ справедливо представление [89, § 58]:

$$f(y) = f(x) + \langle \nabla f(x), y - x \rangle + \frac{1}{2}\langle \nabla^2 f(\tilde{x})(y - x), y - x \rangle,$$

где $\tilde{x} = \tilde{x}(x, y)$ принадлежит отрезку, соединяющему $x$ и $y$. Отсюда можно получить, что для любых $x$ и $y$ выполняется неравенство

$$f(y) \le f(x) + \langle \nabla f(x), y - x \rangle + \frac{L}{2}\|y - x\|_2^2. \tag{1.5}$$

---

[3] Напомним, что стационарной называют такую точку, в которой $\nabla f(x) = 0$.

[4] Строго говоря, из (1.4) выписанное неравенство следует лишь при дополнительном предположении о гладкости оптимизируемой функции. Однако основное неравенство (1.5) может быть получено и непосредственно из (1.4), см. [72, п. 1.2.2].

Из неравенства (1.5) следует, что

$$f(x^{k+1}) \le f(x^k) - h\langle \nabla f(x^k), \nabla f(x^k) \rangle + \frac{Lh^2}{2}\|\nabla f(x^k)\|_2^2 =$$

$$= f(x^k) - h \cdot \left(1 - \frac{Lh}{2}\right)\|\nabla f(x^k)\|_2^2.$$

Выбирая

$$h = \arg\max_{\alpha \ge 0} \alpha \cdot \left(1 - \frac{L\alpha}{2}\right) = \frac{1}{L}, \tag{1.6}$$

получим

$$f(x^{k+1}) \le f(x^k) - \frac{1}{2L}\|\nabla f(x^k)\|_2^2. \tag{1.7}$$

Отсюда, обозначая $x^k = x$ и учитывая, что $f(x^{k+1}) \ge f(x_*)$, где $x_*$ – решение задачи (1.1), получим полезное в дальнейшем неравенство

$$\frac{1}{2L}\|\nabla f(x)\|_2^2 \le f(x) - f(x_*). \tag{1.8}$$

Из неравенства (1.7) следует: для достижения

$$\min_{k=1,\dots,N}\|\nabla f(x^k)\|_2 \le \varepsilon, \tag{1.9}$$

достаточно [82, 447, 579]

$$N = \frac{2L \cdot \left(f(x^0) - f(x^{extr})\right)}{\varepsilon^2} \tag{1.10}$$

итераций метода (1.3) с шагом (1.6). Здесь $\nabla f(x^{\text{extr}}) = 0$ и $x^{\text{extr}}$, вообще говоря, зависит от точки старта $x^0$. Действительно, до момента выполне-ния (1.9) на каждой итерации согласно (1.8) происходит уменьшение значения целевой функции $f(x)$ как минимум на $\varepsilon^2/(2L)$. Таким образом, не более чем после

$$N = \frac{f(x^0) - f(x^{extr})}{\varepsilon^2/(2L)}$$

итераций условие (1.9) должно выполниться первый раз.

Оценка (1.10) на классе функций, удовлетворяющих условию (1.4), с точностью до мультипликативной константы не может быть улучшена (если $n$ не мало [185]) как для метода вида (1.3), так и для любых других методов первого порядка, т. е. использующих только градиент функции [191, 192].

◊ Здесь и далее «с точностью до мультипликативной константы» означает, что в оценке числа итераций можно попробовать улучшить числовой множитель, но не зависимость от параметров задачи. Стоит также отметить, что сделанные оговорки «для класса функций (1.4)» и «для методов первого порядка» – существенные, см., например, [190, 191, 192, 193, 349, 447]. В частности, для достаточно гладких функций $f(x)$ (Липшицев гессиан и т.д.) в классе методов первого порядка (использующих только $f(x)$ и $\nabla f(x)$) можно ожидать улучшение оценки $N \sim \varepsilon^{-2}$ до $N \sim \varepsilon^{-8/5}$ [192], если при этом разрешается использовать в методе старшие производные оптимизируемой функции до порядка $p \geq 1$ включительно, то можно улучшить оценки до $N \sim \varepsilon^{-(p+1)/p}$ [166, 191]. Относительно сложности итерации таких методов при $p = 2$, см. [300, 464, 568] и цитированную там литературу. При $p \geq 3$ такие методы интересны в основном только в теоретическом плане. В упражнении 1.8 3) будет продемонстрировано как можно улучшить оценку (1.10), если $f(x)$ имеет вид суммы функций.

Заметим (см. [70], а также упражнение 1.3), что для задач выпуклой оптимизации дополнительные предположения о Липшицевости старших производных не меняют нижних оценок. Впрочем, если вводить более сильные предположения (самосогласованной) Липшицевости (гладкости) старших производных, то и для задач выпуклой оптимизации в классе методов первого порядка возможно дополнительное ускорение [573].

Отметим также, что в общем случае для функций из класса (1.4) необходимое число итераций (на каждой итерации можно в одной и только одной точке, получить значение функции $f(x)$ и ее старших произ-водных) для поиска такого $x^N$, что $f(x^N) - f(x_*) \le \varepsilon$, где $x \in \mathbb{R}^n$, зависит от $\varepsilon$ существено хуже: $N \sim \varepsilon^{-n/2}$, см. [70, § 6 гл. 1]. Общая идея получения такого типа нижних оценок в наиболее простом виде изложена, например, в [72, теорема 1.1.2]. Допустимое множество разбивается на кубики со стороной $4\sqrt{2\varepsilon/L}$. Под любой метод (алгоритм) подбирается такая функция, которая во всех кубиках, кроме одного, тождественно равна нулю, а в одном кубике, том самом, который данный метод будет просматривать в последнюю очередь, функция «проваливается» (с сохранением условия (1.4)) на глубину $2\varepsilon$. Поиск $x^N$ для так построенной (под рассматриваемый метод) функции потребует «просмотра» всех кубиков, число которых $\sim \left(1/\sqrt{\varepsilon}\right)^n = \varepsilon^{-n/2}$.

Заметим, что далее в пособии будут приводиться примеры универсально плохих функций («худших в мире функций») для рассматриваемых классов задач сразу для всех допустимых алгоритмов, см. упражнения 1.3, 2.1 и приложение. Такую функцию можно построить и в данном случае. Рассмотрим следующую функцию Нестерова–Скокова (другие примеры «плохих» функций см., например, в [76]), обобщающую популярную тестовую функцию Розенброка [33, 37, 82]:

$$f(x) = \frac{1}{4}(x_1 - 1)^2 + \sum_{i=1}^{n-1}\left(x_{i+1} - 2x_i^2 + 1\right)^2 =$$

$$= \frac{1}{4}(x_1 - 1)^2 + \sum_{i=1}^{n-1}(x_{i+1} - \underbrace{P_2(x_i)}_{\text{многочлен Чебышёва}})^2.$$

Эта функция имеет единственный экстремум $x_* = (1,1,\ldots,1)$ ($f(x_*) = 0$), который и есть глобальный минимум. Стартуя с $x^0 = (-1,1,\ldots,1)^T$ ($f(x^0) = 1$, $\left\|\nabla f(x^0)\right\|_2 = 1$) метод наискорейшего спуска (см. замечание 1.4) со временем (на поздних итерациях) обеспечивают малость нормы градиента, что хорошо согласуется с описанной выше теорией, однако, при этом не наблюдается сходимости по функции. Так, при $n = 16$ для градиентного спуска в момент, когда $\left\|\nabla f(x^N)\right\|_2 \approx 10^{-9}$ имеем $f(x^N) - f(x_*) \approx 0.98$. Похожая ситуация имеет место и для многих других популярных на практике методов, например, для метода сопряженных градиентов (см. замечание 1.6) и LBFGS (см. замечание 3 приложения). Особенность рассмотренной функции заключается в экспоненциально осциллирующих оврагах, возникающих в качестве множеств уровня функции. Связано это со следующим свойством многочленов Чебышёва $P_n\left(P_m\left(x\right)\right)=P_{mn}\left(x\right)$ и экспоненциальным ростом числа осцилляций у многочлена $P_{2^i}\left(x\right)$ с ростом $i$ [422]. Фиксируя $x_1=c$ можно заметить, что минимум $f\left(x\right)$ при $x_1=c$ достигается в точке $x_{i+1}\left(c\right)=P_{2^i}\left(x_1\right)=P_{2^i}\left(c\right)$. Даже при близких значениях $c$ старшие координаты $x_{i+1}\left(c\right)$ могут значительно отличаться.

Обратим внимание, что рассмотренные выше функции являются искусственно придуманными. Однако большие сложности, связанные с невыпуклостью (см. ниже) и многоэкстремальностью, часто возникают и в реальных приложениях, например, при обучении нейронных сетей [39] или белковом фолдинге [102]. К счастью, во многих приложениях бывает достаточно найти «хороший» локальный минимум или даже просто локальный минимум, см., например, [282]. ◊

В общем случае, полученный выше результат о сходимости градиентного спуска к экстремальной точке не гарантирует его сходимости даже к локальному минимуму [72, пример 1.2.2]. Впрочем, недавно было показано [396, 397], что метод (1.3) с шагом (1.6) типично сходится именно к локальному минимуму, см. также [127]. Это по-прежнему не означает сходимость к глобальному минимуму.

Если задача (1.1) является *задачей выпуклой оптимизации*, т. е. $f\left(x\right)$ – *выпуклая функция*, что следует из $\lambda_{\min}\left(\nabla^2 f\left(x\right)\right)\geq 0$, то можно гарантировать сходимость метода (1.3) с шагом (1.6) к глобальному минимуму в следующем смысле [72, следствие 2.1.2]:

$$f\left(x^N\right)-f\left(x_*\right)\leq\frac{2LR^2}{N+4}, \tag{1.11}$$

где $x_*$ – решение задачи (1.1), $R^2=\left\|x^0-x_*\right\|_2^2$. Если решение не единственно, то под $x_*$ в (1.11) можно понимать такое решение задачи (1.1), которое наиболее близко в 2-норме к точке старта $x^0$ [56, 72, 82].

◊ В общем случае $f\left(x\right)$ – выпуклая функция, означает, что надграфик $f\left(x\right)$ – выпуклое множество (см. рис. 1). Множество $Q$ – выпуклое, если вместе с любыми двумя своими точками оно содержит отрезок, их соединяющий. Это определение эквивалентно тому, что любая граничная точка множества $Q$ отделима от этого множества, т. е. существует такая разделяющая (опорная) гиперплоскость, касающаяся множества $Q$ в рассматриваемой точке, что множество $Q$ лежит по одну сторо-ну от этой гиперплоскости [62, п. 1.2, 1.3]. В таком виде далее в основном и будет использоваться понятие выпуклости – см. неравенство (1.17). Отметим, что это неравенство верно и для негладких выпуклых функций, если под $\nabla f(x)$ понимать произвольный элемент субдифференциала $\partial f(x)$ [62, п. 1.5]. ◊

Из (1.7) следует, что сходимость в смысле (1.11) влечет сходимость в смысле (1.9). Рассматривая функции скалярного аргумента вида $f_M(x) = x^M$, $M \gg 1$, можно заметить, что, *сходимость по функции*, т. е. в смысле (1.11), не влечет в общем случае *сходимость по аргументу*. Точнее говоря, влечет, но скорость сходимости по аргументу может быть сколь угодно медленной.[5]

Если $f(x)$ – *$\mu$-сильно выпуклая функция* в 2-норме, что следует из $\lambda_{\min}\left(\nabla^2 f(x)\right) \geq \mu$, $\mu > 0$, то для метода (1.3) с шагом (1.6) уже будет иметь место *линейная сходимость* (сходимость со скоростью геометрической прогрессии), причем по аргументу [180, теорема 3.10]:

$$\left\|x^N - x_*\right\|_2^2 \leq R^2 \exp\left(-\frac{\mu}{L}N\right), \tag{1.12}$$

где $x_*$ – решение задачи (1.1), т. е. $\nabla f(x_*) = 0$.

Поясним, каким образом можно прийти к формуле типа (1.12). Для этого заметим, что условие $\lambda_{\min}\left(\nabla^2 f(x)\right) \geq \mu$ влечет (см. вывод (1.5) и [72, п. 2.1.3]) следующее условие для любых $x$ и $y$:

$$f(y) \geq f(x) + \left\langle \nabla f(x), y - x \right\rangle + \frac{\mu}{2}\left\|y - x\right\|_2^2. \tag{1.13}$$

Обычно условие (1.13) и понимают как определение $\mu$-сильной выпуклой функции в 2-норме [72, п. 2.1.3].

Из неравенства (1.13) получается полезное в дальнейшем неравенство

$$\frac{\mu}{2}\left\|x - x_*\right\|_2^2 \leq f(x) - f(x_*). \tag{1.14}$$

[5] Еще точнее, здесь надо ограничить класс используемых методов. При использовании методов типа деления отрезка пополам в условиях абсолютной точности вычислений можно сходиться по аргументу и для таких (вырожденных) примеров (см. упражнение 1.4).

В частности, (1.14) можно использовать для получения неравенств вида (1.12) из неравенств вида (1.16).

◊ Если (1.13) имеет место только для всех $x, y \in Q$ (или $\lambda_{\min}\left(\nabla^2 f(x)\right) \geq \mu$, для всех $x, y \in Q$), где $Q$ – выпуклое множество, а $x_* = \arg\min_{x \in Q} f(x)$, то неравенство (1.14) будет иметь место для всех $x \in Q$. ◊

Из (1.13) следует, что

$$f(x_*) = \min_y f(y) \geq \min_y \left\{ f(x) + \langle \nabla f(x), y - x \rangle + \frac{\mu}{2}\|y - x\|_2^2 \right\} =$$

$$= f(x) - \frac{1}{2\mu}\|\nabla f(x)\|_2^2,$$

т. е.

$$f(x) - f(x_*) \leq \frac{1}{2\mu}\|\nabla f(x)\|_2^2. \tag{1.15}$$

Отсюда, с учетом неравенства (1.7), имеем

$$f(x^{k+1}) - f(x^k) \leq -\frac{1}{2L}\|\nabla f(x^k)\|_2^2 \leq -\frac{\mu}{L} \cdot \left( f(x^k) - f(x_*) \right),$$

т. е.

$$f(x^{k+1}) - f(x_*) \leq \left(1 - \frac{\mu}{L}\right)\left( f(x^k) - f(x_*) \right).$$

Следовательно,

$$f(x^N) - f(x_*) \leq \left(1 - \frac{\mu}{L}\right)^N \left( f(x^0) - f(x_*) \right) \leq \exp\left(-\frac{\mu}{L}N\right)\left( f(x^0) - f(x_*) \right). \tag{1.16}$$

◊ Если вместо настоящего градиента доступен зашумленный градиент в концепции относительной точности[6] $\tilde{\nabla} f(x)$ [83, п. 2 § 1; п. 3 § 2 гл. 4]

$$\|\tilde{\nabla} f(x) - \nabla f(x)\|_2 \leq \alpha \|\nabla f(x)\|_2,$$

[6] Существуют разные концепции шума в градиенте, см., например, [82, гл. 4], [207, 216, 225]. Далее в пособии мы в основном будем работать с концепцией из [225, 227, 228], поскольку с помощью этой конструкции удобно строить универсальные методы, см. § 5.

где $\alpha \in [0,1)$, то выбирая $h$ в (1.3)

$$x^{k+1} = x^k - h\tilde{\nabla} f\left(x^k\right)$$

не по формуле (1.6), а следующим образом

$$h = \frac{1}{L}\frac{1-\alpha}{(1+\alpha)^2},$$

получим вместо (1.7)

$$f\left(x^{k+1}\right) \le f\left(x^k\right) - \frac{1}{2L}\frac{(1-\alpha)^2}{(1+\alpha)^2}\left\|\nabla f\left(x^k\right)\right\|_2^2,$$

что приведет в итоге вместо (1.16) к неравенству

$$f\left(x^N\right) - f\left(x_*\right) \le \left(1 - \frac{\mu}{L}\frac{(1-\alpha)^2}{(1+\alpha)^2}\right)^N \left(f\left(x^0\right) - f\left(x_*\right)\right).$$

Если вместо (1.15) предполагать (1.13), то приведенную оценку можно уточнить [222]

$$\frac{(1-\alpha)^2}{(1+\alpha)^2} \to \mathrm{O}\left(\frac{1-\alpha}{1+\alpha}\right). \ \Diamond$$

**Замечание 1.1 (условие градиентного доминирования).** Формула (1.16) была получена в предположениях (1.4), (1.15). То есть предположение $\lambda_{\min}\left(\nabla^2 f(x)\right) \ge \mu$ на самом деле использовалось лишь в виде своего следствия (1.15). Условие (1.15) называют *условием градиентного доминирования* или *условием Поляка–Лоясиевича* [360, 464]. Приведем пример, когда это условие имеет место, однако нельзя быть уверенным даже в выпуклости $f(x)$ [54], [71, п. 4.3], [492]. Рассмотрим систему нелинейных уравнений $g(x) = 0$, записанную в векторном виде, т. е. $g: \mathbb{R}^n \to \mathbb{R}^m$, $m \le n$. Требуется найти какое-нибудь решение этой системы. Введем матрицу Якоби отображения $g$: $\partial g(x)/\partial x = \left\| \partial g_i(x)/\partial x_j \right\|_{i,j=1}^{m,n}$. Предположим, что существует такое $\mu > 0$, что для всех $x \in \mathbb{R}^n$ имеет место равномерная невырожденность матрицы Якоби:

$$\lambda_{\min}\left(\partial g(x)/\partial x \cdot \left[\partial g(x)/\partial x\right]^T\right) \ge \mu.$$

Тогда функция $f(x)=\|g(x)\|_2^2$ удовлетворяет условию (1.15) для произвольного $x_*$ такого, что $f(x_*)=0$, т.е. $g(x_*)=0$ [464].

Рассмотренная выше тестовая функция Нестерова–Скокова возникает при сведении системы нелинейныйх уравнений вида: $x_1=1$, $x_2=2x_1^2-1$, …, $x_n=2x_{n-1}^2-1$, к задаче оптимизации. Матрица Якоби такой системы получается нижнетреугольной, поэтому условие Поляка–Лоясиевича выполняется. Однако, это не помогает численно минимизировать возникающую функцию, поскольку обусловленность данной задачи оптимизации $L/\mu$ получается экспоненциально большой по $n$.

Отметим также, что минимизировать функционал вида $f(x)=\|g(x)\|_2^2$ при условии Поляка–Лоясиевича можно не только градиентным методом, но и, например, методом Гаусса–Ньютона–Нестерова [71, 452]. Этот метод является методом первого порядка и имеет глобальную скорость сходимости, аналогичную скорости сходимости градиентного спуска, однако локальная скорость сходимости оказывается сверхлинейной (как у метода Ньютона, см. приложение). Метод Гаусса–Ньютона–Нестерова неплохо себя проявил в численных экспериментах с функцией Нестерова–Скокова.

Интересный пример возникновения задач оптимизации с условием Поляка–Лоясиевича недавно был обнаружен в теории управления [269].

Многие результаты, полученные для сильно выпуклой оптимизации, могут быть естественным образом перенесены на случай выполнения условия Поляка–Лоясиевича. В качестве популярного примера, упомянем здесь результаты о скорости сходимости метода стохастического градиентного спуска [112]. ■

Для дальнейшего построения «линейки» основных методов нам будет полезно немного по-другому посмотреть на метод градиентного спуска.

Прежде всего, заметим, что если $f(x)$ – *выпуклая функция* (см. условие (1.13) при $\mu=0$), т. е. для любых $x$ и $y$

$$f(x)+\langle\nabla f(x),y-x\rangle\le f(y), \tag{1.17}$$

то $W(x)=\|x-x_*\|_2^2/2$ также будет функцией Ляпунова системы (1.2).[7] Действительно,

---

[7] Современные приложения аппарата квадратичных функций Ляпунова для получения оценок скорости сходимости численных методов выпуклой оптимизации, см., например, в работах [141, 145, 549, 553].

$$\frac{dW\left(x(t)\right)}{dt} = \left\langle x(t) - x_*, \frac{dx(t)}{dt} \right\rangle =$$

$$= -\left\langle \nabla f\left(x(t)\right), x(t) - x_* \right\rangle \le f\left(x_*\right) - f\left(x(t)\right) \le 0,$$

$$\frac{dW(x)}{dt} = 0 \quad \Leftrightarrow \quad x \in \operatorname{Arg}\min_{x \in \mathbb{R}^n} f(x).$$

Сделанное наблюдение «подсказывает» исследовать поведение последовательности

$$\frac{1}{2}\left\| x^k - x_* \right\|_2^2.$$

Согласно (1.3) имеем,

$$\frac{1}{2}\left\| x^{k+1} - x_* \right\|_2^2 = \frac{1}{2}\left\| \left(x^k - x_*\right) - h\nabla f\left(x^k\right) \right\|_2^2 =$$
$$= \frac{1}{2}\left\| x^k - x_* \right\|_2^2 - h\left\langle \nabla f\left(x^k\right), x^k - x_* \right\rangle + \frac{h^2}{2}\left\| \nabla f\left(x^k\right) \right\|_2^2. \tag{1.18}$$

Следовательно,

$$f\left(x^k\right) - f\left(x_*\right) \overset{1}{\le} \left\langle \nabla f\left(x^k\right), x^k - x_* \right\rangle \overset{2}{\le}$$
$$\overset{2}{\le} \frac{1}{2h}\left\| x^k - x_* \right\|_2^2 - \frac{1}{2h}\left\| x^{k+1} - x_* \right\|_2^2 + Lh \cdot \left( f\left(x^k\right) - f\left(x^{k+1}\right) \right), \tag{1.19}$$

где неравенство 1 вытекает из (1.17), а неравенство 2 – из равенства (1.18) и неравенства (1.7), предполагающего, что $h = 1/L$. Суммируя (1.19) по $k = 0, \ldots, N-1$ и подставляя $h = 1/L$, получим

$$\sum_{k=0}^{N-1}\left( f\left(x^k\right) - f\left(x_*\right) \right) \le \frac{1}{2h}\left\| x^0 - x_* \right\|_2^2 - \frac{1}{2h}\left\| x^N - x_* \right\|_2^2 + Lh \cdot \left( f\left(x^0\right) - f\left(x^N\right) \right) \le$$
$$\le \frac{1}{2h}\left\| x^0 - x_* \right\|_2^2 + Lh \cdot \left( f\left(x^0\right) - f\left(x^N\right) \right) \overset{h=1/L}{=} \frac{LR^2}{2} + f\left(x^0\right) - f\left(x^N\right).$$

Отсюда следует, что

$$\frac{1}{N}\sum_{k=1}^{N}\left( f\left(x^k\right) - f\left(x_*\right) \right) \le \frac{LR^2}{2N}.$$

◊ Альтернативным определением выпуклой функции является *неравенство Иенссена*: для всех $x$ и $y$ и произвольного $\alpha \in [0,1]$:

$$f\left(\alpha x+(1-\alpha)y\right)\le \alpha f(x)+(1-\alpha)f(y),$$

которое по индукции влечет неравенство [176, гл. 3]:

$$f\left(\frac{1}{N}\sum_{k=1}^{N}x^k\right)\le \frac{1}{N}\sum_{k=1}^{N}f\left(x^k\right). \Diamond$$

Таким образом, ввиду выпуклости $f(x)$ имеет место неравенство

$$f\left(\bar{x}^N\right)-f\left(x_*\right)\le \frac{LR^2}{2N}, \qquad (1.20)$$

где

$$\bar{x}^N=\frac{1}{N}\sum_{k=1}^{N}x^k, \qquad (1.21)$$

являющееся аналогом неравенства (1.11).

Резюмируем приведенные выше результаты в немного более точной и симметричной форме [225, 242, 447, 552].

**Теорема 1.1.** *Пусть для численного решения задачи* (1.1)

$$f(x)\to \min_{x\in\mathbb{R}^n},$$

*с функцией* $f(x)$, *удовлетворяющей условию* (1.4), *используется градиентный спуск* (1.3), (1.6):

$$x^{k+1}=x^k-\frac{1}{L}\nabla f\left(x^k\right). \qquad (1.22)$$

*Тогда*

$$\min_{k=1,\ldots,N}\left\|\nabla f\left(x^k\right)\right\|_2\le\sqrt{\frac{2L\cdot\left(f\left(x^0\right)-f\left(x_*\right)\right)}{N}}. \qquad (1.23)$$

*Если дополнительно известно, что* $f(x)$ – $\mu$*-сильно выпуклая функция в* 2*-норме, где* $\mu\ge 0$, *то*[8]

$$\min_{k=1,\ldots,N}f\left(x^k\right)-f\left(x_*\right)\le\frac{LR^2}{2}\min\left\{\frac{1}{N},\exp\left(-\frac{\mu}{L}N\right)\right\}. \qquad (1.24)$$

---

[8] Ввиду (1.7) имеем $\min_{k=1,\ldots,N}f\left(x^k\right)=f\left(x^N\right)$. Однако, начиная со следующего параграфа, в котором допускается наличие неточности (2.3) и более общие способы «проектирования» (2.29) (по сравнению с обычным евклидовым), форма записи (1.24) уже будет по существу.

Приведенные в теореме 1.1 оценки скорости сходимости метода (1.22) точные. Немного могут быть улучшены только числовые множители [191, 192, 242, 552]. При получении оценки (1.23) (впрочем, как и оценки (1.10)) существенным образом использовалось, что рассматривается задача безусловной оптимизации (1.1) [447].

Вместо полного доказательства теоремы 1.1 ниже приводится наглядная интерпретация неравенства (1.7), лежащего в основе доказательства теоремы.

**Замечание 1.2 («геометрия» градиентного спуска).** Если понимать градиентный спуск (1.3) с шагом (1.6) следующим образом:

$$x^{k+1}=x^k-\frac{1}{L}\nabla f\left(x^k\right)=\arg\min_{x\in\mathbb{R}^n}\left\{f\left(x^k\right)+\left\langle\nabla f\left(x^k\right),x-x^k\right\rangle+\frac{L}{2}\left\|x-x^k\right\|_2^2\right\},\quad(1.25)$$

то метод имеет естественную геометрическую интерпретацию. Параболоид вращения

$$\bar{f}_{x^k}(x)=f\left(x^k\right)+\left\langle\nabla f\left(x^k\right),x-x^k\right\rangle+\frac{L}{2}\left\|x-x^k\right\|_2^2$$

касается графика функции $f(x)$ в точке $x^k$ и мажорирует ее на всем пространстве

$$f(x)\le\bar{f}_{x^k}(x)\quad\text{для всех}\quad x\in\mathbb{R}^n.$$

В частности,

$$f\left(x^{k+1}\right)\le\bar{f}_{x^k}\left(x^{k+1}\right)=\min_{x\in\mathbb{R}^n}\bar{f}_{x^k}(x).$$

Но по построению $\bar{f}_{x^k}(x)$ имеем $\bar{f}_{x^k}\left(x^k\right)=f\left(x^k\right)$. Значит, переходя от точки $x^k$ к точке минимума параболоида $x^{k+1}$, мы «выедаем» у функции $f(x)$ не меньше, чем у $\bar{f}_{x^k}(x)$ (см. рис. 2), т. е. не меньше, чем

$$\bar{f}_{x^k}\left(x^k\right)-\bar{f}_{x^k}\left(x^{k+1}\right)=\frac{1}{2L}\left\|\nabla f\left(x^k\right)\right\|_2^2.$$

Таким образом можно прийти к основному соотношению (1.7).

Другая интерпретация имеется, например, в [82, формулы (2), (3) п. 1 § 4, гл. 1], см. также [42, п. 6.4], [53, § 5.2], [209], [472, гл. 4]. В некоторой $r_k$-окрестности точки $x^k$ функция $f(x)$ заменяется линейной функцией (или более сложной моделью)

$$\tilde{f}_{x^k}(x)=f\left(x^k\right)+\left\langle\nabla f\left(x^k\right),x-x^k\right\rangle.$$

Новое положение метода определяется исходя из решения задачи

$$x^{k+1} = \arg\min_{\left\|x-x^k\right\|_2 \le r_k} \tilde{f}_{x^k}(x).$$

С помощью *принципа множителей Лагранжа* [176, гл. 5] данная задача сводится к задаче (1.25), где $L/2$ следует понимать как множитель Лагранжа к ограничению $\left\|x-x^k\right\|_2^2 \le r_k^2$. Такой подход получил название *метода доверительной области* (*trust region*). Вместе с квадратичной (ньютоновской) моделью функции его активно использовали в качестве подхода, *глобализующего сходимость* (см., например, [53, § 5.2]): обеспечивающего попадание исследуемого метода в область сверхлинейной (квадратичной) сходимости [42, п. 6.4], [209], [472, гл. 4]. Однако в последнее десятилетие данный подход в таком «глобализующем» контексте стал частично вытесняться *методом Ньютона с кубической регуляризацией* [444, 464] (см. также приложение), лучше изученным в теоретическом плане. Этот метод можно понимать как перезапись метода доверительной области с квадратичной моделью и ограничением вида $\left\|x-x^k\right\|_2^3 \le r_k^3$, которое заносится в функционал с помощью принципа множителей Лагранжа, что приводит к задаче квадратичной оптимизации с дополнительным кубическим штрафным слагаемым.

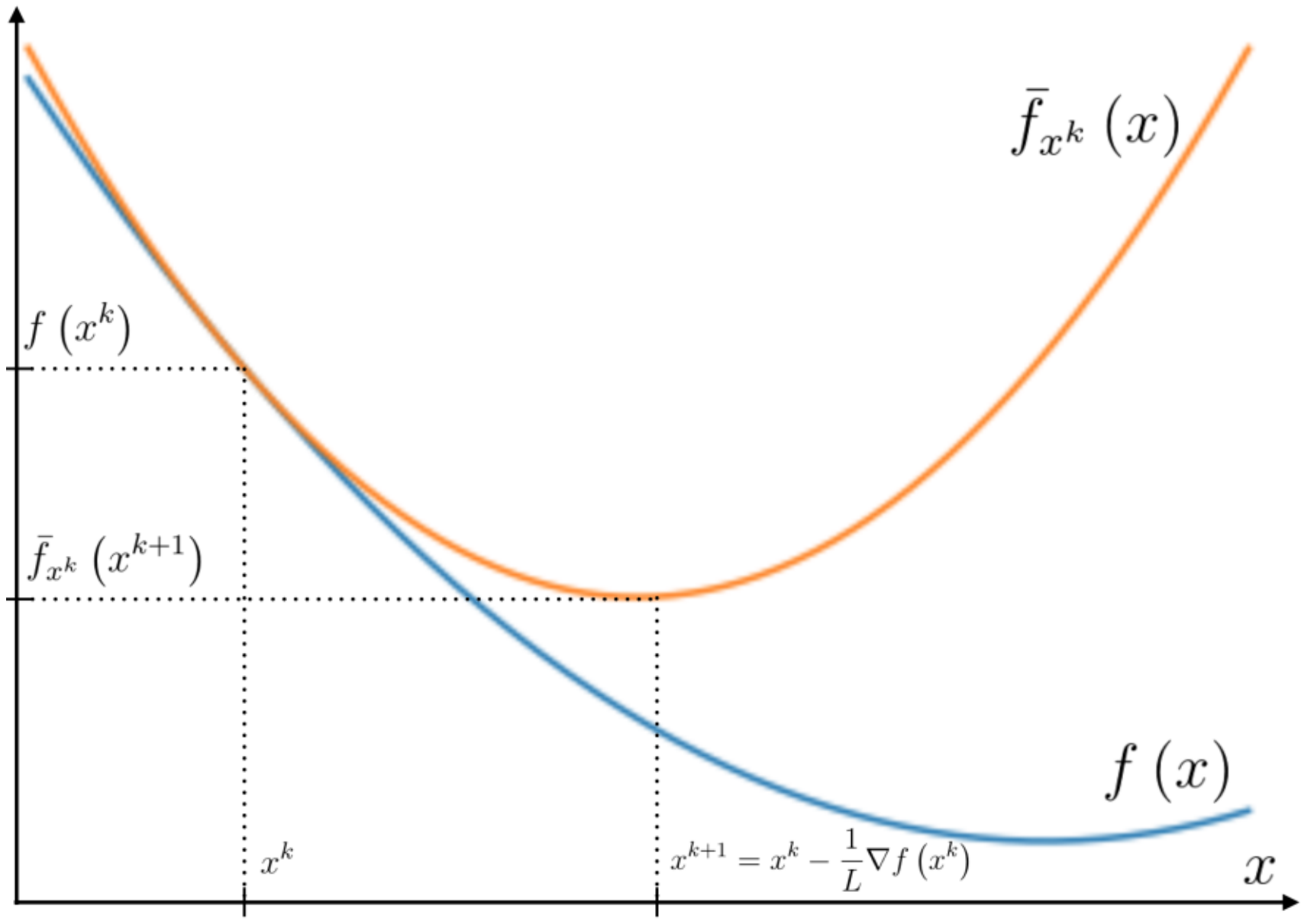

Рис. 2

Отметим также, что если рассматривается задача условной оптимизации на выпуклом множестве $Q$ (см. § 2), то написанные выше интерпретации сохраняются. Причем в случае, когда множество $Q$ компактно, можно выбирать, в частности, $r_k = \infty$. Получившийся в результате метод будет принадлежать к классу *методов условного градиента* [13, 82, 154, 180, 342, 386]. Получившийся метод, как и стандартный метод условного градиента (также используется название *метод Франк–Вульфа*), имеет оценки скорости сходимости на классе гладких выпуклых задач, в целом аналогичные оценкам для обычного градиентного метода [342]. Однако вместо проектирования на $Q$ (см. § 2) на каждой итерации метода необходимо решать вспомогательную задачу минимизации линейного функционала на множестве $Q$. В случае когда $Q$ – симплекс (или шар в 1-норме), на каждой итерации получается разреженное решение вспомогательной задачи (в одной из вершин симплекса), что позволяет существенно уменьшать стоимость итерации, см. упражнение 1.6, а также [3, 18, 180, 211].

Интересный взгляд на метод условного градиента предложил А.С. Немировский [95], [158, п. 5.5.3], см. также [386, chapter 7]. Оказыва-ется, такого типа методы можно также получать, беря за основу быстрые градиентные методы в концепции модели функции (см. § 3) с неточным проектированием (см. упражнение 3.7): вместо задачи (3.3) на каждой итерации решается более простая задача – (3.3) с $1/h \equiv 0$ (в обозначениях упражнения 3.7 $1/\alpha_{k+1} \equiv 0$ ), и решение этой упрощенной задачи интерпретируется как приближенное решение исходной задачи (3.3). ■

**Замечание 1.3 (градиентный спуск в *p*-норме).** Используя приведенную выше схему рассуждений, попробуем распространить метод (1.3) на случай, когда условие (1.4) имеет более общий вид

$$\left\|\nabla f(y)-\nabla f(x)\right\|_* \le L\left\|y-x\right\|, \tag{1.26}$$

где $\left\|y\right\|_* = \max\limits_{\|x\|\le 1}\langle y, x\rangle$ – сопряженная норма к норме $\|\ \|$. В этом случае неравенство (1.5) будет иметь аналогичный вид [122]:

$$f(y) \le f(x) + \langle\nabla f(x), y-x\rangle + \frac{L}{2}\left\|y-x\right\|^2.$$

Тогда естественно заменить метод (1.3) с шагом (1.6) следующим методом:

$$x^{k+1} = \arg\min_{x\in\mathbb{R}^n}\left\{ f\left(x^k\right) + \left\langle\nabla f\left(x^k\right), x-x^k\right\rangle + \frac{L}{2}\left\|x-x^k\right\|^2\right\}. \tag{1.27}$$

Аналог неравенства (1.7) будет иметь вид [122]:

$$f\left(x^{k+1}\right) \le f\left(x^k\right) - \frac{1}{2L}\left\|\nabla f\left(x^k\right)\right\|_*^2.$$

Отсюда можно получить следующую оценку скорости сходимости [122]:

$$f\left(x^N\right) - f\left(x_*\right) \le \frac{2L\tilde{R}^2}{N}, \tag{1.28}$$

где $\tilde{R} = \max\limits_{x:\ f(x)\le f\left(x^0\right)}\left\|x-x_*\right\|$. В случае $\|\ \| = \|\ \|_2$ оценку $\tilde{R}$ можно уточнить:

$$\tilde{R} = R = \left\|x^0 - x_*\right\|_2.$$

Если решение задачи (1.1) не единственно, то можно считать, что $x_*$ в $\tilde{R}^2$ выбирается таким образом, чтобы минимизировать $\tilde{R}^2$. Оценка (1.28) внешне похожа на оценку (1.11). Однако стоит отметить, что константа $L$ в (1.28) определяется согласно (1.26), а не (1.4), и потому при $\|\ \| = \|\ \|_p$, $p \in [1,2)$ можно ожидать, что $L$ в (1.28) меньше, чем в (1.11) [3]. Однако типично, что «выигрыш» в $L$ с запасом нивелируется «проигрышем» в $\tilde{R}^2$, $\tilde{R}^2 \gg R^2$. ■

**Замечание 1.4 (наискорейший спуск).** Будем выбирать в методе (1.3) шаг $h$ не из условия (1.6), а следующим образом [82, § 1, гл. 3]:

$$h^k = \arg\min_{h\geq 0} f\left(x^k - h\nabla f\left(x^k\right)\right). \tag{1.29}$$

Такой метод (*наискорейшего спуска*) является естественным обобщением метода градиентного спуска. Очевидно, что соотношение (1.7) сохраняется. Таким образом, можно ожидать, что *метод наискорейшего спуска* сходится не медленнее градиентного спуска. И, действительно, на практике это часто можно наблюдать. Однако в худшем случае:[9]

$$f(x) = \frac{1}{2}\sum_{i=1}^{n}\lambda_i x_i^2, \quad 0 < \mu = \lambda_1 \leq \ldots \leq \lambda_n = L, \quad x^0 = \left(\frac{1}{\mu}, 0, \ldots, 0, \ldots, 0, \frac{1}{L}\right)^T$$

наискорейший спуск сходится и не быстрее градиентного спуска с постоянным (оптимально выбранным) шагом [82, теорема 3 § 4, гл. 1], вообще говоря, отличным от (1.6) [222, 552], т. е. приведенные выше оценки скорости сходимости градиентного спуска не могут быть принципиально улучшены даже при использовании шага (1.29). Рис. 3, взятый из работы [222], демонстрирует то, как сходится метод наискорейшего спуска в этом случае.

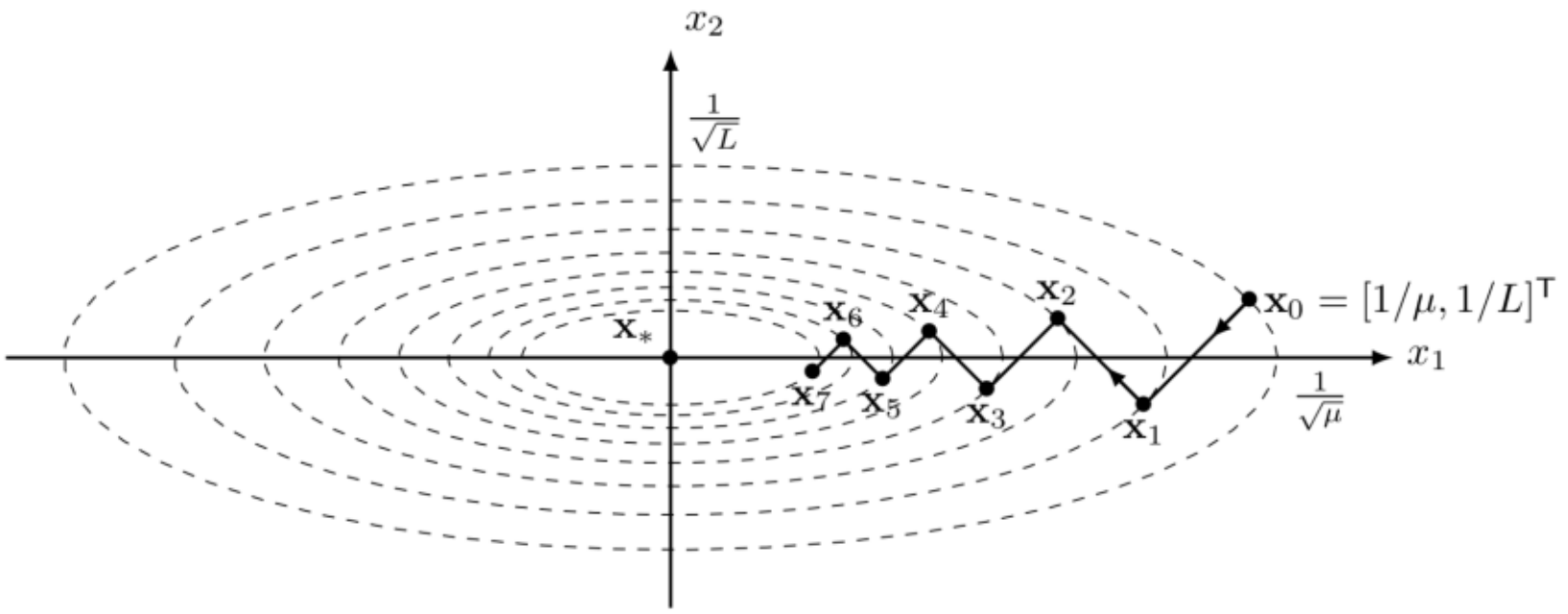

Рис. 3

Детали можно найти, например, в работе [222]. Тем не менее движение в этом направлении (использование вспомогательной одномерной / мало-мерной оптимизации на каждой итерации) может приносить серьезные дивиденды, см. замечания 1.5, 1.6.

---

[9] Не ограничивая общности для задач квадратчиной оптимизации можно считать $\lambda_1 > 0$, поскольку сильную выпуклость квадратичной задачи определяет минимальное отличное от нуля собственное значение матрицы квадратичной формы [294, 434].

Отметим также в этой связи следующий факт [82, § 1, гл. 3]. Если для минимизации положительно определенной квадратичной формы

$$f(x) = \frac{1}{2}\langle Ax, x\rangle - \langle b, x\rangle \to \min_{x\in\mathbb{R}^n} \tag{1.30}$$

использовать градиентный спуск (1.3) с $h^k = 1/\lambda_{k+1}$ где $\lambda_{k+1}$ – $(k+1)$-е собственное значение матрицы $A$ ($0 < \mu = \lambda_1 \le \ldots \le \lambda_n = L$), то независимо от точки старта метод будет конечен: $x^n = x_*$, где $Ax_* = b$. ■

Конструкции замечаний 1.3, 1.4, к сожалению, напрямую не переносятся на обобщения, собранные в последующих параграфах. Во всяком случае, нам о такой возможности не известно. Однако в случае замечания 1.3 существуют «обходные пути», позволяющие за небольшую «дополнительную плату» добиться желаемого обобщения. Подробнее об этом будет написано в следующем параграфе.

Отметим также, что градиентный спуск для класса *гладких* (в смысле (1.4)) выпуклых задач оптимизации не является *оптимальным методом* (см. упражнение 1.3). Однако в следующем параграфе будет отмечено (см. упражнение 2.2), что в условиях шума градиентный спуск может оказаться оптимальным. Для класса невыпуклых гладких задач безусловной оптимизации градиентный спуск является оптимальным методом (по критерию малости 2-нормы градиента) – см. текст после формулы (1.10).

◊ Здесь и далее оптимальность метода на классе задач понимается в смысле Бахвалова–Немировского [70] – число обращений (по ходу работы метода) к оракулу за градиентом (в общем случае старшими производными – см. замечание), т. е. число обращений к подпрограмме расчета градиента, для достижения заданной точности (например, по функции) в зависимости от параметров, характеризующих класс рассматриваемых задач и желаемую точность, может быть уменьшена равномерно на всем рассматриваемом классе только на числовой множитель (в замечании 1.5 и у первого аргумента минимума в оценке (1.45) оптимальность понимается еще более сильно – нельзя улучшить и числовой множитель), не зависящий от этих параметров и размерности пространства. Такой (оракульный) взгляд на сложность задач выпуклой оптимизации оказался очень удобным и популярным [72, 180, 438]. Связано это с тем, что, с одной стороны, существует хорошо разработанная теория оракульной сложности задач выпуклой оптимизации [70], с другой стороны, для большинства методов первого порядка и большинства задач наиболее вычислительно затратной частью итерации является именно расчет градиента. Таким образом, число обращений к оракулу отвечает за число итераций метода, что во многом определяет и общую сложность (время работы) метода.

Книга Немировского–Юдина [70] стала в свое время (с конца 70-х годов XX века) настоящим прорывом. Эта книга во многом определила последующее развитие численных методов выпуклой оптимизации. Запас оригинальных идей, заложенных в данной книге (совсем непростой для чтения), по-прежнему вдохновляет большое число исследователей по всему миру. ◊

В заключение заметим, что естественная попытка перенести метод (1.3) на условные задачи, т. е. задачи с ограничениями $x \in Q$ простой (в смысле проектирования) структуры

$$x^{k+1} = \arg\min_{x\in Q}\left\{\left\langle h\nabla f\left(x^k\right), x - x^k\right\rangle + \frac{1}{2}\left\|x - x^k\right\|_2^2\right\} =$$
$$= \arg\min_{x\in Q}\left\{\frac{1}{2}\left\|x-\left(x^k - h\nabla f\left(x^k\right)\right)\right\|_2^2 - \frac{h^2}{2}\left\|\nabla f\left(x^k\right)\right\|_2^2\right\} = , \tag{1.31}$$
$$= \arg\min_{x\in Q}\left\{\left\|x-\left(x^k - h\nabla f\left(x^k\right)\right)\right\|_2\right\} = \pi_Q\left(x^k - h\nabla f\left(x^k\right)\right)$$

приводит к аналогичному (1.18) выражению

$$\frac{1}{2}\left\|x^{k+1}-x_*\right\|_2^2 = \frac{1}{2}\left\|\pi_Q\left(x^k - h\nabla f\left(x^k\right)\right)-x_*\right\|_2^2 = \frac{1}{2}\left\|\pi_Q\left(x^k - h\nabla f\left(x^k\right)-x_*\right)\right\|_2^2 \le$$
$$\le \frac{1}{2}\left\|x^k - h\nabla f\left(x^k\right)-x_*\right\|_2^2 = \frac{1}{2}\left\|x^k - x_*\right\|_2^2 - h\left\langle\nabla f\left(x^k\right), x^k - x_*\right\rangle + \frac{h^2}{2}\left\|\nabla f\left(x^k\right)\right\|_2^2. \tag{1.32}$$

К сожалению, из этого неравенства уже нельзя получить неравенство (1.19), поскольку используемое при выводе (1.19) неравенство (1.7) уже может быть неверно. В следующем параграфе будет описано, как получить «правильный» аналог (1.18).

**Упражнение 1.1.** Докажите оценку (1.28). Почему в случае $\|\ \| = \|\ \|_2$ имеет место переход $\tilde{R}^2 \to R^2$?

**Упражнение 1.2.** Докажите утверждение из последнего абзаца замечания 1.4.

**Упражнение 1.3 (нижние оценки – гладкий случай / липшицев градиент).** Зафиксируем $N$. Рассмотрим класс методов

$$x^{k+1} \in x^0 + \mathrm{Lin}\left\{\nabla f\left(x^0\right),\dots,\nabla f\left(x^k\right)\right\}. \tag{1.33}$$

Не ограничивая общности, можно считать $x^0 = 0$.

**1)** Покажите, что в этом классе методов для *вырожденной* выпуклой функции

$$f(x) = F_{2N+1}(x) = \frac{L}{8}\left[x_1^2 + \sum_{i=1}^{2N}(x_i - x_{i+1})^2 + x_{2N+1}^2\right] - \frac{L}{4}x_1,$$

удовлетворяющей условию (1.4), при $2N+1 \le n$, где $n = \dim x$, имеют место следующие нижние оценки:

$$\min_{k=1,\dots,N} F_{2N+1}(x^k) - F_{2N+1}(x_*) \ge \frac{3L}{32}\frac{\left\|x^0 - x_*\right\|_2^2}{(N+1)^2},$$

$$\min_{k=1,\dots,N}\left\|x^k - x_*\right\|_2^2 \ge \frac{1}{8}\left\|x^0 - x_*\right\|_2^2,$$

где

$$x_* = \arg\min_{x\in\mathbb{R}^n} F_{2N+1}(x) = \left(1 - \frac{1}{2N+2}, 1 - \frac{2}{2N+2}, \dots, 1 - \frac{2N+1}{2N+2}, 0, \dots, 0\right)^T.$$

**2)** Покажите, что в этом классе методов для следующей $\mu$-сильно выпуклой в 2-норме функции, заданной в пространстве $\mathbb{R}^\infty$ и удовлетворяющей условию (1.4) (число обусловленности $\chi = L/\mu$)

$$f(x) = \frac{\mu\cdot(\chi - 1)}{8}\left[x_1^2 + \sum_{i=1}^{\infty}(x_i - x_{i+1})^2 - 2x_1\right] + \frac{\mu}{2}\|x\|_2^2,$$

при всех $N \ge 1$ имеют место следующие нижние оценки:

$$f(x^N) - f(x_*) \ge \frac{\mu}{2}\left(\frac{\sqrt{\chi}-1}{\sqrt{\chi}+1}\right)^{2N}\left\|x^0 - x_*\right\|_2^2,$$

$$\left\|x^N - x_*\right\|_2^2 \ge \left(\frac{\sqrt{\chi}-1}{\sqrt{\chi}+1}\right)^{2N}\left\|x^0 - x_*\right\|_2^2.$$

◊ Заметим, что [180, теорема 3.1.5]

$$\left(\frac{\sqrt{\chi}-1}{\sqrt{\chi}+1}\right)^{2N} \simeq \exp\left(-\frac{4N}{\sqrt{\chi}}\right). \text{ ◊}$$

**Указание.** См. [72, п. 2.1.2, 2.1.4], [82, § 2, гл. 3, п. 3 § 3, гл. 12], [180, п. 3.5], [239]. Общая идея построения плохих функций под класс методов вида (1.33) следующая:

$$\textit{искать } f(x) \textit{ в виде } f(x) = \sum_{i=1}^{n-1} f_i(x_i, x_{i+1}).$$

Такая структура $f(x)$ позволяет по наперед заданному $x_*$ легко построить $f(x)$ так, чтобы минимум достигался в $x_*$ [82, п. 3 § 3, гл. 12]. Однако основная цель выбора такой структуры – обеспечить при должном выборе точки старта $x^0 = 0$ и $f_i(x_i, x_{i+1})$ условие:

$$x_i^k = 0 \text{ при } i > k$$

для класса методов вида (1.33). В таком случае, как бы не выбирался метод, всегда можно, зная $x_*$, оценивать снизу невязку по функции. Это условие используется также при построении нижних оценок для методов первого порядка для задач выпуклой негладкой оптимизации (см. упражнение 2.1). С другой стороны, чтобы задача была сложной необходимо, чтобы спектр гессиана $\nabla^2 f(x_*)$ был сосредоточен около максимального и минимального собственных значений (см. также замечание 1.6). При этом в выпуклом (но не сильно выпуклом) случае еще желательно потребовать, чтобы обусловленность задачи (отношение максимального собственного значения гессина к минимальному) была не меньше $n^2$ (см. методы типа центров тяжести из указание к упражнению 1.4), и необходимо потребовать, чтобы обусловленность была не меньше величины обратной к относительной точности (по функции), с которой требуется решить задачу. Последнее условие нужно, чтобы задача была *вырожденной* [82, гл. 6]. Достаточно очевидно, что плохих функций можно подобрать много. В данном упражнении подобраны такие функции, для которых все отмеченные выше свойства достигаются на классе трехдиагональных теплицевых форм (а следовательно, достигается на классе квадратичных форм с равномерно ограниченными по $n$ коэффициентами)

$$\underbrace{\begin{pmatrix} 2 & -1 & 0 & 0 & \dots & 0 \\ -1 & 2 & -1 & 0 & \dots & 0 \\ 0 & -1 & 2 & -1 & \dots & 0 \\ \dots & \dots & \dots & \dots & \dots & \dots \\ 0 & 0 & \dots & 0 & -1 & 2 \end{pmatrix}}_{m}.$$

Число обусловленности такой положительно определенной матрицы $\sim m^2$. Это следует из общего факта, что спектр трехдиагональной теплицевой матрицы на главной диагонали которой стоит $a = 2$, а на двух других $b = -1$ и $c = -1$, состоит из чисел [91, 291, 477, 556]

$$a + 2\sqrt{bc}\cos\left(\frac{\pi k}{m+1}\right) = 2\left(1 - \cos\left(\frac{\pi k}{m+1}\right)\right),\ k = 1, \dots, m.$$

Отметим, что рассмотренный здесь способ порождения «плохих» функций для методов первого порядка используется практически в таком же виде и для методов высокого порядка [448] (использующих старшие производные оптимизируемой функции), см. приложение. Заметим, что выписанная выше матрица является также матрицей Лапласа (Кирхгофа), см. пример 4.1.

Приведенные в условиях упражнения нижние оценки с точностью до мультипликативных констант достигаются на классе *ускоренных* (*быстрых, моментных*) *градиентных методов* [72, § 2.2]. За последние десять лет интерес к этому классу методов резко возрос, см., например, [18, 23, 30, 39, 43, 71, 72, 95, 118, 122, 137, 143, 145, 154, 156, 158, 180, 225, 227, 232, 233, 240, 242, 243, 246, 281, 285, 307, 308, 309, 332, 361, 370, 371, 373, 386, 409, 412, 413, 414, 433, 452, 456, 462, 463, 474, 515, 517, 518, 544, 546, 552, 553, 558, 573, 574, 579] и цитированную там литературу. Отметим также, что полученные нижние оценки сохраняют свой вид и для более общего по сравнению с (1.33) класса методов [70, гл. 7].

Приведем в простейшем случае основную идею ускорения градиентного спуска, следуя работам [122, 145]. Ограничимся только выпуклым случаем, отвечающим п. 1) упражнения 1.3. Про перенесение на сильно выпуклый случай см. конец § 5.

Начнем с неформальной идеи [122]. Рассмотрим два режима: 1) на текущей итерации $\left\| f\left(x^k\right)\right\|_2 \geq M$ и 2) на текущей итерации $\left\| f\left(x^k\right)\right\|_2 < M$. Каждый шаг (итерация) обычного градиентного спуска (1.22) в режиме 1 уменьшает невязку по функции согласно (1.7) как минимум на $M^2/(2L)$. Следовательно, верхняя оценка на число таких шагов, необходимых для достижения точности (по функции) $\varepsilon$, будет пропорциональна $L\varepsilon/M^2$. С другой стороны, если все время пребывать в режиме 2 (в этом месте имеется неточность в рассуждениях!), то согласно упражнению 2.1 можно достичь точности (по функции) $\varepsilon$ за число шагов, пропорциональное $M^2/\varepsilon^2$. Если выбрать параметр $M$ так, чтобы сбалансировать обе полученные оценки: $L\varepsilon/M^2 \sim M^2/\varepsilon^2$, то общее число итераций в каждом из режимов $\sim \sqrt{L/\varepsilon}$, что лучше оценки, получаемой при использовании обычного градиентного спуска $\sim L/\varepsilon$.

Перейдем к более формальным выкладкам. Прежде всего заметим, что если в неравенстве (1.19) можно было бы «забыть» про необходимость выбирать шаг согласно правилу (1.6), то выбирая $h = R/\sqrt{2L\Delta f}$, где $\Delta f = f\left(x^0\right) - f\left(x_*\right)$, $R = \left\| x^0 - x_*\right\|_2$, вместо (1.20) получили бы

$$f\left(\bar{x}^N\right)-f\left(x_*\right)\le\frac{\sqrt{2LR^2\Delta f}}{N},$$

где

$$\bar{x}^N=\frac{1}{N}\sum_{k=0}^{N-1}x^k.$$

Следовательно, после $N_1=\sqrt{8LR^2/\Delta f}$ итераций гарантированно бы имели $f\left(\bar{x}^{N_1}\right)-f\left(x_*\right)\le\Delta f/2$. Если *рестартовать* такую процедуру после каждого момента гарантированного уполовинивания невязки по функции, то получится метод, который достигает точности $\varepsilon$ (по функции) после

$$N\simeq\underbrace{\sqrt{\frac{8LR^2}{\Delta f}}}_{N_1}+\underbrace{\sqrt{\frac{8LR^2}{\Delta f/2}}}_{N_2}+\ldots+\sqrt{\frac{8LR^2}{\varepsilon}}= \\ =\mathrm{O}\left(\sqrt{\frac{LR^2}{\varepsilon}}+\sqrt{\frac{LR^2}{2\varepsilon}}+\sqrt{\frac{LR^2}{4\varepsilon}}+\ldots\right)=\mathrm{O}\left(\sqrt{\frac{LR^2}{\varepsilon}}\right) \tag{1.34}$$

итераций. Данная оценка соответствует (с точностью до числового множителя) приведенной в условии 1) упражнения 1.3 нижней оценке. Однако все это было получено при невыполнимом для градиентного спуска предположении, что в неравенстве (1.19) можно выбирать $h=R\big/\sqrt{2L\Delta f}$. Основная проблема связана с тем, что шаг $h$ жестко задается в момент использования неравенства (1.7). Однако, если ввести три последовательности, связанные соотношениями ( $y^0=z^0$ ),

$$y^{k+1}=x^{k+1}-\frac{1}{L}\nabla f\left(x^{k+1}\right), \tag{1.35}$$

$$z^{k+1}=z^k-h\nabla f\left(x^{k+1}\right), \tag{1.36}$$

то из (1.36)

$$\left\langle h\nabla f\left(x^{k+1}\right),z^k-x_*\right\rangle\le\frac{1}{2}\left\|z^k-x_*\right\|_2^2-\frac{1}{2}\left\|z^{k+1}-x_*\right\|_2^2+\frac{\left\|h\nabla f\left(x^{k+1}\right)\right\|_2^2}{2},$$

а из последнего неравенства, подобно (1.19), можно получить неравенство

$$\left\langle\nabla f\left(x^{k+1}\right),z^k-x_*\right\rangle\le\frac{1}{2h}\left\|z^k-x_*\right\|_2^2-\frac{1}{2h}\left\|z^{k+1}-x_*\right\|_2^2+Lh\cdot\left(f\left(x^{k+1}\right)-f\left(y^{k+1}\right)\right),$$

справедливое уже для любого $h>0$. Но решив одну проблему, приобрели другую. Последнее неравенство в отличие от (1.19) не обладает *телескопическим свойством*: при суммировании все слагаемые в правой части неравенства, кроме крайних, взаимосокращаются, а левая часть неравенства мажорирует невязку по функции. Чтобы добиться выполнения телескопического свойства воспользуемся одной, не использованной пока, степенью свободы в определении $x^{k+1}\left(y^k, z^k\right)$. А именно, попробуем так подобрать эту зависимость, чтобы

$$\left\langle \nabla f\left(x^{k+1}\right), x^{k+1}-x_*\right\rangle - Lh\cdot\left(f\left(y^k\right)-f\left(y^{k+1}\right)\right)\le$$
$$\le\left\langle \nabla f\left(x^{k+1}\right), z^{k}-x_*\right\rangle - Lh\cdot\left(f\left(x^{k+1}\right)-f\left(y^{k+1}\right)\right).$$

Это удается сделать, используя выпуклость функции $f(x)$:

$$\left\langle \nabla f\left(x^{k+1}\right), y^{k}-x^{k+1}\right\rangle \le f\left(y^k\right)-f\left(x^{k+1}\right),$$

если

$$x^{k+1}-z^k = Lh\cdot\left(y^k-x^{k+1}\right).$$

Получается следующая простая зависимость (*выпуклая комбинация*)

$$x^{k+1}=\tau z^k+\left(1-\tau\right)y^k, \qquad (1.37)$$

где

$$\tau=\frac{1}{Lh+1}.$$

В результате для описанного здесь метода *линейного каплинга* (1.35)–(1.37) (название взято из работы [122]) можно написать следующую оценку, аналогичную (1.19),

$$\left\langle \nabla f\left(x^{k+1}\right), x^{k+1}-x_*\right\rangle \le \frac{1}{2h}\left\|z^k-x_*\right\|_2^2-\frac{1}{2h}\left\|z^{k+1}-x_*\right\|_2^2+Lh\cdot\left(f\left(y^k\right)-f\left(y^{k+1}\right)\right),$$

обладающую всеми необходимыми свойствами для последующего получения (с помощью рестартов) оптимальной оценки (1.34).

К сожалению, описанный выше подход, во-первых, базируется на знании, как правило, априорно неизвестной величины $R$ (используется при выборе размера шага $h=R\big/\sqrt{2L\Delta f}$ ), во-вторых, для корректности подхода под $R$ вместо $\left\|x^0-x_*\right\|_2$ необходимо понимать заметно бóльшую величину $\max\limits_{x:\, f(x)\le f\left(x^0\right)}\left\|x-x_*\right\|_2$, см. также замечание 1.3, в котором подобная оценка $R$ возникает из-за использования неевклидовой нормы. Обе отмеченные проблемы (явное использование $R$ при выборе шага и переоценка этого параметра) могут быть решены небольшой модификацией описанного подхода [122]. Кратко об этом написано в замечании 1.6 и в конце замечания 2 в приложении.

Впрочем, известно и много других вариантов ускоренных (быстрых, моментных) градиентных методов (см. начало указания), имеющих оценки глобальной скорости сходимости, аналогичные (с точностью до числового множителя) оценке (1.34), которые лишены отмеченных недостатков, см. ниже.

◊ Первым ускоренным градиентным методом с постоянными шагами для не квадратичных задач выпуклой оптимизации был (двухшаговый) *метод тяжелого шарика*:

$$x^{k+1} = x^k - \alpha \nabla f\left(x^k\right) + \beta \cdot \left(x^k - x^{k-1}\right), \qquad (1.38)$$

предложенный Б.Т. Поляком в 1963–1964 гг. [82, п. 1 § 2, гл. 3], [378]. Локальный анализ скорости сходимости метода тяжелого шарика (с помощью *первого метода Ляпунова* [82, § 1, гл. 2]) при специальном выборе параметров шага $\alpha, \beta > 0$ давал правильные порядки (соответствующие нижним оценкам) локальной скорости сходимости в сильно выпуклом случае. Анализ скорости сходимости метода в среднем, давал правильные порядки уже глобальной сходимости [516]. Однако с установлением глобальной сходимости (в худшем случае) были некоторые трудности. В частности, на специально подобранном примере с разрывным гессианом [332] метод может и не сходиться, а в чезаровском смысле траектории метода сходятся медленнее – аналогично (неускоренному) градиентному методу [283]. Несмотря на отмеченные сложности метод тяжелого шарика по-прежнему активно используется и продолжает развиваться [39, 414].

В 1982–1983 гг. Ю.Е. Нестеров в кандидатской диссертации (научным руководителем был Б.Т. Поляк) предложил первый ускоренный (быстрый) градиентный метод с фиксированными шагами (т. е. без вспомогательной маломерной оптимизации на каждой итерации, см. замечание 1.4), для которого удалось доказать глобальную сходимость с оптимальной скоростью (1.34) [73, 75]. Метод был «забыт» почти на 20 лет. Большое внимание этот метод привлек к себе лишь после выхода в 2004 году в издательстве Kluwer на английском языке первого издания книги [72] и работы [456]. Важную роль в привлечении внимания к методу сыграла также статья [156], имеющая большое число цитирований.

Отметим, что метод тяжелого шарика был мотивирован овражным методом Гельфанда–Цетлина [32], [82, п. 3 § 2, гл. 6]. Ускоренный метод

Нестерова также можно понимать, как специальный вариант овражного метода [13, пп. 3–5 § 1, гл. 5]. ◊

Для задач безусловной минимизации наиболее популярным сейчас является следующий (двухшаговый) вариант быстрого градиентного метода [72, 544], см. рис. 4: $x^0 = y^0$,

$$x^{k+1} = y^k - \frac{1}{L}\nabla f\left(y^k\right),$$

$$y^{k+1} = x^{k+1} + \frac{k}{k+3}\left(x^{k+1} - x^k\right),$$

который также можно понимать как *моментный метод* (следует сравнить с формулами (1.38), (1.44))

$$x^1 = x^0 - \frac{1}{L}\nabla f\left(x^0\right),$$

$$x^{k+1} = x^k - \frac{1}{L}\nabla f\left(x^k + \frac{k-1}{k+2}\left(x^k - x^{k-1}\right)\right) + \frac{k-1}{k+2}\left(x^k - x^{k-1}\right). \qquad (1.39)$$

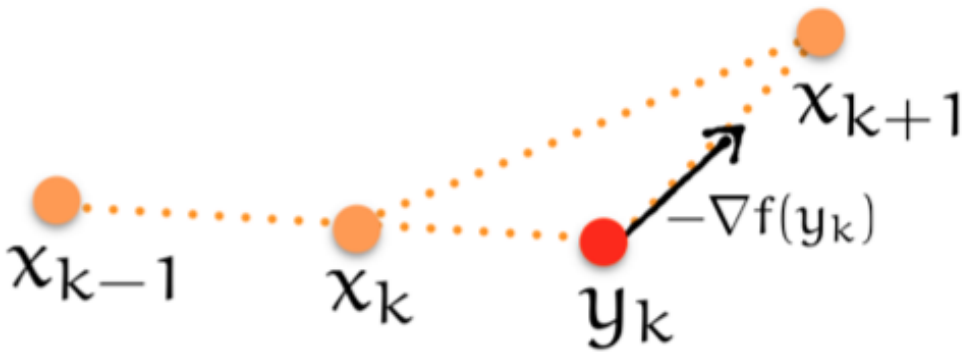

Рис. 4

Идею ускорения поясняет рис. 5, взятый из [508], на котором показаны линии уровня выпуклой функции и поведение траекторий градиентного спуска (слева) и быстрого градиентного метода (справа).

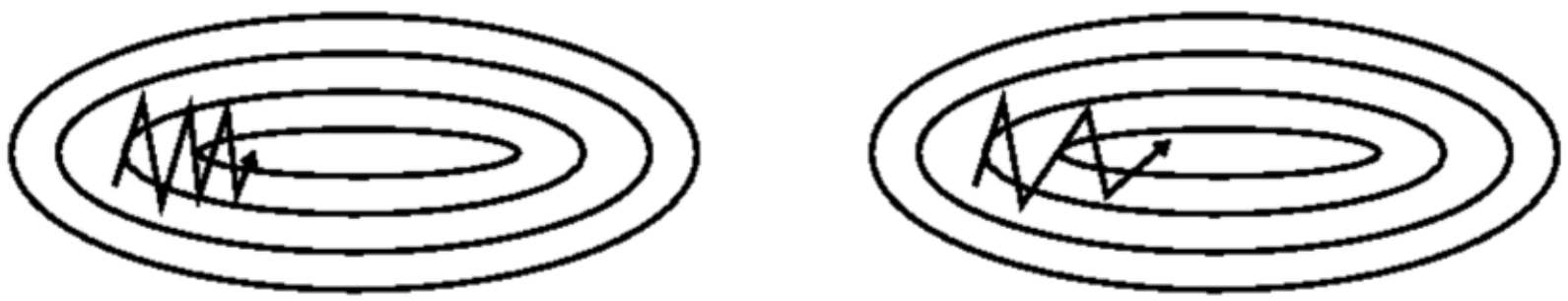

Рис. 5

Для справки приведем также вариант моментного метода для случая, когда оптимизируемая функция является $\mu$-сильно выпуклой в 2-норме [72, (2.38) п. 2.2.1]:

$$x^1 = x^0 - \frac{1}{L}\nabla f\left(x^0\right),$$

$$x^{k+1} = x^k - \frac{1}{L}\nabla f\left(x^k + \frac{\sqrt{L}-\sqrt{\mu}}{\sqrt{L}+\sqrt{\mu}}\left(x^k - x^{k-1}\right)\right) + \frac{\sqrt{L}-\sqrt{\mu}}{\sqrt{L}+\sqrt{\mu}}\left(x^k - x^{k-1}\right). \quad (1.40)$$

Оба метода (1.39), (1.40) сходятся согласно соответствующим нижним оценкам, приведенным в условиях 1), 2) упражнения 1.3, с точностью до числовых множителей при $L$ и $\chi = L/\mu$. Более того, можно даже объединить оба метода (1.39), (1.40) в один, сходящийся (с точностью до числовых множителей при $L$ и $\chi$) по оценке, наилучшей из двух оценок 1) и 2) [30], [72, теорема 2.2.2]. Однако во всех случаях (для метода (1.40) и объединенного метода) требуется априорно знать параметр $\mu$ (см. в этой связи также замечания 1.5, указания к упражнениям 2.3, 4.8 и конец § 5). На данный момент неизвестны такие варианты быстрых градиентных методов для сильно выпуклых задач, которые бы в общем случае обходились без этого знания. Отметим при этом, что незнание параметра $L$ может быть устранено вспомогательной маломерной минимизацией (см. замечания 1.5, 1.6) или адаптивным подбором (см. § 5).

Из сопоставления (1.22) и (1.39), (1.40) легко заметить, что сложность (трудоемкость) итераций у градиентного спуска и у быстрого градиентного метода практически одинаковы, в то время как скорости сходимости отличаются очень существенно. Именно это обстоятельство и обусловило огромную популярность быстрых градиентных методов. ■

**Замечание 1.5.** Среди последних достижений в развитии ускоренных градиентных методов отметим, *оптимизированный* вариант быстрого градиентного метода для задач безусловной оптимизации [239, 240, 370, 371, 373]

$$y^{k+1} = \left(1 - \frac{1}{t_{k+1}}\right)x^k + \frac{1}{t_{k+1}}x^0,$$

$$d^{k+1} = \left(1 - \frac{1}{t_{k+1}}\right)\nabla f\left(x^k\right) + \frac{2}{t_{k+1}}\sum_{j=0}^{k} t_j \nabla f\left(x^j\right), \quad (1.41)$$

$$x^{k+1} = y^{k+1} - \frac{1}{L}d^{k+1},$$

где

$$t_{k+1} = \frac{1+\sqrt{4t_k^2+1}}{2},\ k = 0,\dots,N-2,$$

$$t_N = \theta_N,\ \theta_{k+1} = \frac{1+\sqrt{8\theta_k^2+1}}{2},\ k = 0,\dots,N-1.$$

На классе гладких выпуклых задач метод сходится согласно оценке

$$f\left(x^N\right) - f\left(x_*\right) \le \frac{2LR^2}{\theta_N^2}\left(\le \frac{LR^2}{N^2}\right), \tag{1.42}$$

которая достигается, например, на функции Хьюбера

$$f(x) = \begin{cases} \frac{L}{2}\|x\|_2^2, & \|x\|_2 < \frac{R}{\theta_N^2}, \\ \frac{LR}{\theta_N^2}\|x\|_2 - \frac{LR^2}{\theta_N^4}, & \|x\|_2 \ge \frac{R}{\theta_N^2}. \end{cases} \tag{1.43}$$

Оценка (1.42) является точной оценкой скорости сходимости оптимальных методов вида (1.33) с шагами, не зависящими от оптимизируемой функции на классе гладких выпуклых задач. Другими словами, не существует такого метода вида (1.33) с фиксированными шагами, который бы на всем классе задач гладкой выпуклой оптимизации сходился по оценке лучшей, чем (1.42). Подчеркнем, что нельзя улучшить даже числовой множитель. Из метода (1.41) можно сделать метод, не требующий знания параметра $L$. Для этого в процедуре (1.41) следует заменить

$$x^{k+1} = y^{k+1} - \frac{1}{L}d^{k+1}$$

на

$$x^{k+1} = y^{k+1} - h_{k+1}d^{k+1},$$

где

$$h_{k+1} \in \operatorname{Arg}\min_{h\in\mathbb{R}} f\left(y^{k+1} - hd^{k+1}\right).$$

Такой метод также будет сходиться согласно оценке (1.42) [240]. Отметим, что *жадный метод первого порядка* (следует сравнить с (1.44))

$$x^{k+1} \in \operatorname{Arg}\min_{x\in x^0 + \operatorname{Lin}\left\{\nabla f\left(x^0\right),\dots,\nabla f\left(x^k\right)\right\}} f(x)$$

также будет сходиться согласно оценке (1.42), и для него эта оценка также не может быть улучшена [240].

Для класса гладких сильно выпуклых задач среди методов вида (1.33) не удалось найти простое (аналитическое) описание для оптималь-ного метода [240] (с неулучшаемым числовым множителем в оценке скорости сходимости). Наилучшие известные сейчас методы вида (1.33) с конечной памятью для данного класса задач описаны в [240, 518, 553].

Подобно (1.41) можно предложить (см. [240]) «универсальный» метод со вспомогательной трехмерной оптимизацией, который также не требует на вход никаких параметров. При этом метод оптимально работает (согласно оценке (1.42)) на классе не только гладких выпуклых задач, но и на классе негладких задач. В частности, для негладких задач метод сходится, согласно неулучшаемой (на классе методов вида (1.33) с фиксированными шагами) оценке [241] (следует сравнить с упражнением 2.1):

$$f\left(x^N\right) - f\left(x_*\right) \le \frac{L_0 R}{\sqrt{N+1}},$$

где $L_0$ определяется согласно (2.4). Отметим, что эта оценка также возникает (и является точной) для жадного метода первого порядка, в котором субградиент $\nabla f\left(x^k\right)$ выбирается из субдифференциала $\partial f\left(x^k\right)$ таким образом, чтобы $\left\langle \nabla f\left(x^k\right), \nabla f\left(x^l\right) \right\rangle = 0$, $0 \le l \le k-1$ [240]. К сожалению, пока непонятно, как переносить описанные выше в этом замечании методы на задачи выпуклой минимизации на множествах простой структуры. Впрочем, некоторые шаги в этом направлении уже сделаны [550, 551].

В работе [285] (см. также [463]) был предложен такой универсальный (в смысле § 5) вариант ускоренного метода, который для невыпуклых задач безусловной оптимизации сходится к локальному экстремуму с оптимальной скоростью (с точностью до числового множителя) по критерию малости 2-нормы градиента, а для выпуклых задач – к глобальному минимуму также с равномерно (по классам гладкости задач) оптимальной (с точностью до числового множителя) скоростью по критерию невязки по функции, см. приложение. Впрочем, в глубоком обучении (без особого теоретического обоснования – см., например, работу [497], про недостатки популярного метода Adam [221]) различные варианты быстрого градиентного метода начали использоваться заметно раньше [285], несмотря на невыпуклость возникающих там задач обучения нейронных сетей [39, 508, 546]. Отметим также эффективность быстрых градиентных методов в существенно невыпуклых задачах белкового фолдинга и докинга [102]. ■

◊ Поясним, следуя работам [240, 549, 552, 553], основную идею получения указанных в замечании 1.5 результатов. Сформулируем главную задачу: для заданного метода вида (1.33) с шагами, не зависящими от оптимизируемой функции, или жадного метода первого порядка требуется в заданном классе $F_{\mu,L}$ функций ($\mu$-сильно выпуклых функций с $L$-Липшицевым градиентом – все в 2-норме), при заданной точке старта $x^0$:

$$\alpha\cdot\left(f\left(x^0\right)-f\left(x_*\right)\right)+\beta\left\|\nabla f\left(x^0\right)\right\|_2^2+\gamma\left\|x^0-x_*\right\|_2^2\le C,$$

где $\alpha,\beta,\gamma\ge 0$, подобрать такую функций $f\in F_{\mu,L}$, что при заданном $N\le n-2$ максимально следующее выражение

$$\tilde{\alpha}\cdot\left(f\left(x^N\right)-f\left(x_*\right)\right)+\tilde{\beta}\left\|\nabla f\left(x^N\right)\right\|_2^2+\tilde{\gamma}\left\|x^N-x_*\right\|_2^2,$$

где $\tilde{\alpha},\tilde{\beta},\tilde{\gamma}\ge 0$. Не ограничивая общности можно дополнительно считать, что $x_*=0$, $f\left(x_*\right)=0$.

Сформулированная задача является задачей оптимизации в бесконечномерном пространстве функций из класса $F_{\mu,L}$. Пусть $x^k$ последовательность точек, генерируемых рассматриваемым методом, $f_k=f\left(x^k\right)$, $g^k=\nabla f\left(x^k\right)$. Будем считать, что $k\in I_N=\{0,1,2,\dots,N,*\}$, где $x^*=x_*=0$, $f_*=f\left(x_*\right)=0$, $g^*=\nabla f\left(x_*\right)=0$. Ключевым наблюдением является следующее утверждение: *для набора* $\left\{x^k,g^k,f_k\right\}_{k\in I_N}$ *существует такая функция* $f\in F_{\mu,L}$, *что верно* $f_k=f\left(x^k\right)$, $g^k=\nabla f\left(x^k\right)$, $k\in I_N$ *(в этом случае будем говорить, что набор* $\left\{x^k,g^k,f_k\right\}_{k\in I_N}$ *интерполируется в классе функций* $F_{\mu,L}$*), тогда и только тогда, когда для любых* $i,j\in I_N$ *выполняются условия* $F_{\mu,L}$*-интерполируемости:*

$$f_i-f_j-\left\langle g^j,x^i-x^j\right\rangle\ge$$

$$\ge\frac{1}{2(1-\mu/L)}\left(\frac{1}{L}\left\|g^i-g^j\right\|_2^2+\mu\left\|x^i-x^j\right\|_2^2-2\frac{\mu}{L}\left\langle g^i-g^j,x^i-x^j\right\rangle\right).$$

Доказательство данного утверждения базируется на двух вспомогательных фактах:

1) $\left\{x^k,g^k,f_k\right\}_{k\in I_N}$ *интерполируется в классе функций* $F_{\mu,L}$ *тогда и только тогда, когда* $\left\{x^k,g^k-\frac{\mu}{2}\left\|x^k\right\|_2^2,f_k-\mu x^k\right\}_{k\in I_N}$ *интерполируется в классе функций* $F_{0,L-\mu}$.

Это наблюдение следует из того, что

$$f(x)\in F_{\mu,L} \Leftrightarrow f(x)-\frac{\mu}{2}\|x\|_2^2 \in F_{0,L-\mu}.$$

2) $\{x^k, g^k, f_k\}_{k\in I_N}$ *интерполируется в классе функций* $F_{0,L}$ *тогда и только тогда, когда* $\{g^k, x^k, \langle g^k, x^k\rangle - f_k\}_{k\in I_N}$ *интерполируется в классе функций* $F_{1/L,\infty}$.

Это наблюдение следует из того, что [503, Theorem 23.5]

$$f(x^k)+f^*(g^k)=\langle g^k,x^k\rangle,\ g^k\in\partial f(x^k),\ x^k\in\partial f^*(g^k),$$

где $f^*(g)=\sup_{x\in\mathbb{R}^n}\{\langle g,x\rangle - f(x)\}$ – сопряженная функция к $f(x)$. Причем, по теореме Фенхеля–Моро $f^{**}(x)=f(x)$ [62, пп. 1.4, 2.2].

Сначала критерий $F_{\mu,L}$-интерполируемости устанавливается (несложно понять, что это можно сделать конструктивно в классе кусочно-линейных функций) для класса выпуклых негладких функций $F_{0,\infty}$. В этом случае все довольно наглядно (см. также (1.17)): для любых $i, j\in I_N$

$$f_i - f_j - \langle g^j, x^i - x^j\rangle \ge 0.$$

Затем к случаю $F_{0,\infty}$ последовательно, с помощью двух сделанных наблюдений, сводится и общий случай

$$F_{\mu,L} \xrightarrow{1)} F_{0,L-\mu} \xrightarrow{2)} F_{1/(L-\mu),\infty} \xrightarrow{1)} F_{0,\infty}.$$

Рассмотрим сначала случай, когда выбран метода вида (1.33) с шагами, не зависящими от оптимизируемой функции. В этом случае можно ввести матрицу Грамма $G=P^TP\succ 0$, где $P=[g^0,...,g^N,x^0]$.[10] В терминах матрицы $G$ и $f=(f_0,...,f_N)^T$ задачу можно переписать следующим эквивалентным образом (максимизация осуществляется по $G$ и $f$)

$$\langle f, b_{\tilde\alpha}\rangle + \langle G, \tilde M_{\tilde\beta,\tilde\gamma}\rangle \to \max_{\substack{f_i - f_j + \langle G, A_{ij}\rangle \le 0,\ i,j\in I_N \\ \langle G, M_{\beta,\gamma}\rangle - C\le 0 \\ G\succ 0}},$$

[10] Поскольку векторы $g^k, x^0\in\mathbb{R}^n$, то ранг матрицы $G\succ 0$ не больше чем $n$. Если по $G=P^TP$ нужно восстановить $P$, то для возможности восстановления в общем случае необходимо потребовать, чтобы $N+2\le n$. В этом случае восстановить $P$ всегда можно (причем не единственным образом), например, с помощью разложения Холецкого [91, п. 7.6].

где $\langle G, M \rangle = \mathrm{tr}(GM)$ – скалярное произведение матриц, а матрицы $A_{ij}$ определяются рассматриваемым методом вида (1.33) с шагами, не зависящими от оптимизируемой функции. Важное свойство выписанной задачи – это задача выпуклой оптимизации: в виду линейности функционала эту задачу можно переписать, сохраняя структуру, и как задачу на минимум. Более того, это *задача полуопределенного программирования* (semidefinite programming), см. приложение. С помощью современных пакетов символьных вычислений для ряда конкретных методов (например, метода градиентного спуска) эту задачу, вместе с двойственной к ней, даже удается аналитически решить[11] или хотя бы построить асимптотически (по $N$) точное решение. Более того, в ряде случаев, базируясь на описанной технике, удается найти аналитическое решение и для более сложных задач, в которых метод заранее не задан, и его стоит, в свою очередь, подбирать. В этом случае матрицы $A_{ij}$ сами становятся переменными, по которым необходимо минимизировать, и задача теряет свойство выпуклости. Тем не менее, даже в такой постановке иногда удается получить неожиданные результаты. В частности, именно таким образом был получен результат [372], описанный в замечании 5.3. Еще более удивительной, на первый взгляд, может показаться возможность аналитически решить задачу поиска оптимального метода вида (1.3) с шагами, не зависящими от оптимизируемой функции, для класса $F_{0,L}$, где $L \le \infty$, с коэффициентами в критерии $\tilde{\beta} = 0$, $\tilde{\gamma} = 0$ и в ограничении $\alpha = 0$, $\beta = 0$ [239, 373]. Оказывается, в этом случае результат получается такой же, как если бы в качестве метода использовался жадный метод первого порядка [240].

Рассмотрим этот метод подробнее. Ключевое наблюдение состоит в том, что такой метод генерирует последовательность $\{x^k, g^k\}_{k=0}^N$, которая определяется (однозначно, с точностью до возможного вырождения, $f(x)$ по части аргументов) следующим образом:

---

[11] Отметим, что решение задачи не только дает пример «наихудшей» функции из рассматриваемого класса для исследуемого метода, но и позволяет конструктивно показать, что полученная оценка скорости сходимости является не только нижней, но и верхней и, таким образом, точной. Для этого необходимо рассмотреть двойственную задачу. Решив двойственную задачу, можно использовать двойственные множители как веса, с которыми взвешиваются неравенства в критерии $F_{\mu,L}$-интерполируемости при конструктивном доказательстве того, что рассматриваемый метод сходится на любой функции из рассматриваемого класса $F_{\mu,L}$ не медленнее, чем предписано полученной оценкой – равенство достигается на найденной наихудшей функции.

$$\left\langle g^i, g^j \right\rangle = 0,\ 0 \le j < i = 1,\dots,N;\ \left| \cdot \theta_{ij} \right.$$

$$\left\langle g^i, x^j - x^0 \right\rangle = 0,\ 1 \le j \le i = 1,\dots,N.\ \left| \cdot \tau_{ij} \right.$$

Если ввести матрицу $G = P^T P \succ 0$, где[12]

$$P = [g^0,\dots,g^N, x^1 - x^0,\dots,x^N - x^0, x_* - x^0],$$

то выписанные условия, определяющие метод, вместе с условиями $F_{\mu,L}$-интерполируемости, переписанные в терминах ограничений на элементы матрицы $G$ и значения $f_k$, также как и раньше задают выпуклую (полуопределенную) систему ограничений. Таким образом, исходная задача опять сводится к задаче выпуклой полуопределенной оптимизации. В случае $\mu = 0$, $\tilde{\beta} = 0$, $\tilde{\gamma} = 0$ и $\alpha = 0$, $\beta = 0$ эту задачу и двойственную к ней удалось аналитически решить. При этом, зная решение двойственной задачи, можно явно предъявить оптимальный метод вида (1.33) с шагами, не зависящими от оптимизируемой функции, который в худшем случае сходится не хуже[13], чем в худшем случае сходится жадный метод первого порядка. Действительно, если известно решение двойственной задачи, то известны и $\theta_{ij}$, $\tau_{ij}$ – двойственные множители (множители Лагранжа) к ограничениям, определяющим жадный метод. Система этих ограничений может быть с помощью этих множителей эквивалентным образом переписана как

$$\left\langle g^i, \sum_{j=0}^{i-1} \theta_{ij} g^j + \sum_{j=1}^{i} \tau_{ij} \cdot \left( x^j - x^0 \right) \right\rangle = 0,\ i = 1,\dots,N.$$

Отсюда (в предположении $\tau_{ii} \ne 0$) получается нужное представление искомого оптимального метода вида (1.3):

$$x^k = x^0 - \sum_{j=0}^{i-1} \frac{\theta_{ij}}{\tau_{ii}} \nabla f\left( x^j \right) + \sum_{j=1}^{i-1} \frac{\tau_{ij}}{\tau_{ii}} \cdot \left( x^j - x^0 \right).$$

Обеспечить выполнение выписанного условия можно и другими способами, например, за счет вспомогательной маломерной оптимизации, см., замечание 1.5.

Критерий $F_{\mu,L}$-интерполируемости позволил получить новые интересные результаты в исследовании численных методов первого порядка с конечной памятью и шагами, независящими от номера итерации, в задачах (сильно) выпуклой (стохастической) оптимизации с помощью аппарата квадратичных функций Ляпунова [553] и квадратичных потенциальных функций [549]. Удалось получить необходимые и достаточные условия, когда такие функции (убывающие заданным образом) существуют, и конструктивные способы их построения. В основе лежит следующее наблюдение: условия убывания функции Ляпунова, потенциальной функции на траекториях рассматриваемого метода для всех функций из рассматриваемого класса при заданных матрицах квадратичных форм, записываются как система выпуклых ограничений вида неравенств: полуопределенных ($G \succ 0$) и линейных по $G$, $f$ и матрицам квадратичных форм. При этом необходимо уметь отвечать на вопрос: существует ли такие матрицы квадратичных форм, при которых система совместна? Данная задача сводится к седловой задаче (см. § 4) правильной (выпукло/вогнутой) структуры, которая, в свою очередь, сводится к задаче выпуклого полуопределенного программирования, причем малой размерности, определяемой глубиной памяти метода. Таким образом, например, удается обобщить теорему Ляпунова об устойчивости линейной системы [566] на сильно выпуклые задачи и методы с конечной памятью.

---

[12] В принципе, в определение матрицы $P$ можно не включать $x_* - x^0$.

[13] В рассматриваемом здесь случае можно сказать и точнее: «также сходится».

Заметим также (личное сообщение А. Тейлора), что описанная выше техника может быть перенесена и на задачи, в которых вместо настоящего градиента доступен только зашумленный градиент, например, в концепции аддитивной неточности:

$$\left\|\tilde{\nabla} f(x) - \nabla f(x)\right\|_2 \le \delta$$

или относительной аддитивной неточности:

$$\left\|\tilde{\nabla} f(x) - \nabla f(x)\right\|_2 \le \alpha \left\|\nabla f(x)\right\|_2, \ \alpha \in [0,1).$$

Чтобы воспользоваться описанной выше техникой стоит возвести эти неравенства в квадрат.

Насколько нам известно, это пока еще не сделано. Впрочем, ответ в первом случае ожидается в виде [141, 207, 216, 252, 294] (для неускоренных и ускоренных градиентных спусков): $O\left(\delta \left\|x^0 - x_*\right\|_2\right)$,[14] во втором случае ответ известен только для неускоренных методов – см. текст перед замечанием 1.1. Открытым является вопрос о том, как будет накапливаться неточность для ускоренных методов с относительной аддитивной неточностью в задачах безусловной оптимизации. ◊

**Замечание 1.6 (метод сопряженных градиентов).** Самой характерной задачей выпуклой оптимизации является задача минимизации положительно определенной квадратичной формы (1.30):

---

[14] Это следует из упражнени 3.7 и следующего наблюдения: для любого $L > 0$ выполняется $\left\langle \tilde{\nabla} f(x) - \nabla f(x), y - x \right\rangle \le \frac{1}{2L} \left\|\tilde{\nabla} f(x) - \nabla f(x)\right\|_2^2 + \frac{L}{2}\left\|y - x\right\|_2^2$.

$$f(x) = \frac{1}{2}\langle Ax, x\rangle - \langle b, x\rangle \to \min_{x\in\mathbb{R}^n}.$$

Изучив данный класс задач, можно пытаться понять, как сходятся различные методы хотя бы локально (в окрестности минимума) в задачах выпуклой оптимизации. Кроме того, такие задачи возникают как вспомогательные подзадачи при использовании методов второго порядка (например, метода Ньютона, см. приложение). Известно [36, 67, 70, 72, 82, 91, 137, 472, 556], что для выпуклой задачи квадратичной оптимизации (1.30) *метод сопряженных градиентов* (первая итерация делается согласно (1.29))

$$x^{k+1} \in \operatorname{Arg}\min_{x\in x^0 + \operatorname{Lin}\left\{\nabla f\left(x^0\right),\dots,\nabla f\left(x^k\right)\right\}} f(x) = x^k - \alpha_k \nabla f\left(x^k\right) + \beta_k \cdot \left(x^k - x^{k-1}\right), \quad (1.44)$$

где

$$(\alpha_k, \beta_k) \in \operatorname{Arg}\min_{\alpha,\beta} f\left(x^k - \alpha\nabla f\left(x^k\right) + \beta\cdot\left(x^k - x^{k-1}\right)\right),$$

сходится следующим образом (причем эта оценка в общем случае не может быть улучшена по первому и второму аргументу даже в смысле мультипликативного числового множителя [70, 82, 137, 439])

$$f\left(x^N\right) - f(x_*) \le \min\left\{\frac{LR^2}{2(2N+1)^2}, 2LR^2\left(\frac{\sqrt{\chi}-1}{\sqrt{\chi}+1}\right)^{2N}, \left(\frac{\lambda_{n-N+1}-\lambda_1}{\lambda_{n-N+1}+\lambda_1}\right)^2 R^2\right\}, \quad (1.45)$$

где $N \le n$, $\chi = L/\mu = \lambda_n/\lambda_1$, $R^2 = \left\|x^0 - x_*\right\|_2^2$ и использованы обозначения из замечания 1.4. Второй аргумент в минимуме (1.45) оценивается сверху следующим образом: $2LR^2\exp\left(-2N/\sqrt{\chi}\right)$. Полезно сопоставить эту оценку с соответствующей (сильно выпуклой) частью оценки скорости сходимости обычного градиентного спуска (1.24): $\left(LR^2/2\right)\exp(-N/\chi)$, и с нижней оценкой из п. 2) упражнения 1.3. При $N = n$ метод гарантированно находит точное решение (это свойство является отличительной особенностью методов сопряженных градиентов от их всевозможных обобщений, см., например, [463]), что следует из последней оценки в минимуме (1.45). Сформулированный результат является фундаментальным фактом (жемчужиной) выпуклой оптимизации и вычислительной линейной алгебры одновременно [212, 237], и базируется на наличии рекуррентных формул для *многочленов Чебышёва* [70, 82, 91, 137, 291, 472, 477, 511, 515, 556, 593]. Следуя [72, п. 1.3.2], приведем рассуждения, поясня-ющие правое равенство в (1.44). Введем обозначение (в выкладках воспользовались тем, что по условию $Ax_* = b$ ):

$$\Lambda_k = \mathrm{Lin}\left\{\nabla f\left(x^0\right),...,\nabla f\left(x^{k-1}\right)\right\} = \mathrm{Lin}\left\{Ax^0 - b,...,Ax^{k-1} - b\right\} =$$

$$= \mathrm{Lin}\left\{A\left(x^0 - x_*\right),...,A\left(x^{k-1} - x_*\right)\right\}.$$

Заметим, что

$$\Lambda_k = \mathrm{Lin}\left\{A\left(x^0 - x_*\right),...,A^k\left(x^0 - x_*\right)\right\},$$

т. е. $\Lambda_k$ – есть *линейное подпространство Крылова* [70, 477]. Из определения $x^{k+1}$ (левого равенства в (1.44)) имеем, что:

1) для всех $p \in \Lambda_k$ выполняется: $\left\langle \nabla f\left(x^k\right), p\right\rangle = 0$.

Введем *сопряженные направления* $\delta^k = x^{k+1} - x^k$. Сопряженные направления также порождают $\Lambda_k$ :

2) $\Lambda_k = \mathrm{Lin}\left\{\delta^0,...,\delta^{k-1}\right\}$.

Название «сопряженные направления» обусловлено свойством:

3) для $k \ne i$ выполняется: $\left\langle A\delta^k, \delta^i\right\rangle = 0$.

Действительно, пусть, для определенности, $k > i$, тогда

$$\left\langle A\delta^k, \delta^i\right\rangle = \left\langle A\left(x^{k+1} - x^k\right), \delta^i\right\rangle = \left\langle \nabla f\left(x^{k+1}\right) - \nabla f\left(x^k\right), \delta^i\right\rangle \overset{1)}{=} 0.$$

Из 2) и (1.44) (определения $x^{k+1}$) следует, что

$$x^{k+1} = x^k - h_k \nabla f\left(x^k\right) + \sum_{i=0}^{k-1} \lambda_i \delta^i,$$

т. е.

$$\delta^k = -h_k \nabla f\left(x^k\right) + \sum_{i=0}^{k-1} \lambda_i \delta^i.$$

Взяв скалярное произведение обеих частей этого равенства с вектором $A\delta^j$, по свойству 3) при $j < k-1$ получим

$$0 = \left\langle \delta^k, A\delta^j\right\rangle = -h_k \left\langle \nabla f\left(x^k\right), A\delta^j\right\rangle + \sum_{i=0}^{k-1} \lambda_i \left\langle \delta^i, A\delta^j\right\rangle =$$

$$= -h_k \underbrace{\left\langle \nabla f\left(x^k\right), \nabla f\left(x^{j+1}\right) - \nabla f\left(x^j\right)\right\rangle}_{0} + \lambda_j \left\langle \delta^j, A\delta^j\right\rangle = \lambda_j \left\langle \delta^j, A\delta^j\right\rangle,$$

т. е. $\lambda_j = 0$. Таким образом, доказано правое равенство в (1.44). Оценка (1.45) получается из левого равенства (1.44) с помощью следующего наблюдения (см., например, [82, п. 2 § 2, гл. 3]): $x^N \in x^0 + \Lambda_N$ равносильно тому, что существует такой многочлен $P_N(\lambda) = 1 + a_{1N}\lambda + \ldots + a_{NN}\lambda^N$, где коэффициенты $a_{1N}$, …, $a_{NN}$ могут принимать произвольные действительные значения, что $x^N - x_* = P_N(A)(x^0 - x_*)$. Поэтому для метода (1.44) должно выполняться:

$$f(x^N) - f(x_*) = \frac{1}{2}\langle Ax^N, x^N\rangle - \langle b, x^N\rangle - \underbrace{\left(\frac{1}{2}\langle Ax_*, x_*\rangle - \langle b, x_*\rangle\right)}_{0} =$$

$$= \frac{1}{2}\left\langle A\left(x^N - x_*\right), x^N - x_*\right\rangle =$$

$$= \min_{a_{1N},\ldots,a_{NN}} \left\{\frac{1}{2}\left\langle AP_N(A)^2\left(x^0 - x_*\right), x^0 - x_*\right\rangle\right\} \le$$

$$\le \frac{1}{2} \min_{P_N(\lambda):\, P_N(0)=1} \left\{\max_{\mu=\lambda_1\le\lambda\le\lambda_n=L}\left[\lambda P_N(\lambda)^2\right]\right\}.$$

Таким образом здесь появляются многочлены Чебышёва, наименее уклоняющиеся на рассматриваемом отрезке от нуля [62, п. 6.1, гл. 2]. Если относительно спектра матрицы $A$ делать различные вероятностные предположения, то приведенный выше подход можно уточнять: вместо многочленов Чебышёва наилучшими будут уже другие многочлены [485].

Оценка (1.45) хорошо согласуется с выписанными в условии упражнения 1.3 нижними оценками. При этом оценка (1.45), полученная для класса выпуклых задач квадратичной оптимизации, лучше точной нижней оценки для общего класса задач выпуклой оптимизации (1.42).

Отметим, в виду написанного выше, что на метод сопряженных градиентов можно посмотреть, как на наискорейший спуск (см. замечание 1.4) при выборе вместо антиградиентов сопряженных направлений, порождаемых процедурой ортогонализации Грамма–Шмидта на основе последовательности получаемых градиентов.

◊ Интересные результаты о регуляризующих свойствах метода сопряженных градиентов для вырожденных (некорректных) задач квадратичной оптимизации были получены в цикле работ А.С. Немировского [67, 69, 439]. *Вырожденной* будем называть задачу выпуклой оптимизации, для которой отношение максимального и минимального собственно-го значения функционала (*обусловленность задачи*) не меньше квадрата размерности пространства, в котором происходит оптимизация: $L/\mu \gg n^2$, и не меньше величины обратной к относительной точности, с которой требуется решить задачу, см. также указание к упражнению 1.3. Например, к такому классу задач относятся задача минимизации квадратичной формы, заданной *матрицей Гильберта* [82, п. 5 § 1, гл. 11]

$$\left\| \frac{1}{i+j-1} \right\|_{i,j=1}^{n},$$

см. также приложение. Обусловленность матрицы Гильберта для $n=6$ равна $1.5 \cdot 10^7$, а для $n=10$ – $1.6 \cdot 10^{13}$. Многие задачи, приходящие из реальных приложений, оказываются вырожденными, см., например, [355]. Строить сходящиеся по аргументу алгоритмы для таких задач в общем случае оказывается невозможным. Решение задачи оказывается неустойчивым к неточностям в данных. Для возможности корректного восстановления решения требуются дополнительные предположения (*истокопредставимости*). Здесь мы ограничимся простейшей задачей:

$$Ax = b,$$

в которой не доступны точные значения $A$ и $b$, а доступны только $\tilde{A}$ и $\tilde{b}$, где

$$\left\| \tilde{A} - A \right\|_2 \le \delta_A, \quad \left\| \tilde{b} - b \right\|_2 \le \delta_b,$$

где $\|C\|_2 = \sqrt{\lambda_{\max}\left(C^T C\right)}$. По поставленной задаче можно построить следующие задачи оптимизации

1) $$f_1(x) = \frac{1}{2}\left\| \tilde{A}x - \tilde{b} \right\|_2^2 \to \min_{x \in \mathbb{R}^n},$$

2) $$f_2(x) = \frac{1}{2}\left\langle \tilde{A}x, x \right\rangle - \left\langle \tilde{b}, x \right\rangle \to \min_{x \in \mathbb{R}^n}$$ (если $A^T = A \succ 0$, $\tilde{A}^T = \tilde{A} \succ 0$).

Введем индекс $\tau \in \{1,2\}$, который будет отвечать рассматриваемому случаю. В работе [67] было показано, что в случае, когда выполняется *условие истокопредставимости*

$$x_* = \left(A^T A\right)^{\sigma/2} y_*,\ \|y_*\|_2 \le R_\sigma,\ Ax_* = b,$$

метод сопряженных градиентов с критерием останова вида

$$\left\|\tilde{A}x^N - \tilde{b}\right\|_2 \le 2\left(\delta_A \left\|x^N\right\|_2 + \delta_b\right),$$

стартующий с $x^0 = 0$, сходится для соответствующей задачи $\tau \in \{1,2\}$ следующим образом

$$\omega_N^2 = \left\|\tilde{A}x^N - \tilde{b}\right\|_2^2 = \mathrm{O}\left(\frac{\tilde{L}^{2(1+\sigma)} R_\sigma^2}{N^{2\tau\cdot(1+\sigma)}} + \omega_*^2\right),\ \omega_* = \tilde{L}^\sigma R_\sigma \delta_A + \delta_b,$$

где $\tilde{L} = \max\left\{\|A\|_2, \|\tilde{A}\|_2\right\}$; причем до выполнения *критерия останова*

$$\left\|\tilde{A}x^N - \tilde{b}\right\|_2 \le 2\left(\delta_A \left\|x^N\right\|_2 + \delta_b\right)$$

при $\theta + 2\sigma > 0$, $\theta \in [0,2]$ справедлива следующая оценка

$$v_{\theta,N}^2 = \left\|\left(A^T A\right)^{\theta/4}\left(x^N - x_*\right)\right\|_2^2 = \mathrm{O}\left(R_\sigma^{(2-\theta)/(1+\sigma)} \omega_N^{(\theta+2\sigma)/(1+\sigma)}\right),$$

$$\left\|x^N\right\|_2 = \mathrm{O}\left(\|x_*\|_2\right).$$

Обратим внимание, что в $v_{\theta,N}^2$ стоит настоящая (незашумленная) матрица $A$. Приведенные выше результаты являются точными и не могут быть улучшены за счет использовании других методов. Причем не могут быть улучшены как в части скорости сходимости, так и в части достижимой точности $\mathrm{O}(\omega_*)$. Удивительно здесь, в частности, то, что метод сопря-женных градиентов, безусловно, можно относить к классу ускоренных (оптимальных) методов (см. указание к упражнению 1.3), для которых известно, что в общем случае, неточность в вычислении градиента линейно накапливается с ростом номера итерации [225], см. также упражнение 3.7 и начало приложения. Однако приведенный выше результат свидетельствует об отсутствии накопления неточностей, что соответствует неускоренным методам [225], см. также начало приложения. Отметим при этом, что в общем случае (в отличие от рассмотренного здесь) $\tilde{A}$, $\tilde{b}$ могут зависеть от номера итерации. Но даже в общем случае в численных экспериментах не наблюдалось накопление шума [510]. Объясняется это, по-видимому, тем, что шум в такое модели аддитивен (см. комментарии после замечания 1.5), что отличается от обсуждаемой выше более общей концепции шума, в которой он накапливается. Особенно полезным приведенные выше результаты оказываются при решении некорректных обратных задач в гильбертовых пространствах [281, 355].

Важно заметить, что в описанном выше подходе не производилась регуляризация задачи, см., например, замечание 4.1 и упражнение 4.9. Оптимальный результат был достигнут за счет регуляризирующего свойства самого метода (сопряженных градиентов). Отметим также, что критерий останова (именно с константой 2) может быть не достижим или достижим лишь за очень длительное время, но точно достижим с некоторой другой бóльшей (чем 2) константой. Есть некоторая проблема с тем, чтобы правильно определить эту константу, чтобы критерий оказался достижимым за оптимальное время. Однако, как видно из приведенных выше оценок, не выполнение критерия останова не приводит к тому, что метод плохо работает. Просто, выйдя за оптимальное время на режим $\omega_N = \mathrm{O}(\omega_*)$, метод сопряженных градиентов будет «топтаться на месте» и в итоге так и не сумеет преодолеть нужный порог. Однако остановив его, например, по критерию превышения заданного числа итераций, гарантированно можно получить в этом случае решение с качеством

$$v_{\theta,N}^2 = \mathrm{O}\left(R_\sigma^{(2-\theta)/(1+\sigma)} \omega_*^{(\theta+2\sigma)/(1+\sigma)}\right).$$

Интересно, отметить, что даже в условиях отсутствия шума $\tilde{A} = A$, $\tilde{b} = b$ приведенные выше результаты имеют ряд неожиданных, на первый взгляд, следствий.

Предположим, что $\tau = 1$, $\sigma = 1$, $\theta = 2$ получаем оценку:

$$v_{2,N}^2 = \left\| \left(A^T A\right)^{1/2} \left(x^N - x_*\right) \right\|_2^2 = \left\| Ax^N - b \right\|_2^2 = \mathrm{O}\left( \frac{\tilde{L}^4 R_1^2}{N^4} \right),$$

что соответствует результату, приведенному в упражнении 5.9 в условиях истокопредставимости из упражнения 4.9 с учетом отличия в обозначениях: при $\tau = 1$ имеет место равенство $\tilde{L} \simeq \sqrt{L_{f_1}}$, где $L_{f_1}$ – константа Липшица градиента $f_1(x)$.

Предположим, что $\tau = 2$, $\sigma = 0$, $\theta = 2$ получаем оценку:

$$v_{2,N}^2 = \left\| Ax^N - b \right\|_2^2 = \left\| \nabla f_2\left(x^N\right) \right\|_2^2 = \mathrm{O}\left( \frac{\tilde{L}^2 R_0^2}{N^4} \right) = \mathrm{O}\left( \frac{L_{f_2}^2 \left\| x_* \right\|_2^2}{N^4} \right).$$

Предположим, что $\tau = 2$, $\sigma = -1/2$, $\theta = 2$ получаем оценку:

$$v_{2,N}^2 = \left\| Ax^N - b \right\|_2^2 = \left\| \nabla f_2\left(x^N\right) \right\|_2^2 = \mathrm{O}\left( \frac{\tilde{L} R_{-1/2}^2}{N^2} \right) = \left( \frac{L_{f_2} \cdot \left( f_2(x_0) - f_2(x_*) \right)}{N^2} \right).$$

Оптимальность последней оценки в контексте предпоследней кажется сомнительной, однако в замечании 5.3 будет продемонстрировано, что из последней оценки можно получить предпоследнюю, как следствие.

Отметим также, что последние две оценки удается записать в форме, не учитывающей то, что рассматривается задача квадратичной оптимизации. Это наводит на мысль, что данные оценки должны быть справедливы для оптимальных методов и для более общего класс задач гладкой выпуклой оптимизации. Так оно и есть на самом деле. С точностью до логарифмических множителей этот результат был отмечен еще в работе

[439]. Недавно был предложен метод, который работает точно по приведенным оценкам [372], см. замечание 5.3. ◊

Метод (1.44) ничего не требует на вход (никаких параметров), а работает оптимально на классе гладких выпуклых задач и при этом также оптимально на его подклассе – гладких сильно выпуклых задач. Конечно, хотелось бы, чтобы и для общих задач выпуклой оптимизации метод (1.44) обладал аналогичными свойствами. Однако перенести без изменений (1.44) на весь класс задач выпуклой оптимизации не получилось. Тем не менее в конце 70-х годов XX века А.С. Немировскому удалось предложить две отдельные модификации метода (1.44): для класса гладких выпуклых задач и для класса гладких сильно выпуклых задач, которые доказуемо сходятся (с точностью до числовых множителей при $L$ и $\chi$) по оценкам, соответствующим первому аргументу минимума в (1.45) и второму (в сильно выпуклом случае), см. [70, 308, 433] и цитированную там литературу. Первый метод (для выпуклых задач) также не требует на вход никаких параметров, см., например, вариант метода (1.41) с одномерной минимизацией на каждой итерации. Второй метод (для сильно выпуклых задач) использует процедуру рестартов (см. упражнение 2.3 и конец § 5), и в общем случае требует знания параметра сильной выпуклости. При этом оба метода требуют на каждой итерации решения вспомогательной малоразмерной задачи выпуклой оптимизации. В отличие от задачи (1.44), которая для квадратичных функций решается по явным формулам, см., например, [82, п. 2 § 2, гл. 3], в общем случае на каждой итерации вспомогательную задачу можно решить только приближенно. В [70, § 3, гл. 7] было установлено, что достаточно решать вспомогательную задачу с относительной точностью (по функции) $\delta = O\left(\varepsilon / N(\varepsilon)\right)$, где $\varepsilon$ – желаемая относительная точность (по функции) решения исходной задачи, а $N(\varepsilon)$ – число итераций, которые делает (внешний) метод. Следовательно, вспомогательная задача может быть решена за

$$O\left(\ln\left(N(\varepsilon)/\varepsilon\right)\right) = O\left(\ln\left(\varepsilon^{-1}\right)\right)$$

обращений к оракулу (подпрограмме) за значением оптимизируемой функции (см. указание к упражнению 1.4). Таким образом, оба метода получились вполне практичными. Особенно практичными эти методы оказались для задач гладкой выпуклой оптимизации с функционалом вида $f\left(A^T x\right) + g(x)$, где вычисление $A^T x$ намного дороже по времени, чем вычисление $f(y)$ и $g(x)$, см. [308, 433] и трюк из приложения с быстрым пересчетом $A^T(x+\alpha \mathrm{v})$ для разных $\alpha \in \mathbb{R}$:

$$A^T(x+\alpha \mathrm{v}) = A^T x + \alpha A^T \mathrm{v}.$$

Такие функционалы, например, возникают при решении двойственных задач к задачам минимизации выпуклых сепарабельных функционалов при аффинных ограничениях: $Ay = b$, см. § 4.

Тем не менее до сих пор так и не был найден общий метод типа (1.44) или вариант метода (1.41) с вспомогательной одномерной оптимизацией, не требующий на вход никакой информации (о гладкости / сильно выпуклости оптимизируемого функционала), который сходится оптимально (хотя бы с точностью до числовых множителей) на классе гладких выпуклых задач и при этом также оптимально на его подклассе – гладких сильно выпуклых задач. Мало что известно про методы типа сопряженных градиентов (со вспомогательной маломерной оптимизацией) для задач оптимизации на множествах простой структуры [82, п. 2 § 3, гл. 7], [463]. Также мало что известно про возможные модельные обобщения (см. § 3). Лишь недавно была установлена прямо-двойственность (см. § 4) специального варианта метода линейного каплинга (см. указание к упражнению 1.3), в котором вместо (1.37) осуществляется одномерная оптимизация по параметру $\tau \in [0,1]$ [463] и(или) вместо (1.35) осуществляется одномерный поиск [309, 463]

$$x^{k+1} = \tau_{k+1} z^k + (1-\tau_{k+1}) y^k, \text{ где } \tau_{k+1} \in \operatorname{Arg}\min_{\tau\in[0,1]} f\left(\tau z^k + (1-\tau) y^k\right),$$

$$y^{k+1} = x^{k+1} - h_{k+1}\nabla f\left(x^{k+1}\right),\ h_{k+1} \in \operatorname{Arg}\min_{h\ge 0} f\left(x^{k+1} - h\nabla f\left(x^{k+1}\right)\right),$$

$$z^{k+1} = z^k - \alpha_{k+1}\nabla f\left(x^{k+1}\right),$$

где $\alpha_{k+1}$ определяется из решения квадратного уравнения

$$A_{k+1} f\left(x^{k+1}\right) - \frac{\alpha_{k+1}^2}{2}\left\|\nabla f\left(x^{k+1}\right)\right\|_2^2 = A_{k+1} f\left(y^{k+1}\right),\ \ A_{k+1} = A_k + \alpha_{k+1},\ A_0 = 0.$$

◊ Из (1.7) и того, что

$$y^{k+1} = x^{k+1} - h_{k+1}\nabla f\left(x^{k+1}\right),\ h_{k+1} \in \operatorname{Arg}\min_{h\ge 0} f\left(x^{k+1} - h\nabla f\left(x^{k+1}\right)\right)$$

будет следовать, что[15]

---

[15] В обычных (неадаптивных) ускоренных методах это неравенство следует понимать как равенство. Таким образом, можно понимать, что проблема адаптации

$$\frac{\alpha_{k+1}^2}{A_{k+1}} \ge \frac{1}{L}, \; A_{k+1} = A_k + \alpha_{k+1},$$

т.е.

$$\alpha_{k+1} \ge \frac{1}{2L} + \sqrt{\frac{1}{4L^2} + \alpha_k^2} \ge \frac{1}{2L} + \sqrt{\frac{1}{4L^2} + \frac{A_k}{L}}.$$

При этом

$$f\left(y^N\right) - f\left(x_*\right) \le \frac{R^2}{2A_N}.$$

Далее в тексте пособия мы старались унифицировать обозначения при описании различных вариантов быстрых градиентных методов. За основу взяты обозначения работы [122]. Обратим внимание, что в критерий качества в таких обозначениях следует подставлять точку $y^N$, а не $x^N$.

Важно заметить, что в таком же виде, посредством $A_N$, могут быть представлены оценки скорости сходимости различных ускоренных градиентных методов [18, 71, 72, 452], в том числе, собранных в данном пособии. В частности, в ряде случаев на это явно указано, см. замечания 1.5, 3.3.

Если в описанном методе выбирать

---

ускоренных методов к неизвестному параметру $L$ решается одномерным поиском (см. замечание 1.4) и формулой (1.7), из которой получается оценка неизвестного параметра

$$L = \frac{\left\|\nabla f\left(x^{k+1}\right)\right\|_2^2}{2\left(f\left(x^{k+1}\right) - f\left(y^{k+1}\right)\right)}.$$

Заметим, что для адаптации обычного (неускоренного) градиентного спуска достаточно только одномерного поиска, см. замечание 1.4. Заметим также, что для градиентного спуска существует вариант адаптивной оценки параметра $L$, не требующий вычислений значений функции [420].

$$x^{k+1} = \tau_{k+1} z^k + (1-\tau_{k+1}) y^k,$$

где

$$\tau_{k+1} = \frac{1}{\alpha_{k+1} L} = \frac{\alpha_{k+1}}{A_{k+1}}, \quad A_{k+1} = A_k + \alpha_{k+1}, \quad A_0 = 0,$$

$$y^{k+1} = x^{k+1} - \frac{1}{L} \nabla f\left(x^{k+1}\right),$$

$$z^{k+1} = z^k - \alpha_{k+1} \nabla f\left(x^{k+1}\right),$$

то

$$\alpha_{k+1} = \frac{1}{2L} + \sqrt{\frac{1}{4L^2} + \alpha_k^2} = \frac{1}{2L} + \sqrt{\frac{1}{4L^2} + \frac{A_k}{L}}.$$

Полученный в результате алгоритм является обычным методом линейного каплинга (МЛК) [122], описанным в упрощенной форме в указании к упражнению 1.3. Считая, что $\alpha_k \simeq C \cdot k/L$ при $k \gg 1$, получим уравнение на $C$:

$$C \cdot (k+1) = \frac{1 + \sqrt{1 + 4C^2 k^2}}{2},$$

следовательно,

$$1 + \frac{1}{k} = \frac{1}{2Ck} + \sqrt{1 + \frac{1}{4C^2 k^2}} \simeq \frac{1}{2Ck} + 1 + \frac{1}{8C^2 k^2} \simeq 1 + \frac{1}{2Ck},$$

т.е. $\alpha_k \simeq k/(2L)$, $A_k \simeq k^2/(4L)$. Следует сравнить приведенные рассуждения с немного более точным анализом, например, из [122], см. также конец замечания 2 в приложении. Близкий метод положен в основу проксимального ускоренного метода Монтейро–Свайтера, см. замечание 3.3. Следует также сравнить приведенные варианты МЛК с быстрым градиентным методом из упражнения 3.7 (методом подобных треугольников). ◊

Про другие методы вида (1.33) с вспомогательной маломерной оптимизацией по-прежнему ничего не известно. С учетом того, что по тео-ретическим оценкам использование таких методов не позволяет в общем случае улучшать даже числовые множители (см. также [240], замечание 1.4 и приложение А.С. Немировского в книге [75]), то, кажется, что стоит оставить эти методы в стороне и спокойно двигаться дальше по направлению к методам с постоянными / заданными шагами и конечной памятью, более простым, на первый взгляд, для всестороннего теоретического анализа. В общем-то, далее в пособии реализуется именно такой план.

Однако на практике типично [37], что различные варианты метода сопряженных градиентов (а их насчитывается уже, как минимум, несколько десятков [128, 272]) работают существенно быстрее ускоренных градиентных методов с постоянными шагами [361] и их адаптивных (универсальных) вариантов. На рис. 6, взятом из работы [309], приведен характерный график сходимости одной из наиболее быстрых на практике версий метода сопряженных градиентов [37, 75] (пунктирная линия) и графики сходимости быстрых градиентных методов (в данном случае использовались даже их универсальные варианты, см. § 5). Причина такого различия связана с наличием последнего аргумента минимума в оценке (1.45) для методов типа сопряженных градиентов и следующим наблюдением: методы типа сопряженных градиентов сходятся также как ускоренные методы (с фиксированными шагами) только на специальных (как правило, мало интересных для практики) примерах, см. замечание 1.4 и (1.43) и на локально (в окрестности минимума) квадратчиных функциях с спектром, сконцентрированным около наибольшего и наименьшего собственного значения. Например, для квадратичной формы с равномерно распределенным спектром и случайно (равновероятно) выбранной точкой старта можно ожидать (в среднем [485, 516]) сходимость специальных вариантов методов сопряженых градиентов со скоростью (по функции) $\sim N^{-6}$ вместо ожидаемой скорости $\sim N^{-2}$ (этот результат нам сообщил Ю.Е. Нестеров). В этой связи хотелось бы обратить внимание на важность проблем, затронутых в предыдущем абзаце.

Усилить сказанное в предыдущем абзаце можно следующей цитатой из уже немного устаревшей книги [33, стр. 208]: «Хотя схема сопряженных градиентов далека от идеала, на сегодня она является единственным разумным универсальным средством решения задач безусловной минимизации с очень большим числом переменных». Трудно сейчас всецело согласиться с такой категоричностью. Однако на практике попрежнему именно методы типа сопряженных градиентов и LBFGS (см. замечание 3 приложения) очень часто оказываются среди наилучших при решении выпуклых и невыпуклых задач оптимизации. Подробнее про это написано в более современной книге [37], также посвященной практическим вопросам решения задач оптимизации больших размеров.

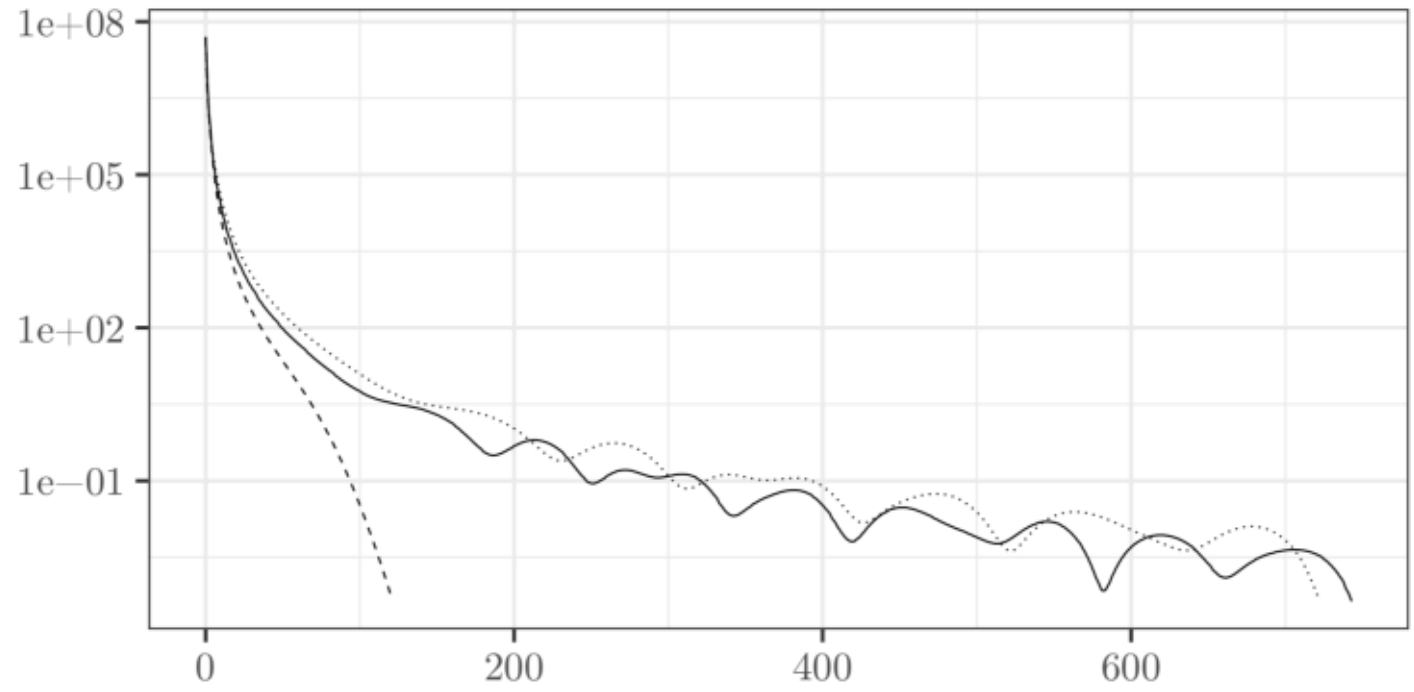

Рис. 6. По оси абсцисс – число итераций, по оси ординат – невязка по функции

Отметим также в связи с рис. 6, что ускоренные градиентные методы с постоянными шагами типично сходятся не монотонно [474] (по-видимому, это объясняется наличием всплесков в устойчивых несимметричных линейных дискретных системах при довольно общих условиях [378, 493]), что порождает различные трудности, см., например, замечание 5.3. Впрочем, добиться монотонности только по функции (не по норме градиента и расстоянию до решения) совсем не сложно [71, п. 1.2.2], см. также упражнение 3.7. ■

◊ Наиболее популярными на практике вариантами метода сопряженных градиентов являются следующие два метода [128], [472, гл. 5]:

$$h_k = \arg\min_{h\in\mathbb{R}} f\left(x^k + hp^k\right),$$

$$x^{k+1} = x^k + h_k p^k,$$

$$p^{k+1} = \nabla f\left(x^{k+1}\right) - \beta_k p^k,\ p^0 = \nabla f\left(x^0\right),$$

$$\beta_k = -\frac{\left\|\nabla f\left(x^{k+1}\right)\right\|_2^2}{\left\|\nabla f\left(x^k\right)\right\|_2^2}, \text{ (формула Флетчера–Ривса)}$$

$$\beta_k = -\frac{\left\langle \nabla f\left(x^{k+1}\right), \nabla f\left(x^{k+1}\right) - \nabla f\left(x^k\right)\right\rangle}{\left\|\nabla f\left(x^k\right)\right\|_2^2} \text{ (формула Полака–Рибьера–Поляка).}$$

Для задач квадратичной оптимизации оба метода идентичны (1.44). Для общих задач выпуклой оптимизации по этим методам не удалось пока получить оптимальные порядки скорости сходимости (установлен только сам факт глобальной сходимости для задач гладкой выпуклой оптимиза-ции). Тем не менее именно эти варианты метода сопряженных градиентов наиболее часто используются при решении практических задач [37] (в том числе не обязательно выпуклых). При этом с некоторой периодичностью (обычно период выбирают пропорционально размерности пространства, в котором происходит оптимизация) требуется перезапускать метод, обнуляя историю: вместо $p^{k+1} = \nabla f\left(x^{k+1}\right) - \beta_k p^k$ в момент рестарта полагают $p^{k+1} = \nabla f\left(x^{k+1}\right)$. По-видимому, необходимость в таких рестартах обусловлена желанием правильно сходиться в случае оптимизации сильно выпуклых функций, см. конец § 5 и [474]. ◊

**Упражнение 1.4.** Предложите способы решения задач выпуклой (но не обязательно сильно выпуклой) одномерной минимизации на отрезке длины $\Delta$ за время $\mathrm{O}\left(\log\left(\Delta/\varepsilon\right)\right)$, где $\varepsilon$ – точность решения задачи по аргументу. Возможно ли такое в двумерном случае? Рассмотрите способы, базирующиеся на вычислениях значения функции и производной. Исследуйте предложенные способы на точность получаемой информации (соответственно, на точность в получаемых значениях функции и в значениях ее производной в разных точках). Покажите, что малейшие шумы могут привести к отсутствию сходимости в требуемую окрестность по аргументу, однако при этом можно сохранить сходимость по функции в требуемую окрестность. Покажите, что оценку можно улучшить до $\mathrm{O}\left(\log\left\lceil\log\left(\Delta/\varepsilon\right)\right\rceil\right)$, если минимум не вырожденный (в точке минимума вторая производная положительная) и точка старта достаточно близка к точке минимума.

◊ Если ориентироваться на сходимость по функции, то для определенного класса методов (например, *метода центров тяжести*) размер области (расстояние о точки старта до решения) уже не будет входить в оценку. Будет входить лишь относительная (по функции) точность. Все это верно лишь для сильно выпуклых задач (см., например, (1.16)), либо для выпуклых задач на выпуклых компактах с оракулом, наделенным дополнительными нетривиальными возможностями, например, находить центр тяжести выпуклого компакта [180, п. 6.7]. На неограниченных множествах даже в одномерном случае в нижнюю оценку скорости сходимости по функции будет входить расстояние от точки старта до решения (в предположении, конечно, что в методе не разрешается бесплатно использовать одномерный поиск на (полу-)прямой). Причем это остается верным даже для выпуклых функций, имеющих ограниченную вариацию на полупрямой [70, упражнение 6 § 3, гл. 4]. К сожалению, для задач выпуклой оптимизации в $\mathbb{R}^n$, $n \gg 1$ в классе методов (1.33) даже вспомогательная маломерная оптимизация не позволяет в общем случае избавиться от вхождения расстояние от точки старта до решения в оценку скорости сходимости метода [70, 75]. Впрочем, в определенных случаях удается и для задач на неограниченных множествах получать оценки, в которые входит только относительная точность по функции [71, гл. 6], [452, гл. 7].

Приведенная оценка $\mathrm{O}\left(\log\left\lceil \log(\Delta/\varepsilon)\right\rceil\right)$, как оценка скорости глобальной сходимости, уже не может быть принципиально улучшена, какие бы не делались дополнительные предположения о гладкости функции и способностях локального оракула, вычисляющего значение функции и ее старшие производные в указанной точке [70, упражнение 2 § 1, гл. 8]. В частности, для класса *чебышёвских методов* [49, п. 2.9], использующих оракулы высокого порядка, можно увеличивать основание логарифмов в рассматриваемой оценке в зависимости от свойств оптимизируемой функции и порядка оракула, тем самым улучшить оценку. Однако при этом структура оценки (повторный логарифм) останется неизменной, см. приложение, а также [138, 464]. ◊

**Указание.** Рассмотрите, например, метод деления отрезка пополам или метод золотого сечения [13, § 3 – § 5, гл. 1]. Для анализа чувствительности методов и возможности ускорения при приближении к минимуму (в случае невырожденного минимума) стоит обратиться к [179, гл. 5].

Отметим, следуя А.С. Немировскому и Д.Б. Юдину, что для задач негладкой выпуклой, негладкой сильно выпуклой и гладкой выпуклой оптимизации на множествах простой структуры $Q$ (см. § 2) в $\mathbb{R}^n$, когда $N \ge n$, нижние оценки (для более общего по сравнению с (1.33) класса методов) числа обращений к оракулу за (суб-)градиентом (что такое субградиент будет пояснено также в § 2; про приложения можно прочитать в [428]) имеют одинаковый (с точностью до числового множителя $C$) вид

$$N \ge Cn\ln\left(\frac{\alpha}{\varepsilon}\right), \tag{1.46}$$

где $\varepsilon$ – относительная точность решения задачи по функции, а $\alpha \le 1$ – отношение длин сторон вписанного и описанного параллелепипеда для $Q$ [70]. Оценка (1.46) с $\alpha = 1$ достигается на *методе центров тяжести* Левина–Ньюмена, представляющего собой естественное обобщение на многомерные пространства метода деление отрезка пополам [70, § 3, гл. 2], [82, теорема 2 § 4, гл. 5], [180, п. 2.1].

◊ В основе метода центров тяжести лежит теорема Грюнбаума–Хаммера [62, п. 3.6], гарантирующая, что любая гиперплоскость, проходящая через центр тяжести ограниченного выпуклого множества, разделяет множество на два выпуклых подмножества, объем каждого из кото-рых не меньше $1/e$ от объема исходного множества. Считая, что у нас есть возможность эффективно находить центр тяжести ограниченных выпуклых множеств на основе этой теоремы можно предложить очень простой метод: проводить через центр тяжести текущего множества, на котором происходит оптимизация, гиперплоскость, перпендикулярную вектору (суб-)градиента, посчитанному в центре тяжести. Далее отсечь от множества ту часть, в которую смотрит антиградиент. На полученном в результате множестве можно повторить описанную процедуру и т.д. После $N$ итераций объем множества уменьшится не менее чем в $(1-1/e)^N$ раз. Выберем $N$ так, чтобы $(1-1/e)^N < \varepsilon^n$, т.е. $N = O\left(n \ln \varepsilon^{-1}\right)$. После такого числа итераций можно быть уверенным, что найдется хотя бы одно направление $r$, в проекции на которое размер множества уменьшился строго более чем в $1/\varepsilon$ раз. Отсюда следует, что значение функции в точках множества, оставшегося после $N$ итераций, будет не больше, чем значение в точке $x_* + \varepsilon r$, не попавшей в это оставшееся множество. Поскольку оптимизируемая функция выпуклая, то

$$f(x_* + \varepsilon r) - f(x_*) \le \varepsilon \cdot \left(f(x_* + r) - f(x_*)\right).$$

Отсюда следует оценка (1.46). ◊

Поскольку искать центр тяжести выпуклого множества в общем случае проблематично[16] [180, item 6.7], то у метода центров тяжести дорогая итерация [180, п. 6.7], поэтому на практике часто используют *метод эллипсоидов* [70, § 5, гл. 2], работающий по оценке[17]

$$N = O\left(n^2 \ln\left(\varepsilon^{-1}\right)\right),$$

[16] Одним из лучших известных сейчас способов решения этой задачи является рандомизированный алгоритм Hit and Run, в основу которого положена идея Markov Chain Monte Carlo. Алгоритм предполагает наличие граничного оракула, который по уравнению прямой выдает точки пересечения этой прямой с границей множества. Время работы алгоритма оценивается как $\tilde{O}\left(n^6\right)$.

[17] С точностью до логарифмического по $n$ множителя в случае оптимизации на евклидово-ассиметричных множествах, типа шара в 1-норме пространства $\mathbb{R}^n$.

но с относительно дешевой стоимостью (трудоемкости) итерации[18] $O(n^2)$.

◊ В основу метода эллипсоидов положено простое наблюдение: вокруг половины шара (полученной гиперплоскостью, проходящей через центр шара) в $n$-мерном пространстве за $O(n^2)$ арифметических операций можно описать (явно задав уравнением) эллипсоид (содержащий выбранную половину), объем которого не будет превышать

$$\left(1-\frac{1}{(n+1)^2}\right)^{n/2}$$

объема исходного шара [158, 180, 452]. Таким образом, если изначально решение находится в явно заданном эллипсоиде, то за $O(n^2)$ арифметических операций можно так линейно преобразовать пространство, что эллипсоид преобразуется в шар. Далее через центр шара проводится гиперплоскость согласно выбранному вектору из субдифференциала целевой функции в центре шара. Затем отбрасывается та половина шара, в которую смотрит выбранный субградиент. Вокруг оставшейся половины описывается эллипсоид (см. выше) и рассуждения повторяются (по индукции). При этом, если в какой-то момент центр шара окажется точкой, которая не принадлежит изначальному эллипсоиду (в момент старта), то вместо гиперплоскости, определяемой субградиентом целевой функции в этой точке за $O(n^2)$ арифметических операций строится отделяющая гиперплоскость (эта гиперплоскость отделяет центр шара от изначального эллипсоида), и отбрасывается та половина шара, которая не пересекается с изначальным эллипсоидом. ◊

В 2015 г. был предложен метод, который с точностью до логарифмического по $n$ множителя работает по оценке (1.46) в смысле требуемого числа итераций и по оценке[19] $\tilde{O}(n^2)$ в смысле стоимости итерации

---

[18] В оценку этой стоимости изначально не входит расчет (суб-)градиента. Однако обычно (суб-)градиент можно посчитать за $O(n^2)$ (см., например, начало § 2), поэтому можно считать, что выписанная оценка – есть оценка общей сложности итерации.

[19] К сожалению, с достаточно большим (полиномиально) логарифмическим множителем.

[400]. Отметим также *метод вписанных эллипсоидов*, предложенный в 1986 г. [97, стр. 253–259] и метод Вайды, предложенный в 1989 г. [180, п. 2.3], [563], проигрывающие методу 2015 года лишь в оценке стоимости итерации. Стоимость одной итерации этих методов оценивается, соответсвенно, как $\tilde{O}(n^{3.5})$ и $\tilde{O}(n^{2.37})$ [400]. Отметим, что оценка $\tilde{O}(n^{2.37})$ связана с шагом метода Ньютона[20], т.е. с решением системы линейных уравнений с симметричной матрицей, см. также указание к упражнению 5.9 и конец приложения. Продвижение работы [400] во многом базируется на построении нового гибридного барьера и возможности эффективно осуществлять шаг метода Ньютона для минимизации такого барьера, т.е. на возможности быстро перерешивать возникающую систему линейных уравнений с симметричной матрицей: $\tilde{O}(n^{2.37}) \to \tilde{O}(n^2)$, см. также указание к упражнению 5.9.

За дополнительную мультипликативную плату, пропорциональную с точностью до логарифмических множителей размерности пространства $n$, написанное в предыдущем абзаце переносится и на безградиентные методы [84, 399] – вместо (суб-)градиента оракул выдает значение функции. В гладком случае это довольно очевидное утверждение, поскольку градиент можно восстановить по частным производным, каждая из которых требует расчета значения функции в двух точках, причем одна из этих точек общая для всех компонент.

◊ Упоминая метод эллипсоидов, нельзя не отметить изящное обоснование Л.Г. Хачияном в 1978 г. с его помощью полиномиальной сложности задачи линейного программирования в битовой сложности [97, С. 453–461].[21]

Предположим необходимо точно ответить на вопрос: совместна ли система линейных неравенств $Ax \le b$, с размерами $n = \dim x$, $m = \dim b$?

---

[20] Для вспомогательной подзадачи поиска центра *объемного барьера*, построенного для оставшегося (к текущей итерации) от исходного множества политопа после всех отсечений (на предыдущих итерациях) гиперплоскостями полупространств. Стоит также отметить, что эта оценка представляет больше теоретический интерес, чем практический. Наилучшую практическую производительность демонстрируют методы, которые недалеко ушли в теоретическом плане от метода Штрассена с оценкой $O(n^{\log_2 7}) \gg \tilde{O}(n^{2.37})$ [333] и разложения Холецкого $O(n^3)$ [597].

[21] Отметим, что в общем случае поиск точного решения задачи оптимизации – NP-полная задача. Например, в упражнении 1.7 показывается, что задача поиска точного минимума выпуклого многочлена четвертой степени на единичной сфере (не выпуклом множестве) – NP-полная задача.

Предположим, что все элементы $A$ и $b$ являются целыми числами. Случай рациональных чисел очевидным образом сводится к целым. Если система совместна, требуется предъявить хотя бы одно точное решение $x_*$. Решение этой задачи (с точностью до логарифмического множителя в оценке сложности) эквивалентно точному решению задачи линейного программирования (ЛП)

$$\langle c, x \rangle \to \min_{Ax \le b}$$

с целочисленными $A$, $b$ и $c$. Напомним, для того чтобы найти точное решение совместной системы линейных уравнений $Ax = b$ достаточно воспользоваться алгоритмом Гаусса со сложностью $O(n^3)$ арифметических операций. Аналогичный вопрос относительно системы линейных неравенств $Ax \le b$ оставался открытым до конца 70-х годов прошлого века. Первый полиномиальный алгоритм был построен Л. Г. Хачияном в 1978 г. (тогда аспирантом ФУПМ МФТИ на базовой кафедре ВЦ РАН), как реакция на доклад А. С. Немировского с рассказом о скорости сходимости метода эллипсоидов. В конце 70-х годов прошлого века А. С. Немировский читал в ЦЭМИ РАН курс лекций по вышедшей впоследствии книге [70]. В 1982 г. Л. Г. Хачиян, А. С. Немировский и Д. Б. Юдин за цикл исследований, включающий, в том числе, и разработку метода эллипсоидов с доказательством полиномиальности задач ЛП в битовой сложности, были удостоены премии Фалкерсона. Опишем вкратце основной результат, полученный Л. Г. Хачияном. Положим

$$\Lambda = \sum_{i,j=1,1}^{m,n} \log_2 |a_{ij}| + \sum_{i=1}^{m} \log_2 |b_i| + \log_2 (mn) + 1.$$

Первое простое, но важное наблюдение: если система $Ax \le b$ совместна, тогда существует такой $x_*$, что $\|x_*\|_\infty \le 2^\Lambda$ и $Ax_* \le b$, иначе для всех $x$

$$\left\| (Ax - b)_+ \right\|_\infty \ge 2^{-(\Lambda - 1)},$$

где

$$\left[ (z)_+ \right]_i = \begin{cases} z_i, & z_i \ge 0 \\ 0, & z_i < 0 \end{cases}.$$

Второе простое, но также важное наблюдение: можно переформулировать исходную задачу совместности системы $Ax \le b$, как задачу негладкой выпуклой оптимизации

$$\left\| (Ax-b)_+ \right\|_\infty \to \min_{\|x_*\|_\infty \le 2^\Lambda},$$

которую достаточно решить с точностью $\varepsilon = 2^{-\Lambda}$. Собственно, метод эллипсоидов и был выбран в качестве алгоритма решения данной задачи в виду геометрической скорости сходимости, т.е. зависимости числа итераций от относительной точности вида $N(\varepsilon) \sim \ln\left(\tilde{\varepsilon}^{-1}\right)$. Было показано [97, С. 453–461], что, работая в $\mathrm{O}(n\Lambda)$-битной арифметике, с затратами (оперативной) памяти $\tilde{\mathrm{O}}\left(mn+n^2\right)$ описанным выше образом можно найти $x_*$ за $\tilde{\mathrm{O}}\left(n^3\left(n^2+m\right)\Lambda\right)$ арифметических операций.

Вопрос о полиномиальной сложности задачи ЛП в конце 70-х годов стоял достаточно остро. Связано это было, в первую очередь, с примером Кли–Минти (1972), в котором самый популярный на тот момент метод решения задач ЛП – симплекс-метод прежде чем найти решение последовательно проходит через все $2^m$ вершины гиперкуба, отвечающего $m$ специально сконструированным аффинным ограничениям [440]. Только в первой половине 80-х годов удалось объяснить хорошую работу симплекс-метода на практике. Оказалось, что математическое ожидание времени работы симплекс-метода при некоторых вполне естественных способах задания распределения вероятностей на параметрах задачи ЛП можно оценить следующим образом: $\tilde{\mathrm{O}}\left(m^3\right)$ [15, 170, 531].

В начале этого столетия Д. Спилманом и А. Тенгом был показан более тонкий результат [536]: если матрица аффинных ограничений в задача ЛП имеет вид $A + \|A\| G$, где матрица $G = \left\| g_{ij} \right\|_{i,j=1,1}^{m,n}$ состоит из независимых одинаково распределенных нормальных случайных величин $g_{ij} \in N\left(0,\sigma^2\right)$ и $T_\sigma(A)$ – время работы специальной версии симплекс-метода, необходимое, для нахождения точного решения, то [214]

$$E_G\left[T_\sigma(A)\right] = \mathrm{O}\left(n^2\sqrt{\ln m}\cdot\sigma^{-2} + n^3\ln^{3/2} m\right).$$

Можно ли в этой формуле (или ее аналогах) для математического ожидания времени работы, заменить $\sigma^{-2}$ на $\log^r\left(\sigma^{-1}\right)$, насколько нам известно, до сих пор открытый вопрос.

Рекордные результаты по сложности решения описанной задачи ЛП таковы: число итераций специального варианта метода внутренней точки $\tilde{O}\left(\sqrt{\operatorname{rank} A}\right)$, см. [398], указание к упражнению 5.9 и замечание 4 приложения. На каждой итерации необходимо решать $\tilde{O}(1)$ систем линейных уравнений, что может быть сделано за время $\tilde{O}(nnz(A))$, см. [400], указание к упражнению 5.9 и замечание 4 приложени. Однако на данный момент эти результаты представляют в основном только теоретический инетерес (методы не практичные).

Отметим, что вопрос о полиномиальной сложности задач ЛП в идеальной арифметике (в которой допустимы любые числа, и все арифметические операции, в частности, умножение $\pi \cdot e$, выполняется за $O(1)$) по-прежнему остается открытым, и был выделен С. Смейлом, в качестве одной из главных математических проблем этого столетия [168].

Отметим, что на практике для решения задач небольшого размера ($n \sim 10^2$) часто используют *bundle method*[22] (см. [158, item 5.4], [402]) и вариации метода эллипсоидов (*методы с процедурой растяжения пространства*), восходящие к работам Н.З. Шора [82, § 4, гл. 5], [86], [100]. Стоит также отметить большую роль, которые сыграли работы Н.З. Шора в появлении *субградиентного метода* (см. § 2). ◊

В случае, когда $N \le n$ (обычно это соответствует задачам оптимизации в пространстве большой размерности) и (при определенных условиях) в гладком сильно выпуклом случае оценка (1.46) перестает быть оптимальной (точной нижней границей сложности) [70]. Оптимальные оценки в этом случае будут достигаться на методах типа градиентного спуска, см., например, [72, 180], указание к упражнению 1.3 и замечания 1.5, 1.6. ∎

**Упражнение 1.5 (Нестеров–Пасечнюк–Стонякин, 2018 [484]).** Рассмотрим следующий метод решения задачи минимизации выпуклой липшицевой функции (с константой Липшица $L_0$) на квадрате в $\mathbb{R}^2$ со стороной $R$. Через центр квадрата проводится горизонтальная прямая. На отрезке, высекаемом из квадрата этой прямой, с точностью $\sim \varepsilon / \log\left(L_0 R / \varepsilon\right)$ (по функции) решается задача одномерной оптимизации. В найденной точке вычисляется вектор (суб-)градиента функции и определяется, в какой из двух прямоугольников он «смотрит», этот прямоугольник «отбрасывается». Через центр оставшегося прямоугольника проводится вертикальная прямая, на отрезке, высекаемом этой прямой в прямоугольнике, также с точностью $\sim \varepsilon/\log(L_0 R/\varepsilon)$ (по функции) решается задача одномерной оптимизации. В найденной точке вычисляется вектор (суб-)градиента функции и определяется, в какой из двух квадратов он «смотрит», этот квадрат «отбрасывается». В результате такой процедуры линейный размер исходного квадрата уменьшается вдвое.

[22] Названный в работе [72, п. 3.3.3] методом уровней. Однако в данном пособии под методом уровней, понимается немного другой метод, см. пример 3.2.

1) Покажите, что если оптимизируемая функция гладкая (дифференцируемая), то после $\sim \log(L_0 R/\varepsilon)$ повторений такой процедуры можно найти с точностью $\varepsilon$ (по функции) решение исходной задачи?

2) Покажите, что для негладкой выпуклой функции такой метод может не сходиться к решению задачи даже по функции.

**Указание.** 2) Рассмотрите функцию $|x-y|+0.9x$ на $[0,1]^2$.

**Упражнение 1.6 (метод условного градиента для задач квадратичной оптимизации на симплексе [18, п. 4.2.2, 4.3.3]).** Рассмотрим задачу квадратичной выпуклой оптимизации:

$$f(x)=\frac{1}{2}\langle Ax,x\rangle \to \min_{x\in S_n(1)},$$

где все элементы матрицы $A\succ 0$ по модулю не больше $M$, число ненулевых элементов в каждом столбце (строке) матрицы $A$ не больше $s \ll n$. Для решения этой задачи будем использовать метод условного градиента (см. замечание 1.2). Выберем одну из вершин симплекса и возьмем точку старта $x^0$ в этой вершине. Далее действуем по индукции, шаг которой имеет следующий вид. Решаем задачу

$$\left\langle \nabla f\left(x^k\right), y\right\rangle = \left\langle Ax^k, y\right\rangle \to \min_{y\in S_n(1)}.$$

Обозначим решение этой задачи через

$$y^k=(0,\ldots,0,1,0,\ldots,0),$$

где 1 стоит на позиции

$$i_k \in \operatorname{Arg}\min_{i=1,\ldots,n} \partial f\left(x^k\right)/\partial x_i .$$

Несложно показать, что решение такого вида всегда есть. Далее положим

$$x^{k+1}=(1-\gamma_k)x^k+\gamma_k y^k, \quad \gamma_k=\frac{2}{k+2}.$$

Заметим, что в такой метод не входят никакие параметры!

Имеет место следующая оценка скорости сходимости описанного метода:

$$f\left(x^N\right)-f\left(x_*\right)\le\frac{2L^pR_p^2}{N},$$

где $R_p^2=\max\limits_{x,y\in S_n(1)}\|y-x\|_p^2$, $L^p=\max\limits_{\|h\|_p\le 1}\langle h,Ah\rangle$, $1\le p\le\infty$, причем $p$ тут можно выбирать произвольно. С учетом того, что оптимизация происходит на симплексе, выберем $p=1$. Несложно показать, что этот выбор оптимален. В результате получим, что $R_1^2=2$,

$$L^1=\max_{i,j=1,\ldots,n}\left|A_{ij}\right|\le M.$$

Покажите, что после предварительных приготовлений (*препроцессинга*), имеющих сложность $\mathrm{O}(n)$, каждую итерацию можно осуществлять с трудоемкостью $\mathrm{O}(s\log_2 n)$. Таким образом, общая трудоемкость метода будет

$$\mathrm{O}\left(n+\frac{M}{\varepsilon}s\log_2 n\right),$$

что может быть значительно лучше оценки

$$\mathrm{O}\left(sn\sqrt{\frac{M\ln n}{\varepsilon}}\right)=\underbrace{\mathrm{O}(sn)}_{\substack{\text{сложность}\\\text{итерации}}}\underbrace{\mathrm{O}\left(\sqrt{\frac{M\ln n}{\varepsilon}}\right)}_{\text{число итераций}},$$

которая получается при использовании быстрого градиентного метода (оптимального по числу обращений к оракулу за градиентом) с наилучшей для данной задачи прокс-структурой (из известных) – энтропийной (см. конец § 2 и упражнение 3.7).

◊ Методы условного градиента можно при должной модификации применять к задачам сильно выпуклой оптимизации (см. [188, 366, 381, 386] и цитированную там литературу) и к выпукло-вогнутым седловым задачам (см. [287, 386] и цитированную там литературу). Недавно появились тензорные методы условного градиента [234]. Используя нижние оценки для методов типа условного градиента (методов с линейным минимизационным оракулом) в не сильно выпуклом случае [386, item 7.1.3] и конструкцию каталист (см. замечание 3.3 и приложение), можно получить нижние оценки в сильно выпуклом случае [104]. Эти оценки показывают, что в отличие от выпуклого случая, в сильно выпуклом случае методы условного градиента сходятся существенно медленнее. Попытки исправить ситуацию требуют дополнительных предположений относительно возможностей линейного минимизационного оракула [386, item 7.1.1.4]. С другими возможностями и обобщениями методов условного градиента можно познакомиться по работе [386, chapter 7]. ◊

**Упражнение 1.7 (Ю. Е. Нестеров [74]).** Пусть задан набор неотрицательных целых чисел $\{a_i\}_{i=1}^n$ полиномиально по $n$ неограниченных. Покажите, что следующие задачи NP-полны [280]:

**1)** Минимизация (точный поиск минимума) многочлена четвертой степени

$$f(x)=\sum_{i=1}^n x_i^4-\frac{1}{n}\left(\sum_{i=1}^n x_i^2\right)^2+\left(\sum_{i=1}^n a_i x_i\right)^4+(1-x_1)^4 .$$

Эта задача эквивалентна задаче минимизации выпуклого многочлена четвертой степени

$$P_4(x)=\sum_{i=1}^n x_i^4+\left(\sum_{i=1}^n a_i x_i\right)^4+(1-x_1)^4$$

на сфере.

**2)** Задача оптимального управления

$$P_4(x(1))\to\min_{u(\,)}$$

$$\frac{dx}{dt}=\frac{\langle x,u\rangle}{\|x\|_2^2}x-u\,,\quad 0\le t\le 1$$

где $x(0)=x_0\in\mathbb{R}^n$ – задан, причем $\|x_0\|_2=1$. Требуется найти точное решение.

**3)** Поиск направления убывания невыпуклой негладкой функции

$$f(x)=\left(1-\frac{1}{\gamma}\right)\max_{i=1,\dots,n}|x_i|-\min_{i=1,\dots,n}|x_i|+|\langle a,x\rangle|,$$

где $\gamma = \sum_{i=1}^{n} a_i > 1$, в точке $x = 0$, т.е. нужно найти хотя бы один такой $x \in \mathbb{R}^n$, что $f(x) < f(0) = 0$.

**Указание. 1)** Заметим, что

$$\sum_{i=1}^{n} x_i^4 - \frac{1}{n}\left(\sum_{i=1}^{n} x_i^2\right)^2 = \left\langle B[x]^2, [x]^2 \right\rangle \geq 0,$$

где $\left[[x]^2\right]_i = x_i^2$, $B = I - \frac{1}{n}1_n 1_n^T \succ 0$, причем $B[x]^2 = 0$ тогда и только тогда, когда $[x]^2 = \text{const} \cdot 1_n$. Значит, поиск такого $x \in \mathbb{R}^n$, что

$$f(x) = \min_{x \in \mathbb{R}^n} f(x) = 0$$

равносилен поиску такого $x$, что

$$(1 - x_1)^4 = 0, \; [x]^2 = \text{const} \cdot 1_n, \; \left(\sum_{i=1}^{n} a_i x_i\right)^4 = 0,$$

что эквивалентно решению *задачи о рюкзаке* (*ранце*)

$$a_1 + \sum_{i=2}^{n} a_i x_i = 0, \; x_i = \pm 1.$$

То есть если бы можно было эффективно точно решить исходную задачу оптимизации, то тогда можно было бы эффективно решить и задачу о рюкзаке. С помощью *быстрого преобразования Фурье* можно решить задачу о рюкзаке за время $\mathrm{O}\left(\ln n \cdot \sum_{i=1}^{n} |a_i|\right)$. Однако по условию задачи мы не можем считать эту оценку полиномиальной, т.е. мы не предполагаем, что выполняются условия, позволяющие убрать NP-полную задачу о рюкзаке из класса NP-полных задач.[23]

**2)** Заметим, что

$$\frac{d\|x\|_2^2}{dt} = \left\langle 2x, \frac{\langle x, u \rangle}{\|x\|_2^2} x - u \right\rangle = 2\left\langle x, \left(\frac{xx^T}{\|x\|_2^2} - I\right) u \right\rangle \equiv 0,$$

[23] Другой пример, подобный рассмотренному, см. в докладе S. Sra [607].

т.е. $x(1)$ принадлежит единичной сфере, как и $x(0)$. Более того, произвола в выборе $\{u(t)\}_{t\in[0,1]}$ достаточно, чтобы получить в качестве $x(1)$ произвольную точку на единичной сфере. Таким образом, исходная постановка задачи сводится к задаче, рассмотренной в п. 1.

**3)** Сначала заметим, что если набор $\sigma_i = \pm 1$ удовлетворяет условию

$$\langle a, \sigma \rangle = 0, \text{ то } f(\sigma) = -1/\gamma < 0.$$

Пусть $f(x) < 0$. Ввиду линейной однородности $f(x)$ можно считать, не ограничивая общности, что $\max\limits_{i=1,\dots,n} |x_i| = 1$. Обозначим $\delta = \langle a, \sigma \rangle$. Тогда $|x_i| > 1 - 1/\gamma + \delta$, $i = 1, \dots, n$. Вводя $\sigma_i = \operatorname{sign} x_i$, получим $\sigma_i x_i > 1 - 1/\gamma + \delta$, следовательно $|\sigma_i - x_i| = 1 - \sigma_i x_i < 1/\gamma - \delta$. Поэтому

$$|\langle a, \sigma \rangle| \le |\langle a, x \rangle| + |\langle a, \sigma - x \rangle| \le \delta + \gamma \max_{i=1,\dots,n} |\sigma_i - x_i| < (1-\gamma)\delta + 1 < 1.$$

Поскольку вектор $a$ имеет целочисленные компоненты по предположению, то последнее неравенство означает просто, что $\langle a, \sigma \rangle = 0$. ■

**Упражнение 1.8. 1)** Рассмотрим задачу невыпуклой оптимизации

$$f(x) \to \min_{x \in \mathbb{R}^n}$$

в условиях $L$-Липшицева градиента (1.4), и метод градиентного спуска с переменным шагом $h_k \le 1/(2L)$ и неточным градиентом $\tilde{\nabla} f(x)$:

$$x^{k+1} = x^k - h_k \tilde{\nabla} f(x^k), \quad \text{если } \left\| \tilde{\nabla} f(x^k) \right\|_2^2 > 2\varepsilon^2/5;$$

$$\text{STOP}, \quad \text{иначе},$$

где

$$\left\| \tilde{\nabla} f(x) - \nabla f(x) \right\|_2^2 \le \varepsilon^2/10.$$

Покажите, что такой метод гарантированно остановится не более чем за

$$N = 4\left( f(x^0) - f(x_*) \right) / \left( \varepsilon^2 h_{\min} \right)$$

итераций, см. (1.10), где $\varsigma = \min\limits_{k=1,\dots,\text{STOP}} \left\| \tilde{\nabla} f(x^k) \right\|_2^2 h_k \Big/ \varepsilon^2$. Таким образом, можно гарантировать, что в момент остановки

$$\left\| \nabla f(x^k) \right\|_2^2 = \left\| \nabla f(x^k) - \tilde{\nabla} f(x^k) + \tilde{\nabla} f(x^k) \right\|_2^2 \le$$

$$\leq 2\left\|\tilde{\nabla} f(x)-\nabla f(x)\right\|_2^2+2\left\|\tilde{\nabla} f\left(x^k\right)\right\|_2^2 \leq \varepsilon^2 .$$

◊ В следующих двух пунктах для наглядности вместо неравенства $\left\|\tilde{\nabla} f(x)-\nabla f(x)\right\|_2^2 \leq \varepsilon^2/10$ используется его ослабленный вариант: $E_\xi\left[\left\|\tilde{\nabla} f(x)-\nabla f(x)\right\|_2^2\right] \leq \varepsilon^2/10$. Использование неравенств концентрации меры (см. приложение) позволяет проработать этот момент аккуратнее. На самом деле, приведенные ниже результаты (о сходимости в среднем) можно получить и без предположений, которые требуются для возможности применения данных неравенств (субгауссовские хвосты). ◊

**2) (Stochastic Gradient Descent (SGD))** Покажите, что если вместо градиента $\nabla f(x)$ доступен стохастический градиент $\nabla_x f(x,\xi)$:

$$E_\xi\left[\nabla_x f(x,\xi)\right]=\nabla f(x), \; E_\xi\left[\left\|\nabla_x f(x,\xi)-\nabla f(x)\right\|_2^2\right] \leq D,$$

то для *пробатченного градиента*

$$\tilde{\nabla} f(x)=\frac{1}{r}\sum_{l=1}^{r}\nabla_x f\left(x,\xi^l\right),$$

где случайные величины $\left\{\xi^l\right\}_{l=1}^r$ – независимы и одинаково распределены, также как $\xi$, имеет место оценка

$$E_{\left\{\xi^l\right\}_{l=1}^r}\left[\left\|\tilde{\nabla} f(x)-\nabla f(x)\right\|_2^2\right] \leq \frac{D}{r}.$$

Выбирая $h_k \equiv 1/(2L)$ и $r \approx 10D/\varepsilon^2$, с помощью установленной оценки покажите, что не более чем за

$$N=\mathrm{O}\left(\frac{L\cdot\left(f\left(x^0\right)-f\left(x_*\right)\right)}{\varepsilon^2}\right)$$

итераций метод сможет найти такую точку $x^k$, что $E\left[\left\|\nabla f\left(x^k\right)\right\|_2^2\right] \leq \varepsilon^2$.

**3) (метод редукции дисперсии и SPIDER [268])** Предположим дополнительно, что

$$f(x)=\frac{1}{m}\sum_{l=1}^{m} f_l(x),$$

где все функции $f_l(x)$ имеют $L$-Липшицев градиент в 2-норме. Введем стохастический градиент (здесь все $\xi^k, \xi^{k-1}, \xi^{k,j}, \ldots$ – независимы, одинаково равномерно распределены среди чисел $1, \ldots, m$)

$$\nabla_x f\left(x^k, \left\{\xi^k\right\}\right) = \begin{cases} \dfrac{1}{r}\sum_{j=1}^{r} \nabla f_{\xi^{k,j}}\left(x^k\right), \text{ если } r < m \text{ и } k \text{ делится на } q, \\ \nabla f\left(x^k\right), \text{ если } m \le r \text{ и } k \text{ делится на } q, \\ \nabla f_{\xi^k}\left(x^k\right) - \nabla f_{\xi^k}\left(x^{k-1}\right) + \nabla_x f\left(x^{k-1}, \left\{\xi^{k-1}\right\}\right), \text{ иначе,} \end{cases}$$

где $r \approx 20D/\varepsilon^2$; $q = r$, если $r < m$ и $q = m$, если $m \le r$. Покажите, что если

$$h_k = \begin{cases} \varepsilon^2 \Big/ \left(20L\sqrt{D}\left\|\nabla_x f\left(x^k, \left\{\xi^k\right\}\right)\right\|_2\right), \text{ если } r < m, \\ \varepsilon \Big/ \left(L\sqrt{10m}\left\|\nabla_x f\left(x^k, \left\{\xi^k\right\}\right)\right\|_2\right), \text{ если } m \le r, \end{cases}$$

то для любого $k$ до момента остановки

$$E_{\left\{\xi^k\right\}}\left[\left\|\nabla_x f\left(x^k, \left\{\xi^k\right\}\right) - \nabla f(x)\right\|_2^2\right] \le \varepsilon^2/10.$$

Получите отсюда, что не более чем за

$$N = O\left(\min\left\{\frac{L\cdot\left(f\left(x^0\right) - f\left(x_*\right)\right)\sqrt{D}}{\varepsilon^3} + \frac{D}{\varepsilon^2}, \frac{L\cdot\left(f\left(x^0\right) - f\left(x_*\right)\right)\sqrt{m}}{\varepsilon^2} + m\right\}\right) =$$

$$= L\cdot\left(f\left(x^0\right) - f\left(x_*\right)\right) O\left(\min\left\{\frac{\sqrt{D}}{\varepsilon^3}, \frac{\sqrt{m}}{\varepsilon^2}\right\}\right)$$

итераций (вычислений $\nabla f_\xi(x)$) метод сможет найти такую точку $x^k$, что

$$E\left[\left\|\nabla f\left(x^k\right)\right\|_2\right] \le \varepsilon.$$

**Указание. 1)** Из (1.5) и сноски 14, в которой выбираем

$$L := 1/(2h_k),$$

имеем

$$f\left(x^{k+1}\right) \le f\left(x^k\right) + \left\langle \nabla f\left(x^k\right), x^{k+1} - x^k \right\rangle + \frac{L}{2}\left\|x^{k+1} - x^k\right\|_2^2 \le$$

$$\le f\left(x^k\right) + \left\langle \tilde{\nabla} f\left(x^k\right), x^{k+1} - x^k \right\rangle +$$

$$+\left(\frac{L}{2}+\frac{1}{4h_k}\right)\left\|x^{k+1}-x^k\right\|_2^2+h_k\left\|\tilde{\nabla} f\left(x^k\right)-\nabla f\left(x^k\right)\right\|_2^2 \le$$

$$\le f\left(x^k\right)-\frac{h_k}{2}\left\|\tilde{\nabla} f\left(x^k\right)\right\|_2^2+\frac{h_k}{10}\varepsilon^2 .$$

Если

$$\left\|\tilde{\nabla} f\left(x^k\right)\right\|_2^2 > 2\varepsilon^2/5 ,\ \varepsilon^2\varsigma \le \left\|\tilde{\nabla} f\left(x^k\right)\right\|_2^2 h_k ,$$

то из последнего неравенства следует, что

$$f\left(x^{k+1}\right) < f\left(x^k\right)-\frac{h_k}{2}\left\|\tilde{\nabla} f\left(x^k\right)\right\|_2^2+\frac{h_k}{4}\left\|\tilde{\nabla} f\left(x^k\right)\right\|_2^2 \le f\left(x^k\right)-\frac{\varsigma}{4}\varepsilon^2 .$$

**2)** См. начало приложения и [133, 284]. Отметим, что оценка на общее число вычислений $\nabla_x f\left(x,\xi\right)$ не может быть улучшена (с точностью до числового множителя) никаким другим методом при сделанных предположениях [133].

**3)** См. начало приложения и [268, 590]. Описанный способ формирования стохастического градиента (и близкие к нему способы, см. приложение) в оптимизационном сообществе принято называть *методом редукции дисперсии* [274, section 7]. Смысл данного названия поясняет следующее наблюдение (SPIDER [268]): для

$$k = a\cdot q + p ,\ 0 \le p < q \text{ и } \left\|\nabla_x f\left(x^k,\left\{\xi^k\right\}\right)\right\|_2 h_k \le \tilde{h}$$

справедливо

$$E_{\left\{\xi^k\right\}}\left[\left\|\nabla_x f\left(x^k,\left\{\xi^k\right\}\right)-\nabla f\left(x\right)\right\|_2^2\right] \le$$

$$\le \begin{cases} E_{\left\{\xi^{a\cdot q}\right\}}\left[\left\|\nabla_x f\left(x^{a\cdot q},\left\{\xi^{a\cdot q}\right\}\right)-\nabla f\left(x\right)\right\|_2^2\right]+L^2\tilde{h}^2 p, & \text{если } r<m, \\ L^2\tilde{h}^2 p, & \text{если } m\le r \end{cases} \le$$

$$\le \begin{cases} \varepsilon^2/20+\varepsilon^4 p/\left(400D\right), & \text{если } r<m, \\ \varepsilon^4 p/\left(10m\right), & \text{если } m\le r \end{cases} \le \frac{\varepsilon^2}{10}$$

То есть делая $q$ итераций, мы вычисляем $\nabla f_\xi\left(x\right)$ всего $3q$ раз. При этом дисперсию стохастического градиента удается уменьшить в $\sim D/\varepsilon^2$ . Однако, за это есть и «плата» в виде «маленьких» шагов, что приводит к увеличению необходимого числа итераций. Тем не менее, в итоге получается выиграть на числе вычислений $\nabla f_\xi\left(x\right)$ по сравнению с подходом, изложенном в предыдущем пункте.

Отметим, что приведенная оценка на число вычислений $\nabla f_{\xi}(x)$ не может быть улучшена (с точностью до числового множителя) никаким другим методом при сделанных предположениях [133]. Более того, полученная оценка (в варианте $r < m$) не требует, чтобы

$$f(x) = \frac{1}{m}\sum_{l=1}^{m} f_l(x).$$

Достаточно, чтобы [133, 268]

$$E_{\xi}\left[\left\|\nabla_x f(y,\xi) - \nabla_x f(x,\xi)\right\|_2^2\right] \le L^2 \left\|y - x\right\|_2^2. \qquad (1.47)$$

При этом в варианте $m \le r$ полученная оценка не требует, чтобы $D < \infty$.

Отметим также, что в предположении ограниченности старших производных оптимизируемой функции (достаточно Липшицевости гессиана и возможности умножения оценки гессиана на вектор; привлечение (оценок) старших производных ничего не дает) можно получить оценку $N \sim \varepsilon^{-3}$ без предположения (1.47) [134]. ■

# § 2. Метод проекции градиента

Рассмотрим задачу выпуклой оптимизации

$$f(x) \to \min_{x \in Q}. \tag{2.1}$$

Это значит, что $f(x)$ – выпуклая функция, а $Q \subseteq \mathbb{R}^n$ – выпуклое множество, которые мы считаем достаточно *простым* в том смысле, что решение вспомогательной задачи проектирования на это множество (1.31) занимает существенно меньше времени, чем расчет градиента $\nabla f(x)$. В качестве наглядного примера можно рассмотреть задачу минимизации квадратичной функции на параллелепипеде, задав, например,

$$f(x) = \frac{1}{2}\|Ax\|_2^2, \; Q = \prod_{k=1}^{n}[a_k, b_k].$$

В случае плотной матрицы $A$ расчет $\nabla f(x) = A^T \cdot (Ax)$ будет стоить $O(n^2)$ *арифметических операций*[24], а проектирование на $Q$ согласно (1.31) делается по явным формулам за $O(n)$. Совсем необязательно, чтобы проектирование осуществлялось по явным формулам. Однако в подавляющем большинстве рассматриваемых в приложениях случаев простых множеств $Q$ проектирование может быть осуществлено за (см. упражнение 4.6):

$$O\left(n \ln^2\left(\frac{n}{\varepsilon}\right)\right), \tag{2.2}$$

где $\varepsilon$ – относительная точность проектирования (в смысле сходимости по аргументу или по функции – в данном случае неважно). В противном случае множество уже, как правило, не считают простым, и его стараются описывать с помощью функциональных ограничений. Тогда становится правильнее говорить уже о задаче *условной оптимизации* [82], см. также пример 3.2, замечание 4.3 и упражнение 5.5.

---

[24] Операций типа сложения, умножения, деления двух чисел типа *float* [42, п. 1.3] – все эти операции сопоставимы (с точностью до логарифмического множителя от длины операндов) по сложности [59, глава 29], [63]. Число арифметических операций определяет время работы программы (длительность вычислений, трудоемкость), поэтому далее вместо числа арифметических операций также будет использоваться словосочетание *время работы*.

Существенным недостатком подхода из § 1 является предположение о том, что неравенство (1.4) имеет место на всем пространстве $\mathbb{R}^n$. Легко понять, что это довольно обременительное условие. Например, простая выпуклая функция скалярного аргумента $f(x) = x^4$ не удовлетворяет этому условию. Далее в этом параграфе путем специальной компактификации, вообще говоря, неограниченного множества $Q$, мы избавимся от отмеченной проблемы.

Другим недостатком подхода § 1 является невозможность его использования для выпуклых, но негладких функций $f(x)$. Действительно, следуя [82, § 3, гл. 5], рассмотрим

$$f(x_1, x_2) = |x_1 - x_2| + 0.2|x_1 + x_2|.$$

Естественно пытаться заменять градиент в подходе § 1 субградиентом (произвольным элементом субдифференциала) в точках, в которых $f(x)$ не является гладкой.

◊ Напомним, что *субдифференциал* – это в общем случае выпуклый компакт [62, п. 1.5]. Например, для функции скалярного аргумента $f(x) = |x|$ субдифференциал будет иметь вид

$$\partial f(x) = -1,\ x < 0;\ \partial f(x) = 1,\ x > 0 \text{ и } \partial f(x) = [-1,1],\ x = 0.$$

Напомним также, что мера (Лебега) точек негладкости выпуклой функции равна нулю по теореме Радемахера [101, 427], однако часто решение негладких задач достигается в одной из таких точек, и получается, что градиентный спуск может проводить заметную долю времени в окрестности таких точек. ◊

Рассмотрим одну из таких точек $(1, 1)$. Используя, например, соотношение (1.17) несложно проверить, что вектор (субградиент) $\nabla f(1,1) = (1.2, -0.8)$ будет принадлежать субдифференциалу $\partial f(1, 1)$. Однако при любом выборе шага $h > 0$ в методе (1.3) функция $f(x)$ из точки $(1,1)$ может только возрастать по направлению $\nabla f(1,1)$. Таким образом, в негладком случае рассчитывать на основное неравенство (1.7) не приходится, что и не удивительно, поскольку в это неравенство входит константа Липшица градиента $L$, предполагающая гладкость $f(x)$.

Более того, рассматривая простейшую негладкую выпуклую функцию скалярного аргумента с острым минимумом $f(x) = |x|$, имеем

$\left|\partial f\left(x\right)\right| = 1$, если только, случайно, мы не оказались в точке $x = 0$. Поэтому для метода (1.3) с шагом $h$ для почти всех точек старта $x^0$ имеем $\left|x^{k+1} - x^k\right| = h$, что влечет для любого $k$:

$$\max\left\{ f\left(x^k\right),\, f\left(x^{k+1}\right)\right\} \geq h/2\,.$$

Значит, необходимо выбирать $h$ пропорционально желаемой точности решения задачи $\varepsilon$, либо считать, что $h_k \to 0$ при $k \to \infty$, чтобы оказаться в нужной окрестности решения. Это существенно отличается от способа выбора шага (1.6) в гладком случае.

Несмотря на отмеченные сложности далее, следуя Ю.Е. Нестерову [462], мы постараемся единообразно посмотреть на гладкий и негладкий случаи.

Определим множество

$$B_{R,Q}\left(x_*\right) = \left\{x \in Q:\ \left\|x - x_*\right\|_2 \leq R\right\},$$

где $x_*$ – решение задачи (2.1), $R = \left\|x^0 - x_*\right\|_2$. Если решение не единственно, то под $x_*$ будем понимать такое решение задачи (2.1), которое наиболее близко в 2-норме к точке старта $x^0$. Предположим, что для любых $x, y \in B_{R,Q}\left(x_*\right)$

$$f\left(y\right) \leq f\left(x\right) + \left\langle \nabla f\left(x\right), y - x\right\rangle + \frac{L}{2}\left\|y - x\right\|_2^2 + \delta, \qquad (2.3)$$

где $\delta > 0$. В частности, если градиент $f\left(x\right)$ удовлетворяет *условию Гёльдера*, точнее для любых $x, y \in B_{R,Q}\left(x_*\right)$ имеет место неравенство

$$\left\|\nabla f\left(y\right) - \nabla f\left(x\right)\right\|_2 \leq L_\nu \left\|y - x\right\|_2^\nu,\ \nu \in \left[0,1\right],\ L_0 < \infty, \qquad (2.4)$$

то (2.3) имеет место с

$$L = L_\nu \cdot \left[\frac{L_\nu}{2\delta}\frac{1-\nu}{1+\nu}\right]^{\frac{1-\nu}{1+\nu}}. \qquad (2.5)$$

Детали см. в работе [227]. См. также рис. 7, соответствующий $\nu = 0$.

◊ В случае $\nu = 0$ условие (2.4) (аналогично (2.3)) выполняется для любых элементов соответствующих субдифференциалов $\partial f\left(x\right)$ и $\partial f\left(y\right)$. Фактически, условие (2.4) отвечает тому, что у функции $f\left(x\right)$ равномерно ограниченны все элементы субдифференциалов во всех точках, т. е.

функция $f(x)$ имеет равномерно ограниченную константу Липшица. Заметим, что с небольшими оговорками это условие отвечает самому общему классу всех собственных выпуклых функций на открытом множестве, содержащем $B_{R,Q}(x_*)$, см. также [152, 304].

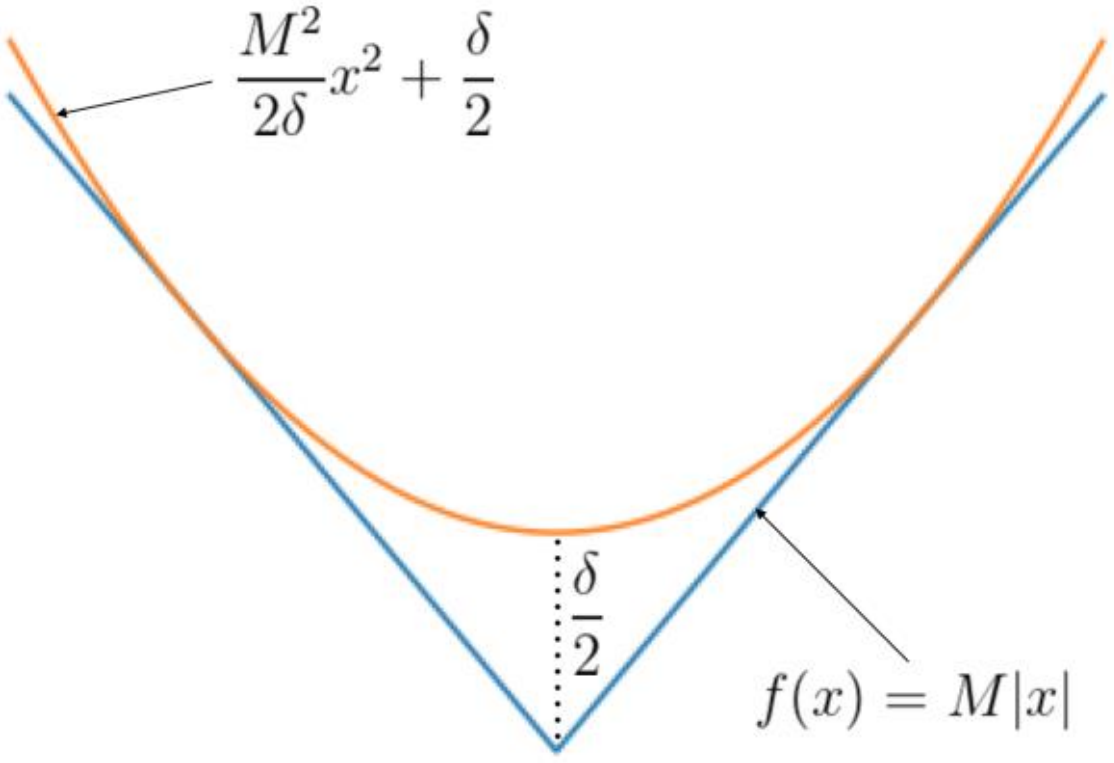

Рис. 7

Отметим также, что с точки зрения формальной логики формула (2.4) некорректна, потому что незамкнута относительно параметра $\nu$ [79]. Однако это было сделано вполне осмысленно. Дело в том, что в настоящем параграфе на формулу (2.4) будем смотреть только с точки зрения конкретного $\nu$. В § 5 уже будем «играть» на выборе параметра $\nu\in[0,1]$, считая, что (2.4) имеет место для любого $\nu\in[0,1]$, но при этом допуская, возможность того, что $L_\nu=\infty$ начиная с некоторого $\nu\in(0,1]$.

Во всех последующих формулах, если не оговорено противного (см. (3.4), (3.5)), использование без уточнений $\nabla f(x)$ подразумевает, что формулы справедливы при любом выборе $\nabla f(x)\in\partial f(x)$. ◊

Рассмотрим (см. также формулу (1.31)) простейший метод проекции градиента с шагом $h\le 1/L$

$$x^{k+1}=\pi_Q\left(x^k-h\nabla f\left(x^k\right)\right)=\arg\min_{x\in Q}\left\{\left\langle h\nabla f\left(x^k\right),x-x^k\right\rangle+\frac{1}{2}\left\|x-x^k\right\|_2^2\right\}=$$
$$=\arg\min_{x\in Q}\left\{f\left(x^k\right)+\left\langle\nabla f\left(x^k\right),x-x^k\right\rangle+\frac{1}{2h}\left\|x-x^k\right\|_2^2\right\}. \quad (2.6)$$

Следствием (2.6) является условие, которое получается из (1.17) при $x = x^{k+1}$, $y = x$: для всех $x \in Q$:

$$\left\langle \nabla_x \left( \left\langle h\nabla f\left(x^k\right), x - x^k \right\rangle + \frac{1}{2}\left\| x - x^k \right\|_2^2 \right)\Bigg|_{x = x^{k+1}}, x - x^{k+1} \right\rangle = \\ = \left\langle h\nabla f\left(x^k\right) + x^{k+1} - x^k, x - x^{k+1} \right\rangle \geq 0, \tag{2.7}$$

т. е. для всех $x \in Q$:

$$\left\langle h\nabla f\left(x^k\right), x^{k+1} - x \right\rangle \leq \left\langle x^{k+1} - x^k, x - x^{k+1} \right\rangle = \frac{1}{2}\left\| x - x^k \right\|_2^2 - \frac{1}{2}\left\| x - x^{k+1} \right\|_2^2 - \frac{1}{2}\left\| x^{k+1} - x^k \right\|_2^2. \tag{2.8}$$

Следуя [122], введем

$$\mathrm{Prog}^h\left(x^k\right) \overset{\mathrm{def}}{=} -\min_{x \in Q}\left\{ \left\langle \nabla f\left(x^k\right), x - x^k \right\rangle + \frac{1}{2h}\left\| x - x^k \right\|_2^2 \right\} \geq 0. \tag{2.9}$$

По формуле (2.6) для всех $x \in Q$ имеет место следующее «правильное» обобщение равенства (1.18) (следует сравнить с «неправильным» / «грубым» вариантом (1.32))

$$\left\langle h\nabla f\left(x^k\right), x^k - x \right\rangle = \left\langle h\nabla f\left(x^k\right), x^k - x^{k+1} \right\rangle + \left\langle h\nabla f\left(x^k\right), x^{k+1} - x \right\rangle \leq \\ \leq \left\langle h\nabla f\left(x^k\right), x^k - x^{k+1} \right\rangle + \frac{1}{2}\left\| x - x^k \right\|_2^2 - \frac{1}{2}\left\| x - x^{k+1} \right\|_2^2 - \frac{1}{2}\left\| x^{k+1} - x^k \right\|_2^2 = \\ = -h\left\{ \left\langle \nabla f\left(x^k\right), x^{k+1} - x^k \right\rangle + \frac{1}{2h}\left\| x^{k+1} - x^k \right\|_2^2 \right\} + \frac{1}{2}\left\| x - x^k \right\|_2^2 - \frac{1}{2}\left\| x - x^{k+1} \right\|_2^2 \leq \\ \leq h\mathrm{Prog}^h\left(x^k\right) + \frac{1}{2}\left\| x - x^k \right\|_2^2 - \frac{1}{2}\left\| x - x^{k+1} \right\|_2^2. \tag{2.10}$$

Поскольку $h \leq 1/L$, то, если $x^k, x^{k+1} \in B_{R,Q}\left(x_*\right)$, из (2.3) имеем

$$\mathrm{Prog}^h\left(x^k\right) = -\min_{x \in Q}\left\{ \left\langle \nabla f\left(x^k\right), x - x^k \right\rangle + \frac{1}{2h}\left\| x - x^k \right\|_2^2 \right\} = \\ = -\left( \left\langle \nabla f\left(x^k\right), x^{k+1} - x^k \right\rangle + \frac{1}{2h}\left\| x^{k+1} - x^k \right\|_2^2 \right) = \\ = f\left(x^k\right) - \underbrace{\left( f\left(x^k\right) + \left\langle \nabla f\left(x^k\right), x^{k+1} - x^k \right\rangle + \frac{1}{2h}\left\| x^{k+1} - x^k \right\|_2^2 + \delta \right)}_{\geq f\left(x^{k+1}\right)} + \delta \leq \\ \leq f\left(x^k\right) - f\left(x^{k+1}\right) + \delta. \tag{2.11}$$

Подставляя (2.11) в (2.10), в предположении $x^k, x^{k+1} \in B_{R,Q}(x_*)$, аналогично (1.19), получим

$$h\left\langle \nabla f\left(x^k\right), x^k - x\right\rangle \le h\cdot\left(f\left(x^k\right) - f\left(x^{k+1}\right) + \delta\right) + \frac{1}{2}\left\|x - x^k\right\|_2^2 - \frac{1}{2}\left\|x - x^{k+1}\right\|_2^2. \quad (2.12)$$

По выпуклости $f(x)$ имеем (см. (1.17)):

$$f\left(x^k\right) - f(x) \le \left\langle \nabla f\left(x^k\right), x^k - x\right\rangle, \quad (2.13)$$

также по выпуклости $f(x)$ имеем

$$f\left(\bar{x}^m\right) \le \frac{1}{m}\sum_{k=1}^{m} f\left(x^k\right), \quad (2.14)$$

где (см. (1.21))

$$\bar{x}^m = \frac{1}{m}\sum_{k=1}^{m} x^k .$$

Положим в (2.12) $x = x_*$, если решение не единственно, то выберем то $x_*$, для которого $\left\|x^0 - x_*\right\|_2^2$ минимально.

Суммируя (2.12) с учетом (2.13):

$$h\cdot\left(f\left(x^k\right) - f\left(x_*\right)\right) \le h\cdot\left(f\left(x^k\right) - f\left(x^{k+1}\right) + \delta\right) + \frac{1}{2}\left\|x_* - x^k\right\|_2^2 - \frac{1}{2}\left\|x_* - x^{k+1}\right\|_2^2,$$

т. е.

$$h\cdot\left(f\left(x^{k+1}\right) - f\left(x_*\right)\right) \le h\delta + \frac{1}{2}\left\|x_* - x^k\right\|_2^2 - \frac{1}{2}\left\|x_* - x^{k+1}\right\|_2^2 \quad (2.15)$$

по $k = 0,...,m-1$, получим с учетом (2.14):

$$mh\cdot\left(f\left(\bar{x}^m\right) - f\left(x_*\right)\right) \le mh\delta + \frac{1}{2}\left\|x_* - x^0\right\|_2^2 - \frac{1}{2}\left\|x_* - x^m\right\|_2^2, \quad (2.16)$$

т. е.

$$\frac{1}{2}\left\|x_* - x^m\right\|_2^2 \le mh\cdot\left(\delta - \left(f\left(\bar{x}^m\right) - f\left(x_*\right)\right)\right) + \frac{1}{2}\left\|x_* - x^0\right\|_2^2. \quad (2.17)$$

Вполне естественно (см. § 4, § 5) рассчитывать на то, что метод останавливается когда

$$\varepsilon = f\left(\bar{x}^N\right) - f\left(x_*\right) \approx 2\delta > \delta, \quad (2.18)$$

где

$$\bar{x}^N = \frac{1}{N}\sum_{k=1}^{N} x^k .$$

Как мы увидим далее (см. (2.22) и упражнение 2.2), получить точность $\varepsilon < \delta$ в общем случае не представляется возможным. Поэтому из (2.17) и (2.18) имеем

$$\frac{1}{2}\left\| x_* - x^k \right\|_2^2 \le \frac{1}{2}\left\| x_* - x^0 \right\|_2^2, \; k = 0,...,N \,. \tag{2.19}$$

Другими словами, если $x^0 \in B_{R,Q}\left(x_*\right)$ (а это выполняется по построению $B_{R,Q}\left(x_*\right)$), то для любого $k = 0,..,N$ также верно, что $x^k \in B_{R,Q}\left(x_*\right)$. Таким образом, оговорку о том, что $x^k, x^{k+1} \in B_{R,Q}\left(x_*\right)$ можно опустить.

Строго говоря, мы вывели этот факт, как бы опираясь на него самого (см. оговорку «в предположении $x^k, x^{k+1} \in B_{R,Q}\left(x_*\right)$» около формулы (2.12)). Однако несложно понять, что предполагая условие (2.3) выполненным и вне множества $B_{R,Q}\left(x_*\right)$, то есть на всем $Q$, с теми же параметрами $\left(L,\delta\right)$ мы уже без всяких оговорок получаем, что $x^k \in B_{R,Q}\left(x_*\right)$. Но это означает, что последовательность $\left\{x^k\right\}_{k=0}^{N}$ никогда не выйдет за пределы множества $B_{R,Q}\left(x_*\right)$ и поэтому от того, что именно мы предполагали о выпуклой на всем множестве $Q$ функции $f\left(x\right)$ вне множества $B_{R,Q}\left(x_*\right)$, ничего не зависит – мы никогда не окажемся вне $B_{R,Q}\left(x_*\right)$.

Вернемся к формуле (2.16) при $m = N$, которую перепишем следующим образом:

$$f\left(\bar{x}^N\right) - f\left(x_*\right) \le \frac{1}{2hN}\left\| x_* - x^0 \right\|_2^2 + \delta \,. \tag{2.20}$$

Вспоминая, что на $h$ было условие $h \le 1/L$, аналогично (1.6) выберем

$$h = \frac{1}{L} \,. \tag{2.21}$$

Подставляя (2.21) в (2.20), аналогично (1.20), получим

$$f\left(\bar{x}^N\right) - f\left(x_*\right) \le \frac{LR^2}{2N} + \delta \,. \tag{2.22}$$

**Замечание 2.1 (условие слабой квази-выпуклости).** Вместо $\bar{x}^N$ в приведенных выше формулах можно использовать

$$\hat{x}^N = \arg\min_{k=1,\dots,N} f\left(x^k\right). \tag{2.23}$$

Несложно заметить, что в случае подхода с $\hat{x}^N$ приведенные выше рассуждения используют лишь свойство (2.13) с $x = x_*$:

$$\left(f\left(x^k\right) - f\left(x_*\right)\right) \le \left\langle \nabla f\left(x^k\right), x^k - x_* \right\rangle,$$

т. е. «полноценная» выпуклость $f\left(x\right)$ не требуется. По-видимому, это наблюдение восходит к Н.З. Шору [82, п. 5.4].

Отметим, что условие (2.13) также ослабляют следующим образом (следует сравнить с однородными относительно $x_*$ функциями [82, п. 4 § 3, гл. 3] и звездной выпуклостью [329, 464]):

$$\alpha \cdot \left(f\left(x^k\right) - f\left(x_*\right)\right) \le \left\langle \nabla f\left(x^k\right), x^k - x_* \right\rangle, \ \alpha \in \left(0,1\right]. \tag{2.24}$$

Условие (2.24) иногда называют *условием α-слабой квази-выпуклости* функции $f\left(x\right)$. В последнее время оно стало достаточно популярно в связи с приложениями, возникающими в *Глубоком Обучении* (*Deep Learning*) [316]. Несложно показать, что приведенные выше рассуждения переносятся и на этот случай. При этом оценка (2.22) «портится» следующим образом [308]:

$$f\left(\hat{x}^N\right) - f\left(x_*\right) \le \frac{LR^2}{2\alpha N} + \delta,$$

где $\hat{x}^N$ определяется формулой (2.23).

В приложениях часто используют также такое неравенство

$$\begin{aligned} f\left(\bar{x}^N\right) - f\left(x_*\right) &\le \underbrace{\sup_{x\in Q} \frac{1}{N} \sum_{k=0}^{N-1} \left\langle \nabla f\left(x^k\right), x^k - x \right\rangle}_{\text{сертификат точности}} \le \\ &\le \frac{f\left(x^0\right) - f\left(x^N\right) + L \sup_{x\in Q} \left\| x - x^0 \right\|_2^2}{2N} + \delta, \end{aligned} \tag{2.25}$$

где (следует сопоставить с (1.21))

$$\bar{x}^N = \frac{1}{N} \sum_{k=0}^{N-1} x^k.$$

Введенный в (2.25) *сертификат точности* (accuracy certificate) играет ключевую роль в обосновании прямодвойственности исследуемого метода [443] (см. также § 4). Сертификат точности и его аналоги из § 4 являются вычислимыми (не требуют знания, как правило, неизвестных значений $f(x_*)$ или $R^2$) и потому могут использоваться в качестве критерия останова методов. ∎

Предположим, что неравенство (2.3) имеет вид (см. также замечание 1.3):

$$f(y) \le f(x) + \langle \nabla f(x), y - x \rangle + \frac{L}{2}\|y - x\|^2 + \delta . \tag{2.26}$$

Для этого, например, достаточно, чтобы имело место неравенство (см. также (1.26)), аналогичное (2.4):

$$\|\nabla f(y) - \nabla f(x)\|_* \le L_\nu \|y - x\|^\nu , \; \nu \in [0,1], \; L_0 < \infty . \tag{2.27}$$

Тогда (2.26) имеет место с константой $L$, рассчитываемой по формуле (2.5). Попробуем, следуя А.С. Немировскому [158], распространить метод градиентного спуска на этот случай. Как уже отмечалось в предыдущем параграфе, сделать это за счет такого обобщения метода не получается (см. формулу (1.27)):

$$x^{k+1} = \arg\min_{x \in Q} \left\{ \left\langle h \nabla f\left(x^k\right), x - x^k \right\rangle + \frac{1}{2}\left\|x - x^k\right\|^2 \right\}. \tag{2.28}$$

Причина прежде всего в том, что, например, $\left\|x - x^k\right\|_1^2$ не есть сильно выпуклая функция в 2-норме и тем более в 1-норме. Отсутствие этого свойства, как мы увидим чуть ниже, и не позволяет сделать необходимое обобщение. Рассмотрим, однако, близкий к (2.28) метод

$$x^{k+1} = \arg\min_{x \in Q} \left\{ \left\langle h \nabla f\left(x^k\right), x - x^k \right\rangle + V\left(x, x^k\right) \right\}. \tag{2.29}$$

Получим условия на выпуклую по $x$ функцию $V\left(x, x^k\right)$, при которых приведенная в § 2 *конструкция* вывода основных оценок *сохраняется*. Первым ключевым местом, в котором использовались свойства функции $V(x,y) = \|x - y\|_2^2 / 2$, было неравенство (2.8). В случае (2.29) неравенство (2.8) должно было бы принять вид

$$\begin{gathered} \left\langle h \nabla f\left(x^k\right), x^{k+1} - x \right\rangle \le \left\langle \nabla_{x^{k+1}} V\left(x^{k+1}, x^k\right), x - x^{k+1} \right\rangle \overset{1}{=} \\ \overset{1}{=} V\left(x, x^k\right) - V\left(x, x^{k+1}\right) - V\left(x^{k+1}, x^k\right). \end{gathered} \tag{2.30}$$

Таким образом, достаточно потребовать выполнение равенства 1[25] тождественно по $x$. Например, если считать, что имеет место следующее представление:

$$V(x,y) = d(x) - d(y) - \langle \nabla d(y), x - y \rangle \tag{2.31}$$

с выпуклой функцией $d(x)$, то тождество 1 также имеет место. Вторым и заключительным ключевым местом было неравенство (2.11), которое в нашем случае останется верным, если

$$V\left(x^{k+1}, x^k\right) \geq \frac{1}{2}\left\|x^{k+1} - x^k\right\|^2. \tag{2.32}$$

Для этого достаточно, чтобы в представлении (2.31) функция $d(x)$ была 1-сильно выпукла в выбранной норме $\|\ \|$. Функцию $d(x)$ называют *прокс-функцией*, а функцию $V(x,y)$ – порожденным ею *расхождением* или *дивергенцией Брэгмана* (Bregman divergence) [12, 23, 158, 180]. Отметим, что для метода (2.29) в оценку (2.22) будет входить $R^2 = 2V\left(x_*, x^0\right)$. Если решение не единственно, то оценка (2.22) будет верна в том числе и для того решения $x_*$, которое доставляет минимум $R^2$. Отметим также, что для «сохранения конструкции» достаточно, чтобы условия (2.26), (2.27) имели место только при

$$x, y \in \left\{x \in Q:\ V(x_*, x) \leq R^2\right\}.$$

Рассуждения здесь аналогичны рассуждениям, использованным при выводе (2.19).

Примеры прокс-функций для множеств $Q$ вида шаров в различных нормах собраны в табл. 1 [23, 158]. Параметр

$$a = \frac{2\ln n}{2\ln n - 1} \simeq 1 + \frac{1}{2\ln n}.$$

Дополнительно к тому, что приведено в табл. 1, особо отметим «Spectrahedron setup» [180, п. 4.3]. Приведенные в табл. 1 прокс-функции можно распространить и на прямые произведения шаров [158, п. 5.3.3].

---

[25] Условие на функцию $V(x,y)$, накладываемое равенством 1, можно понимать и как неравенство в сторону «$\leq$». Однако, как будет видно в дальнейшем, удается подобрать функцию $V(x,y)$ и с равенством, что дополнительно привносит меньше грубости в рассуждения.

Таблица 1

| $Q=B_p^n(1)$ | $1\le p\le a$ | $a\le p\le 2$ | $2\le p\le\infty$ |
|---|---|---|---|
| $\\|\ \\|$ | $\\|\ \\|_1$ | $\\|\ \\|_p$ | $\\|\ \\|_2$ |
| $d(x)$ | $\frac{1}{2(a-1)}\\|x\\|_a^2$ | $\frac{1}{2(p-1)}\\|x\\|_p^2$ | $\frac{1}{2}\\|x\\|_2^2$ |
| $R^2$ | $O(\ln n)$ | $O\left((p-1)^{-1}\right)$ | $O\left(n^{1/2-1/p}\right)$ |

По-видимому, в табл. 1 в общем случае нельзя избавиться от дополнительного $\ln n$ фактора (множителя) в оценке $R^2$ с помощью дивергенции Брэгмана по сравнению с оценкой $R^2$, равной квадрату соответствующей нормы. Однако эта плата позволяет переносить все основные свойства, присущие работе в евклидовом случае, на множества типа шаров в 1-норме. Кроме того, во всех использующихся примерах прокс-структур эта мультипликативная плата по порядку не превышает $\ln n$. Такой порядок у этой константы будет, например, для единичного симплекса (см. также текст ниже) при

$$d(x)=\sum_{i=1}^{n}x_i\ln x_i,\ V(x,y)=\sum_{i=1}^{n}x_i\ln\left(x_i/y_i\right)\ \text{и}\ x^0=\left(n^{-1},\dots,n^{-1}\right).$$

Все приведенные в табл. 1 прокс-функции 1-сильно выпуклы в указанных нормах на всем пространстве. Поэтому их можно использовать и в том случае, когда мы заранее не знаем, где локализовано решение [3]. Например, если стартовать с $x^0=0$ и заранее знать, что решение разреженно (большинство компонент равно нулю), то, естественно, выбирать 1-норму и соответствующую прокс-функцию (см. табл. 1). Действительно, в этом случае можно ожидать, что $R^2$ в оценке (2.22) от выбора нормы не будет сильно зависеть, в то время как для константы

$$L:=L^p=\sup_{x\in Q}\max_{\\|h\\|_p\le 1}\left\langle h,\nabla^2 f(x)h\right\rangle$$

отличие может быть в $n$ раз, поскольку $L^2/n\le L^1\le L^2$. Скажем, для функции (1.30)

$$f(x)=\frac{1}{2}\langle Ax,x\rangle-\langle b,x\rangle$$

несложно получить, что $L^1=\max_{i,j=1,\dots,n}\left|A_{ij}\right|$, а $L^2=\lambda_{\max}(A)$.

◊ Аналогично можно определить и коснтанту сильной выпуклости относительно $p$-нормы

$$\mu^p = \inf_{x\in Q} \min_{\|h\|_p \le 1} \left\langle h, \nabla^2 f(x) h \right\rangle .$$

Строго говоря, приведенные здесь определения констант $L^p$ и $\mu^p$ справедливы только при дополнительном предположении о существовании и конечности матрицы Гессе $\nabla^2 f(x)$ при $x \in Q$. Последнее предположение не всегда имеет место. В частности, функция Хьюбера (1.43) имеет Липшицев градиент, однако гессиан определен не везде. ◊

Как мы увидим в дальнейшем (см. упражнение 4.6), задача (2.29) при известном векторе $\nabla f\left(x^k\right)$ для примеров из табл. 1 решается за время $\mathrm{O}\left(n \ln^2 (n/\varepsilon)\right)$, где $\varepsilon$ – точность решения в смысле сходимости по аргументу. Выделим особо один частный случай, когда эту оценку можно улучшить до $\mathrm{O}(n)$:

$$Q = S_n(1) = \left\{ x \in \mathbb{R}^n_+ : \sum_{i=1}^n x_i = 1 \right\},$$

т. е. $Q$ – единичный симплекс в $\mathbb{R}^n$. В этом случае, выбирая

$$d(x) = \sum_{i=1}^n x_i \ln x_i , \tag{2.33}$$

получим

$$x_i^{k+1} = \frac{x_i^k \exp\left(-h \partial f\left(x^k\right) / \partial x_i\right)}{\sum_{j=1}^n x_j^k \exp\left(-h \partial f\left(x^k\right) / \partial x_j\right)}, \quad i = 1,\ldots,n .$$

Отметим, что метод (2.29) при $h \le 1/L$ ввиду $V(x,y) \ge \|x-y\|^2 / 2$ имеет геометрическую интерпретацию, аналогичную обычному градиентному спуску (см. замечания 1.2, 1.3).

Конструкция (2.7) – (2.12) с учетом (2.29), (2.31) может быть распространена на более общий класс задач. В следующем параграфе приводится основная схема такого распространения, лежащая в основе получения результатов в наиболее общем виде.

**Упражнение 2.1 (нижние оценки – негладкий случай).** Покажите, что в условии (2.4) с $\nu = 0$ оценка (2.22) для метода (2.6), (2.21) с

$$h = \frac{1}{L} = \frac{2\delta}{L_0^2} = \frac{\varepsilon}{L_0^2}, \text{ где } \varepsilon = \frac{LR^2}{2N} + \delta = \frac{LR^2}{2N} + \frac{\varepsilon}{2},$$

т. е. с[26]

$$h = \frac{\varepsilon}{L_0^2} = \frac{R}{L_0\sqrt{N}}, \tag{2.34}$$

будет иметь вид

$$f\left(\bar{x}^N\right) - f\left(x_*\right) \le \frac{L_0 R}{\sqrt{N}}. \tag{2.35}$$

Покажите, что в классе методов

$$x^{k+1} \in x^0 + \text{Lin}\left\{\partial f\left(x^0\right), ..., \partial f\left(x^k\right)\right\} \tag{2.36}$$

оценка (2.35) при $N \le n - 1$, где $n = \dim x$, не может быть улучшена с точностью до числового множителя.[27] Предполагается, что в рассматриваемом классе методов (2.36) нельзя выбирать субградиент из субдифференциала. Это означает, что при оценивании скорости сходимости необходимо исходить из того, что субградиент из субдифференциала может выбираться наиболее неудачным образом.

**Указание.** Следует сравнить это упражнение с упражнением 1.3. Согласно (2.5) $L = L_0^2/(2\delta)$. С учетом этого найдите минимум правой части неравенства (2.22) по $\delta > 0$. Получите отсюда оценку (2.35).

Для получения нижней оценки воспользуйтесь, например, [72, п. 3.2.1], [180, п. 3.5]. По заданному $N \le n - 1$ определите

$$f(x) = F_{N+1}(x) = L_0 \max_{1 \le i \le N+1} x_i + \frac{\mu}{2}\|x\|_2^2, \quad \mu = \frac{L_0}{R\sqrt{N+1}}.$$

Из решения задачи

---

[26] То, что для негладких задач шаг градиентного метода $h = \text{const} \cdot \varepsilon / L_0^2$, а для гладких – $h = \text{const}/L_1$, можно было понять и из П-теоремы теории размерностей [4, 47]. Отметим также, что в негладком случае еще возможен вариант $h = \text{const} \cdot R/L_0$. Подробнее об этом см. в упражнении 2.6.

[27] На самом деле если ограничиться только классом методов вида (2.36) с фиксированными шагами (не зависящими от оптимизируемой функции, см. замечание 1.5), то в (2.35) можно немного улучшить только знаменатель в правой части неравенства $\sqrt{N} \to \sqrt{N+1}$, см. [241].

$$L_0\tau + \frac{\mu\cdot(N+1)}{2}\tau^2 \to \min_\tau$$

определите $\tau_* = -R\big/\sqrt{N+1}$, $x_* = (\underbrace{\tau_*,\ldots,\tau_*}_{N+1},0,\ldots,0)$. Тогда

$$\|x_*\|_2^2 = (N+1)\tau_*^2 = R^2,\; f(x_*) = \min_{x\in\mathbb{R}^n} F_{N+1}(x) = -L_0R\Big/\left(2\sqrt{N+1}\right) = F_{N+1}^*.$$

Если $x^0 = 0$, тогда для метода вида (2.36) после $k \le N$ итераций при специальном (наиболее неблагоприятном) выборе субградиентов из субдифференциалов имеет место условие: $x_i^k = 0$ при $i > k$. Действительно, если это условие верно для шага $k-1$, то для того, чтобы появились новые ненулевые компоненты у вектора $x^k$ (на шаге $k$) по сравнению с $x^{k-1}$ необходимо, чтобы $\max\limits_{1\le i\le N} x_i^{k-1} \le 0$ в этом случае субдифференциал $\partial \max\limits_{1\le i\le N} x_i^{k-1}$ будет определяться выпуклой комбинацией таких единичных ортов $e_i$, для которых $x_i^{k-1} = 0$. В частности, можно взять $\partial \max\limits_{1\le i\le N} x_i^{k-1} = e_k$. Такой выбор (в случае, если $\max\limits_{1\le i\le N} x_i^{k-1} < 0$) обусловлен желанием обеспечить как можно более медленную скорость сходимость метода вида (2.36), что в конечном итоге должно приводить к наиболее точным нижним оценкам. Таким образом, поскольку $x_{N+1}^N = 0$, то $\max\limits_{1\le i\le N+1} x_i^N = 0$. Значит $F_{N+1}\left(x^N\right) \ge 0$. Следовательно,

$$F_{N+1}\left(x^N\right) - F_{N+1}^* \ge -F_{N+1}^* = \frac{L_0R}{2\sqrt{N+1}} = \frac{L_0^2}{2\mu\cdot(N+1)}. \qquad (2.37)$$

Заметим, что оценка (2.37) одновременно является нижней оценкой в классе методов (2.36) для $\mu$-сильно выпуклых задач в 2-норме. Оценка вида (2.37) будет нижней и для более общего класса методов [70, гл. 4].

Стоит сделать несколько замечаний. Исходная задача – задача безусловной оптимизации, таким образом под $R$ в (2.37) стоит понимать расстояние от точки старта до решения, собственно, как и в (2.35). Не существует сильно выпуклой функции заданной на всем пространстве, у которой была бы равномерно ограничена константа Липшица (см. также конец § 5). Однако, константы, которые входят в оценку (2.37), характеризуют функцию в шаре $B_R\left(x^0\right)$, в котором в виду (2.19) естественно было бы ожидать пребывание всей траектории метода вида (2.36). В действительности, константа $L_0$, входящая в (2.37), немного меньше настоящей константы Липшица функции $F_N(x)$ в шаре $B_R(x^0)$ в виду наличия у $F_N(x)$ композитного (сильно выпуклого) слагаемого $\mu\|x\|_2^2/2$. Несложно учесть это слагаемое и должным образом скорректировать оценку (2.37). Общий вывод при это сохранится [72, п. 3.2.1], [180, п. 3.5]. Однако, в виду примера 3.1 полезно заметить, что и исходном виде оценка (2.37) представляет ценность, поскольку правильно отражает сложность класса задач композитной выпуклой оптимизации. ■

**Упражнение 2.2**. Покажите, что при $\delta > 0$ оценку (2.22) нельзя принципиально улучшить в части зависимости от $N$, не ухудшив в части зависимости ее от $\delta > 0$.

**Указание.** Проведем рассуждения, следуя [227]. Допустите противное, т. е. что можно получить следующую оценку:

$$f(\bar{x}^N) - f(x_*) \le C_1 \frac{LR^2}{N^{1+\gamma}} + C_2\delta, \ \gamma > 0. \tag{2.38}$$

Рассмотрите негладкий случай, т. е. используйте (2.4) с $\nu = 0$. Согласно нижним оценкам (см. упражнение 2.1), для любого $N \le n$, и для любого метода вида (2.36) существует такая выпуклая функция из гёльдерова класса с $\nu = 0$ и константой $L_0$, что

$$f(\bar{x}^N) - f(x_*) \ge \frac{L_0 R}{2\sqrt{N}}.$$

Покажите, что при достаточно большом $N$ (а следовательно, и $n$) это противоречит (2.38). Для этого, согласно (2.5), подставьте в (2.38) $L = L_0^2/(2\delta)$ и специально подберите $\delta = L_0 R\sqrt{C_1/(2C_2 N^{1+\gamma})}$. Тогда

$$f(\bar{x}^N) - f(x_*) \le C_1 \frac{L_0^2 R^2}{2\delta N^{1+\gamma}} + C_2\delta = \sqrt{2C_1C_2}\,\frac{L_0 R}{\sqrt{N^{1+\gamma}}}. \ ■$$

**Упражнение 2.3 (техника рестартов). 1)** Как из метода, работающего по оценке, аналогичной (2.35),

$$f(\bar{x}^N) - f(x_*) \le \frac{L_0\left\|x_* - x^0\right\|_2}{\sqrt{N}} + \delta,$$

где $\delta > 0$ достаточно мало, получить метод, который для $\mu$-сильно выпуклых задач в 2-норме работает по оценке

$$f(\tilde{x}^N) - f(x_*) \le \frac{128 L_0^2}{\mu N} + 2\delta,$$

точнее, по оценке

$$f\left(\tilde{x}^N\right)-f\left(x_*\right)\le\frac{256L_0^2}{\mu N},\quad N\le\frac{128L_0^2}{\mu\delta}? \tag{2.39}$$

◊ При $\delta=0$ существует много способов уменьшения константы 512 в оценке (2.39) на два порядка [321, 351, 354, 382, 496]. Однако при этом следует отметить, что согласно упражнению 2.2, оценка (2.39) неулучшаема с точностью до мультипликативной константы. Это общее свойство техники рестартов, описанной в указании к этому упражнению: *из оптимального метода рестарты получают оптимальный метод для сильно выпуклых задач*. Во всяком случае пока не удалось придумать контрпримера, равно как не удалось придумать ситуацию с перенесением метода на сильно выпуклые задачи, в которой бы не было своего варианта рестарт метода. В основе техники рестартов лежит простая идея: рестартовать метод, т. е. запускать по-новому, с возможно новыми значениями параметров, каждый раз в момент, когда есть гарантия, что генерируемая последовательность оказалась на «расстоянии» (или невязке по функции) в два раза ближе к решению по сравнению с моментом последнего рестарта. Важно подчеркнуть, что основное свойство описываемой техники (в плане ее обоснования) существенно завязано именно на рестарты по расстоянию от текущей точки до решения (или невязке по функции). Использование более удобных критериев для рестарта, например для гладких выпуклых задач безусловной оптимизации можно было бы использовать легко вычислимую величину нормы градиента функционала, к сожалению, не позволяет строго обосновать сохранение свойства оптимальности метода в общем случае [474], см. также замечание 5.3. ◊

**2)** Попробуйте обобщить этот результат на случай, когда используется произвольная норма $\|\ \|$ и не евклидова прокс-структура, однако имеет место следующее свойство: $d\left(x\right)\le C_n\left\|x\right\|^2$. Покажите, что при таком предположении оценка (2.39) немного ухудшится: $\mu\to\mu/\left(2C_n\right)$. Заметим, что для прокс-функций из табл. 1 $C_n=\mathrm{O}\left(\ln n\right)$.[28]

**Указание. 1)** См., например, [68, 70, 351, 354, 438, 446, 507]. Покажите, что при $N_1\le L_0^2R_0^2/\delta^2$

[28] По-видимому, такой оценки $C_n$ можно добиться во всех интересных для практики случаях при правильном выборе прокс-функции. Заметим, что выбор прокс-функции (2.33) для $Q=S_n\left(1\right)$ в этом смысле не будет правильным.

$$\frac{\mu}{2}\underbrace{\left\|\bar{x}^{N_1}-x_*\right\|_2^2}_{R_1^2}\overset{1}{\le} f\left(\bar{x}^{N_1}\right)-f\left(x_*\right)\overset{2}{\le}\frac{L_0\overbrace{\left\|x^0-x_*\right\|_2}^{R_0}}{\sqrt{N_1}}+\delta\overset{3}{\le}\frac{2L_0\left\|x^0-x_*\right\|_2}{\sqrt{N_1}}. \quad (2.40)$$

Неравенство 1 имеет место ввиду $\mu$-сильной выпуклости $f(x)$ (см. (1.14)), неравенство 2 – ввиду (2.35), а неравенство 3 – ввиду $N_1 \le L_0^2 R_0^2/\delta^2$. Выберите $N_1$ из условия $R_1 = R_0/2$. Из (2.40) получите $N_1 \simeq 64L_0^2/\left(\mu^2 R_1^2\right)$.

◊ В этом месте появляется проблема с практической реализацией схемы рестартов. Дело в том, что в такой реализации предписано сделать число итераций, зависящее явно от параметра $\mu$, который, как правило, либо просто неизвестен, либо грубо оценен, не говоря уже о возможности локальной настройки метода на значение этого параметра, отвечающего текущему положению метода. В отличие от адаптивной настройки на гладкость задачи (см. § 5), на данный момент не известны общие способы настройки на параметр сильной выпуклости лучше, чем рестарты по этому неизвестному параметру [446, 472]: решаем задачу с $\mu = \mu_0$, метод не сходится, полагаем $\mu := \mu/2$, снова решаем и т. д., пока не получим сходимость. При таком подходе есть некоторые тонкости с детектированием сходимости. Несложно показать, что число дополнительных вычислений при этом увеличится не более чем в 8 раз [20]. Впрочем, в некоторых случаях можно более изящно решать отмеченную проблему [265, 270, 370, 474, 499, 507, 544]. Мы вернемся к этому вопросу в конце § 5. ◊

Далее, после $N_1$ итераций, рестартуйте исходный метод и положите $x^0 := \bar{x}^{N_1}$. Определите $N_2 \le L_0^2 R_1^2/\delta^2$ из условия $R_2 = \left\|\bar{x}^{N_2}-x_*\right\|_2 = R_1/2$. Получите, что $N_2 \simeq 64L_0^2/\left(\mu^2 R_2^2\right)$ и т. д. … После $k$ таких рестартов общее число итераций будет

$$N = N_1 + \ldots + N_k \simeq \frac{256L_0^2}{\mu^2 R_0^2}\left(1+4^1+\ldots+4^{k-1}\right) \approx \frac{4^{k+4}L_0^2}{\mu^2 R_0^2}. \quad (2.41)$$

Обозначьте через $\tilde{x}^N$ то, что получается после $N$ описанных итераций рестартованным методом. Обозначьте через $\varepsilon = f\left(\tilde{x}^N\right)-f\left(x_*\right)$. Из (2.40) получите

$$\varepsilon = \frac{\mu R_k^2}{2} = \frac{2L_0 R_{k-1}}{\sqrt{N_k}}. \quad (2.42)$$

Покажите, что из $N_k \le L_0^2 R_{k-1}^2 / \delta^2$, с учетом (2.42), следует $\varepsilon \ge 2\delta$. Покажите, что из (2.41) следует

$$\frac{\mu R_k^2}{2} = \frac{128 L_0^2}{\mu N}. \tag{2.43}$$

Объединяя (2.42), (2.43) и $\varepsilon \ge 2\delta$, получите оценку (2.39). ■

**Упражнение 2.4**. Покажите, что прокс-функции, собранные в табл. 1, действительно 1-сильно выпуклы в соответствующих нормах.

**Указание.** См. [158, п. 5.6]. ■

**Упражнение 2.5.** Предложите норму и прокс-функцию для задачи оптимизации на прямом произведении симплексов.

**Указание.** См. [29]. ■

**Упражнение 2.6 (адаптивные субградиентные методы [72, 152, 248, 455, 554]).** Покажите, что для задачи негладкой выпуклой оптимизации (2.1) метод

$$x^{k+1} = x^k - h_k \nabla f\left(x^k\right)$$

при

$$h_k \equiv \frac{\varepsilon}{L_0^2},\ h_k \equiv \frac{R}{L_0\sqrt{N}},\ h_k \equiv \frac{\varepsilon}{\left\|\nabla f\left(x^k\right)\right\|_2^2},\ h_k \equiv \frac{R}{\left\|\nabla f\left(x^k\right)\right\|_2 \sqrt{N}},$$

где $R = \left\|x_* - x^0\right\|_2$, будет сходиться согласно оценке (2.35) с

$$\bar{x}^N = \hat{x}^N \in \operatorname{Arg}\min_{k=0,\dots,N-1} f\left(x^k\right).$$

Покажите, что в общем случае решения задачи оптимизации на множестве простой структуры и использовании неевклидовой нормы результат останется верным с заменой метода (2.1) *на метод зеркального спуска (2.29)* [152, 158, 180] и[29]

$$R = \sqrt{2V\left(x_*, x^0\right)},\ \left\|\nabla f\left(x^k\right)\right\|_2 \to \left\|\nabla f\left(x^k\right)\right\|_*.$$

Используя работы [152, 351, 543], предложите обобщение алгоритма (2.29) на задачи негладкой выпуклой оптимизации с негладкими выпуклыми функциональными ограничениями.

**Указание.** См. [4]. Для простоты рассмотрим здесь задачу безусловной оптимизации и ограничимся сначала случаем

---

[29] При аккуратном анализе, необходимом, например, для получения оптимальных оценок в приложении рассматриваемых методов к задаче о многоруких бандитах [31], [181, items 5.2, 7.1], вместо $\left\|\nabla f\left(x^k\right)\right\|_*$ можно писать более точное выражение.

$$h_k \equiv h = \frac{R}{L_0\sqrt{N}}.$$

Из структуры метода следует, что

$$\left\|x - x^{k+1}\right\|_2^2 = \left\|x - x^k + h\nabla f\left(x^k\right)\right\|_2^2 = \left\|x - x^k\right\|_2^2 + 2h\left\langle \nabla f\left(x^k\right), x - x^k\right\rangle +$$

$$+h^2\left\|\nabla f\left(x^k\right)\right\|_2^2 \le \left\|x - x^k\right\|_2^2 + 2h\left\langle \nabla f\left(x^k\right), x - x^k\right\rangle + h^2 L_0^2.$$

Отсюда (при $x = x_*$) следует, что

$$f\left(\hat{x}^N\right) - f\left(x_*\right) = \min_{k=0,\ldots,N-1} f\left(x^k\right) - f\left(x_*\right) \le$$

$$\le \frac{1}{N}\sum_{k=0}^{N-1} f\left(x^k\right) - f\left(x_*\right) \le \frac{1}{N}\sum_{k=0}^{N-1}\left\langle \nabla f\left(x^k\right), x^k - x_*\right\rangle \le$$

$$\le \frac{1}{2hN}\sum_{k=0}^{N-1}\left\{\left\|x_* - x^k\right\|_2^2 - \left\|x_* - x^{k+1}\right\|_2^2\right\} + \frac{hL_0^2}{2} =$$

$$= \frac{1}{2hN}\left(\left\|x_* - x^0\right\|_2^2 - \left\|x_* - x^N\right\|_2^2\right) + \frac{hL_0^2}{2}.$$

Полагая $h = R\big/\left(L_0\sqrt{N}\right)$, получите

$$f\left(\hat{x}^N\right) - f\left(x_*\right) \le \frac{L_0 R}{\sqrt{N}}.$$

Заметим, что если $f\left(x\right)$ имеет $L_1$ липшицев градиент, то вместо оценки $\left\|\nabla f\left(x^k\right)\right\|_2^2 \le L_0^2$, которая для негладких задач в общем случае не может быть улучшена по мере приближения к решению, можно использовать оценку (1.8) $\left\|\nabla f\left(x^k\right)\right\|_2^2 \le 2L\cdot\left(f\left(x^k\right) - f\left(x_*\right)\right)$, которая уже отражает уменьшение нормы градиента по мере приближения к решению.

В результате при выборе шага $h = 1/(2L_1)$ получается следующая оценка (следует сравнить с оценкой (2.22))

$$f(\hat{x}^N) - f(x_*) \le \frac{2L_1R^2}{N}.$$

Для негладких задач выпуклой оптимизации на множествах простой структуры рассуждения будут аналогичными тем, что были приведены в § 2 при выводе формулы (2.10). Только для оценки $\mathrm{Prog}^h(x^k)$ придется ограничиться неравенством Гёльдера (в виду отсутствия липшицевости градиента)

$$\mathrm{Prog}^h(x^k) = \max_{x\in Q}\left\{\langle \nabla f(x^k), x^k - x\rangle - \frac{1}{2h}\|x^k - x\|_2^2\right\} \le$$

$$\le \max_z\left\{\langle \nabla f(x^k), z\rangle - \frac{1}{2h}\|z\|_2^2\right\} = \frac{h}{2}\|\nabla f(x^k)\|_2^2.$$

Простая структура приведенных выше рассуждений позволяет аналогичным образом провести рассуждения и в случае адаптивного выбора шага $h_k$. Ограничимся случаем

$$h_k \equiv \frac{R}{\|\nabla f(x^k)\|_2\sqrt{N}}.$$

В этом случае

$$\frac{R}{\sqrt{N}}\sum_{k=0}^{N-1}\left\langle \frac{\nabla f(x^k)}{\|\nabla f(x^k)\|_2}, x^k - x_*\right\rangle \le \frac{1}{2}\sum_{k=0}^{N-1}\left\{\|x_* - x^k\|_2^2 - \|x_* - x^{k+1}\|_2^2\right\} + \frac{R^2}{2} \le R^2.$$

Откуда следует (см. упражнение 2.7), что

$$f(\hat{x}^N) - f(x_*) \le \frac{L_0}{N}\min_{k=0,\dots,N-1}\left\langle \frac{\nabla f(x^k)}{\|\nabla f(x^k)\|_2}, x^k - x_*\right\rangle \le \frac{L_0R}{\sqrt{N}}.$$

Отметим, что адаптивность (суб-)градиентных методов можно получить за счет других соображений, см. § 5. Причем подход § 5 является более универсальным, в частности потому, что позволяет работать с общей концепцией модели функции, рассмотренной в § 3. Отметим, что приведенные в данном упражнении адаптивные методы непонятно как обобщать даже на частный случай общей модельной концепции: на задачи композитной оптимизации [19], см. пример 3.1. ■

**Упражнение 2.7 (квазивыпуклые функции [72, пп. 3.2.2–3.2.4], [75, п. 1.5], [322, 543]).** Функция $f(x)$ называется квазивыпуклой на выпуклом множестве $Q$, если для любого $\alpha \in [0,1]$ и любых $x, y \in Q$ выполняется

$$f(\alpha x + (1-\alpha) y) \le \max\{f(x), f(y)\}.$$

Покажите, что если функция $f(x)$ – квазивыпуклая, то множества Лебега этой функции (множества вида $\{x \in Q:\ f(x) < C\}$) будут выпуклыми. Любая выпуклая функция будет квазивыпуклой, обратное в общем случае не верно.

Покажите (можно ограничиться случаем евклидовой нормы), что если $\nabla f(x) \ne 0$, то

$$f(x) - f(x_*) \le \omega(\mathrm{v}(x)),$$

где модуль непрерывности

$$\omega(t) = \max\{f(x) - f(x_*):\ \|x - x_*\| \le t\},$$

$$\mathrm{v}(x) = \frac{\langle \nabla f(x), x - x_* \rangle}{\|\nabla f(x)\|_*}.$$ [30]

Используя упражнение 2.6 предложите численный метод решения задач оптимизации с квазивыпуклым функционалом.

**Указание.** Ограничимся установлением неравенства

$$f(x) - f(x_*) \le \omega(\mathrm{v}(x))$$

в евклидовым случаем $\|\ \| = \|\ \|_2$ (в неевклидовом случае см. [152, 459]).

На рис. 8 кривой линией изображена линия уровня функции $f(x)$. Обозначим через $y$ проекцию точки $x_*$ на касательную к этой линии уровня в точке $x$. По предположению квазивыпуклости $f(x)$ (множества Лебега выпуклые) точки $x_*$ и $y$ лежат по разные стороны от этой линии, поэтому

$$f(x) - f(x_*) \le f(y) - f(x_*),$$

но последнюю величину в виду определения функции $\omega(t)$ и того простого наблюдения, что в евклидовом случае расстояние между точками $x_*$ и $y$ равно $\text{v}(x)$, можно оценить сверху $\omega(\text{v}(x))$.

---

[30] Здесь не приводится аккуратное определение $\nabla f(x)$. Чтобы лучше понять, какие на этом пути возникают сложности, можно вспомнить пример невыпуклой функции из § 1, с помощью которого получалась нижняя оценка сложности класса задач невыпуклой гладкой оптимизации. Описанная там функция была плохая с точки зрения локального оракула, выдающего локальную информацию о функции в запрошенной точке. В то же время, данная функция является квазивыпуклой, а значит (согласно упражнению 2.7) вполне эффективно может быть прооптимизированна. На первый взгляд кажется, что получается противоречие. Однако противоречия тут нет, поскольку эта функция во всех кубиках кроме одного тождественно равна нулю и, следовательно, $\nabla f(x) \equiv 0$ во всех таких кубиках. Таким образом, если понимать под $\nabla f(x)$, например, некоторый «перпендикуляр» к множеству Лебега, то это сразу приводит к нарушению условия локальности.

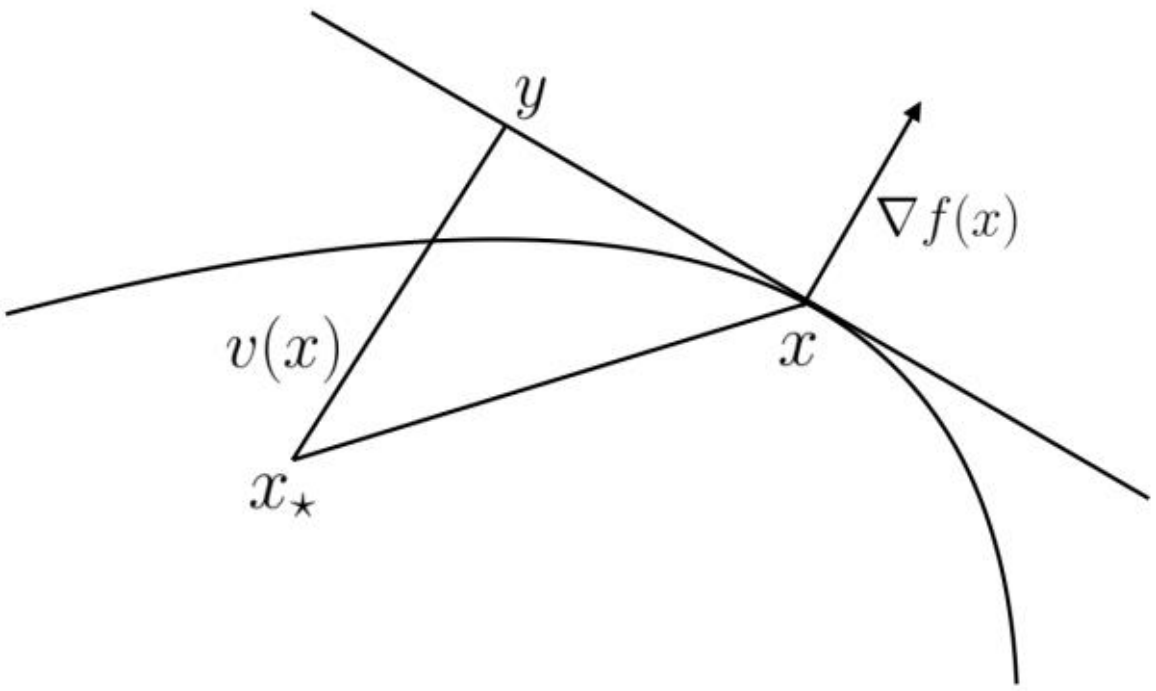

Рис. 8

В связи с полученным результатом, интересно заметить, что

$$\text{v}(x) \le \|x - x_*\|,$$

поэтому неравенство

$$f(x) - f(x_*) \le \omega(\text{v}(x))$$

является более тонким, чем просто условие Липшицевой непрерывности функции $f(x)$.

Численный метод можно получить, если специальным образом выбирать шаг в субградиентном методе из упражнении 2.6

$$h_k = \frac{R}{\|\nabla f(x^k)\|_* \sqrt{N}}. \blacksquare$$

# § 3. Общая схема получения оценок скорости сходимости. Структурная оптимизация

Как и в § 2, рассмотрим задачу выпуклой оптимизации (2.1):

$$f(x) \to \min_{x \in Q}.$$

Сначала обобщим условие (2.26) (см. также (2.3)). Будем говорить, что имеем $(\delta, L)$-*модель функции* $f(x)$ в точке $x$ (относительно нормы $\| \ \|$), и обозначать эту модель $\left(f_\delta(x); \psi_\delta(y,x)\right)$, если для любого $y \in Q$ справедливо неравенство [95]:

$$0 \le f(y) - \left(f_\delta(x) + \psi_\delta(y,x)\right) \le \frac{L}{2}\|y-x\|^2 + \delta, \tag{3.1}$$

где $\psi_\delta(y,x)$ – выпуклая функция по $y$, $\psi_\delta(x,x) = 0$, $\delta > 0$.

◊ Из (3.1) при $y = x$ следует, что $0 \le f(x) - f_\delta(x) \le \delta$, поэтому под $(\delta, L)$-моделью функции $f(x)$ в точке $x$ можно понимать только такую выпуклую по $y \in Q$ функцию $\psi_\delta(y,x)$, что для всех $y \in Q$

$$f(x) + \psi_\delta(y,x) - \delta \le f(y) \le f(x) + \psi_\delta(y,x) + \frac{L}{2}\|y-x\|^2 + \delta. \ ◊$$

Частным случаем, отвечающим условиям

$$f_\delta(x) = f(x),\ \psi_\delta(y,x) = \langle \nabla f(x), y - x \rangle, \tag{3.2}$$

такого определения является условие (2.26). Если не налагать условия $f_\delta(x) = f(x)$ в (3.2), то концепция (3.1), (3.2) совпадает с концепцией $(\delta, L)$-оракула из работы [227], см. также [82, гл. 4], [143, 216, 527]. Близкие концепции модели функции также имеются в работах [418, 473], см. также [139, 245]. Дальнейшее развитие данной концепции отражено в работах [253, 541, 542, 554].

◊ Из дальнейшего будет ясно, что в правой части неравенства (3.1) можно заменить $\|y-x\|^2$ на $2V(y,x)$. При этом правое неравенство в (3.1) интерпретируют уже не как условие гладкости $f(x)$ (липшицевости градиента), а как *условие относительной гладкости* [149, 415]. Важное преимущество, которое приобретается в случае такой замены, – отсутствие условия (2.32) на дивергенцию Брэгмана $V(y,x)$ (правда, и некоторые сложности приобретаются, например, задача (3.21) может стать более сложной в таком случае). Это наблюдение позволяет другим способом, отличным от описанного ранее в пособии, бороться с возможной неограниченностью параметра $L$, определяемого условиями (2.3) или (2.26), в случае неограниченного множества $Q$. Тут можно вспомнить пример $f(x)=x^4$, $Q=\mathbb{R}$ из § 2. С другой стороны, здесь, так же как и ранее в условии (2.3) (см. также (2.26)), достаточно потребовать, чтобы условие (3.1) выполнялось только для всех

$$x,y\in\left\{x\in Q:\ V(x_*,x)\le R^2\right\},$$

где $R^2=V(x_*,x^0)$, см. вторую половину § 2 и (3.16) ниже. Если решение не единственно, то в определении $R^2$ выбирается такое решения $x_*$, которое доставляет минимум $R^2$. ◊

Заметим, что (3.1) включает в себя намного больше свободы (см. [227]) по сравнению с (2.3). В частности, (3.1) включает возможность неточного вычисления (суб-)градиента и значения функции, а не только игру на гладкости (см. (2.4), (2.5)). Мы вернемся к более подробному обсуждению вопросов, связанных с концепцией (3.1), ниже, см. примеры 1, 2 § 3 и упражнения 3.2, 3.3, 4.3.

Подобно (2.29), рассмотрим следующий метод (пояснение записи (3.3) приведено ниже (3.4)):

$$x^{k+1}=\arg_{\tilde{\delta}}\min_{x\in Q}\underbrace{\left\{\psi_\delta\left(x,x^k\right)+\frac{1}{h}V\left(x,x^k\right)\right\}}_{\Psi\left(x,x^k\right)},\qquad(3.3)$$

где $V(x,x^k)$ – дивергенция Брэгмана, определенная в конце предыдущего параграфа. Если задача (3.3) точно решена, то существует такой

$$\nabla_{x^{k+1}}\Psi\left(x^{k+1},x^k\right)\in\partial_x\Psi\left(x,x^k\right)\Big|_{x=x^{k+1}},$$

что для любого $x\in Q$

$$\left\langle \nabla_{x^{k+1}} \Psi\left(x^{k+1}, x^k\right), x - x^{k+1} \right\rangle \geq 0 .$$

Однако мы будем допускать, что задача (3.3) решается лишь в следующем смысле:

$$\left\langle \nabla_{x^{k+1}} \Psi\left(x^{k+1}, x^k\right), x_* - x^{k+1} \right\rangle \geq -\tilde{\delta} ,$$

т. е. (следует сравнить с [143] и [158, п. 5.5.1.2])

$$\left\langle \nabla_{x^{k+1}} \Psi\left(x^{k+1}, x^k\right), x^{k+1} - x_* \right\rangle \leq \tilde{\delta} , \tag{3.4}$$

где $\tilde{\delta} > 0$. Добиться выполнения (3.4) можно по-разному, в зависимости от сложности задачи (3.3) (см. упражнение 3.1). Для возможности перенесения описанного в этом параграфе подхода на следующий параграф, другими словами, для обоснования прямодвойственности метода (3.3), необходимо отказаться от того, что $x = x_*$ (3.4). В этом случае нужно предполагать, что существует такой

$$\nabla_{x^{k+1}} \Psi\left(x^{k+1}, x^k\right) \in \partial_x \Psi\left(x, x^k\right)\Big|_{x=x^{k+1}} ,$$

что

$$\max_{x\in Q} \left\langle \nabla_{x^{k+1}} \Psi\left(x^{k+1}, x^k\right), x^{k+1} - x \right\rangle \leq \tilde{\delta} . \tag{3.5}$$

Введем $\mathrm{Prog}^h_{\psi,V}\left(x^k\right)$, см. также (2.9), (2.11):

$$\mathrm{Prog}^h_{\psi,V}\left(x^k\right) = -\left( \psi_\delta\left(x^{k+1}, x^k\right) + \frac{1}{h} V\left(x^{k+1}, x^k\right) \right). \tag{3.6}$$

Из выпуклости $\psi_\delta\left(x, x^k\right)$ по $x$, определения $x^{k+1}$ (формула (3.3)) и тождества 1 в (2.30), подобно выводу (2.10), из (3.5) (или из (3.4), в этом случае можно сразу положить в последующих выкладках $x = x_*$) получим

$$\begin{aligned} -\tilde{\delta} &\leq \left\langle \nabla_{x^{k+1}} \Psi\left(x^{k+1}, x^k\right), x - x^{k+1} \right\rangle = \left\langle \nabla_{x^{k+1}} \psi_\delta\left(x^{k+1}, x^k\right) + \frac{1}{h} \nabla_{x^{k+1}} V\left(x^{k+1}, x^k\right), x - x^{k+1} \right\rangle = \\ &= \left\langle \nabla_{x^{k+1}} \psi_\delta\left(x^{k+1}, x^k\right), x - x^{k+1} \right\rangle + \frac{1}{h} V\left(x, x^k\right) - \frac{1}{h} V\left(x, x^{k+1}\right) - \frac{1}{h} V\left(x^{k+1}, x^k\right) \leq \\ &\leq \psi_\delta\left(x, x^k\right) - \psi_\delta\left(x^{k+1}, x^k\right) + \frac{1}{h} V\left(x, x^k\right) - \frac{1}{h} V\left(x, x^{k+1}\right) - \frac{1}{h} V\left(x^{k+1}, x^k\right). \end{aligned} \tag{3.7}$$

Отсюда следует, что

$$-\psi_\delta\left(x, x^k\right) \leq \mathrm{Prog}^h_{\psi,V}\left(x^k\right) + \tilde{\delta} + \frac{1}{h} V\left(x, x^k\right) - \frac{1}{h} V\left(x, x^{k+1}\right). \tag{3.8}$$

Согласно неравенству (3.1) при $y = x = x^k$:

$$0 \le f\left(x^k\right) - f_\delta\left(x^k\right) \le \delta . \qquad (3.9)$$

Отсюда по левому неравенству (3.1) при $y = x$, $x = x^k$:

$$f\left(x^k\right) - f\left(x\right) - \delta \le f_\delta\left(x^k\right) - f\left(x\right) \le -\psi_\delta\left(x, x^k\right). \qquad (3.10)$$

При $h \le 1/L$ из (3.6) имеем

$$\begin{gathered}\operatorname{Prog}^h_{\psi,V}\left(x^k\right) = -\left(\psi_\delta\left(x^{k+1}, x^k\right) + \frac{1}{h}V\left(x^{k+1}, x^k\right)\right) \overset{1}{\le} \\ \overset{1}{\le} f_\delta\left(x^k\right) - f\left(x^{k+1}\right) + \delta \overset{2}{\le} f\left(x^k\right) - f\left(x^{k+1}\right) + \delta.\end{gathered} \qquad (3.11)$$

Неравенство 1 следует из (2.32) и правого неравенства (3.1) при $x = x^k$, $y = x^{k+1}$, неравенство 2 следует из (3.9). Подставляя неравенство (3.11) в (3.8), получим

$$-\psi_\delta\left(x, x^k\right) \le f\left(x^k\right) - f\left(x^{k+1}\right) + \tilde{\delta} + \delta + \frac{1}{h}V\left(x, x^k\right) - \frac{1}{h}V\left(x, x^{k+1}\right). \qquad (3.12)$$

Подставляя (3.10) в (3.12), получим аналог неравенства (2.15):

$$f\left(x^k\right) - f\left(x\right) \le f\left(x^k\right) - f\left(x^{k+1}\right) + \tilde{\delta} + 2\delta + \frac{1}{h}V\left(x, x^k\right) - \frac{1}{h}V\left(x, x^{k+1}\right),$$

т. е. при $h \le 1/L$ имеет место основное неравенство

$$f\left(x^{k+1}\right) - f\left(x\right) \le \frac{1}{h}V\left(x, x^k\right) - \frac{1}{h}V\left(x, x^{k+1}\right) + \tilde{\delta} + 2\delta . \qquad (3.13)$$

Мы остановимся на этой формуле, поскольку все дальнейшие рассуждения в точности совпадают с аналогичными рассуждениями из предыдущего параграфа. Общий вывод, который можно сделать из (3.13), сформулируем следующим образом.

**Теорема 3.1.** *Пусть нужно решить задачу* (2.1). *Для метода* (3.3), (2.21):

$$x^{k+1} = \arg_{\tilde{\delta}} \min_{x \in Q}\left\{\psi_\delta\left(x, x^k\right) + LV\left(x, x^k\right)\right\}, \qquad (3.14)$$

*в условиях* (3.1), (3.4) *имеют место оценки, аналогичные оценкам*[31] (2.22), (2.19):

[31] С оговоркой, аналогичной (2.18), в случае оценки (3.16).

$$f\left(\bar{x}^N\right) - f\left(x_*\right) \le \frac{LR^2}{N} + \tilde{\delta} + 2\delta, \tag{3.15}$$

*где*

$$\bar{x}^N = \frac{1}{N}\sum_{k=1}^{N} x^k ,$$

$$V\left(x_*, x^k\right) \le V\left(x_*, x^0\right), \tag{3.16}$$

$R^2 = V\left(x_*, x^0\right)$. *Если решение* $x_*$ *не единственно, то оценки* (3.15), (3.16) *будут верны для того решения* $x_*$, *которое доставляет минимум* $R^2$.

Рассмотрим пару примеров задач *структурной оптимизации* [71], демонстрирующих полезность рассмотрения более общих ситуаций, чем (3.2).

**Пример 3.1 (композитная оптимизация).** Рассмотрим задачу *композитной оптимизации* (composite optimization) [156, 446]:

$$f\left(x\right) = F\left(x\right) + g\left(x\right) \to \min_{x \in Q} \tag{3.17}$$

с выпуклой функцией $F\left(x\right)$, удовлетворяющей условию (2.3), и, вообще говоря, негладкой выпуклой функцией $g\left(x\right)$ простой структуры. Последнее означает, что множества Лебега

$$\Lambda_y = \left\{x \in Q:\ g\left(x\right) < y\right\} \tag{3.18}$$

имеют простую структуру. К такой задаче, например, можно отнести задачу LASSO:

$$\frac{1}{2}\left\|Ax - b\right\|_2^2 + \lambda\left\|x\right\|_1 \to \min_{x \in \mathbb{R}^n}.$$

Естественным обобщением метода (2.29) для задачи (3.17) будет

$$x^{k+1} = \arg\min_{x \in Q}\left\{\left\langle \nabla F\left(x^k\right), x - x^k\right\rangle + g\left(x\right) + LV\left(x, x^k\right)\right\}. \tag{3.19}$$

Метод (3.19) в точности соответствует методу (3.14) с

$$\psi_\delta\left(y, x\right) = \left\langle \nabla F\left(x\right), y - x\right\rangle + g\left(y\right) - g\left(x\right). \tag{3.20}$$

Таким образом, все приведенные выше результаты удается полностью перенести на задачи композитной оптимизации (3.17). В частности, имеет место оценка скорости сходимости (2.22). Стоит особо подчеркнуть, что в полученную оценку скорости сходимости никак не вошла информация о композите $g\left(x\right)$. Это может показаться странным, однако все становится на места, если заметить, что по принципу множителей Лагранжа [62, п. 2.1], использованном в «обратном направлении», при весьма общих условиях существует такое $y$, что задача (3.17) эквивалентна задаче

$$F(x) \to \min_{x \in \Lambda_y}$$

с множеством $\Lambda_y$ (см. (3.18)) простой структуры.[32]

Беря в качестве композитного члена индикаторные функции выпуклых множеств простой структуры, можно получить результаты § 2 из композитного подхода с $Q = \mathbb{R}^n$.

Беря в качестве композитного члена линейные функции, несложно понять, что скорость сходимости метода (3.19) в негладком случае (см. (2.3) – (2.5) с $\nu = 0$) зависит от константы $L_0$, а не от константы Липшица оптимизируемого функционала [156, 454]. Это можно понять непосредственно из самой оценки (2.22), т. е. без композитной оптимизации. Однако с композитной оптимизацией это свойство становится более ясным. ■

В связи с примером 3.1 можно заметить, что если для обычной (некомпозитной) задачи, вообще говоря, негладкой выпуклой оптимизации (2.1)

$$f(x) \to \min_{x \in Q}$$

взять в описанном выше подходе (для простоты считаем $\delta = 0$)

$$\psi_\delta(y, x) = f(y) - f(x)$$

и выбрать произвольное $L$ в условии (3.1), то полученный по формуле (3.3) метод

$$x^{k+1} = \arg\min_{x \in Q} \left\{ f(x) + LV\left(x, x^k\right) \right\} \tag{3.21}$$

становится известным *прокс-методом* решения задачи (2.1) [82, § 1, гл. 6], [150, 174, 200, 233, 339, 358, 451], см. также замечание 3.2 в случае евклидовой прокс-структуры. Метод (3.21) будет сходиться согласно оценке (3.15) из теоремы 3.1, т. е. быстрее, чем следует исходя из нижней оценки, см. упражнение 2.1. Проблема, однако, в том, что в оценку (3.15) входит $\tilde{\delta}$ – «точность» решения вспомогательной задачи. Согласно упражнени-ям 2.3, 3.1 сложность решения вспомогательной задачи, которую можно понимать уже как задачу композитной оптимизации с $L$-сильно выпуклым композитом[33] $LV\left(x,x^k\right)$, будет не меньше, чем $\tilde{O}\left(L_0^2/\left(L\tilde{\delta}\right)\right)$, где $L_0$ определяется по формуле (2.4).[34] Здесь под сложностью понимается число вычислений $\nabla f\left(x\right)$ и число решений уже стандартных вспомогательных подзадач вида (2.29). Комбинируя оценку (3.15) с оценкой $O\left(L_0^2/\left(L\tilde{\delta}\right)\right)$, выбирая $\tilde{\delta}\sim\varepsilon$, где $\varepsilon$ – желаемая точность (по функции) решения исходной задачи (2.1), получим оценку вида (2.35), что уже соответствует нижней оценке (2.37). Действительно, выбирая $N$ в (3.15) из условия $LR^2/N\sim\varepsilon$, получим для итоговой сложности

$$N\frac{L_0^2}{L\tilde{\delta}}\sim\frac{LR^2}{\varepsilon}\frac{L_0^2}{L\varepsilon}=\frac{L_0^2R^2}{\varepsilon^2}, \tag{3.22}$$

что соответствует оценке (2.35), если приравнять правую часть (3.22) $N$ из (2.35) и выразить $\varepsilon\left(N\right)$.

[32] Есть и другой способ, объясняющий факт отсутствия в оценке скорости сходимости информации о $g(x)$ [23, замечание 6].

◊ На первый взгляд, кажется, что нет никакой выгоды от описанного в предыдущем абзаце подхода. Однако выгода получается в случае, когда приведенную конструкцию используют для задач с более сложной структурой, например, для задач композитной оптимизации (3.17), но уже без предположения о простой структуре негладкой выпуклой функции $g\left(x\right)$, как в примере 3.1. В случае $V\left(x,y\right)=\frac{1}{2}\left\|x-y\right\|_2^2$ можно показать, что исходную задачу (3.17) можно решить с точностью по функции $\varepsilon$ за $O\left(L_{0,g}^2R^2/\varepsilon^2\right)$ обращений к оракулу за субградиентом $g\left(x\right)$, где $L_{0,g}$

---

[33] То, что композит сильно выпуклый нужно только для последующего обоснования градиентного слайдинга на базе прокс-метода. В теореме 3.1, описывающей работу, в том числе, прокс-метода, достаточно только, чтобы $V\left(y,x\right)$ было дивергенцией Брэгмана, т.е. имело место представление

$$V\left(y,x\right)=d\left(y\right)-d\left(x\right)-\left\langle\nabla d\left(x\right),y-x\right\rangle.$$

[34] Строго говоря, в упражнении 2.3 рассматривается не композитная постановка, однако ввиду примера 3.1 несложно перенести результаты данного упражнения и на композитную постановку. Необходимые рассуждения, дословно повторяющие написанное в указании к упражнению 2.3, было решено здесь опустить, детали см., например, в [18, п. 2.3]. Также важно отметить, что задача (3.21) должна решаться с точностью $\tilde{\delta}$ в более сильном смысле (3.4), чем «по функции».

определяется согласно (2.4) с $\nu = 0$ и $f \equiv g$, и $O\left(L_{1,F}R^2/\varepsilon\right)$ обращений к оракулу за градиентом $F(x)$, где $L_{1,F}$ определяется согласно (2.4) с $\nu = 1$ и $f \equiv F$. Удалось как бы «расщепить» задачу (3.17) на две задачи, отвечающие отдельным слагаемым, и организовать процедуру решения исходной задачи таким образом, чтобы сложность этой процедуры соответствовала суммарной сложности решения отдельных подзадач. В случае, когда вычисление градиента $F(x)$ занимает намного больше времени, чем вычисление субградиента $g(x)$, такое расщепление дает очевидные преимущества [353, 385]. Прием, с помощью которого удается достичь описанного результата, называется *градиентный слайдинг* [43, 44, 337, 359, 385, 386]. По-видимому, впервые он был предложен в работе [353]. В последние годы этот прием стал достаточно популярным в связи с большим числом приложений в задачах анализа изображений. За многочисленными обобщениями и приложениями данного приема можно следить, например, по работам Дж. Лана [383]. Отметим, в частности, что схема слайдинга переносится и на функционалы, состоящие из суммы двух гладких композитов (см. [386, chapter 8] и упражнение 3.8), что также находит многочисленные приложения [2].

К сожалению, в общем случае техника слайдинга требует довольно тонких и весьма громоздких рассуждений для своего обоснования. Пожалуй, это единственная известная нам и достаточно широко использующаяся конструкция в современной выпуклой оптимизации, суть которой пока так и не удалось раскрыть (в случае, когда целевая функция представляется в виде суммы гладкого и негладкого слагаемого) с помощью элементарных соображений. Впрочем, популярное изложение более частного результата все же имеется [539]. Для всех остальных основных конструкций в данном пособии предпринимается попытка представить их достаточно простым (естественным) способом. ◊

Заметим также, что в работах [244, 409, 410] на базе *проксимального подхода* (описанного выше, см. (3.21)) был предложен новый общий способ ускорения различных неускоренных методов, получивший название *Каталист* (Catalyst). Немного подробнее об этом будет написано в замечании 3.3 и приложении.

**Пример 3.2 (метод уровней).** На пример 3.1 можно посмотреть и с немного другой точки зрения. Как уже отмечалось в самом начале § 2, типично имеется большой зазор между сложностью выполнения итерации $\tilde{O}(n)$ (сложностью проектирования) и сложностью вычисления градиента $O(n^2)$.[35] Можно заметить, что простота множества $\Lambda_y$ (композитной функции) на самом деле в рассуждениях примера 3.1 никак не использовалась. Она была нужна, чтобы не задумываться о сложности проектирования. Поэтому можно понимать пример 3.1 как способ (аддитивного) перенесения части сложности задачи в итерацию, благо для этого имеется хороший запас. Ведь все равно, чтобы сделать шаг метода, нужно посчитать градиент, поэтому сложность «проектирования» вполне можно «утяжелять», например, за счет отмеченной идеи композитной оптимизации, до сложности расчета градиента. Общая сложность итерации по порядку сохранится, но зато число итераций может существенно уменьшиться. Продолжая движение в намеченном направлении, приведем другой пример задачи «со структурой», которая также позволяет заносить часть сложности задачи в «проектирование», сохраняя общую конструкцию [68], [70, § 4, гл. 7], [72, п. 4.3], [384, п. 4], [438]:

$$f(x) = F\left(f_1(x), \ldots, f_m(x)\right) \to \min_{x \in Q}, \qquad (3.23)$$

где все функции выпуклые, причем функция $F(y)$ еще и неубывающая по каждому из своих аргументов. Также предполагаем, что все функции $f_j(x)$, $j = 1, \ldots, m$ удовлетворяют условию (2.3) с $L = L_j$, $j = 1, \ldots, m$, а функция $F(y)$ удовлетворяет условиям (2.3) – (2.5) с $\nu = 0$ ($L = L_0$) и $\| \; \| = \| \; \|_1$. Подобно (3.2), (3.20), положим

$$\psi_\delta(y,x) = F\left(f_1(x) + \langle \nabla f_1(x), y - x \rangle, \ldots, f_m(x) + \langle \nabla f_m(x), y - x \rangle\right) - f(x). \qquad (3.24)$$

Сделанные предположения позволяют утверждать, что условие (3.1) выполняется при $L = L_0 \sum_{j=1}^{m} L_j$. Получаемый при таком выборе $\psi_\delta(y,x)$ (см. (3.24)) метод (3.3) называют *методом уровней* (level method).

Достаточно популярным частным случаем задачи (3.23) является задача, в которой $F(y) = \max_{j=1,\ldots,m} y_j$ [72, п. 2.3]. К такой задаче с помощью *метода нагруженного функционала* [13, § 19, гл. 5], [82, § 3, гл. 9] сводятся и задачи *условной оптимизации* (задачи с функциональными ограничениями) вида

---

[35] Такой большой зазор ($n \leftrightarrow n^2$) имеет место не всегда. Однако в типичных ситуациях вычисление градиента занимает значительно больше времени, чем последующее проектирование.

$$f_0(x) \to \min_{\substack{f_1(x)\le 0,\dots,f_m(x)\le 0,\\ x\in Q}}. \qquad (3.25)$$

Действительно, задачу (3.25) можно переписать следующим образом: найти такой $t = t_*$ и соответствующий $x(t_*)$, доставляющий решение вспомогательной задачи минимизации в (3.26), что $G(t) > 0$ при $t < t_*$ и $G(t_*) = 0$, где

$$G(t) = \min_{x\in Q} \max\left\{ f_0(x) - t, f_1(x), \dots, f_m(x) \right\}. \qquad (3.26)$$

Очевидно, что $G(t)$ невозрастающая функция. Чуть посложнее показывается, что $G(t)$ – выпуклая функция.

**Замечание 3.1.** Это следует из двух общих фактов выпуклого анализа [176, гл. 3]:

1) пусть $\tilde{F}(x,y)$ – выпуклая функция, как функция $x$, тогда функция

$$f(x) = \max_{y\in\tilde{Q}} \tilde{F}(x,y)$$

также выпуклая. Хорошей иллюстрацией тут является представление $|x| = \max\{-x, x\}$, $x \in \mathbb{R}$.

2) пусть $\bar{F}(x,y)$ – выпуклая функция, как функция $(x,y)$, а $\bar{Q}$ – выпуклое множество, тогда функция

$$f(x) = \min_{y:\,(x,y)\in\bar{Q}} \bar{F}(x,y)$$

также выпуклая. Это следует из того, что пересечение надграфика выпуклой функции с выпуклым цилиндром с основанием $Q$ также будет выпуклым множеством и его проекция вдоль $y$ также будет выпуклым множеством. ■

Из общих результатов о поиске корня скалярного нелинейного уравнения [179, гл. 4], [452, Appendix A1] можно пытаться найти $t_*$ с относительной точностью $\varepsilon$ за $O\left(\ln\left(\varepsilon^{-1}\right)\right)$ вычислений значения $G(t)$. Каждое такое вычисление приводит к необходимости решения задачи вида (3.23). Поскольку задачу (3.23) можно решить в общем случае только приближенно, то и посчитать[36] $x(t)$ можно только приближенно. Это обстоятельство приводит к необходимости более тонкого анализа. Детали см., например, в [72, п. 2.3]. Однако сохраняется общий вывод о возрастании сложности решения задачи (3.25) по сравнению с (3.23) в $O\left(\ln\left(\varepsilon^{-1}\right)\right)$ раз при рассматриваемом подходе.

Заметим, что в случае, когда задача (3.25) – негладкая, существуют и другие эффективные численные способы ее решения, базирующиеся на *методе зеркального спуска* (см. упражнение 2.6) и замечании 4.3 и на методе зеркального спуска с переключениями, см., например, [152]. ■

Все последующие рассуждения могут проводиться в общности, выбранной в данном параграфе. Однако в методических целях далее мы намеренно не будем «гнаться за общностью» и стараться формулировать результаты таким образом, чтобы подчеркнуть в первую очередь обсуждаемую идею.

**Упражнение 3.1.** Пусть для задачи выпуклой оптимизации

$$f(x) \to \min_{x \in Q}$$

найдено $\varepsilon$-приближенное по функции решение $x_\varepsilon \in Q$, т. е.

$$f(x_\varepsilon) - f(x_*) \le \varepsilon .$$

**1)** Пусть функция $f(x)$ удовлетворяет условию (2.27) при $\nu = 1$ с константой $L_1$ и $\nabla f(x_*) = 0$. Покажите, что тогда для всех $x \in Q$ имеет место следующая оценка:

$$\left\langle \nabla f(x_\varepsilon), x_\varepsilon - x \right\rangle \le \|x_\varepsilon - x\| \sqrt{2L_1\varepsilon} .$$

Пусть $R = \max_{x,y \in Q} \|y - x\| < \infty$. Покажите, что тогда имеет место следующая оценка:

$$\max_{x \in Q} \left\langle \nabla f(x_\varepsilon), x_\varepsilon - x \right\rangle \le R\sqrt{2L_1\varepsilon} .$$

**2)** Пусть функция $f(x)$ удовлетворяет на отрезке, соединяющим точки $x_\varepsilon$ и $x_*$, условию (2.27) при $\nu = 1$ с константой $L_1$, является $\mu$-сильно выпуклой в норме $\| \ \|$. Покажите, что тогда для всех $x \in Q$ имеет место следующая оценка:

---

[36] Здесь $x(t)$ – решение задачи вспомогательной задачи минимизации по $x \in Q$ в (3.26).

$$\langle \nabla f(x_\varepsilon), x_\varepsilon - x\rangle \le \left(L_1\|x_\varepsilon - x\| + \|\nabla f(x)\|_*\right)\sqrt{2\varepsilon/\mu},$$

и, как следствие,

$$\max_{x\in Q}\langle \nabla f(x_\varepsilon), x_\varepsilon - x\rangle \le \left(L_1 R + \|\nabla f(x_*)\|_*\right)\sqrt{2\varepsilon/\mu}.$$

**Упражнение 3.2.** Пусть $f(x) = \min_{y\in\tilde{Q}} \bar{F}(y,x)$, где $\tilde{Q}$ – ограниченное выпуклое множество, а $\bar{F}(y,x)$ – такая достаточно гладкая, выпуклая по совокупности переменных функция, что при $y, y' \in \tilde{Q}$, $x, x' \in \mathbb{R}^n$:

$$\|\nabla\bar{F}(y',x') - \nabla\bar{F}(y,x)\|_2 \le L\|(y',x') - (y,x)\|_2.$$

Пусть для произвольного $x$ можно найти такой $\tilde{y}_\delta(x) \in \tilde{Q}$, что (следует сравнить с (3.4)):

$$\max_{y\in\tilde{Q}}\langle \nabla_y \bar{F}(\tilde{y}_\delta(x), x), \tilde{y}_\delta(x) - y\rangle \le \delta.$$

Покажите, что для любых $x, x' \in \mathbb{R}^n$

$$\bar{F}(\tilde{y}_\delta(x), x) - f(x) \le \delta,\ \|\nabla f(x') - \nabla f(x)\|_2 \le L\|x' - x\|_2$$

и

$$\left(\bar{F}(\tilde{y}_\delta(x), x) - 2\delta; \langle \nabla_y \bar{F}(\tilde{y}_\delta(x), x), y - x\rangle\right)$$

будет $(6\delta, 2L)$-моделью для функции $f(x)$ в точке $x$ относительно 2-нормы.

**Указание.** См. [23]. Интересно сопоставить это упражнение с леммой 13 из п. 5 § 1, гл. 5 [82]. ■

**Упражнение 3.3 (прокс-метод с неточным решением задачи минимизации на итерации).** Рассмотрим функцию (следует сравнить с задачей из (3.21)):

$$f_L(x) = \min_{y\in Q}\underbrace{\left\{f(y) + \frac{L}{2}\|y - x\|_2^2\right\}}_{\Psi(y,x)}.$$

Предположим, что $f(y)$ – выпуклая функция и

$$\max_{y\in Q}\left\{\Psi(y(x), x) - \Psi(y, x) + \frac{L}{2}\|y - y(x)\|_2^2\right\} \le \delta.$$

Покажите, что тогда

$$\left( f\left(y(x)\right)+\frac{L}{2}\left\|y(x)-x\right\|_2^2-\delta;\left\langle L\cdot\left(x-y(x)\right),y-x\right\rangle \right)$$

будет $(\delta, L)$-моделью функции $f_L(x)$ в точке $x$ относительно 2-нормы.

**Указание.** См. [227].

**Замечание 3.2 (прокс-метод и *сглаживание по Моро–Иосиде* [403]).** Введем функции

$$F_{L,x}(y)=f(y)+\frac{L}{2}\|y-x\|_2^2,$$

$$f_L(x)=\min_{y\in Q}F_{L,x}(y)=F_{L,x}\left(y_L(x)\right).$$

Для любого $L\geq 0$ имеет место неравенство

$$f_L(x)\leq f(x),$$

причем выпуклая функция $f_L(x)$ будет иметь $L$-липшицев градиент. Кроме того, согласно [82, теорема 5 п. 2 § 1, гл. 6],

$$x_*\in \operatorname{Arg}\min_x f_L(x)\left(=\operatorname{Arg}\min_{x\in Q} f_L(x)\right)\Rightarrow x_*\in\operatorname{Arg}\min_{x\in Q} f(x),\ \ f_L(x_*)=f(x_*).$$

Таким образом, вместо исходной задачи можно рассматривать (сглаженную по Моро–Иосиде) задачу

$$f_L(x)\to\min_{x\in\mathbb{R}^n}\left(\min_{x\in Q}\right).$$

На эту задачу можно смотреть как на обычную задачу гладкой выпуклой оптимизации. Согласно упражнению 1.3 сложность решения этой задачи с точностью $\varepsilon$ по функции (число вычислений градиента $\nabla f_L(x)$, т.е. число раз, которое необходимо решить вспомогательную задачу) быстрым градиентным методом можно оценить следующим образом

$$\mathrm{O}\left(\sqrt{\frac{LR^2}{\varepsilon}}\right).$$

Чем меньше выбирается параметр $L$, тем эта оценка будет лучше, но при этом тем сложнее на каждой итерации решать вспомогательную задачу. Заметим, что

$$\nabla f_L(x) = L \cdot (x - y_L(x)).$$

Поэтому обычный градиентный метод будет иметь вид

$$x^{k+1} = x^k - \frac{1}{L} L \cdot \left(x^k - y_L\left(x^k\right)\right) = y_L\left(x^k\right).$$

Однако, согласно упражнению 3.3 внешнюю задачу можно решать и быстрым градиентным методом, работающим с концепцией $(\delta, L)$-модели функции. Например, методом подобных треугольников из упражнения 3.7, см. также замечание 3.3. ■

**Упражнение 3.4 (градиентное отображение).** Изложенный в этом параграфе подход является далеко не единственным способом получения части описанных в § 3 результатов. Удобным инструментом также является использование *градиентного отображении*, см., например, [71], [72, п. 2.3]. С помощью градиентного отображения обобщается (путем замены градиента на градиентное отображение) основной набор базовых формул, из которых выводятся все последующие оценки,[37] см., например, [72, п. 2.2.3, 2.3.2]. Попробуйте получить собранные в § 3 результаты с помощью градиентного отображения.

**Упражнение 3.5 (модель для невыпуклой функции).** Предложите обобщение концепции модели функции (3.1), пригодное для работы с невыпуклыми функциями.

---

[37] В связи с этим тезисом отметим, что сложность задачи оптимизации гладкой/негладкой выпуклой/[сильно выпуклой] функции для рассматриваемого класса численных методов (1.33) равносильна сложности задачи оптимизации функции, удовлетворяющей лишь определенному (явно выписываемому и конечному!) набору условий, связывающих значения функции и её (суб-)градиента в генерируемых методом (1.33) точках [552]. Это наблюдение позволяет получать точные минимаксные оценки скорости сходимости различных итерационных процедур вида (1.33) [222, 239–242, 371, 373, 550–553], см. также замечание 1.5.

**Указание.** См. [26, 169, 255, 418, 473, 541]. ■

**Упражнение 3.6.** В работе [457] в связи с изучением процессов, происходящих в ходе избирательных компаний, и в связи с изучение быстрых способов кластеризации многомерных данных предлагается искать решение следующей задачи выпуклой оптимизации:

$$f_\mu\left(x=(z,p)\right)=g\left(z,p\right)+\mu\sum_{k=1}^{n} z_k\ln z_k+\frac{\mu}{2}\|p\|_2^2\to\min_{z\in S_n(1),\ p\in\mathbb{R}_+^m}.$$

$$x$$

Введем норму $\|x\|^2=\|(z,p)\|^2=\|z\|_1^2+\|p\|_2^2$. Убедитесь, что $\|\ \|$ действительно норма. Предположим, что

$$\left\|\nabla g\left(x_2\right)-\nabla g\left(x_1\right)\right\|_*\le L\left\|x_2-x_1\right\|,$$

где $L\le\mu$. Покажите, что если

$$\psi_\delta\left(y=\left(y^z,y^p\right),x=\left(x^z,x^p\right)\right)=\left\langle\nabla g\left(x\right),y-x\right\rangle+\mu\sum_{k=1}^{n}y_k^z\ln y_k^z+\frac{\mu}{2}\left\|y^p\right\|_2^2-$$

$$-\left(\mu\sum_{k=1}^{n}x_k^z\ln x_k^z+\frac{\mu}{2}\left\|x^p\right\|_2^2\right)-\left(L\sum_{k=1}^{n}y_k^z\ln\left(y_k^z/x_k^z\right)+\frac{L}{2}\left\|y^p-x^p\right\|_2^2\right),$$

то при $\mu\ge L$, имеем: $\psi_\delta\left(y,x\right)$ – выпуклая по $y$ функция, $\psi_\delta\left(x,x\right)=0$ и для любых $y,x\in S_n\left(1\right)\otimes\mathbb{R}_+^m$

$$f_\mu\left(x\right)+\psi_\delta\left(y,x\right)\le f_\mu\left(y\right)\le$$

$$\le f_\mu\left(x\right)+\psi_\delta\left(y,x\right)+2L\sum_{k=1}^{n}y_k^z\ln\left(y_k^z/x_k^z\right)+\frac{2L}{2}\left\|y^p-x^p\right\|_2^2.$$

Заметим, что выпуклость или простота функции $g\left(x\right)$ здесь не требуется! Используя концепцию модели функции в условиях относительной гладкости предложите численный способ решения исходной задачи.

**Указание.** Идея такой модели была заимствована из работ [119, 122, 540]. ■

**Упражнение 3.7 (*метод подобных треугольников* [30, 95, 110, 452]).** Для задачи (2.1) рассмотрите следующий вариант быстрого (ускоренного) градиентного спуска с одной проекцией, работающий с моделью функции (3.1):

$$y^0=z^0,\ A_0=\alpha_0=0,$$

$$\alpha_{k+1} = \frac{1}{2L} + \sqrt{\frac{1}{4L^2} + \alpha_k^2},$$

$$A_{k+1} = A_k + \alpha_{k+1},$$

$$x^{k+1} = \frac{\alpha_{k+1} z^k + A_k y^k}{A_{k+1}},$$

$$z^{k+1} = \arg_{\tilde{\delta}} \min_{x \in Q} \left\{ \alpha_{k+1} \psi_\delta \left( x, x^{k+1} \right) + V\left( x, z^k \right) \right\},$$

$$y^{k+1} = \frac{\alpha_{k+1} z^{k+1} + A_k y^k}{A_{k+1}}.$$

Покажите, что

$$f\left(y^N\right) - f\left(x_*\right) = O\left( \frac{LR^2}{N^2} + \frac{L\tilde{\delta}}{N} + N\delta \right).$$

Полезно сравнить эту формулу с (3.15).

Покажите, что оценка скорости сходимости описанного метода не ухудшится, если на каждой итерации делать в конце дополнительное присваивание: в качестве $x^{k+1}$ выбирать ту точку среди $\left\{ y^{k+1}, u^{k+1}, x^{k+1} \right\}$, которая доставляет наименьшее значение целевой (минимизиремой) функции. Для задач безусловной оптимизации с простейшей моделью функции (3.2) покажите, что получившийся в результате быстрый градиентный метод будет *релаксационным*, т.е. на генерируемой таким метод последовательности точек $\left\{ x^k \right\}_k$ целевая функция будет монотонно убывать, см. также [157].

По аналогии с (2.19) и (3.16) покажите, что (похожую оценку, см., например, в [30], [452, Remark 6.1.1])

$$\max\left\{ V\left( x_*, x^k \right), V\left( x_*, y^k \right), V\left( x_*, z^k \right) \right\} \le V\left( x_*, x^0 \right).$$

**Замечание 3.3 (каталист и оптимальные тензорные методы [21, 43, 233, 339, 410, 426]).** Упражнение 3.7 позволяет строить ускоренный метод с моделью функции

$$\psi_\delta\left(y, x\right) = f\left(y\right) - f\left(x\right).$$

Однако получающиеся в результате вспомогательные задачи за счет роста $\alpha_k \sim k$ с ростом номера итерации будут все хуже и хуже обусловленны-ми. Более эффективным представляется способ, базирующийся на упражнении 3.3 и замечании 3.2, см. также [233].[38] В таком случае вспомогательные задачи будут намного проще: их обусловленность не меняется с ростом номера итерации. Платой за это является: 1) евклидова прокс-структура (впрочем, см. [233]), 2) необходимость достаточно точно решать вспомогательные задачи и 3) итоговый критерий качества работы метода $f_L\left(x^N\right)-f_L\left(x_*\right)$ вместо желаемого $f\left(x^N\right)-f\left(x_*\right)$. Заметим, что

$$f_L\left(x^N\right)-f_L\left(x_*\right)=f_L\left(x^N\right)-f\left(x_*\right)\le f\left(x^N\right)-f\left(x_*\right).$$

Если функция $f\left(x\right)$ имеет $L_f$-липшицев градиент, то проблема 2 отсутствует, поскольку вспомогательные задачи решаются в нужном смысле за линейное время. Если дополнительно известно, что $f\left(x\right)$ еще и $\mu_f$-сильно выпуклая функция, то $f_L\left(x\right)$ также будет сильно выпуклой функцией с константой [403]

$$\tilde{\mu}_f=\mu_f\frac{L}{\mu_f+L}\simeq\mu_f \quad (\mu_f\ll L),$$

поэтому исчезает и проблема 3. В действительности, проблему 3 можно решить и без предположения о сильной выпуклости за счет выбора специального варианта ускоренного внешнего метода [410]. Более того, проблемы 2, 3 просто и не возникают, если «правильно» выбрать ускоренный внешний метод. Ниже в обозначениях замечания 3.2 при $Q=\mathbb{R}^n$ приводится вариант «правильно» ускоренного градиентного метода (следует сравнить с методом линейного каплинга (МЛК) из указания к упражнению 1.3 и замечания 1.6 в части выбора $y^{k+1}$):

**Инициализация (метод Монтейро–Свайтера)**

*Задаем* $z^0$, $y^0$, $A_0=0$;

[38] Далее в этом замечании 3.3 используются обозначения, введенные в замечании 3.2.

**Основной цикл**

*Подбираем* $L_{k+1}$ *и* $y^{k+1}$ *так, что*

$$a_{k+1} = \frac{1/L_{k+1} + \sqrt{1/L_{k+1}^2 + 4A_k/L_{k+1}}}{2}, \quad A_{k+1} = A_k + a_{k+1},$$

$$x^{k+1} = \frac{A_k}{A_{k+1}} y^k + \frac{a_{k+1}}{A_{k+1}} z^k,$$

$$\left\|\nabla F_{L_{k+1},x^{k+1}}\left(y^{k+1}\right)\right\|_2 \le \frac{L_{k+1}}{2}\left\|y^{k+1} - x^{k+1}\right\|_2 \quad \text{// условие Монтейро–Свайтера}$$

$$z^{k+1} = z^k - a_{k+1}\nabla f\left(y^{k+1}\right).$$

Для последовательности $\left\{x^k, y^k, z^k\right\}_{k=1}^N$, генерируемой методом Монтейро–Свайтера, справедливы следующие неравенства [426]

$$\frac{1}{2}\left\|z^N - x_*\right\|_2^2 + A_N \cdot \left(f\left(y^N\right) - f\left(x_*\right)\right) + \frac{1}{4}\sum_{k=1}^N A_k L_k \left\|y^k - x^k\right\|_2^2 \le \frac{1}{2}\left\|x^0 - x_*\right\|_2^2 = \frac{R^2}{2},$$

$$f\left(y^N\right) - f\left(x_*\right) \le \frac{R^2}{2A_N}, \quad \left\|z^N - x_*\right\|_2 \le R,$$

$$\sum_{k=1}^N A_k L_k \left\|y^k - x^k\right\|_2^2 \le 2R^2.$$

Последние два неравенства являются следствием первого (см. также замечание 1.6). Можно также получить оценку

$$A_N \ge \frac{1}{4}\left(\sum_{k=1}^N \frac{1}{\sqrt{L_k}}\right)^2.$$

Условие Монтейро–Свайтера позволяет вместо точного решения $y_{L_{k+1}}\left(x^{k+1}\right)$ вспомогательной задачи, для которого (напомним, см. замечание 3.2, что $F_{L,x}\left(y\right)=f\left(y\right)+\frac{L}{2}\left\|y-x\right\|_2^2$)

$$\left\|\nabla F_{L_{k+1},x^{k+1}}\left(y_{L_{k+1}}\left(x^{k+1}\right)\right)\right\|_2=0,$$

искать только такое решение $y^{k+1}$, что (см. также [43, 147])

$$\left\|y^{k+1}-y_{L_{k+1}}\left(x^{k+1}\right)\right\|_2\le\frac{L_{k+1}}{3L_{k+1}+L_f}\left\|x^{k+1}-y_{L_{k+1}}\left(x^{k+1}\right)\right\|_2,$$

следовательно [339],

$$\left\|\nabla F_{L_{k+1},x^{k+1}}\left(y^{k+1}\right)\right\|_2\le\frac{L_{k+1}}{2}\left\|y^{k+1}-x^{k+1}\right\|_2,$$

что влечет неравенство (неулучшаемое в общем случае с точностью до числового множителя)

$$\left\|\nabla f_{L_{k+1}}\left(x^{k+1}\right)-L_{k+1}\cdot\left(x^{k+1}-y^{k+1}\right)\right\|_2\le\left\|\nabla f_{L_{k+1}}\left(x^{k+1}\right)\right\|_2,$$

где

$$\nabla f_{L_{k+1}}\left(x^{k+1}\right)=L_{k+1}\cdot\left(x^{k+1}-y_{L_{k+1}}\left(x^{k+1}\right)\right).$$

Последнее неравенство можно понимать так, что для задачи

$$f_{L_{k+1}}\left(x\right)\to\min_{x\in\mathbb{R}^n}$$

доступен зашумленный градиент $L_{k+1}\cdot\left(x^{k+1}-y^{k+1}\right)$, с относительными детерминированными помехами. В работе [82, п. 3 § 2, гл. 4] отмечается, что такие помехи (с точностью до числового множителя, меньшего единицы, в правой части неравенства) не меняют по порядку картину сходимости обычного градиентного спуска, см. также текст перед замечанием 1.1. Впрочем, в данном подходе используется ускоренный градиент-ный метод, поэтому приведенные рассуждения не следует воспринимать, как доказательство.

В описанном подходе параметр $L$ можно выбирать по-разному на разных итерациях. Ограничимся сначала случаем, когда $L_k \equiv L$, и попробуем подобрать значение этого параметра $L$. Для этого предположим, что у нас есть некоторый метод, позволяющий решать задачи выпуклой оптимизации с целевым функционалом $f(x)$, обладающим $L_f$-липшицевым градиентом и являющимся $\mu_f$-сильно выпуклым, со сложностью $\tilde{O}\left(\Xi\left(L_f/\mu_f\right)\right)$ – число вычислений $\nabla f(x)$. Тогда «стоимость» каждой итерации метода Монтейро–Свайтера будет

$$\tilde{O}\left(\Xi\left(\frac{L_f + L}{\mu_f + L}\right)\right),$$

а число итераций $\tilde{O}\left(\sqrt{L/\tilde{\mu}_f}\right) = \tilde{O}\left(\sqrt{L/\mu_f}\right)$. Таким образом, разумно подбирать значение параметра $L$ из условия

$$\Xi\left(\left(L_f + L\right)/\left(\mu_f + L\right)\right)\sqrt{L/\mu_f} \to \min_{\mu_f \le L \le L_f}.$$

В частности, если $\Xi\left(L_f/\mu_f\right) = O\left(L_f/\mu_f\right)$, то следует выбрать $L \simeq L_f$. Тогда

$$\tilde{O}\left(\Xi\left(\left(L_f + L\right)/\left(\mu_f + L\right)\right)\sqrt{L/\mu_f}\right) = \tilde{O}\left(\sqrt{L_f/\mu_f}\right),$$

что соответствует ускоренному градиентному методу, см. § 1. Таким образом, с помощью указанной выше конструкции на базе неускоренного градиентного спуска можно построить ускоренный. Описанный общий способ ускорения неускоренных методов первого и нулевого порядка, т.е. использующих производные оптимизируемой функции и ее значения, получил название *каталист* (*catalyst*) [409, 410]. Конструкция каталист переносится и на вариационные неравенства, см. замечание 5.1 и [189, item 3], [483]. Имеется естественное обобщение данной конструкции на задачи стохастической оптимизации [379]. Содержательные примеры ускорения неускоренных рандомизированных методов каталистом будут приведены в упражнении 3.8 и приложении.

Поскольку внешний метод в конструкции Монтейро–Свайтера первого порядка, то может показаться, что описанная выше конструкция не в состоянии ускорять методы более высокого порядка. Так оно и есть, если считать параметр $L$ фиксированным на итерациях. Однако, чем меньше выбирается параметр $L$, тем оценка скорости сходимости внешнего метода

$$O\left(\left(\frac{LR^2}{\varepsilon}\right)^{1/2}\right)$$

будет лучше, но при этом тем сложнее на каждой итерации решать вспомогательную задачу, чтобы посчитать градиент $\nabla f_L(x)$ с нужной точностью. Идея, позволяющая использовать подход Монтейро–Свайтера для ускорения методов высокого порядка, состоит в следующем:

1) Вместо задачи $F_{L,x^{k+1}}(y)\to\min\limits_{y\in\mathbb{R}^n}$ с фиксированным $L$ рассмотреть параметрическое семейство задач $F_{L_{k+1},x^{k+1}}(y)\to\min\limits_{y\in\mathbb{R}^n}$ со специальным образом убывающей (на внешних итерациях) последовательностью $\{L_{k+1}\}_k$. Все эти задачи оптимизации имеют одинаковое решение $x_*$, которое необходимо найти.

2) При этом считать точно $\nabla f_{L_{k+1}}\left(x^{k+1}\right)$ нет возможности, поэтому для решения вспомогательной задачи $F_{L_{k+1},x^{k+1}}(y)\to\min\limits_{y\in\mathbb{R}^n}$ используется неускоренный (тензорный) метод $p$-го порядка: $y^{k+1}=T_{p,pM_p}^{F_{L_{k+1},x^{k+1}}}\left(x^{k+1}\right)$, где

$$T_{p,pM_p}^{F_{L,x}}(x)=\arg\min_{y\in\mathbb{R}^n}\left\{\sum_{r=0}^{p}\frac{1}{r!}\left[\nabla_z^r F_{L,x}(z)\right]_{z=x}\underbrace{\left[y-x,\ldots,y-x\right]}_{r}+\frac{pM_p}{(p+1)!}\left\|y-x\right\|_2^{p+1}\right\}$$

и предполагается, что (см. приложение)

$$\left\|\nabla^p f(y)-\nabla^p f(x)\right\|_2\le M_p\left\|y-x\right\|_2,\ x,y\in\mathbb{R}^n,\ M_p\le\infty .$$

Если для такого метода выбирать $L_{k+1}$ так, чтобы выполнялось условие Монтейро–Свайтера и условие [21, 426]

$$\frac{2(p+1)}{p!}\frac{M_p}{L_{k+1}}\left\|y^{k+1}-x^{k+1}\right\|_2^{p-1}\geq\frac{1}{2},$$

то число внешних итераций (число вычислений оператора $T_{p,pM_p}^{F_{L,x}}(x)$, а, следовательно, и $\left\{\nabla^r f(x)\right\}_{r=1}^p$) будет определяться оценкой

$$\mathrm{O}\left(\left(\frac{M_p R^{p+1}}{\varepsilon}\right)^{2/(3p+1)}\right),$$

неулучшаемой для данного класса задач на классе методов $p$-го порядка, см. приложение.

Из работы [448] следует, что для всех $x\in\mathbb{R}^n$ имеет место следующее неравенство

$$\left\|\nabla F\left(T_{p,pM_p}^{F}(x)\right)\right\|_2\leq\frac{(p+1)M_p}{p!}\left\|T_{p,pM_p}^{F}(x)-x\right\|_2^p$$

Таким образом, $L_{k+1}$ нужно подбирать из условий

$$\frac{1}{2}\leq\frac{2(p+1)}{p!}\frac{M_p}{L_{k+1}}\left\|T_{p,pM_p}^{F_{L_{k+1},x^{k+1}}}\left(x^{k+1}\right)-x^{k+1}\right\|_2^{p-1}\leq 1.$$

Далее заметим, что при $x^{k+1}\neq x_*$ найдется такое, вообще говоря, достаточно маленькое значение $\breve{L}_{k+1}>0$, что

$$\frac{2(p+1)}{p!}\frac{M_p}{\breve{L}_{k+1}}\left\|T_{p,pM_p}^{F_{\breve{L}_{k+1},x^{k+1}}}\left(x^{k+1}\right)-x^{k+1}\right\|_2^{p-1}\geq 1$$

и достаточно большое значение $\bar{L}_{k+1}>0$, что

$$\frac{2(p+1)}{p!}\frac{M_p}{\bar{L}_{k+1}}\left\|T_{p,pM_p}^{F_{\bar{L}_{k+1},x^{k+1}}}\left(x^{k+1}\right)-x^{k+1}\right\|_2^{p-1}\le\frac{1}{2}.$$

Отсюда, ввиду непрерывной зависимости $T_{p,pM_p}^{F_{L_{k+1},x^{k+1}}}\left(x^{k+1}\right)$ от $L_{k+1}$, имеем, что подобрать $L_{k+1}$, можно с помощью процедуры вида: $L_{k+1}=2L_k$, $L_{k+1}:=L_{k+1}/\sqrt{2}$. Конечно, есть риск «проскочить» нужный диапазон. В таком случае можно предусмотреть процедуру «возврата» вида $L_{k+1}:=\sqrt[4]{2}L_{k+1}$ и т.д. В типичных ситуация можно ожидать, что число вызовов оператора $T_{p,pM_p}^{F_{L,x}}(x)$ на одной итерации внешнего метода будет $\mathrm{O}(1)$, см. также [184, 346]. При этом каждый вызов такого оператора порождает свою выпуклую задачу (то что задача получается выпуклой – нетривиальный факт, который был обнаружен совсем недавно [448, 450]). Сложность решения такой задачи (т.е. вычисление $T_{p,pM_p}^{F_{L,x}}(x)$) с нужной точностью [359] сопоставима при $p=2,3$ по объему вычислений со сложностью итерации метода Ньютона, т.е. оценивается как $\tilde{\mathrm{O}}\left(n^{2.37}\right)$. Отметим, что при решении возникающей задачи при $p=3$ используется концепция относительной гладкости и относительной сильной выпуклости [448, 450, 542], см. также начало § 3. В оценке $\tilde{\mathrm{O}}\left(n^{2.37}\right)$ не учитывается время расчета $\left[\nabla_z^r F_{L,x}(z)\right]_{z=x}\underbrace{\left[y-x,\ldots,y-x\right]}_{r}$. Приведенное выражение во многих интересных на практике случаях может быть эффективно посчитано с помощью автоматического дифференцирования [448]. Дополнительную информацию о тензорных методах можно найти, например, в приложении. ■

◊ Близкие к описанному выше подходу рассуждения имеются в работе [233]. Различие в том, что в [233] в качестве внешнего ускоренного метода используется метод без вспомогательного одномерного поиска, что в конечном итоге приводит к чуть более плохим теоретическим оценкам оракульной сложности предложенных в [233] тензорных методов (методов высокого порядка). Данная ситуация была исправлена в работе [449], см. также [451]. Отметим, что ускоренный проксимальный метод из работы [233] может работать с неевклидовым проксом вида дивергенции Брэгмана (см. § 2) в предположении, что прокс-функция имеет ограниченную константу Липшица градиента на допустимом множестве. ◊

**Упражнение 3.8 (гладкий / ускоренный слайдинг [386, item 8.2]).** Рассмотрим следующую задачу

$$f(x)+g(x)\to\min_x ,$$

где $f(x)$ и $g(x)$ имеют $L_f$ и $L_g$-Липшицевы градиенты в 2-норме, причем $L_f \le L_g$, а функция $g(x)$ – $\mu$-сильно выпуклая в 2-норме, причем $\mu \le L_f$. Покажите, что для решения рассмотренной задачи с заданной точностью[39] достаточно $\tilde{O}\left(\sqrt{L_f/\mu}\right)$ вычислений $\nabla f(x)$ и $\tilde{O}\left(\sqrt{L_g/\mu}\right)$ вычислений $\nabla g(x)$. Попробуйте получить аналогичный результат в модельной общности.

**Указание.** Применим к рассмотренной задаче технику каталист[40]. Тогда вместо исходной задачи потребуется $\tilde{O}\left(\sqrt{L/\mu}\right)$ раз решать задачу вида

$$f(x)+g(x)+\frac{L}{2}\left\|x-x^k\right\|_2^2\to\min_x .$$

---

[39] Не важно, с какой именно точностью. Эта точность будет входить под логарифмами в приведенные далее оценки, а для наглядности логарифмические сомножители было решено опустить в данном упражнении. Далее в указании к этому упражнению оговорки о точности решения возникающих подзадач также опускаются, поскольку все это влияет только на логарифмические сомножители в итоговых оценках, которые опущены.

[40] Обойтись без этой техники не получается по тем же причинам (см. начало замечания 3.3), по которым из ускоренного метода, описанного в упражнении 3.7, не получается с помощью модельной общности получить ускоренный проксимальный метод с оптимальной оценкой скорости сходимости.

Последнюю задачу можно решать неускоренным композитным градиентным методом (см. пример 3.1), считая $g(x)+\frac{L}{2}\|x-x^k\|_2^2$ композитом. Число итераций такого метода будет совпадать с числом вычислений $\nabla f(x)$ и равно $\tilde{O}(L_f/(L+\mu))$, считаем $\mu \le L \le L_f$. Но в условиях задачи не предполагалась проксимальная дружественность функции $g(x)$, поэтому возникающую на каждой итерации неускоренного композитного градиентного метода задачу вида

$$\left\langle \nabla f(\tilde{x}^l), x-\tilde{x}^l \right\rangle + \frac{L_f}{2}\|x-\tilde{x}^l\|_2^2 + g(x) + \frac{L}{2}\|x-x^k\|_2^2 \to \min_x ,$$

в свою очередь, необходимо будет решать. Для решения данной задачи можно использовать ускоренный композитный градиентный метод для задач сильно выпуклой оптимизации (см. [30, 446], а также упражнение 3.7 и конец § 5), считая $\frac{L_f}{2}\|x-\tilde{x}^l\|_2^2 + \frac{L}{2}\|x-x^k\|_2^2$ композитом. Число итераций такого метода будет $\tilde{O}\left(\sqrt{L_g/(L_f+L+\mu)}\right)$. Таким образом, общее число вычислений $\nabla g(x)$ будет

$$\tilde{O}\left(\sqrt{L/\mu}\right)\cdot\left[\tilde{O}(L_f/(L+\mu))\cdot\tilde{O}\left(\sqrt{L_g/(L_f+L+\mu)}\right)+1\right].$$

Выбирая параметр $L \in [\mu, L_f]$ так, чтобы последнее выражение было минимальным, получим (с учетом сделанных предположений $L_f \le L_g$ и $\mu \le L_f$), что $L \simeq L_f$. Следовательно, общее число вычислений $\nabla g(x)$ будет равно $\tilde{O}\left(\sqrt{L_g/\mu}\right)$. ■

◊ Заметим, что данное упражнение можно обобщить на случай, когда вместо $\nabla g(x)$ доступно только $g(x)$ [337]. Чтобы понять как все это можно сделать, рекомендуется ознакомиться с безградиентными методами в начале приложения и цитированной там литературой. Приведенные в упражнении результаты можно улучшить, если дополнительно известно, что, функция $g(x)$ имеет представление в виде суммы функций [44, 337]. Тогда вместо ускоренного композитного градиентного метода для задач сильно выпуклой оптимизации можно использовать ускоренный композитный метод редукции дисперсии для задач сильно выпуклой оптимизации, см., например, [388] и замечание 1 в приложении. Описанный в упражнении 3.8 слайдинг переносится и на тензорные методы [359].

Отметим также, что в описанной в указании к упражнению 3.8 трехуровневой схеме слайдинга можно убрать промежуточный уровень, связанный с использованием неускоренного композитного градиентного спуска. Для этого следует использовать композитный вариант ускоренной проксимальной оболочки (с гладким композитом), см., например, [43, 358, 359]. ◊

**Упражнение 3.9 (проксимальный метод Синхорна [260, 581]).** В последнее время в различных приложениях часто встречается расстояние Монжа–Канторовича–Васерштейна [488] между двумя вероятностными мерами. Вычисление такого расстояния для дискретных мер сводится к классической транспортной задаче линейного программирования (ЛП)

$$\sum_{i,j=1}^{n} c_{ij}x_{ij} \to \min_{\substack{\sum_{j=1}^{n} x_{ij}=l_i,\ i=1,\dots,n \\ \sum_{i=1}^{n} x_{ij}=w_j,\ j=1,\dots,n \\ x_{ij}\ge 0,\ i,j=1,\dots,n}},$$

где $\sum_{i=1}^{n} l_i = \sum_{j=1}^{n} w_i = 1$. Используя неускоренный прокс-метод (3.21) с $V(x,y)=\sum_{i.j=1}^{n} x_{ij}\ln\left(x_{ij}/y_{ij}\right)$:

$$x^{k+1} = \arg\min_{\substack{\sum_{j=1}^{n} x_{ij}=l_i,\ i=1,\dots,n \\ \sum_{i=1}^{n} x_{ij}=w_j,\ j=1,\dots,n \\ x_{ij}\ge 0,\ i,j=1,\dots,n}} \left\{\sum_{i,j=1}^{n} c_{ij}x_{ij} + L\sum_{i.j=1}^{n} x_{ij}\ln\left(x_{ij}/x_{ij}^{k}\right)\right\},$$

предложите способ решения транспортной задачи.

Предложите адаптивный способ подбора параметра $L$.

**Указание.** Следует учесть, что возникающую на каждой итерации прокс-метода задачу можно приближенно решать путем перехода к двойственной задаче и использования метода альтернативных направлений: двойственную функцию можно явно прооптимизировать по группе (двойственных) множителей Лагранжа, отвечающих ограничениям

$$\sum_{j=1}^{n} x_{ij} = l_i, \ \ i = 1,\ldots,n$$

и при «замороженных» остальных множителях, аналогичное можно проделать и по группе оставшихся множителей, отвечающих ограничениям

$$\sum_{i=1}^{n} x_{ij} = w_j, \quad j = 1,\ldots,n\,.$$

Чередуя такие оптимизации, получим метод Синхорна–(Брэгмана–Шелейховского), который также называют методом балансировки, представляющий собой метод альтернативных направлений[41] для двойственной задачи [155]. Про этот метод (в приложении к данной задаче) известно, что его трудоемкость имеет вид [260, 275, 540]

$$n^2 \min\left\{ \tilde{O}\left(\frac{1}{L\tilde{\varepsilon}}\right), \tilde{O}\left(\exp\left(\frac{C(n)}{L}\right)\right)\right\}.$$

Внешний прокс-метод согласно теореме 3.1 с $R^2 = O\left(\ln n^2\right)$ сойдется за

$$O\left(LR^2/\varepsilon\right) = \tilde{O}\left(L/\varepsilon\right)$$

итераций, где точность решения внутренней задачи $\tilde{\varepsilon}$ должна быть существенно выше точности решения исходной задачи: $\tilde{\varepsilon} \ll \varepsilon$. Таким образом, описанный выше проксимальный метод при оптимальном выборе $L$ будет иметь трудоемкость[42]

---

[41] Метод альтернативных направлений в худшем случае сходится как градиентный спуск (в наилучшей норме, см. замечание 1.3 и работу [559]) с константой Липшица градиента, отвечающей наименьшей из констант Липшица градиента целевой функции по соответствующей группе переменных [155]. Этот результат недавно был перенесен и на (ускоренные) блочно-покомпонентные методы (см., например, приложение), в которых вместо шага типа градиентного спуска в одном из блоков осуществляется явная оптимизация, по соответствующим этому блоку переменным [231]. Также недавно были предложены (прямо-двойственные) ускоренные варианты метода альтернативных направлений с $m$ блоками [231, 307]. В теоретическом плане метод из работы [307] требует в $m$ раз больше итераций, чем обычный ускоренный метод, но на практике работает заметно быстрее последнего.

[42] Заметим, что если рассматриваемую транспортную задачу решать с помощью энтропийной регуляризации, то согласно замечанию 4.1 нужно выбирать $L = \varepsilon/\left(2\ln n^2\right)$ (более тонкий анализ конкретного случая: энтропийной

$$n^2\tilde{O}\left(C(n)/\varepsilon\right),$$

где $C(n) \gg n$. В действительности, на практике описанный метод работает заметно лучше [540, 541]. Отметим в этой связи, что наилучшие с точки зрения теоретических оценок способы решения исходной транспортной задачи [167, 347, 398, 482]:

$$n^2 \cdot \min\left\{\tilde{O}(1/\varepsilon), \tilde{O}\left(\sqrt{n}\right)\right\}$$

и (рассмотренного в этом упражнении) ее $L$-энтропийно регуляризованного варианта [121, 208]:

$$\tilde{O}\left(n^2/L\right),$$

пока далеки от практической эффективности. В частности, солверы, решающие транспортную задачу, как задачу ЛП с помощью методов внутренней точки (см. текст после замечания 4 приложения) имеют практическую сложность $\tilde{O}\left(n^3\right)$ [486, 488], такую же сложность имеет для этой задачи симплкекс метод [277, 548] и венгреский алгоритм для специальной подзадачи «о назначениях» [377], см. также указание к упражнению 1.4.

Способ подбора параметра $L$ можно построить, например, на базе следующей идеи. На первой итерации прокс-метода стартуем с завышенной оценки $L$, решаем задачу, затем полагаем $L := L/2$ и перерешиваем задачу, и так до тех пор, пока не детектируем существенное увеличение (например, в 10 раз) сложности решения вспомогательной задачи энтропийно-линейного программирования по сравнению со стартовой сложностью. Найденное значение $L$ можно использовать и на последующих итерациях проксимального метода. В качестве точки старта новой итерации такого метода можно выбирать решение вспомогательной задачи с предыдущей итерации. ■

**Упражнение 3.10 (имплементируемость тензорных методов третьего порядка; Ю.Е. Нестеров, 2018 [448]).** Рассмотри шаг тензорного метода из замечания 3.3 при $p=3$

$$T_{3,3M_3}^{F_{L,x}}(x) = \arg\min_{y\in\mathbb{R}^n} \underbrace{\left\{ \sum_{r=0}^{3} \frac{1}{r!} \left[\nabla_z^r F_{L,x}(z)\right]_{z=x} \underbrace{\left[y-x,\ldots,y-x\right]}_{r} + \frac{3M_3}{4!}\|y-x\|_2^4 \right\}}_{f(x)} =$$

---

регуляризации транспортной задачи, имеется в работе [572]), что приведет в итоге к оценке трудоемкости метода Синхорна $\tilde{O}\left(n^2/\varepsilon^2\right)$.

$$= \arg\min_{y\in\mathbb{R}^n}\left\{F_{L,x}(x)+\langle\nabla F_{L,x}(x),y-x\rangle+\frac{1}{2}\langle\nabla^2F_{L,x}(x)(y-x),y-x\rangle+\right.$$

$$\left.+\frac{1}{6}\nabla^3F_{L,x}(x)[y-x,y-x,y-x]+\frac{M_3}{8}\|y-x\|_2^4\right\}.$$

Используя $\nabla F_{L,x}(x)$, $\nabla^2F_{L,x}(x)$ предложите способ приближенного[43] поиска $T_{3,3M_3}^{F_{L,x}}(x)$.

**Указание.** Прежде всего, заметим, что вычислять $\nabla^3F_{L,x}(x)$ нет необходимости, потому что в описываемом далее подходе используется только

$$\nabla^3F_{L,x}(x)[y-x,y-x]\approx\frac{\nabla^2F_{L,x}(x+\tau(y-x))(y-x)-\nabla^2F_{L,x}(x)(y-x)}{\tau}.$$

Далее введем условие относительной сильной выпуклости. Будем говорить, что $f(x)$ $\mu$-сильно выпуклая и $L$-гладкая функция относительно дивергенции Брэгмана $V(y,x)=d(y)-d(x)-\langle\nabla d(x),y-x\rangle$ если для всех $x,y\in\mathbb{R}^n$ (можно обобщить данное определение и основные выводы и на случай неточной модели $f(x)$ [542])

$$f(x)+\langle\nabla f(x),y-x\rangle+\mu V(y,x)\le f(y)\le f(x)+\langle\nabla f(x),y-x\rangle+LV(y,x).$$

Достаточным условием для этого будет

$$\mu\nabla^2 d(y)\prec\nabla^2 f(y)\prec L\nabla^2 d(y).$$

Описанный в формуле (3.21) градиентный спуск

$$x^{k+1}=\arg_{\tilde{\delta}}\min_{x\in\mathbb{R}^n}\left\{\langle\nabla f(x^k),x-x^k\rangle+LV(x,x^k)\right\}$$

будет сходиться следующим образом для такой функции согласно формуле [448, 451, 542]

$$\min_{k=0,\dots,N} f(x^k)-f(x_*)\le(L+\mu)V(x_*,x^0)\min\left\{\frac{1}{N},\left(\frac{L-\mu}{L+\mu}\right)^N\right\}+\tilde{\delta}.$$

Для рассматриваемой в условиях упражнения задачи оптимизации (для простоты обозначений будем считать, что ищется $T_{3,3M_3}^{F_{L,x}}(0)$) в [448] была предложена следующая прокс-функция

$$d(y)=\frac{1}{2}\langle\nabla^2 f(0)y,y\rangle+\frac{M_3}{8}\|y\|_2^4.$$

[43] Для получения условий на точность решения см. [235, 358, 359, 449, 450, 451].

Для такой прокс-функции было установлено [451], что

$$\mu = 1 - \frac{1}{\sqrt{2}}, \quad L = 1 + \frac{1}{\sqrt{2}}.$$

Отсюда следует, что решение задачи поиска $T_{3,3M_3}^{F_{L,x}}(0)$ по сложности эквивалентно (с точностью до логарифмического множителя от желаемой точности поиска $T_{3,3M_3}^{F_{L,x}}(0)$) задаче вида

$$\langle b, y\rangle + \frac{1}{2}\langle Ay, y\rangle + \lambda\|y\|_2^4 \to \min_{y\in\mathbb{R}^n},$$

что, в свою очередь, эквивалентно задаче

$$\langle b, y\rangle + \frac{1}{2}\langle Ay, y\rangle \to \min_{\|y\|_2^4 \le C(\lambda)},$$

что эквивалентно

$$\langle b, y\rangle + \frac{1}{2}\langle Ay, y\rangle \to \min_{\|y\|_2^2 \le \sqrt{C(\lambda)}},$$

что в итоге эквивалентно

$$\langle b, y\rangle + \frac{1}{2}\langle Ay, y\rangle + \tilde{\lambda}(\lambda)\|y\|_2^2 \to \min_{y\in\mathbb{R}^n},$$

где поиск $\tilde{\lambda}(\lambda)$ может быть осуществлен за логарифмическое время от желаемой точности. Последняя задача имеет такую же сложность (не больше), как и итерация метода Ньютона (см. приложение). Таким образом, с точностью до квадрата логарифмического по желаемой точности множителя задача поиска $T_{3,3M_3}^{F_{L,x}}(x)$ по сложности эквивалентна выполнению итерации метода Ньютона. Более того, также как и в методе Ньютона о функции используется информация не выше второго порядка $\nabla F_{L,x}(x)$, $\nabla^2 F_{L,x}(x)$. Другими словами, с описанным здесь внутренним методом второго порядка, получается, что тензорный метод третьего порядка в такой реализации становится, по сути, методом второго порядка. Если все описанные здесь действия осуществить более аккуратно (с контролем точности), то всему этому можно придать вполне законченный вид [358, 449, 451]. Полученные в результате такого подхода методы были названы Ю.Е. Нестеровым *супербыстрыми тензорными методами второго порядка* [358, 449, 451]. Об этом также вкратце будет рассказано в приложении. ■

**Упражнение 3.11 (статистический предобуславливатель для задач минимизации суммы выпуклых функций [327]).** Рассмотрим задачу

$$f(x) = \frac{1}{m}\sum_{l=1}^{m} f_l(x) \to \min_{x\in\mathbb{R}^n},$$

где $f_l(x)$ – $\mu$-сильно выпуклые в 2-норме гладкие функции. Будем считать, что $m$ – большое число. Предположим, что есть централизованная архитектура с $r \ll m$ узлами, см. текст после упражнения 4.7. Первый узел – центральный, т.е. связан со всеми остальными. Поместим в первый узел случайно отобранные $\tilde{m}$ ($\tilde{m} \ll m$) слагаемых из суммы. Остальные слагаемые распределим по остальным узлам. Обозначим соответствующую первому узлу нормированную подсумму через $\tilde{f}(x)$. Рассмотрим градиентный спуск (3.21), выполняемый на центральном (первом) узле

$$x^{k+1} = \arg_{\tilde{\delta}}\min_{x\in\mathbb{R}^n}\left\{\left\langle \nabla f\left(x^k\right), x - x^k\right\rangle + LV\left(x, x^k\right)\right\},$$

где $V(y,x) = d(y) - d(x) - \left\langle \nabla d(x), y - x\right\rangle$, $\quad d(x) = \tilde{f}(x) + \frac{\gamma}{2}\|x\|_2^2$, $\quad \gamma > 0$. Каждая итерация такого градиентного спуска отвечает коммуникации центрального узла с остальными, чтобы получить $\nabla f\left(x^k\right)$.

Покажите, что если $\left\|\nabla^2 \tilde{f}(x) - \nabla^2 f(x)\right\|_2 \le \gamma$ (определение матричной 2-нормы, согласованной с векторной 2-нормой, см. в приложении), то

$$\frac{\mu}{\mu + 2\gamma}\nabla^2 d(x) \prec \nabla^2 f(x) \prec \nabla^2 d(x).$$

Из матричного неравенства Хёфдинга с вероятностью $\ge 1 - \sigma$ имеем, что $\gamma \sim \sqrt{\ln(n/\sigma)/\tilde{m}}$ [557]. Используя это наблюдение и указание к упражнению 3.10, оцените скорость сходимости градиентного спуска (3.21).

◊ Заметим, что в работе [327] предлагается и ускоренный адаптивный вариант градиентного метода в условиях относительной сильной выпуклости и гладкости. Как известно, в общем случае, ускорения такой метод может не давать [238]. Однако в данном контексте ускорение удается наблюдать.

Также отметим, что адаптивность ускоренного метода работы [327] хорошо согласуется с общей идей адаптивного подбора параметров, характеризующих гладкость целевой функции, продемонстрированной также в § 5. А именно, из доказательства выделяется условие, следующее из гладкости, и явно зашивается в сам алгоритм. Итерация осуществляется только в том случае, когда не известный параметр подберется так, что это условие выполнится. В отличие от параметра, характеризующего сильную выпуклость, в такие условия для параметра гладкости не входят неизвест-ные величины (кроме самого этого параметра), как, например, само решение задачи $x_*$.

Несмотря на то, что в упражнениях 3.10 и 3.11 применяется один и тот же метод (3.21), стоит отметить, что применяется он существенно по-разному. В упражнении 3.10 это внутренний метод. В упражнении 3.11, наоборот, это внешний метод – оболочка. Тем не менее, для решения вспомогательной задачи, возникающей в упражнении 3.11, можно использовать тензорные методы с шагом, описанным в упражнении 3.10. Тогда получится, что метод (3.21) будет возникать сразу в двух местах. Заметим, что идея применения методов второго порядка (типа Ньютона) к минимизации целевой функции вида суммы может быть мотивирована тем, что вычисление гессиана суммы (если сумма не очень большого числа слагаемых) может существенно не влиять на сложность итерации, определяющуюся стоимостью обращения (регуляризованного) гессиана. Другими словами, метод как бы не замечает до определенного момента (числа слагаемых), что минимизируется сумма. Время работы метода будет практически такое же, как если бы сумма была из одного слагаемого. Собственно, поэтому описанный в упражнении 3.11 прием, с помощью которого можно сводить задачу оптимизации большого числа слагаемых к задаче минимизации меньшего числа слагаемых такого же типа, дает дополнительную мотивацию для изучения тензорных методов, см. [21, 143, 358, 448, 451] и цитированную в этих работах литературу. ◊

# § 4. Прямодвойственная структура градиентного спуска

Как и в § 2, 3 рассмотрим сначала общую задачу выпуклой оптимизации (2.1):

$$f(x) \to \min_{x \in Q}.$$

Под *прямодвойственным методом* решения задачи (2.1) будем понимать такой метод, сходимость которого может быть сформулирована (представлена) в терминах *сертификата точности* (2.25) [443] (по А.С. Немировскому) или, в общем случае, в терминах неравенств типа (4.2) [18, 71, 72, 143, 452, 455] (по Ю.Е. Нестерову).

В данном параграфе будет продемонстрирована прямодвойственная природа градиентного спуска (2.6) [4]. Сначала на примере задачи минимизации выпуклой функции при аффинных ограничениях демонстрируется как градиентный спуск применяется к двойственной задаче (решение прямой задачи удается восстановить за счет прямодвойственности метода), а затем (в конце параграфа) градиентный спуск будет применен к исходной (прямой) задаче при дополнительном предположении, что аффинные ограничения имеют достаточно простую структуру (решение двойственной задачи также удается восстановить за счет прямодвойственности метода). Отмеченные возможности прямодвойственных методов, рассмотренные далее на примере только градиентного спуска, отчасти и объясняют их название [455].

◊ В действительности, при правильном взгляде [443], практически любой численный метод оптимизации с фиксированными шагами (см. указание к упражнению 1.3 и замечание 1.6) является прямодвойственным. Нетривиальный пример – метод эллипсоидов. Как уже отмечалось ранее, к прямодвойственным методам относят методы, в которых имеются оценки на *сертификат точности* (2.25) [443]. Отметим, что в данном параграфе мы явно не используем сертификат точности, поскольку его использование приводит к наличию дополнительного слагаемого в правой части оценки (2.25), от которого на самом деле можно избавиться. Однако стоит отметить, что в идейном плане в § 4 используется, по сути, тот же самый подход, что и в работах [443, 455]. ◊

Пусть сначала планируется решать двойственную задачу. Для двойственных задач множество $Q$ – либо все пространство, либо неотрицательный ортант, либо прямое произведение пространства на неотрицательный ортант. В любом из этих случаев имеет смысл выбирать 2-норму и евклидову прокс-структуру (см. § 2).

Итак, вернемся к формуле (2.12) с $h = 1/L$ (2.21). Перепишем её следующим образом:

$$f\left(x^{k+1}\right) \le \left\{ f\left(x^k\right) + \left\langle \nabla f\left(x^k\right), x - x^k \right\rangle \right\} + \frac{L}{2}\left\|x - x^k\right\|_2^2 - \frac{L}{2}\left\|x - x^{k+1}\right\|_2^2 + \delta . \quad (4.1)$$

Суммируя (4.1) по $k = 0, \dots, N-1$, учитывая выпуклость функции $f(x)$ и произвол в выборе $x \in Q$, получим

$$f\left(\bar{x}^N\right) \le \frac{1}{N} \min_{x \in Q} \left\{ \sum_{k=0}^{N-1} \left[ f\left(x^k\right) + \left\langle \nabla f\left(x^k\right), x - x^k \right\rangle \right] + \frac{L}{2}\left\|x - x^0\right\|_2^2 \right\} + \delta , \quad (4.2)$$

где

$$\bar{x}^N = \frac{1}{N}\sum_{k=1}^{N} x^k .$$

Данная формула является обоснованием *прямодвойственности* метода градиентного спуска (2.6), (2.21) [18, гл. 3], [4, 443, 455]. Как будет продемонстрировано ниже, сходимость метода в таком смысле (более сильном, чем просто по функции) позволяет строить сходящуюся с такой же скоростью последовательность и для сопряженной (двойственной) задачи.

Будем считать, что $\delta \le \varepsilon/2$ (см. (2.18)). Рассмотрим следующий (вычислимый! – ввиду простоты множества $Q$, см. § 2) критерий останова метода:

$$f\left(\bar{x}^N\right) - \frac{1}{N} \min_{x \in B_{R,Q}\left(x^0\right)} \left\{ \sum_{k=0}^{N-1} \left[ f\left(x^k\right) + \left\langle \nabla f\left(x^k\right), x - x^k \right\rangle \right] \right\} \le \varepsilon . \quad (4.3)$$

Обратим внимание, что минимум в (4.3) берется по множеству $B_{R,Q}\left(x^0\right)$, где $R = \left\|x_* - x^0\right\|_2$, а не $B_{R,Q}\left(x_*\right)$, поскольку $x_*$ нам не известно.[44] При

[44] Строго говоря, и $R = \left\|x_* - x^0\right\|_2$ также неизвестен. Однако при использовании $B_{R,Q}\left(x^0\right)$ удается ограничиться лишь одним неизвестным $R$, по которому можно делать рестарты подобно указанию к упражнению 2.3, а в ряде случаев достаточно здесь исходить из размера множества $Q$.

этом ранее в § 2 мы показали, что $x_* \in B_{R,Q}\left(x^0\right)$. Значит по выпуклости $f\left(x\right)$ (см. (2.13)) имеем нужное нам неравенство

$$\begin{aligned} f\left(\bar{x}^N\right) - f\left(x_*\right) &\le f\left(\bar{x}^N\right) - \frac{1}{N}\sum_{k=0}^{N-1}\left[ f\left(x^k\right) + \left\langle \nabla f\left(x^k\right), x_* - x^k \right\rangle \right] \le \\ &\le f\left(\bar{x}^N\right) - \frac{1}{N}\min_{x\in B_{R,Q}\left(x^0\right)}\left\{\sum_{k=0}^{N-1}\left[ f\left(x^k\right) + \left\langle \nabla f\left(x^k\right), x - x^k \right\rangle \right]\right\} \le \varepsilon. \end{aligned} \tag{4.4}$$

С другой стороны, из (4.2) имеем

$$\begin{aligned} f\left(\bar{x}^N\right) &\le \frac{1}{N}\min_{x\in Q}\left\{\sum_{k=0}^{N-1}\left[ f\left(x^k\right) + \left\langle \nabla f\left(x^k\right), x - x^k \right\rangle \right] + \frac{L}{2}\left\|x - x^0\right\|_2^2\right\} + \frac{\varepsilon}{2} \le \\ &\le \frac{1}{N}\min_{x\in B_{R,Q}\left(x^0\right)}\left\{\sum_{k=0}^{N-1}\left[ f\left(x^k\right) + \left\langle \nabla f\left(x^k\right), x - x^k \right\rangle \right]\right\} + \frac{LR^2}{2N} + \frac{\varepsilon}{2}. \end{aligned} \tag{4.5}$$

Значит, с учетом (2.5), (2.19) метод (2.6), (2.21) при условии (2.4) гарантированно остановится по критерию (4.3), сделав не более

$$N = \frac{LR^2}{\varepsilon} \le \left(\frac{L_\nu R^{1+\nu}}{\varepsilon}\right)^{\frac{2}{1+\nu}} \tag{4.6}$$

итераций (вычислений $\nabla f\left(x^k\right)$).

Рассмотрим конкретный пример использования оценки типа (4.2) [4, 226, 260, 555]. Пусть необходимо решить задачу (в данном случае численно решать планируется двойственную задачу к (4.7), поэтому переменные в прямой задаче (4.7) обозначили через $y$):

$$\varphi\left(y\right) \to \min_{Ay=b,\ y\in\tilde{Q}}, \tag{4.7}$$

где функция $\varphi\left(y\right)$ – $\mu$-сильно выпуклая в $p$-норме на $\tilde{Q}$ $\left(1 \le p \le 2\right)$. Решение задачи (4.7) обозначим через $y_*$, а оптимальное значение функционала – через $\varphi_*$ ($\varphi_* = \varphi\left(y_*\right)$).

Построим (с точностью до знака) *двойственную задачу* к задаче (4.7):

$$f\left(x\right) = \max_{y\in\tilde{Q}}\left\{\left\langle x, b - Ay\right\rangle - \varphi\left(y\right)\right\} \to \min_{x\in\mathbb{R}^n}. \tag{4.8}$$

◊ Опишем общий принцип построения двойственных задач [176, гл. 5]. Итак, пусть исходная задача выпуклой оптимизации имеет вид

$$\varphi\left(y\right) \to \min_{h\left(y\right)\le 0,\ Ay=b,\ y\in\tilde{Q}}.$$

Тогда

$$\min_{h(y)\le 0,\ Ay=b,\ y\in\tilde{Q}} \varphi(y) = \min_{y\in\tilde{Q}} \left\{ \varphi(y) + \max_{z\ge 0} \langle z, h(y) \rangle + \max_{x} \langle x, Ay-b \rangle \right\} \stackrel{?}{=}$$

$$\stackrel{?}{=} \max_{z\ge 0,\ x} \min_{y\in\tilde{Q}} \left\{ \varphi(y) + \langle z, h(y) \rangle + \langle x, Ay-b \rangle \right\}.$$

Равенство со знаком вопроса обосновывается с помощью теорем типа *фон Неймана* или *Сиона–Какутани* [158, приложение D.4]. К сожалению, при таком подходе требуется компактность множества $Q$ или возможность компактифицировать двойственные переменные $(z,x)$ (см. [452, item 3.1.8], а также замечание 4.2 и упражнение 4.1). В любом случае в реальных задачах, как правило, удается обосновать это равенство [176], которое также называют *сильной двойственностью* [176, гл. 5]. Таким образом, решение исходной задачи сводится к двойственной задаче (с точностью до знака):

$$\max_{y\in\tilde{Q}} \left\{ \langle x, b-Ay \rangle - \langle z, h(y) \rangle - \varphi(y) \right\} \to \min_{z\ge 0,\ x}. \ \Diamond$$

Точное решение вспомогательной max-задачи в (4.8) будем обозначать через $y(x)$. Во многих важных приложениях основной вклад в сложность расчета $y(x)$ дает умножение $Ay$. Это так, например, для сепарабельных функционалов

$$\varphi(y) = \sum_{i=1}^{m} \varphi_i(y_i)$$

и параллелепипедных ограничений $\tilde{Q}$. В таких случаях задача (4.8) сводится к $n$ задачам одномерной оптимизации, которые с запасом могут быть решены за время (2.2) (см. упражнение 1.4) при условии, что $Ay$ уже было посчитано.

Для двойственного функционала $f(x)$, определяемого согласно (4.8), выполняется условие (2.4) с $\nu = 1$ и $L_1 = L = \frac{1}{\mu} \max_{\|y\|_p \le 1} \|Ay\|_2^2$ [4, 456]. В частности, для $p = 1$

$$L = \frac{1}{\mu} \max_{j=1,\dots,m} \left\| A^j \right\|_2^2,$$

где $A^j$ – $j$-й столбец матрицы $A$. Для $p = 2$

$$L = \frac{1}{\mu}\lambda_{\max}\left(A^T A\right) \overset{def}{=} \frac{1}{\mu}\sigma_{\max}\left(A\right).$$

**Замечание 4.1 (метод регуляризации и техника двойственного сглаживания).** Добиться сильной выпуклости $\varphi(y)$ всегда можно с помощью *регуляризации* задачи. Опишем, в чем состоит *техника регуляризации* (см., например, [6], [14, гл. 9] и цитированную там литературу), восходящая к работам трех отечественных научных школ: А.Н. Тихонова [90] (Москва), М.М. Лаврентьева [61] (Новосибирск), В.К. Иванова [52] (Свердловск), занимавшихся изучением некорректных задач. Рассмотрим новую задачу:

$$\varphi^{\mu}\left(y\right) = \varphi\left(y\right) + \mu V\left(y, y^0\right) \to \min_{Ay=b,\ y\in\tilde{Q}}, \tag{4.9}$$

где $V\left(y, y^0\right)$ – 1-сильно выпуклая в $p$-норме функция $y$. Обозначим через $\varphi_*^{\mu}$ оптимальное значение функционала в задаче (4.9). Пусть[45]

$$\mu \le \frac{\varepsilon}{2V\left(y_*, y^0\right)}, \tag{4.10}$$

и удалось найти $\varepsilon/2$-решение задачи (4.9), т. е. нашелся такой $y_{\varepsilon/2}$, что $Ay_{\varepsilon/2} = b$, $y_{\varepsilon/2} \in \tilde{Q}$:

$$\varphi^{\mu}\left(y_{\varepsilon/2}\right) - \varphi_*^{\mu} \le \varepsilon/2.$$

Тогда

$$\varphi\left(y_{\varepsilon/2}\right) - \varphi_* \le \varepsilon.$$

Действительно,

$$\varphi\left(y_{\varepsilon/2}\right) - \varphi_* \le \varphi^{\mu}\left(y_{\varepsilon/2}\right) - \varphi_* \le \varphi^{\mu}\left(y_{\varepsilon/2}\right) - \varphi_*^{\mu} + \varepsilon/2 \le \varepsilon.$$

Здесь использовались определение $\varphi_*^{\mu}$ и формула (4.10):

$$\varphi_*^{\mu} = \min_{Ay=b,\ y\in\tilde{Q}}\left\{\varphi\left(y\right) + \mu V\left(y, y^0\right)\right\} \le \varphi\left(y_*\right) + \mu V\left(y_*, y^0\right) \le \varphi_* + \varepsilon/2.$$

[45] Как правило, величина $V\left(y_*, y^0\right)$ неизвестна, поэтому на практике используются рестарты по параметру $\mu$, приводящие к увеличению общего числа итераций в несколько раз [20, 118], см. также указание к упражнению 2.3.

Стоит отметить, что если изначально рассматривалась задача вида (4.8), то говорят, что функционал $f(x)$ представим в *форме Лежандра*. Пусть $\tilde{Q}$ – выпуклое компактное множество простой структуры. В этом случае описанная выше техника регуляризации $\varphi(y) \to \varphi^{\mu}(y)$, в которой вместо $R^2 = V(y_*, y^0)$ используется $\tilde{R}^2 = \max_{y \in \tilde{Q}} V(y, y^0)$ с $\mu \le \varepsilon / (2\tilde{R}^2)$, приводит к сглаживанию функции

$$\begin{gathered} f(x) \to f_{\mu}(x) = \max_{y \in \tilde{Q}} \left\{ \langle x, b - Ay \rangle - \varphi(y) - \mu V(y, y^0) \right\}, \\ 0 \le f(x) - f_{\mu}(x) \le \varepsilon/2. \end{gathered} \tag{4.11}$$

При этом $f_{\mu}(x)$ будет иметь константу Липшица градиента в 2-норме:

$$L_{\varepsilon} = \frac{2\tilde{R}^2}{\varepsilon} \max_{\|y\|_p \le 1} \|Ay\|_2^2 .$$

Простейший пример такого сглаживания:

$$\begin{gathered} f(x) = \max_{l=1,\dots,m} \langle c^l, x \rangle = \max_{y \in S_m(1)} \sum_{l=1}^{m} y_l \langle c^l, x \rangle \to \\ \to \max_{y \in S_m(1)} \left\{ \sum_{l=1}^{m} y_l \langle c^l, x \rangle - \mu \sum_{l=1}^{m} y_l \ln\left( \frac{y_l}{1/m} \right) \right\} = \\ = \mu \ln\left( \sum_{l=1}^{m} \exp\left( \langle c^l, x \rangle / \mu \right) \right) - \mu \ln m = f_{\mu}(x), \end{gathered}$$

где $\mu = \varepsilon / (2 \ln m)$. Описанную выше конструкцию (4.11) обычно называют *двойственным сглаживанием* или *техникой сглаживания по Нестерову* [71, гл. 5], [456]. В классе рассматриваемых в этой главе неускоренных методов данная техника по оценкам не дает преимуществ: задача выпуклой оптимизации с негладким функционалом для решения с точностью по функции $\varepsilon$ требует $\sim \varepsilon^{-2}$ вычислений (суб-)градиента (см. упражнение 2.1), и сглаженная задача также требует $\sim L_{\varepsilon}\varepsilon^{-1} \sim \varepsilon^{-2}$ вычислений (суб-)градиента. Однако для ускоренных методов техника двойственного сглаживания приводит к лучшим оценкам [71, гл. 5], [452, 456]:

$$\sim \sqrt{L_{\varepsilon}\varepsilon^{-1}} \sim \sqrt{\varepsilon^{-2}} \sim \varepsilon^{-1} .$$

Разумеется, имеет смысл говорить о двойственном сглаживании только в случае, когда задача максимизации в (4.11) является относительно простой. Как следствие, описанная техника сглаживания применима к намного более узкому классу задач, чем регуляризация. Более того, конструкция, описанная в замечании 5.1 в части решения седловых задач, позволяет получать аналогичные результаты при более общих условиях. Тем не менее стоит отметить, что в определенных (композитных) случаях описанная техника позволяет получать новые результаты, недостижимые с помощью техники замечания 5.1, см., например, [2, 202, 247, 386]. Отметим также, что есть и другие способы сглаживания (см., например, [118]), впрочем, также имеющие весьма ограниченную область применимости.

Хорошо известный пример использования регуляризации – способ вычисления (понимания) *псевдообратной матрицы* [82, 177]:

$$A^{+} = \lim_{\mu \to 0+} \left( A^{T} A + \mu I \right)^{-1} A^{T}.$$

Такое понимание эквивалентно тому, что

$$x_{*} = A^{+} b = \lim_{\mu \to 0+} \arg\min_{x} \left\{ \frac{1}{2} \left\| Ax - b \right\|_2^2 + \frac{\mu}{2} \left\| x \right\|_2^2 \right\}$$

является решением задачи

$$\frac{1}{2} \left\| Ax - b \right\|_2^2 \to \min_{x}$$

с наименьшим значением 2-нормы, если решение не единственно, см. также упражнение 5.9. В анализе данных описанная регуляризация имеет простую содержательную интерпретацию – байесовская регуляризация для задачи нормальной регрессии с нормальным (гауссовским) априорным распределением параметров [99, лекции 13, 14]. На рассмотренном примере также удобно демонстрировать связь метода регуляризации и метода *штрафных функций*, см. замечание 4.3 и [130].

Менее известный, но не менее интересный пример *итеративной регуляризации/сглаживания* имеется в работе [288], в которой решение седловой билинейной задачи сводится к последовательности задач выпуклой оптимизации в условиях острого минимума [82, 507]. ■

◊ Так же как и в упражнении 2.3, здесь можно отметить, что *из оптимального метода для сильно выпуклой задачи можно получить с помощью регуляризации оптимальный метод для просто выпуклой задачи*. Во всяком случае пока не удалось придумать ни одного контрпримера, когда бы это было не так. ◊

Вернемся к задаче (4.8), в которой для больше наглядности будем считать, что $\tilde{Q} = \mathbb{R}^{n_y}$, $n_y = \dim y$,

$$f(x)=\max_{y\in\mathbb{R}^{n_y}}\left\{\langle x,b-Ay\rangle-\varphi(y)\right\}\to\min_{x\in\mathbb{R}^n}.$$

Положим[46] $x^0=0$, $N=2LR^2/\varepsilon$, где $R=\|x_*\|_2$. Рассмотрим метод градиентного спуска (1.22):

$$x^{k+1}=x^k-\frac{1}{L}\nabla f\left(x^k\right).$$

Подобно (4.5) можно написать:

$$f\left(\bar{x}^N\right)\le\frac{1}{N}\min_{x\in\mathbb{R}^n}\left\{\sum_{k=0}^{N-1}\left[f\left(x^k\right)+\left\langle\nabla f\left(x^k\right),x-x^k\right\rangle\right]+\frac{L}{2}\left\|x-x^0\right\|_2^2\right\}\le$$
$$\le\frac{1}{N}\min_{x\in B_{2R}(0)}\left\{\sum_{k=0}^{N-1}\left[f\left(x^k\right)+\left\langle\nabla f\left(x^k\right),x-x^k\right\rangle\right]\right\}+\underbrace{\frac{2LR^2}{N}}_{\varepsilon}.$$

Здесь выбирается шар радиуса $2R$, чтобы впоследствии можно было получить оценку (4.14), для чего нужна оценка (4.13) именно с $2R$, что и требует выбора радиуса шара равным $2R$. Ввиду (4.8) отсюда *по формуле Демьянова–Данскина* [40, 41, 162]

$$\nabla f(x)=b-Ay(x)$$

имеем

$$\begin{aligned}&f\left(\bar{x}^N\right)-\frac{1}{N}\sum_{k=0}^{N-1}\left\langle x^k,b-Ay\left(x^k\right)\right\rangle+\frac{1}{N}\sum_{k=0}^{N-1}\varphi\left(y\left(x^k\right)\right)-\\&-\frac{1}{N}\min_{x\in B_{2R}(0)}\left\{\sum_{k=0}^{N-1}\left\langle b-Ay\left(x^k\right),x-x^k\right\rangle\right\}\le\varepsilon.\end{aligned}\tag{4.12}$$

По выпуклости $\varphi(y)$ (см. (4.7)) из (4.12) имеем

$$f\left(\bar{x}^N\right)+\varphi\underbrace{\left(\frac{1}{N}\sum_{k=0}^{N-1}y\left(x^k\right)\right)}_{\bar{y}^N}+\max_{x\in B_{2R}(0)}\left\{\left\langle A\frac{1}{N}\sum_{k=0}^{N-1}y\left(x^k\right)-b,x\right\rangle\right\}\le\varepsilon,$$

т. е.

$$f\left(\bar{x}^N\right)+\varphi\left(\bar{y}^N\right)+2R\left\|A\bar{y}^N-b\right\|_2\le\varepsilon.\tag{4.13}$$

[46] Если выбирать точку старта $x^0\ne0$, то оценка $R$ в приводимых далее выкладках (результатах) ухудшится $R=\left\|x^0\right\|_2+\left\|x_*-x^0\right\|_2$. Ухудшатся и числовые множители. Детали см., например, в [338].

Из (4.13) и *слабой двойственности* $-\varphi(y_*) \le f(x_*)$ имеем

$$\varphi(\bar{y}^N) - \varphi(y_*) \le \varphi(\bar{y}^N) - \varphi(y_*) + 2R\|A\bar{y}^N - b\|_2 \le$$
$$\le \varphi(\bar{y}^N) + f(x_*) + 2R\|A\bar{y}^N - b\|_2 \le$$
$$\le \varphi(\bar{y}^N) + f(\bar{x}^N) + 2R\|A\bar{y}^N - b\|_2 \le \varepsilon.$$

◊ Фактически слабая двойственность – это отражение простого факта, что всегда имеет место неравенство

$$\max_y \min_x L(x, y) \le \min_x \max_y L(x, y).$$

На самом деле во всех естественных ситуациях, когда рассматриваются невырожденные (совместные) выпуклые задачи, в этом неравенстве имеет место равенство, т. е. имеет место сильная двойственность [176, гл. 5]. ◊

Поскольку $x_*$ одновременно является решением двойственной задачи (4.8) и множителем Лагранжа к ограничению $Ay = b$ в задаче (4.7) (см. представление (4.8)), то

$$f(x_*) = \underbrace{\langle x_*, b - Ay_* \rangle}_{0} - \varphi(y_*) = \max_y \{\langle x_*, b - Ay \rangle - \varphi(y)\} \ge$$
$$\ge \langle x_*, b - Ay^N \rangle - \varphi(y^N),$$

т.е.

$$\varphi(y_*) - \langle x_*, A\bar{y}^N - b \rangle \le \varphi(\bar{y}^N).$$

Объединяя это неравенство с полученным ранее неравенством

$$\varphi(\bar{y}^N) - \varphi(y_*) + 2R\|A\bar{y}^N - b\|_2 \le \varepsilon,$$

получим

$$R\|A\bar{y}^N - b\|_2 \le \varepsilon.$$

Таким образом, установлен следующий результат.

**Теорема 4.1.** *Пусть нужно решить задачу* (4.7) *в следующем смысле*

$$\varphi(\bar{y}^N) - \varphi(y_*) \le \varepsilon, \quad \|A\bar{y}^N - b\|_2 \le \tilde{\varepsilon}. \tag{4.14}$$

*Для этого рассмотрим двойственную задачу* (4.8), *которую будем решать градиентным спуском* (1.22):

$$x^{k+1} = x^k - \frac{1}{L}\nabla f(x^k)$$

*с* $x^0 = 0$. *Выберем в качестве критерия останова метода следующие условия* (*зазор двойственности и невязку в ограничении*):

$$f\left(\bar{x}^N\right)+\varphi\left(\bar{y}^N\right)\le\varepsilon\,,\ \left\|A\bar{y}^N-b\right\|_2\le\tilde{\varepsilon}\,,$$

*где*

$$\bar{x}^N=\frac{1}{N}\sum_{k=1}^{N}x^k\,,\ \bar{y}^N=\frac{1}{N}\sum_{k=0}^{N-1}y\left(x^k\right),$$

*из которых вытекает* (4.14). *Тогда метод гарантированно остановится, сделав не более чем*

$$\max\left\{\frac{2LR^2}{\varepsilon},\frac{2LR}{\tilde{\varepsilon}}\right\} \tag{4.15}$$

*итераций, где* $L=\frac{1}{\mu}\max_{\|y\|_p\le 1}\|Ay\|_2^2$, $R=\|x_*\|_2$. *Если решение задачи* (4.8) $x_*$ *не единственно, то в оценке* $R$ *в* (4.15) *выбирается то решение, которое имеет наименьшую* 2-*норму.*

**Замечание 4.2 (оценка размера решения двойственной задачи).** В оценку (4.15) входит неизвестный размер решения двойственной задачи $R=\|x_*\|_2$ (4.8). Если решение $x_*$ не единственно, то выбирается наименьшее по 2-норме (см. § 1, 2). Это $R$ можно оценить следующим образом [71, п. 4.3.4], [387]:

$$R^2=\|x_*\|_2^2\le\left\|\nabla\varphi\left(y_*\right)\right\|_2^2\Big/\tilde{\sigma}_{\min}\left(A\right), \tag{4.16}$$

где

$$\tilde{\sigma}_{\min}\left(A\right)=\min\left\{\lambda>0:\ \exists\ x\ne 0: AA^Tx=\lambda x\right\}.$$

Действительно, исходя из определения $x_*$ и $y_*$, имеем для любого $y\in\tilde{Q}$

$$-\varphi\left(y_*\right)=\left\langle x_*,b-Ay_*\right\rangle-\varphi\left(y_*\right)=\max_{y\in\tilde{Q}}\left\{\left\langle x_*,b-Ay\right\rangle-\varphi\left(y\right)\right\}\ge$$

$$\ge\left\langle x_*,b-Ay\right\rangle-\varphi\left(y\right)=\left\langle x_*,Ay_*-Ay\right\rangle-\varphi\left(y\right)=\left\langle A^Tx_*,y_*-y\right\rangle-\varphi\left(y\right).$$

Следовательно, для любого $y\in\tilde{Q}$

$$\varphi\left(y\right)\ge\varphi\left(y_*\right)+\left\langle -A^Tx_*,y-y_*\right\rangle.$$

По выпуклости $\varphi\left(y\right)$ отсюда следует, что

$$-A^Tx_*=\nabla\varphi\left(y_*\right).$$

Точнее, $-A^T x_* \in \partial\varphi(y_*)$. Осталось только заметить, для любого $x \in \left(\operatorname{Ker} A^T\right)^{\perp}$ имеет место неравенство

$$\left\|-A^T x\right\|_2^2 = \left\langle -A^T x, -A^T x\right\rangle = \left\langle x, AA^T x\right\rangle \ge \tilde{\sigma}_{\min}(A)\|x\|_2^2 . \blacksquare$$

**Пример 4.1 (децентрализованная распределенная оптимизация [252, 386, 387, 436, 513, 514, 562]).** Пусть необходимо решать задачу выпуклой оптимизации

$$\sum_{i=1}^{n} \varphi_i(y) \to \min_{y\in\mathbb{R}} . \qquad (4.17)$$

Для большей наглядности считаем $y$ скалярной величиной. Заметим, однако, что от этого упрощения легко отказаться. Будем считать, что $\varphi_i''(y) \ge \mu$, $i = 1,...,n$, $y \in \mathbb{R}$. Предположим, что есть связанная сеть (коммуникационный граф) $G = \langle V, E\rangle$ из $n$ узлов. В $i$-м узле хранится функция $\varphi_i(y)$. Зададим матрицу инцидентности

$$\mathrm{I} = \left\|\mathrm{I}_{ij}\right\|_{i,j=1}^{n} : \mathrm{I}_{ij} = 1,\ (i,j) \in E;\ \mathrm{I}_{ij} = 0,\ (i,j) \ne E .$$

По матрице $\mathrm{I}$ построим симметричную неотрицательно определенную *матрицу Лапласа* (*Кирхгофа*) $W \succ 0$:

$$W_{ij} = \begin{cases} -\mathrm{I}_{ij}, i \ne j, \\ \sum_{j=1}^{n} \mathrm{I}_{ij}, i = j. \end{cases}$$

По *теореме Фробениуса–Перрона* [77, § 7, 8, гл. 2]:

$$Wy = 0 \iff y_1 = ... = y_n . \qquad (4.18)$$

Ввиду (4.18) перепишем задачу (4.17) следующим образом:

$$\varphi(y) = \sum_{i=1}^{n} \varphi_i(y_i) \to \min_{Wy=0} . \qquad (4.19)$$

Построим (с точностью до знака) двойственную задачу к задаче (4.19) (см. (4.8)):

$$\begin{aligned} f(x) = \varphi^*(-Wx) &= \max_{y\in\mathbb{R}^n}\left\{-\langle x, Wy\rangle - \varphi(y)\right\} = \\ &= \sum_{i=1}^{n} \max_{y_i\in\mathbb{R}}\left\{[-Wx]_i\, y_i - \varphi_i(y_i)\right\} \to \min_{x\in\mathbb{R}^n}. \end{aligned} \qquad (4.20)$$

Считаем, что $i$-я вспомогательная задача максимизации в (4.20) может эффективно решаться в $i$-м узле. Заметим, что для функции $f(x)$ константу Липшица градиента можно оценить как $L = \sigma_{\max}(W)/\mu$ [4, 456], а на размер двойственного решения есть оценка (см. замечание 4.2) $R^2 = \|x_*\|_2^2 \le \|\nabla\varphi(y_*)\|_2^2 / \tilde{\sigma}_{\min}(W)$. Согласно написанному ранее в этом параграфе решать задачу (4.20) можно методом

$$x_i^{k+1} = x_i^k - \frac{1}{L}\left[\nabla f\left(x^k\right)\right]_i = x_i^k + \frac{1}{L}\left[W\tilde{y}\left(-Wx^k\right)\right]_i, \qquad (4.21)$$

где через $\tilde{y}(-Wx)$ обозначается решение задачи (4.20). Итак, пусть в каждом узле хранятся $\left\{x_i^k, \tilde{y}_i\left(\left[-Wx^k\right]_i\right)\right\}_k$. Ключевое наблюдение: чтобы вычислить $\tilde{y}_i\left(\left[-Wx^k\right]_i\right)$, $i$-у узлу необходимо обратиться только к своим непосредственным соседям за соответствующими компонентами вектора $x^k$ (см. (4.20)), а чтобы вычислить $x_i^{k+1}$, $i$-му узлу также необходимо обратиться только к своим непосредственным соседям за соответствующими компонентами вектора $\tilde{y}(-Wx)$ (см. (4.21)). Таким образом, один шаг градиентного спуска для двойственной задачи приводит к коммуникации каждого узла со своими соседями два раза (передается два числа). Поскольку вычислительные возможности узлов, как правило, на несколько порядков выше скорости передачи информации по сети, то полученный дисбаланс (решать вспомогательную задачу поиска $\tilde{y}_i\left(\left[-Wx^k\right]_i\right)$ заметно труднее, чем послать и принять несколько чисел) хорошо способствует эффективному решению задачи.

Несложно заметить (см. (4.15)), что время работы алгоритма будет прямо пропорционально числу обусловленности матрицы $W^T W = W^2$, т. е. $\sigma_{\max}(W)/\tilde{\sigma}_{\min}(W)$. В действительности, можно улучшить описанный выше подход, если сделать замену $W \to \sqrt{W}$ [514]:

$$\sqrt{W}y = 0 \;\Leftrightarrow\; y_1 = \ldots = y_n,$$

$$\tilde{y}\left(-\sqrt{W}x\right) = \arg\max_{y\in\mathbb{R}^n}\left\{-\left\langle\sqrt{W}x, y\right\rangle - \varphi(y)\right\},$$

$$\sqrt{W}x^{k+1} = \sqrt{W}x^k + \frac{1}{L}W\tilde{y}\left(-\sqrt{W}x^k\right).$$

Обозначая $z = \sqrt{W}x$, запишем метод в новых переменных:[47]

$$\tilde{y}(z) = \arg\max_{y\in\mathbb{R}^n}\left\{\langle z, y\rangle - \varphi(y)\right\},$$

$$z^{k+1} = z^k + \frac{1}{L}W\tilde{y}\left(-z^k\right).$$

Легко понять, что такой метод также может работать распределено [387, 436, 514, 562], причем один шаг такого варианта градиентного спуска для двойственной задачи приводит к коммуникации каждого узла со своими соседями всего один раз. Таким образом, можно редуцировать

$$\sigma_{\max}(W)/\tilde{\sigma}_{\min}(W) \text{ к } \sigma_{\max}\left(\sqrt{W}\right)/\tilde{\sigma}_{\min}\left(\sqrt{W}\right) = \sqrt{\sigma_{\max}(W)/\tilde{\sigma}_{\min}(W)}.$$

На основе ускоренных градиентных методов можно построить более быстрые децентрализованные распределенные алгоритмы решения задачи (4.17), см. [252, 294, 514, 562].

К сожалению, во всех случаях (ускоренном и неускоренном) не удается построить адаптивные / универсальные (см. § 5) варианты таких методов (в децентрализованном случае). Равно как и не удается предложить эффективный (практический) критерий останова методов (см. замечание 2.1 и § 4).

Отметим также, что между рассмотренной в этом примере задачей и задачами типа распределения ресурсов (см. упражнение 4.7) имеется связь – двойственная задача к задаче о распределении ресурсов имеет вид (4.17). Таким образом, появляется возможность решать задачу распределения ресурсов не централизованным образом, как предлагается в указании к упражнению 4.7, а децентрализованным образом, как в разобранном примере. Детали см., например, в [199] (см. также [175, п. 7.3.1]). ■

◊ После работы [175] распределенная оптимизация (Grid-технологии) прочно закрепилась в современном анализе данных, см., например, концепцию Google: Federated Learning [363, 375, 394, 423].

Хорошей практической демонстрацией, описанной в примере 4.1 техники, является распределенный способ вычисления барицентра Васерштейна, учитывающий наличие явного *представления Лежандра* (сопряженного представления) для расстояния Монжа–Канторовича–Васерштейна [213, 256, 260, 376, 488, 561], см. также замечание 4.4.

---

[47] Заметим, что если сделать такую замену в самой двойственной задаче, то придется вместо условия $z \in \mathbb{R}^n$ писать $z \in \operatorname{Im}\sqrt{W} = \left(\ker\sqrt{W}\right)^{\perp}$, что порождает сложности с интерпретацией метода. Описанный в этом примере подход (замена переменных в самом алгоритме) является, пожалуй, наиболее простым способом преодоления этих сложностей.

Подобно теореме 1 можно распространить градиентный спуск и на невыпуклые задачи распределенной оптимизации [545]. ◊

При описанном выше подходе, к сожалению, возникает невязка в ограничении $Ay = b$ в задаче (4.7). Эту невязку можно полностью устранить, изменив подход. Далее мы будем в основном следовать работе [455]. Немного обобщим постановку задачи (4.7):

$$\varphi(y) = F(y) + g(y) \to \min_{Ay \le b,\, y \in \tilde{Q}}. \tag{4.22}$$

Вместо равенства $Ay = b$ в (4.7) стали рассматривать неравенство $Ay \le b$ в (4.22) и добавили простой выпуклый композитный член $g(y)$. Будем предполагать, что выпуклая функция $F(y)$ удовлетворяет (только) условию (2.3) (см. также (2.26)), которое в данном случае будет иметь вид

$$F(y) \le F(z) + \langle \nabla F(z), y - z \rangle + \frac{L}{2}\|y - z\|^2 + \delta.$$

Рассмотрим метод вида (3.19) с шагом $h = 1/L$ (2.21) для задачи (4.22):

$$y^{k+1} = \arg\min_{Ay \le b,\, y \in \tilde{Q}} \left\{ \left\langle \nabla F\left(y^k\right), y - y^k \right\rangle + g(y) + LV\left(y, y^k\right) \right\}. \tag{4.23}$$

Для наглядности будем считать, что задача (4.23) решается на каждой итерации $k$ явно (точно). В приложениях задача (4.23) может быть простой, например, когда $Ay \le b$ имеет вид $y \le \bar{y}$ или $y \ge \bar{y}$ [18, гл. 1, 3].

Повторяя рассуждения примера 3.1 (см. формулы (3.12), (3.20)), из (4.23) получим подобно оценке (4.2):

$$\varphi\left(\bar{y}^N\right) \le$$

$$\le \min_{Ay \le b,\, y \in \tilde{Q}} \left\{ \frac{1}{N} \sum_{k=0}^{N-1} \left[ F\left(y^k\right) + \left\langle \nabla F\left(y^k\right), y - y^k \right\rangle + g(y) \right] + \frac{LV\left(y, y^0\right)}{N} \right\} + \delta, \tag{4.24}$$

где

$$\bar{y}^N = \frac{1}{N} \sum_{k=1}^{N} y^k .$$

Обозначим множитель Лагранжа к ограничению $Ay \le b$ в (4.24) через $\tilde{x}^N \ge 0$. Тогда (4.24) можно переписать следующим образом: для любого $\tilde{y} \in \tilde{Q}$

$$\varphi\left(\bar{y}^N\right) \le \min_{y\in\tilde{Q}}\left\{\frac{1}{N}\sum_{k=0}^{N-1}\underbrace{\left[F\left(y^k\right)+\left\langle\nabla F\left(y^k\right),y-y^k\right\rangle+g\left(y\right)\right]}_{\le\varphi(y)}+\right.$$

$$\left.+\left\langle\tilde{x}^N,Ay-b\right\rangle+\frac{LV\left(y,y^0\right)}{N}\right\}+\delta\le\varphi\left(\tilde{y}\right)+\left\langle\tilde{x}^N,A\tilde{y}-b\right\rangle+\frac{LV\left(\tilde{y},y^0\right)}{N}+\delta. \tag{4.25}$$

Введем, подобно (4.8), двойственную (с точностью до знака) функцию

$$f\left(x\right)=\max_{y\in\tilde{Q}}\left\{\left\langle x,b-Ay\right\rangle-\varphi\left(y\right)\right\}. \tag{4.26}$$

Обозначим, как и раньше, через $y\left(x\right)$ решение задачи максимизации в (4.26). Тогда, выбирая в (4.25) $\tilde{y}=y\left(\tilde{x}^N\right)$ и обозначая через $R^2=V\left(y\left(\tilde{x}^N\right),y^0\right)$, получим

$$0\le\varphi\left(\bar{y}^N\right)-\varphi\left(y_*\right)\le\varphi\left(\bar{y}^N\right)+f\left(\tilde{x}^N\right)\le\frac{LR^2}{N}+\delta\,. \tag{4.27}$$

◊ В отличие от (4.14) в (4.27) используется допустимая точка $\bar{y}^N$: $A\bar{y}^N\le b$, поэтому в (4.27) имеет место оценка снизу $0\le\varphi\left(\bar{y}^N\right)-\varphi\left(y_*\right)$. Из слабой двойственности имеем

$$\varphi\left(\bar{y}^N\right)-\varphi\left(y_*\right)+f\left(\bar{x}^N\right)-f\left(x_*\right)\le\varphi\left(\bar{y}^N\right)+f\left(\bar{x}^N\right).$$

С учетом этих неравенств из (4.27) имеем

$$0\le f\left(\bar{x}^N\right)-f\left(x_*\right)\le\varphi\left(\bar{y}^N\right)+f\left(\bar{x}^N\right)\le\frac{LR^2}{N}+\delta\,. \ ◊$$

Из формулы (4.27) (с $\delta=\varepsilon/2$) при условии, что функция $F\left(y\right)$ удовлетворяет (только) условию (2.27), следует, что метод вида (3.19) гарантированно остановится по критерию

$$\varphi\left(\bar{y}^N\right)+f\left(\bar{x}^N\right)\le\varepsilon\,,$$

сделав не более

$$N=\frac{2LR^2}{\varepsilon}\le\left(\frac{2L_\nu R^{1+\nu}}{\varepsilon}\right)^{\frac{2}{1+\nu}} \tag{4.28}$$

итераций (вычислений $\nabla F\left(y^k\right)$).

Описанный способ решения задачи композитной оптимизации (4.22) может быть осуществлен в модельной общности (см. § 3) [94].

В заключение подчеркнем, в чем различие в описанных в этом параграфе подходах. В подходе (4.23) решается исходная задача (4.22), в то время как в подходе, описанном в первой половине параграфа, решается двойственная задача (4.8). Отметим, что оба описанных подхода можно сочетать друг с другом [19].

**Упражнение 4.1 (условие Слейтера).** Рассматривается задача выпуклой оптимизации

$$f(x) \to \min_{h(x)\le 0,\, x\in Q},$$

где $h\colon \mathbb{R}^n \to \mathbb{R}^m$. Двойственная задача (с точностью до знака) имеет вид

$$\varphi(y) = \max_{x\in Q}\left\{-\langle y, h(x)\rangle - f(x)\right\} \to \min_{y\in\mathbb{R}^m_+}. \qquad (4.29)$$

Обозначим через $y^*$ решение двойственной задачи. Предположим, что выполняется *условие Слейтера*:

*существует такая точка* $\bar{x}\in Q$*, что* $h(\bar{x})<0$.

Пусть $\gamma = \min_{i=1,\dots,m}\left\{-h_i(\bar{x})\right\}$. Покажите, что

$$\left\|y^*\right\|_1 \le \frac{1}{\gamma}\left(f(\bar{x})+\varphi(0)\right) = \frac{1}{\gamma}\left(f(\bar{x}) - \min_{x\in Q} f(x)\right).$$

**Указание.** См. [163]. Ключевое неравенство:

$$\varphi(0) \ge \varphi(y^*) = \max_{x\in Q}\left\{-\sum_{i=1}^m y_i^* h_i(x) - f(x)\right\} \ge -\sum_{i=1}^m y_i^* h_i(\bar{x}) - f(\bar{x}).$$

Аналогичным образом можно получать оценки на размер решения двойственной задачи и в более общих случаях (см., например, [27]). ■

**Упражнение 4.2.** Пусть задача из упражнения 4.1 решена в следующем смысле: найден такой $\tilde{y}\in\mathbb{R}^m_+$, что

$$\left\langle \tilde{y}, \nabla\varphi(\tilde{y})\right\rangle \le \varepsilon, \quad \left\|\left[-\nabla\varphi(\tilde{y})\right]_+\right\|_2 \le \tilde{\varepsilon},$$

где по формуле Демьянова–Данскина [176, гл. 3] $\nabla\varphi(\tilde{y}) = -h(x(\tilde{y}))$, $x(\tilde{y})$ – решение вспомогательной задачи максимизации в (4.29). Тогда

$$f(x(\tilde{y})) - f(x_*) \le \varepsilon, \quad \left\|\left[h(x(\tilde{y}))\right]_+\right\|_2 \le \tilde{\varepsilon}.$$

**Указание.** Ключевая выкладка:

$$-\left\langle \tilde{y}, h\left(x\left(\tilde{y}\right)\right)\right\rangle - f\left(x\left(\tilde{y}\right)\right) \ge -\langle \underset{\tilde{y}\ge 0}{\tilde{y}}, \underset{h(x_*)\le 0}{h\left(x_*\right)} \rangle - f\left(x_*\right) \ge -f\left(x_*\right).$$

Заметим также, что если здесь вместо ограничения в виде неравенства $h(x) \le 0$ имели бы аффинное ограничение в виде равенства $Ax - b = 0$, то тогда в проведенных рассуждения следовало бы сделать следующую корректировку: $\nabla\varphi(\tilde{y}) = b - Ax(\tilde{y})$, как следствие, условие $\left\|\nabla\varphi(\tilde{y})\right\|_2 \le \tilde{\varepsilon}$ обеспечивает выполнение условия $\left\|Ax(\tilde{y}) - b\right\|_2 \le \tilde{\varepsilon}$. ■

**Упражнение 4.3. 1)** Рассматривается задача поиска *седловой точки* вида (4.8)

$$f(x) = \max_{y\in\tilde{Q}}\left\{\langle x, b - Ay\rangle - \varphi(y)\right\} \to \min_{x\in\mathbb{R}^n},$$

где функция $\varphi(y)$ – $\mu$-сильно выпуклая относительно $p$-нормы ($1 \le p \le 2$). Покажите, что функция $f(x)$ будет гладкой, с константой Липшица градиента в 2-норме:

$$L = \frac{1}{\mu}\max_{\|z\|_p \le 1}\|Az\|_2^2.$$

Более того, если $y_\delta(x)$ – решение вспомогательной задачи максимизации с точностью по функции $\delta$, то

$$\left(\langle x, b - Ay_\delta(x)\rangle - \varphi\left(y_\delta(x)\right); b - Ay_\delta(x)\right)$$

будет $(\delta, 2L)$-моделью функции $f(x)$ в точке $x$ относительно 2-нормы (см. начало § 3).

Предложите численный метод решения поставленной (седловой) задачи, и сравните его трудоемкость с нижними оценками [480]. Если $\varphi(y)$ дополнительно имеет Липшицев градиент, можно ли это как-то использовать (см. текст, следующий после упражнения 4.8 до замечания 4.4 включительно)?

Покажите, что если вместо $\langle x, b - Ay\rangle$ в определении $f(x)$ использовать функцию $F(x, y)$, которая выпуклая по $x$, вогнутая по $y$, а также достаточно гладкая:

$$\left\|\nabla_x F(x_2, y) - \nabla_x F(x_1, y)\right\|_2 \le L_{xx}\left\|x_2 - x_1\right\|_2,$$

$$\left\|\nabla_x F(x, y_2) - \nabla_x F(x, y_1)\right\|_2 \le L_{xy}\left\|y_2 - y_1\right\|_2,$$

то

$$\left(F\left(x, y_\delta(x)\right) - \varphi\left(y_\delta(x)\right); \nabla_x F\left(x, y_\delta(x)\right)\right)$$

будет $(2\delta, 2L)$-моделью, где $L = L_{xx} + 2L_{xy}^2 / \mu$, функции $f(x)$ в точке $x$ относительно 2-нормы.

**2)** Задачу вида (4.7) можно решать с помощью *модифицированной функции Лагранжа* (*augmented Lagrangians*) [161]. В основе подхода лежит идея переписывания задачи (4.7) следующим образом ( $\mu \ge 0$ – выбираемый параметр):

$$\varphi(y) + \frac{\mu}{2}\|Ay - b\|_2^2 \to \min_{Ay=b,\ y\in\tilde{Q}},$$

и стандартный переход к двойственной задаче:

$$f(x) = \max_{y\in\tilde{Q}} \underbrace{\left\{\langle x, b - Ay\rangle - \varphi(y) - \frac{\mu}{2}\|Ay - b\|_2^2\right\}}_{\Psi(y,x)} \to \min_{x\in\mathbb{R}^n}.$$

Покажите, что если $y_\delta(x)$ – решение вспомогательной задачи максимизации в смысле (3.5):

$$\max_{y\in\tilde{Q}} \left\langle \nabla_y \Psi\left(y_\delta(x), x\right), y - y_\delta(x)\right\rangle \le \delta,$$

то

$$\left(\left\langle x, b - Ay_\delta(x)\right\rangle - \varphi\left(y_\delta(x)\right) - \frac{\mu}{2}\|Ay_\delta(x) - b\|_2^2; b - Ay_\delta(x)\right)$$

будет $(\delta, \mu^{-1})$-моделью функции $f(x)$ в точке $x$ относительно 2-нормы (см. начало § 3).

**Указание.** См. [2, 43, 202, 227, 328, 456, 567, 586]. ■

◊ Описанная в первом пункте упражнения 4.3 конструкция позволяет, в частности, решать гладкие $\mu_x$-сильно выпуклые и $\mu_y$-сильно вогнутые седловые задачи (см. замечание 5.1) за[48] $\tilde{O}\left(1/\sqrt{\mu_x\mu_y}\right)$ вычислений $\nabla_x$ и $\tilde{O}\left(1/\sqrt{\mu_x^2\mu_y}\right)$ вычислений $\nabla_y$ или, наоборот (в зависимости от того, что нам выгоднее), за $\tilde{O}\left(1/\sqrt{\mu_x\mu_y}\right)$ вычислений $\nabla_y$ и $\tilde{O}\left(1/\sqrt{\mu_x\mu_y^2}\right)$ вычислений $\nabla_x$, что улучшает наилучшую известную сейчас оценку

[48] Здесь и далее в этом абзаце для наглядности в оценках приводится зависимость только от $\mu_x$, $\mu_y$ и $\varepsilon$.

$\tilde{O}\left(1/\min\left\{\mu_x,\mu_y\right\}\right)$ (нижняя оценка[49] $\tilde{O}\left(1/\sqrt{\mu_x\mu_y}\right)$ [586]) на число вычислений $\nabla_x$ и $\nabla_y$ (см. указание к упражнению 5.4), например, по числу вычислений $\nabla_x$ за счет ухудшения оценки на число вычислений $\nabla_y$. Если же, скажем, по $x$ есть только выпуклость, то за счет регуляризации можно обеспечить $\mu_x \simeq \varepsilon/R^2$, см. замечание 4.1. В этом случае седловую задачу можно решить за $\tilde{O}\left(1/\sqrt{\mu_y\varepsilon}\right)$ вычислений $\nabla_x$ и $\tilde{O}\left(1/\sqrt{\mu_y^2\varepsilon}\right)$ вычислений $\nabla_y$. Отметим также, что для седловых задач вида суммы [328] и ускоренных рандомизированных алгоритмов из замечания 1 приложения приведенные здесь оценки можно дополнительно улучшить.

Описанный во втором пункте упражнения 4.3 метод модифицированной функции Лагранжа лежит в основе одного из самых популярных алгоритмов распределенной оптимизации ADMM [175, 386, 479]. ◊

**Замечание 4.3 (метод штрафных функций).** Метод модифицированной функции Лагранжа тесно связан с методом *штрафных функций* [13, § 16, гл. 5]. Для полноты картины приведем здесь соответствующие идеи. Вместо исходной, вообще говоря, невыпуклой задачи условной оптимизации

$$f(x) \to \min_{g(x)=0} \tag{4.30}$$

рассматривается задача безусловной оптимизации:

$$f(x) + \frac{K}{2}\|g(x)\|_2^2 \to \min_x . \tag{4.31}$$

К задаче (4.31) можно прийти, например, *релаксировав* исходную постановку следующим образом [130]:

---

[49] В работе [411] было показано как каталист (см. также замечание 3.3 и приложение) позволяет получить (с точностью до логарифмических множителей) оптимальную оценку на базе описанного подхода. Исходная седловая задача переписывается как задача оптимизации $f(x) = \max_y F(x,y) \to \min_x$. Эту задачу оптимизации предлагается решать ускоренным проксимальным методом (см. замечание 3.3), на каждой итерации которого необходимо решать задачу вида $\min_x \max_y \left\{F(x,y) + \left\|x - x^k\right\|_2^2\right\}$. Число таких итераций будет $\tilde{O}\left(1/\sqrt{\mu_x}\right)$. Вспомогательная седловая задача на каждой итерации может быть решена за $\tilde{O}\left(1/\sqrt{\mu_y}\right)$ вычислений $\nabla_x$ и $\nabla_y$. Таким образом, в итоге и получается желаемая оценка [43, 567].

$$f(x) \to \min_{\frac{1}{2}\|g(x)\|_2^2 \le \frac{1}{2}\varepsilon^2}.$$

В таком случае $K := K(\varepsilon)$ можно понимать как множитель Лагранжа к ограничению

$$\frac{1}{2}\|g(x)\|_2^2 \le \frac{1}{2}\varepsilon^2.$$

Обозначим решение задачи (4.31) через $x^K$, а решение исходной задачи (4.30) через $x_*$. Тогда также имеет место следующая связь метода множителей Лагранжа и метода штрафных функций (см., например, [13, § 17, гл. 5], [82, п. 5 § 2, гл. 8])

$$Kg(x^K) \xrightarrow[K\to\infty]{} \lambda, \text{ т. е. } g(x^K) \simeq \frac{\lambda}{K},$$

$$f(x^K) - f(x_*) = \mathrm{O}\left(\frac{\|\lambda\|_2^2}{K}\right),$$

где $\lambda$ – множитель Лагранжа к ограничению $g(x) = 0$. Метод модифицированной функции Лагранжа является промежуточным методом между отмеченными двумя и может быть проинтерпретирован как их комбинация (сочетание). Метод штрафных функций является одним из наиболее простых и универсальных способов сведения задач условной оптимизации к задачам безусловной оптимизации [37].

Для задач выпуклой оптимизации с ограничениями (см. упражнение 4.1)

$$f(x) \to \min_{h(x)\le 0,\, x\in Q},$$

где $h:\ \mathbb{R}^n \to \mathbb{R}^m$ в случае возможности ограничить решение двойственной задачи $\lambda_* \in Q_\lambda \subseteq \mathbb{R}_+^m$ (например, из соображений упражнения 4.1) можно построить выпуклую негладкую точную штрафную функцию [452, item 3.1.7]:

$$x_* \in \operatorname{Arg}\min_{h(x)\le 0,\, x\in Q} f(x) = \operatorname{Arg}\min_{x\in Q}\max_{\lambda\in\mathbb{R}_+^m}\left\{f(x) + \langle\lambda, h(x)\rangle\right\} =$$

$$= \operatorname{Arg}\min_{x\in Q}\left\{f(x) + \max_{\lambda\in\mathbb{R}_+^m}\langle\lambda, h(x)\rangle\right\} = \operatorname{Arg}\min_{x\in Q}\left\{f(x) + \max_{\lambda\in Q_\lambda}\langle\lambda, h(x)\rangle\right\}. \quad \blacksquare$$

◊ На метод штрафных функций (с немного другой функцией штрафа – см. (1.43)) в случае аффинных ограничений $g(x) = Ax - b = 0$ можно посмотреть и с точки зрения двойственного сглаживания, см. замечание 4.1. Действительно, задачу

$$f(x) \to \min_{Ax=b}$$

можно переписать следующим образом

$$f(x) + \sup_{\lambda} \langle \lambda, Ax - b \rangle \to \min_{x}.$$

Существует такой $\lambda_*$ (множитель Лагранжа), что последняя задача равносильна

$$f(x) + \langle \lambda_*, Ax - b \rangle \to \min_{x}.$$

Будем считать, что $\|\lambda_*\|_2 \le R_\lambda$. Тогда исходная задача равносильна следующей

$$f(x) + \max_{\|\lambda\|_2 \le R_\lambda} \langle \lambda, Ax - b \rangle \to \min_{x}.$$

Отсюда имеем, что $\varepsilon/2$ -решение следующей задачи

$$f(x) + \max_{\|\lambda\|_2 \le R_\lambda} \left\{ \langle \lambda, Ax - b \rangle - \frac{\varepsilon}{2R_\lambda^2} \|\lambda\|_2^2 \right\} =$$

$$= f(x) + \begin{cases} \dfrac{R_\lambda^2}{2\varepsilon} \|Ax - b\|_2^2, & \|Ax - b\|_2 \le \dfrac{\varepsilon}{R_\lambda} \\ R_\lambda \|Ax - b\|_2 - \dfrac{\varepsilon}{2}, & \|Ax - b\|_2 > \dfrac{\varepsilon}{R_\lambda} \end{cases} \to \min_{x}$$

будет $\varepsilon$ -решением исходной.

Близкий способ формирования штрафа [130] (следует сравнить с задачей (4.31) и выбором $K(\varepsilon)$)[50]

---

[50] Данный критерий также можно понимать как (байесовскую) свертку двух критериев

$$f(x) \to \min_{x}, \quad \|Ax - b\|_2^2 \to \min_{x}.$$

Чтобы определить с какими весами сворачивать эти критерии (не ограничивая общности, достаточно из двух весов оставить только вес при квадратичном функционале, полагая вес при $f(x)$ равным 1), нужно вспомнить, что хочется достичь (см. (4.14)): $f(x^N) - f(x_*) \le \varepsilon$, $\|Ax^N - b\|_2 \le \varepsilon/R_\lambda$. Переписывая последнее условие, как $\dfrac{R_\lambda^2}{\varepsilon}\|Ax^N - b\|_2^2 \le \varepsilon$, получим, что критерии

$$F_\varepsilon(x) = f(x) + \frac{R_\lambda^2}{\varepsilon}\|Ax - b\|_2^2 \to \min_x$$

также позволяет восстанавливать решение исходной задачи. Действительно, из

$$F_\varepsilon(x^N) - \min_x F_\varepsilon(x) \le \varepsilon ,$$

следует

$$f(x^N) - f(x_*) + \frac{R_\lambda^2}{\varepsilon}\|Ax^N - b\|_2^2 \le \varepsilon .$$

В частности,

$$f(x^N) - f(x_*) \le \varepsilon .$$

В виду неравенства (см. вывод теоремы 4.1 в обозначениях $f \to \varphi$, $x \to y$, $y \to -\lambda$)

$$-R_\lambda\|Ax - b\|_2 \le \langle \lambda_*, Ax - b\rangle \le f(x) - f(x_*),$$

отсюда имеем, что

$$-R_\lambda\|Ax^N - b\|_2 + \frac{R_\lambda^2}{\varepsilon}\|Ax^N - b\|_2^2 \le \varepsilon .$$

Откуда следует, что

$$R_\lambda\|Ax^N - b\|_2 \le \frac{1+\sqrt{5}}{2}\varepsilon < 2\varepsilon .$$

---

$$f(x) \to \min_x , \quad \frac{R_\lambda^2}{\varepsilon}\|Ax - b\|_2^2 \to \min_x$$

следует сворачивать с одинаковыми весами (равными 1), что и было сделано.

Заметим, что если $f(x)$ – негладкая выпуклая функция с константой Липшица $M$, то с помощью слайдинга Дж. Лана [385] можно решить исходную задачу с указанной выше точностью

$$F_\varepsilon\left(x^N\right)-\min_x F_\varepsilon(x)\le\varepsilon,$$

используя

$$O\left(\frac{M^2R_x^2}{\varepsilon^2}\right)$$

вычислений $\nabla f(x)$ и

$$O\left(\sqrt{\frac{\lambda_{\max}\left(A^TA\right)R_\lambda^2R_x^2}{\varepsilon^2}}\right)=O\left(\frac{\lambda_{\max}\left(\sqrt{A^TA}\right)R_\lambda R_x}{\varepsilon}\right)$$

умножений $A^TAx$. Этот результат при $A=\sqrt{W}$ и $b=0$ (см. пример 4.1) позволяет довольно просто объяснить различие между числом коммуникационных шагов и числом вызовов оракула (выдающего $\nabla f_k(x)$) в задачах негладкой децентрализованной распределенной оптимизации, а также построить оптимальные методы децентрализованной распределенной оптимизации для гладких задач, см. [252, 387] и замечание 4.4. ◊

**Упражнение 4.4.** Технику регуляризации, описанную в замечании 4.1, можно применять не только к задаче (4.7), но и к задаче (4.8):

$$f^\mu(x)=f(x)+\frac{\mu}{2}\|x\|_2^2\to\min_{x\in\mathbb{R}^n}. \tag{4.32}$$

Обозначим через $x_*^\mu$ – решение задачи (4.32). Покажите, что

$$\|\nabla f(x)\|_2\le\|\nabla f^\mu(x)\|_2+\mu\|x\|_2,$$

$$\langle x,\nabla f(x)\rangle\le\frac{L_\mu}{\mu}\left(f^\mu(x)-f^\mu\left(x_*^\mu\right)\right),\ L_\mu=L+\mu,$$

$$\|x_*^\mu\|_2^2\le\frac{2}{\mu}\left(f(0)-f(x_*)\right),$$

где для всех $x,y\in\mathbb{R}^n$ по постановке имеет место неравенство

$$\|\nabla f(y)-\nabla f(x)\|_2\le L\|y-x\|_2.$$

Используя эти оценки и упражнение 4.2, исследуйте скорость сходимости метода (1.3) с шагом $h = 1/L_\mu$ и $x^0 = 0$, сходящегося по оценке (1.24), на задаче (4.32) с критерием останова:

$$\left\langle x^N, \nabla f\left(x^N\right)\right\rangle \left( \le \frac{L_\mu}{\mu}\left( f^\mu\left(x^N\right) - f^\mu\left(x_*^\mu\right)\right)\right) \le \varepsilon,$$

$$\left\|\nabla f\left(x^N\right)\right\|_2 \left( \le \left\|\nabla f^\mu\left(x^N\right)\right\|_2 + \mu\left\|x^N\right\|_2\right) \le \tilde{\varepsilon}.$$

Учитывая, что $\left\|x^N\right\|_2 \le 2\left\|x_*^\mu\right\|_2$ (см. (1.12)), предложите способ выбора параметра регуляризации $\mu$. Сопоставьте полученную таким образом оценку скорости сходимости с оценкой (4.16), учитывая, что $\left\|x_*^\mu\right\|_2 \le \left\|x_*\right\|_2$.

Обобщите полученные результаты на случай, когда в исходной и регуляризованной постановке задачи (4.32) вместо $x \in \mathbb{R}^n$ стоит $x \in \mathbb{R}^{n_1} \otimes \mathbb{R}_+^{n_2}$.

**Указание.** Детали см., например, в работе [20]. Из неравенства (1.8) имеем

$$f^\mu(x) - f^\mu\left(x_*^\mu\right) \ge \frac{\left\|\nabla f^\mu(x)\right\|_2^2}{2L_\mu} = \frac{\left\|\nabla f(x) + \mu x\right\|_2^2}{2L_\mu} \ge \frac{\mu\left\langle \nabla f(x), x\right\rangle}{L_\mu}.$$

Из неравенства (1.14) имеем

$$\frac{\mu}{2}\left\|x_\mu^*\right\|_2^2 = \frac{\mu}{2}\left\|0 - x_\mu^*\right\|_2^2 \le f^\mu(0) - f^\mu\left(x_*^\mu\right) \le f(0) - f(x_*).$$

Последнее неравенство имеет место ввиду $f^\mu(0) = f(0)$ и

$$f^\mu\left(x_*^\mu\right) = f\left(x_*^\mu\right) + \frac{\mu}{2}\left\|x_*^\mu\right\|_2^2 \ge f(x_*) + \frac{\mu}{2}\left\|x_*^\mu\right\|_2^2 \ge f(x_*).$$

В случае наличия у $f(x)$ представления (4.8) имеет место также оценка

$$f(0) - f(x_*) \le \min_{Ay=b,\ y\in\tilde{Q}} \varphi(y) - \min_{y\in\tilde{Q}} \varphi(y).$$

Заметим, что для оценки $\left\|x_*^\mu\right\|_2$ сверху можно было бы также использовать неравенство $\left\|x_*^\mu\right\|_2 \le \left\|x_*\right\|_2$ и (4.16). ■

**Упражнение 4.5.** Обобщите рассуждения § 4 в случае, когда в задаче (4.22) вместо неравенства $Ay \le b$ рассматриваются общие выпуклые ограничения: $Ay = b$, $h(y) \le 0$.

**Указание.** См. [455]. ■

**Упражнение 4.6 (сложность проектирования).** В самом начале § 2 была приведена оценка (2.2). Покажите, что если множество $Q$ есть шар в $p$-норме и(или) задается небольшим числом сепарабельных выпуклых неравенств вида $\sum_{i=1}^{n} h_i^j(x_i) \le 0$, $j = 1,\ldots,m$, то задачи вида (2.6), (2.29) и при определенных условиях (3.3) могут быть решены (в смысле (3.4)) за время $\mathrm{O}\left(nm^2 \ln^2(n/\varepsilon)\right)$.

**Указание.** Характерный пример получения такой оценки разбирается в [27] (см. также [158, п. 5.3.3]). В основе подхода – решение малоразмерной двойственной задачи каким-нибудь прямодвойственным быстро (линейно) сходящимся методом типа метода эллипсоидов [211, 443], см. также упражнения 1.4, 5.5. Предварительно двойственные переменные компактифицируются (см. упражнение 4.1), а на заключительном этапе при рассмотрении исходной (прямой) задачи уже используется оценка из упражнения 3.1. Для расчета градиента двойственного функционала необходимо решить $n$ одномерных задач с точным оракулом, что может быть сделано за линейное время (см. упражнение 1.4), и поскольку двойственную задачу мы также можем решать за линейное время, то все «огрубления», накопленные по ходу описанных рассуждений, соберутся под логарифмом и испортят лишь мультипликативную константу в итоговой оценке. Отметим также, что в (2.29) в функционале присутствует не сепарабельное слагаемое вида (см. табл. 1 в § 2):

$$\|x\|_p^2 = \left(\sum_{i=1}^{n} |x_i|^p\right)^{2/p}.$$

Однако, введя новую переменную $y = \|x\|_p^2$, можно занести это слагаемое в ограничение, заменив в функционале $\|x\|_p^2$ на $y$ и добавив сепарабельное выпуклое ограничение вида неравенства $\|x\|_p^p \le y^{p/2}$, где $p/2 < 1$. ■

**Упражнение 4.7 («нащупывание» цен по Вальрасу и централизованная распределенная оптимизация [17, 276, 336, 338, 340, 341]).** Пусть руководство города владеет $n$ пекарнями. Затраты $i$-й пекарни на выпечку $x_i$ тонн хлеба в день равны $f_i(x_i)$ – сильно выпуклые возраста-

ющие функции. Задача руководства: производить не меньше $C$ тонн хлеба в день ($C$ – объём спроса на хлеб в день со стороны населения города) так, чтобы суммарные затраты всех пекарен были бы минимальны. Формально задача может быть поставлена следующим образом:

$$\sum_{i=1}^{n} f_i(x_i) \to \min_{\substack{\sum_{i=1}^{n} x_i \geq C \\ x_i \geq 0,\, i=1,\ldots,n}} . \qquad (4.33)$$

Обозначим решение этой задачи $x^* = \{x_i^*\}_{i=1}^{n}$.

**1)** Предположим теперь, что у пекарен есть собственники, которые продают хлеб руководству города (распределяющего этот хлеб среди населения) по цене $p^k$ в $k$-й день. Таким образом, собственники решают задачи:

$$x_i(p^k) = \arg\max_{x_i \geq 0} \overbrace{\{\underbrace{p^k x_i}_{\text{выручка}} - \underbrace{f_i(x_i)}_{\text{затраты}}\}}^{\text{прибыль}}, \quad i = 1,\ldots,n . \qquad (4.34)$$

Руководство города действует по следующему правилу: в каждый день $k$ у руководства есть представление о том, в каком отрезке лежит равновесная цена $\left[p_{\min}^k, p_{\max}^k\right]$. Выставив цену $p^k = \frac{1}{2}\left(p_{\min}^k + p_{\max}^k\right)$, руководство собирает следующую информацию с пекарен: $\sum_{i=1}^{n} x_i(p^k)$. Далее

$$\left[p_{\min}^{k+1}, p_{\max}^{k+1}\right] = \left[p_{\min}^k, \frac{1}{2}\left(p_{\min}^k + p_{\max}^k\right)\right], \text{ если } \sum_{i=1}^{n} x_i(p^k) > C ;$$

$$\left[p_{\min}^{k+1}, p_{\max}^{k+1}\right] = \left[\frac{1}{2}\left(p_{\min}^k + p_{\max}^k\right), p_{\max}^k\right], \text{ если } \sum_{i=1}^{n} x_i(p^k) \leq C .$$

Покажите, что

$$x_i(p^k) \xrightarrow[k\to\infty]{} x_i^* .$$

Попробуйте оценить скорость сходимости. Предложите способ оценки $\left[p_{\min}^0, p_{\max}^0\right]$.

**2)** Переписав задачу (4.33) следующим (равносильным) образом:

$$\sum_{i=1}^{n} f_i(x_i) \to \min_{\substack{x_i \ge y_i,\, i=1,\dots,n \\ \sum_{i=1}^{n} y_i \ge C \\ x_i, y_i \ge 0,\, i=1,\dots,n}},$$

попробуйте предложить алгоритм нащупывания равновесной цены, когда каждый день пекарни выставляют свою цену на хлеб, по которой готовы (хотят) продавать хлеб руководству, а руководство закупает хлеб у пекарни, выставившей самую низкую цену. Если таких пекарен, выставивших наименьшую (одинаковую) цену, несколько, то руководство города каким-то (произвольным) образом может распределять закупки среди таких пекарен (и только таких – даже если суммарно эти пекарни производят меньше хлеба, чем нужно руководству).

**Указание.** Заметим, что пекарни не имеют информации о производственных процессах друг друга, а руководство не имеет представление о производственных процессах на всех пекарнях. С одной стороны, описанный выше процесс можно понимать как процесс нащупывания равновесной цены [81, гл. 10]. С другой стороны, этот процесс можно понимать как распределенный централизованный алгоритм решения задачи выпуклой оптимизации (4.33) [436]: задача хранится на $n$ рабочих узлах/пекарнях (slave nodes), взаимодействие которых осуществляется через центр/руководство (master node), см. также рис. 10 (с точностью до $n-1 \to n$). На каждом узле осуществляется работа только со своей частью задачи (решаются вспомогательные задачи (4.34)), коммуникация осуществляется, как показано на рис. 9.

◊ Централизованная распределенная оптимизация очень похожа на параллельные вычисления. Однако принципиально важным отличием от параллельных вычислений является наличие распределенной (оперативной) памяти, в которой можно хранить описание большой задачи, не пытаясь собрать все задачу целиком где-то в одном месте [164]. ◊

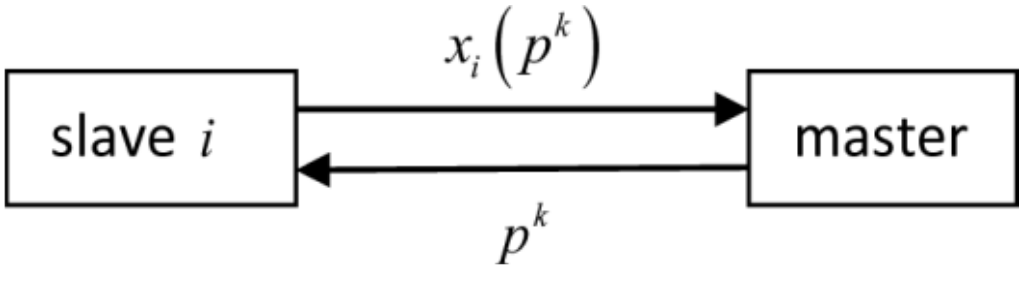

Рис. 9

1) Для решения задачи нужно построить двойственную задачу (с точностью до знака):

$$\psi(p) = \sum_{i=1}^{n} \left( p x_i(p) - f_i(x_i(p)) \right) - Cp \to \min_{p \ge 0}$$

и заметить, что

$$\psi'(p) \overset{def}{=} \frac{d\psi(p)}{dp} = \sum_{i=1}^{n} x_i(p) - C .$$

Из условий задачи $p_{\min}^0 \geq 0$, а $p_{\max}^0$ можно оценить с помощью упражнения 4.1. Далее для решения двойственной задачи можно использовать метод деления отрезка пополам, см. упражнение 1.4.

2) В данном случае двойственная задача (с точностью до знака) будет иметь вид

$$\psi(p) := \psi(p_1, ..., p_n) =$$
$$= \sum_{i=1}^{n} \left( p_i x_i(p_i) - f_i(x_i(p_i)) \right) - C \min_{i=1,...,n} p_i \to \min_{p=(p_1,...,p_n) \in \mathbb{R}_+^n}, \qquad (4.35)$$

где $p_i$ – множитель Лагранжа к ограничению $x_i \geq y_i$. Тогда

$$\partial \psi(p) = \begin{pmatrix} x_1(p) \\ \vdots \\ x_n(p) \end{pmatrix} - C \begin{pmatrix} \lambda_1 \\ \vdots \\ \lambda_n \end{pmatrix}, \begin{pmatrix} \lambda_1 \\ \vdots \\ \lambda_n \end{pmatrix} \in S_n(1), \ \lambda_i > 0 \ \Rightarrow \ i \in \operatorname{Arg} \min_{l=1,...,n} p_l .$$

Двойственная задача (4.35) получилась негладкой, потому что функционал исходной (прямой) задачи не зависел от $\{y_i\}_{i=1}^n$, т. е. не был сильно выпуклым по всей совокупности прямых переменных, см. также указание к упражнению 4.8. Для решения задачи (4.35) можно использовать, например, субградиентный метод [17, 455] или *сходящийся субградиентный метод* Нестерова–Шихмана, который будет иметь в данном контексте вполне естественную интерпретацию [465, 466]. Отличие этого метода от близкого ему субградиентного метода [17], описанного также в упражнении 2.1, в том, что сходимость по функции теперь будет иметь место в обычном (не чезаровском) смысле, поэтому в название метода и вошло слово *сходящийся*. Подумайте, можно ли ускорить процедуры нащупывания равновесия [17, 338, 465, 466], сохранив возможность содержательной интерпретации, если смотреть на слагаемое $C \min_{i=1,...,n} p_i$ в задаче (4.35) как на композитный член (см. пример 3.1). ■

◊ В развитие примера 4.1 и упражнения 4.7 отметим, что для, так называемой, *звездной топологии* коммуникационного графа (см. рис. 10), в которой выделяется главный узел (master node), связанный со всеми остальными $n-1$ узлами (slave nodes), диаметр графа равен 2, а $\sigma_{\max}(W)/\tilde{\sigma}_{\min}(W) = n^2$.

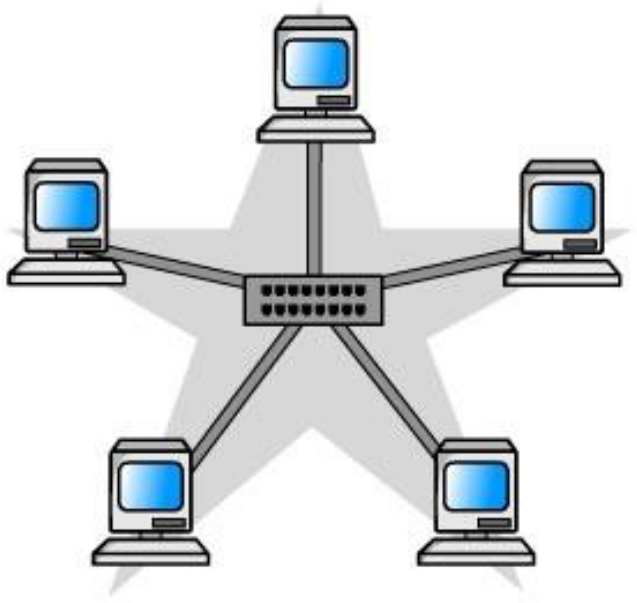

Рис. 10

При этом, если (симметричную) матрицу $W$, которая в данном случае соответствует системе линейных уравнений

$$\begin{aligned} (n-1)y_1 - y_2 - \ldots - y_n &= 0 \\ -y_1 + y_2 &= 0 \\ \ldots\ldots\ldots\ldots\ldots\ldots\ldots\ldots \\ -y_1 + y_n &= 0 \end{aligned},$$

определять из немного другой системы

$$\begin{aligned} \frac{n-1}{n}y_1 - \frac{1}{n}y_2 - \ldots - \frac{1}{n}y_n &= 0 \\ -y_1 + y_2 &= 0 \\ \ldots\ldots\ldots\ldots\ldots\ldots\ldots\ldots \\ -y_1 + y_n &= 0 \end{aligned},$$

то для получившейся в результате (несимметричной) матрицы $W$ выполняется: $\sigma_{\max}(W)/\tilde{\sigma}_{\min}(W) \approx n$. Это утверждение следует из того, что спектр новой матрицы $W$ будет иметь вид $[0;1;1;\ldots;1;2-1/n]$, жорданова нормальная форма матрицы $W$ – диагональная, собственный вектор, отвечающий собственному значению 0 равен $\mathrm{v}_0 = (1,1,\ldots,1)^T$, а максимальному собственному значению – $\mathrm{v}_{\max} = \mathrm{v}_0 - (2-1/n)\underbrace{(1,0,0,..,0)^T}_{e_1}$, собственное подпространство, отвечающее собственному значению 1 ортогонально $\mathrm{v}_0$ и $e_1$, поэтому $\tilde{\sigma}_{\min}(W) = 1$, а $\sigma_{\max}(W) = \max_{\|h\|_2=1}\|Wh\|_2^2 = \|We_1\|_2^2 \approx n$ [91, 291, 556]. Таким образом, отказавшись от симметричности матрицы в рассматриваемом случае можно уменьшить $\sigma_{\max}(W)/\tilde{\sigma}_{\min}(W)$ в $\sim n$ раз, однако это не позволяет ускорить сходимость градиентного спуска, по-скольку в случае несимметричной матрицы $W$ уже нельзя работать, как описано выше, с $\sqrt{W}$ .

Заметим, что в общем случае задача наилучшего подбора системы линейных уравнений, равносильной $y_1 = \ldots = y_n$ (и соответствующей этой системе матрицы $W$ ), тесно связана с задачами выпуклого полуопределенного программирования, рассмотренными в [519]. Другой способ подбора, см. в комментарии к упражнению 5.9 (конструкция предобуславливателя).

Отметим также популярную в последнее время тенденцию (см., например, [240, 241, 417, 515, 553] и замечание 1.5): сведение задачи поиска наилучшего численного метода для рассматриваемого класса задач оптимизации, в свою очередь, к задаче оптимизации (седловой задаче). Если удается решить такую (мета-)задачу, то удается найти наилучший метод (в минимаксном смысле). Однако, как правило, точное и удобное для практического применения решение удается найти в редких случаях (но все же удается, см. замечание 1.5), поэтому интересны и различные релаксации (упрощения) сложной метазадачи. ◊

**Упражнение 4.8 (децентрализованная распределенная оптимизация в сильно выпуклом случае [294, 514]).** Покажите, что если в примере 4.1 дополнительно предположить, что $\varphi_i''(y) \le L_\varphi$, то «двойственная» функция $f(x)$ (4.20) будет сильно выпуклой в 2-норме с константой $\mu_f = \tilde{\sigma}_{\min}(W)/L_\varphi$ на $\left(\mathrm{Ker}\, W^T\right)^\perp = \left(\mathrm{Ker}\, W\right)^\perp$ (при редуцированном подходе $\mu_f = \tilde{\sigma}_{\min}\left(\sqrt{W}\right)/L_\varphi$ и на $\left(\mathrm{Ker}\, \sqrt{W}\right)^\perp$).

**Указание**. См. [247, утверждение 2.1], [456, теорема 1], [357, теорема 6], [413, Lemma 3.1], [434], [503, Theorem 23.5] . Следует обратить внимание на неполную симметричность следующих связей: 1) сильная выпуклость прямой задачи порождает гладкость (липшицевость градиента) двойственной и 2) гладкость в прямой задаче порождает сильную выпуклость двойственной. Во втором случае требуется, чтобы при переходе к двойственной задаче все ограничения с помощью множителей Лагранжа переносились в функционал. В первом случае этого не требуется. Эта несимметричность вместе с теоремой Фенхеля–Моро [62, п. 1.4, 2.2] отчасти объясняет отмеченную ранее асимметрию в возможности адаптивной настройки методов на неизвестную константу Липшица градиента и на отсутствие такой возможности для константы сильной выпуклости. Во всяком случае, пока не придумали, как можно было в общем случае адаптивно осуществлять такую настройку без серьезных дополнительных усилий, см. указание к упражнению 1.3. ■

◊ Отметим, что сильная выпуклость двойственного функционала (4.20) $f(x)=\varphi^*(-Wx)$ имеет место только на $(\mathrm{Ker}\,W)^\perp$, см. упражнение 4.8. Точнее говоря, сильная выпуклость имеет место на любой гиперплоскости вида $x^0+(\mathrm{Ker}\,W)^\perp$. При этом важно заметить, что

$$\nabla f(x)=-W^T\,\nabla\varphi^*(z)\big|_{z=-Wx}\in \mathrm{Im}\left(-W^T\right)=\mathrm{Im}\,W^T=(\mathrm{Ker}\,W)^\perp,$$

т.е. $x^k\in x^0+(\mathrm{Ker}\,W)^\perp$ для любого $k$. Таким образом, траектория градиентного спуска все время будет находиться в той же самой гиперплоскости $x^0+(\mathrm{Ker}\,W)^\perp$, в которой функция $f(x)$ $\mu_f$-сильно выпуклая. Последнее означает, что для градиентного спуска (аналогично и для быстрого градиентного спуска) двойственную задачу можно считать $\mu_f$-сильно выпуклой, «забывая» про вырождение на подпространстве $\mathrm{Ker}\,W$.

Изменяющиеся со временем коммуникационные графы стали достаточно популярными в последнее время в связи с ростом интереса к изучению различных беспроводных систем мобильных объектов с ограниченным радиусом действия антенн. ◊

**Замечание 4.4 (распределенная оптимизация и модели консенсуса).** Заметим, что распределенное решение задач оптимизации можно осуществлять, не переходя к двойственной задаче. Поясним идею прямого подхода. Пусть матрица Лапласа неориентированного связного графа рассматриваемой коммуникационной сети равна $W$, см. пример 4.1. В $i$-м узле графа хранится число $x_i^0$, $i=1,\dots,n$. Требуется, так организовать коммуникацию,[51] чтобы система сошлась к консенсусу [1, 593] за наименьшее число коммуникационных шагов. Поставим в соответствие этой задаче следующую задачу квадратичной оптимизации

$$\frac{1}{2}\langle x,Wx\rangle\to\min_{x\in\mathbb{R}^n}.$$

[51] По условию на одном шаге коммуникации каждый узел может обмениваться информацией только со своими прямыми соседями.

Очевидно, что решением этой задачи будет любой вектор вида $x_* = \mathrm{const}\cdot(1,...,1)^T$. Рассмотрим быстрый градиентный метод из указания к упражнению 1.3, который, подобно примеру 4.1, имеет естественную интерпретацию с точки зрения коммуникации узлов с прямыми соседями:[52]

$$x^1 = x^0 - \frac{1}{\lambda_{\max}(W)} W x^0;$$

$$x^{k+1} = x^k - \frac{1}{\lambda_{\max}(W)} W \cdot \left( x^k + \frac{\sqrt{\lambda_{\max}(W)} - \sqrt{\tilde{\lambda}_{\min}(W)}}{\sqrt{\lambda_{\max}(W)} + \sqrt{\tilde{\lambda}_{\min}(W)}} \left( x^k - x^{k-1} \right) \right) +$$

$$+ \frac{\sqrt{\lambda_{\max}(W)} - \sqrt{\tilde{\lambda}_{\min}(W)}}{\sqrt{\lambda_{\max}(W)} + \sqrt{\tilde{\lambda}_{\min}(W)}} \left( x^k - x^{k-1} \right),$$

где $\tilde{\lambda}_{\min}(W) = \tilde{\sigma}_{\min}\left(\sqrt{W}\right)$ – наименьшее из положительных собственных значений матрицы $W$. Аналогично комментарию к упражнению 4.8 можно показать, что этот метод сходится к такому $x_*$, что $x_* \in x^0 + (\mathrm{Ker}\, W)^{\perp}$, т.е. к *консенсусу* (состоянию, когда все узлы «знают» среднее арифметическое начальных состояний узлов):

$$x_* = \frac{1}{n} \sum_{i=1}^{n} x_i^0 \cdot (1,...,1)^T.$$

Причем для того, чтобы

$$\left\| x^N - x_* \right\|_2 \le \varepsilon \left\| x^0 - x_* \right\|_2,$$

достаточно

---

[52] Заметим, что стандартный метод градиентного спуска (см. § 1) для рассматриваемой задачи можно понимать как метод простой итерации [91, п. 3.4]. Этот метод будет иметь еще более простую интерпретацию.

$$N = \mathrm{O}\left(\sqrt{\chi(W)}\ln\left(\varepsilon^{-1}\right)\right), \text{ где } \chi(W) = \frac{\lambda_{\max}(W)}{\tilde{\lambda}_{\min}(W)}$$

итераций (коммуникационных шагов).[53] Если теперь считать, что вместо $x_i^0$, $i = 1,\ldots,n$ в узлах хранятся (вычисляются) $\nabla\varphi_i(y)$, $i = 1,\ldots,n$, то аналогичным образом за дополнительную мультипликативную плату $\mathrm{O}\left(\sqrt{\chi(W)}\ln\left(\varepsilon^{-1}\right)\right)$, можно добиться того, что все узлы «знают» градиент всего функционала

$$\varphi(y) = \sum_{i=1}^{n}\varphi_i(y) \to \min_y .$$

Таким образом,[54] если каждая функция $\varphi_i(y)$ в исходной задаче оптимизации имеет $L_\varphi$-липшицев градиент и константу сильной выпуклости $\mu_\varphi$, а вместо базового градиентного метода используется (и в прямом и в двойственном подходе) быстрый градиентный метод, то общее число коммуникационных шагов в прямом подходе будет

$$\mathrm{O}\left(\sqrt{nL_\varphi/(n\mu_\varphi)}\ln\left(\varepsilon^{-1}\right)\cdot\sqrt{\chi(W)}\ln\left(\varepsilon^{-1}\right)\right) = \mathrm{O}\left(\sqrt{L_\varphi/\mu_\varphi}\cdot\sqrt{\chi(W)}\ln^2\left(\varepsilon^{-1}\right)\right),$$

где $\varepsilon$ – желаемая относительная точность решения задачи по аргументу, при двойственном подходе:

$$\mathrm{O}\left(\sqrt{L_\varphi/\mu_\varphi\cdot\chi(W)}\ln\left(\varepsilon^{-1}\right)\right),$$

---

[53] Считаем, что $\chi(W) \ll \varepsilon^{-1}$.

[54] На самом деле, эти рассуждения не точны, и приведены здесь лишь для наглядности восприятия материала. В действительности, здесь требуется дополнительно пенализировать исходную постановку задачи квадратичным штрафом, описанным в конце замечания 4.2, но общая связь с консенсусным алгоритмом сохраняется [252, 294, 406, 407]. Впрочем, имеются и прямые рассуждения [506, 584], максимально приближенные к описанной процедуре. Отметим, что проверить корректность доказательств из работы [584] пока не получилось.

см. упражнения 1.3, 4.8. В прямом подходе число вызовов градиентного оракула в каждом узле будет

$$O\left(\sqrt{nL_\varphi/(n\mu_\varphi)}\ln\left(\varepsilon^{-1}\right)\right)=O\left(\sqrt{L_\varphi/\mu_\varphi}\ln\left(\varepsilon^{-1}\right)\right),$$

для двойственного подхода (см. формулу (4.20) с $W\to\sqrt{W}$, см. также ниже в этом замечании) число вызовов двойственного градиентного оракула:

$$O\left(\sqrt{L_\varphi/\mu_\varphi\cdot\chi(W)}\ln\left(\varepsilon^{-1}\right)\right).$$

Стоит отметить, что в прямом подходе возможно дополнительное ускорение за счет следующей простой конструкции [164, 514]. Строим остовное дерево для исходной коммуникационной сети. В принципе, это можно сделать и распределенным образом [279, 514]. Объявляем одну из вершин корнем. Выделение корня порождает иерархию родитель $\to$ потомок. Цель та же самая, что и раньше. На первом шаге коммуникации каждый узел, у которого нет потомков (лист дерева) посылает свое значение (градиент) своему родителю. Родители прибавляют полученные от потомков значения к своему. Система переходит на следующий шаг коммуникации, в котором роль листьев будут играть родители, получившие на прошлом шаге коммуникации значения от своих предков. Процесс повторяется до тех пор, пока вся сумма значений (во всех узлах) не соберется в корневом узле. Корневой узел, базируясь на градиенте всего функционала, может сделать шаг выбранного итерационного метода, и полученное в результате значение(-я) распространить среди потомков. Запустив процесс в обратном порядке. При самом неблагоприятном выборе корня дерева число коммуникационных шагов при таком подходе будет не больше удвоенного диаметра исходного коммуникационного графа. При этом во всех случаях этот диаметр будет не больше, чем корень из числа обусловленности матрицы Лапласа этого графа $\sqrt{\chi(W)}$. Часто это величины одного порядка [514]. Таким образом, появляется возможность «отыграть» потерю в логарифмическом множители в оценке числа коммуникационных шагов при прямом подходе. Более того, для рассмотренной в комментарии к упражнению 4.7 звездной топологии диаметр графа оказывается в $\sim\sqrt{n}$ меньше корня из числа обусловленности. То есть дополнительный выигрыш от второго подхода может быть и существенным.[55]

Отсюда, однако, не стоит делать вывод о том, что прямой подход всегда лучше двойственного. Во-первых, стоимость прямого и двойственного градиентного оракула может быть существенно разной [256, 561], и выбор в пользу одного из подходов следует осуществлять, учитывая и стоимость вызова соответствующего оракула. Действительно, в задаче вычисления барицентра Монжа–Канторовича–Васерштейна (МКВ) [488] $n$ вероятностных дискретных мер (с носителем на $d$ точках) вместо настоящего расстояния МКВ считают $\mu$-энтропийно регуляризованное расстояние, где $\mu=\varepsilon/\left(2n\ln d^2\right)$, см. замечание 4.1, поэтому каждое слагаемое в сумме (4.17) будет иметь вид:

$$\varphi_l(y)=\varphi_{w^l}(y)=\min_{\substack{\sum_{j=1}^{d}\pi_{ij}=y_i,\ \sum_{i=1}^{d}\pi_{ij}=w_j^l\\ \pi_{ij}\ge 0,\ i,j=1,\dots,d}}\left\{\sum_{i,j=1}^{d}c_{ij}\pi_{ij}+\mu\sum_{i,j=1}^{d}\pi_{ij}\ln\pi_{ij}\right\}=$$

$$=\max_{x\in\mathbb{R}^d}\min_{\substack{\sum_{i=1}^{d}\pi_{ij}=w_j^l\\ \pi_{ij}\ge 0,\ i,j=1,\dots,d}}\left\{\sum_{i,j=1}^{d}c_{ij}\pi_{ij}+\sum_{i=1}^{d}x_i\cdot\left(y_i-\sum_{j=1}^{d}\pi_{ij}\right)+\mu\sum_{i,j=1}^{d}\pi_{ij}\ln\pi_{ij}\right\}=$$

$$=\max_{x\in\mathbb{R}^d}\left\{\langle y,x\rangle-\underbrace{\mu\sum_{j=1}^{d}w_j^l\ln\left(\frac{1}{w_j^l}\sum_{i=1}^{d}\exp\left(\frac{-c_{ij}+x_i}{\mu}\right)\right)}_{\varphi_l^*(x)=\varphi_{w^l}^*(x)}\right\},$$

где $l=1,\dots,n$, $y\in S_d(1)$. Достаточно точный подсчет $\varphi_{w^l}(y)$ и $\nabla\varphi_{w^l}(y)$ при малых значениях $\varepsilon$ требует $\tilde{O}\left(d^3\right)$ арифметических операций [125, 277, 486, 548] (см. также указание к упражнению 1.4 и указание к упражнению 3.9, в котором отмечается, что при не малых значениях $\varepsilon$ суще-ствуют более эффективные способы), в то время как подсчет $\varphi^*_{w^l}(y)$ и $\nabla\varphi^*_{w^l}(y)$ – всего лишь $O(d^2)$. Отметим также, что двойственная задача – задача безусловной оптимизации, а прямая задача – условной, поскольку есть ограничение $y \in S_d(1)$. Однако наличие ограничения простой структуры не скажется на возможности применения прямого подхода.

[55] К сожалению, второй подход, требующий вычисления остовного дерева, плохо подходит для работы на изменяющихся со временем графах.

Во-вторых, для негладких, но сильно выпуклых прямых задач[56] прямой подход требует

$$\tilde{O}\left(M^2/(n\mu\cdot\varepsilon)\cdot\sqrt{\chi(W)}\right)$$

коммуникационных шагов и

$$O\left(M^2/(n\mu\cdot\varepsilon)\right)$$

вызов прямого градиентного оракула в каждом узле (см. упражнение 2.3 и табл. 2 в приложении), где $M$ – константа Липшица прямого функционала $\varphi(y)=\sum_{l=1}^{n}\varphi_{w^l}(y)$, $y\in\tilde{Q}\subseteq\mathbb{R}^d$. В двойственном подходе к задаче

$$\tilde{\varphi}\left(Y=(y_1,\ldots,y_n)\right)=\sum_{l=1}^{n}\varphi_{w^l}(y_l)\to\min_{\substack{\sqrt{W}Y=0\\ y_l\in\tilde{Q},\, l=1,\ldots,n}}$$

строится (с точностью до знака) двойственный функционал, который необходимо минимизировать (см. (4.20) с $W\to\sqrt{W}$ )

$$f(x)=\sum_{l=1}^{n}\varphi_l^*\left(\left[-\sqrt{W}x\right]_l\right)\to\min_{x\in\mathbb{R}^{nd}},$$

---

[56] Такой задачей, как раз, будет задача вычисления $\mu$-энтропийно регуляризованного барицентра МКВ – $\varphi_l(y)$ будет $\mu$-сильно выпуклой функцией в 2-норме.

где в выписанных формулах $\sqrt{W} := \sqrt{W} \otimes 1_d 1_d^T$, а $\otimes$ – кронекерово произведение матриц [387]. Число коммуникационных шагов и число вызовов оракула, выдающего градиент соответствующей (рассматриваемому узлу) двойственной функции, будут совпадать (это общее положение для двойственного оракула, см., например, оценки выше в гладком сильно выпуклом случае) и равны

$$O\left(\sqrt{\frac{L_f R_x^2}{\varepsilon}}\right) = O\left(\sqrt{\frac{\left(\sigma_{\max}\left(\sqrt{W}\right)/\mu\right)\left(\left\|\nabla\tilde{\varphi}(Y_*)\right\|_2^2/\tilde{\sigma}_{\min}\left(\sqrt{W}\right)\right)}{\varepsilon}}\right) =$$

$$= O\left(\sqrt{\frac{\left\|\nabla\tilde{\varphi}(Y_*)\right\|_2^2}{\mu\varepsilon}\cdot\chi(W)}\right),$$

см. указание к упражнению 1.3, замечание 4.2, пример 4.1 и табл. 2 в приложении. Для того чтобы сравнить приведенные оценки заметим, что

$$M^2 = \max_{y\in\tilde{Q}\cap B_{2R}\left(y^0\right)}\left\|\nabla\varphi(y)\right\|_2^2 \text{ и } \left\|\nabla\varphi(y_*)\right\|_2^2 \le n\left\|\nabla\tilde{\varphi}(Y_*)\right\|_2^2.$$

Из этих оценок можно сделать, в общем случае, довольно грубый вывод о том, что

$$M^2/n \approx \left\|\nabla\tilde{\varphi}(Y_*)\right\|_2^2.$$

В таком приближении имеем: в двойственном подходе оценки на число коммуникаций и вызовов оракула (в каждом узле) получаются извлечением корня из оценок в прямом подходе.

Отметим при этом, что на основе двойственного подхода можно предложить такой способ решения исходной задачи распределенным образом с прямым оракулом, т.е. выдающим градиент прямых функций (в соответствующих узлах), при котором число коммуникационных шагов будет как в двойственном подходе. Для этого необходимо использовать двойственный подход, а для вычисления градиента двойственной функ-ции (соответствующей рассматриваемому узлу) использовать прямой оракул. А именно, чтобы с нужной точностью вычислить

$$\nabla\varphi_l^*\left(\left[-\sqrt{W}x\right]_l\right)=-\sqrt{W}y_l\left(x\right)$$

нужно с помощью прямого оракула достаточно точно решить задачу

$$y_l\left(x\right)=\arg\max_{y_l\in\tilde{Q}}\left\{\left\langle\left[-\sqrt{W}x\right]_l,y_l\right\rangle-\varphi_l\left(y_l\right)\right\}.$$

К сожалению, при таком подходе без дополнительных предположений[57] не удается получить ожидаемые (из прямого подхода) оценки на число вызовов прямого оракула.[58] Оказывается, что можно предложить такой распределенный способ решения исходной задачи [294, 387], который обеспечивает оптимальное число коммуникаций и, одновременно, ожидаемое (из прямого подхода) число вызовов оракула в общем случае, причем как в случае выпуклого, так и в, рассмотренном выше, случае сильно выпуклого (прямого) функционала. Для этого стоит обратиться к конструкции, изложенной в комментарии к замечанию 4.3 (стоит иметь в виду, что обозначения сильно отличаются). В этой конструкции в выпуклом случае описывается способ получения оптимальных оценок на число коммуникационных шагов и «ожидаемых» оценок на число вызовов пря-мого оракула. С помощью рестартов (см. указание к упражнении 2.3 и конец § 5) можно перенести эти результаты и на сильно выпуклый случай, получив анонсированный результат.

[57] Например, таких, что для регуляризованной двойственной функции доступен градиент сопряженной (двойственной) к ней функции, см. замечание 4.1 в части двойственного сглаживания. То есть доступен не прямой оракул, а двойственно-сглаженный прямой оракул.

[58] Здесь не говорится об оптимальном числе вызовов прямого оракула (на узле), потому что функционал имеет определенную структуру – вида суммы. В замечании 1 приложения (см. также [387]) обсуждаются оценки сложности оптимальных (нераспределенных) методов для данного класса задач. По-видимому, при заданном (оптимальном) числе коммуникационных шагов, обсуждаемые «ожидаемые» оценки также являются оптимальными [252, 584]. Однако на данный момент, для задач распределенной выпуклой оптимизации строго установлена оптимальность приводимых здесь оценок только в части оценок на число коммуникационных шагов [136, 513, 514], да и то в предположении «враждебного» распределения частей функции по имеющимся узлам коммуникационного графа. Отметим работы [323–326, 374, 394, 575, 580], в которых обсуждаются вопросы возможности распространения ускорения из замечания 1 приложения на распределенные (параллельные) алгоритмы.

Заметим, что оптимальные оценки на число коммуникационных шагов можно получать из разобранных выше двух случаев просто за счет регуляризации исходной (прямой) задачи, см. замечание 4.1. При этом двойственный оракул должен будет выдавать градиент не сопряженной функции к прямой, а сопряженной функции к регуляризованной прямой. Отметим также, что в гладком случае «ожидаемые» оценки (с точностью до логарифмического множителя) на число вызовов прямого оракула получаются после регуляризации из разобранного в самом начале случая (гладкого и сильно выпуклого), см. также [294, 406]. ■

◊ Отметим, что приведенные выше результаты частично переносятся и на меняющиеся со временем графы [506]. Можно отказаться и от неориентированности коммуникационного графа / симметричности матрицы $W$, сохраняя при этом его связность. В этом случае, используя другую технику, также можно получить результаты, частично похожие на те, что были приведены в пособии, см., например, [435] и цитированную там литературу. К сожалению, даже если рассматривать только ориентированные, не изменяющиеся со временем коммуникационные графы, то на данный момент неизвестно, можно ли (а если можно, то каким образом) ускорить сходимость, как это было сделано в случае неориентированных графов [294, 514]. Впрочем, существует достаточно общий способ ускорения, в том числе, и распределенных централизованных алгоритмов – каталист (см. замечание 3.3 и приложение), в котором предлагается решать (внутреннюю) вспомогательную (сильно выпуклую) задачу рассматриваемым неускоренным распределенным централизованным методом. В частности, таким образом, можно например, ускорить асинхронный централизованный алгоритм из работы [424]. ◊

**Упражнение 4.9 (схема регуляризации А. Н. Тихонова [6, лекции 5, 7]).** Пусть необходимо решить систему линейных уравнений $Ax = b$ *в условиях истокопредставимости*:

$$x_* = A^T y_*, \ \|y_*\|_2 \le R.$$

Как будет видно из указания условие истокопредставимости можно понимать также как (см. замечание 1.6):

$$x_* = \left(A^T A\right)^{1/2} \tilde{y}_*, \ \left\|\tilde{y}_*\right\|_2 \le R.$$

Результат не изменится.

Предположим, что точное значение правой части $b$ неизвестно. Зато известно значение $\tilde{b}$ такое, что $\left\|\tilde{b} - b\right\|_2 \le \delta$. Схема А.Н. Тихонова предлагает искать решение системы $Ax = b$ путем решения задачи оптимизации

$$\frac{1}{2}\left\|Ax - \tilde{b}\right\|_2^2 + \frac{\mu}{2}\left\|x\right\|_2^2 \to \min_{x \in \mathbb{R}^n}.$$

Считая, что известно точное решение этой задачи

$$x_*^\mu\left(\tilde{b}\right) = \left(A^T A + \mu I\right)^{-1} A^T \tilde{b},$$

предложите способ выбора параметра $\mu(\delta)$ так чтобы в оценке

$$\left\|x_* - x_*^\mu\left(\tilde{b}\right)\right\|_2 \le \sigma(\delta),$$

функция $\sigma(\delta) \ge 0$ была как можно меньше.

**Указание.** По неравенству треугольника

$$\left\|x_* - x_*^\mu\left(\tilde{b}\right)\right\|_2 \le \left\|x_* - x_*^\mu(b)\right\|_2 + \left\|x_*^\mu(b) - x_*^\mu\left(\tilde{b}\right)\right\|_2.$$

Заметим, что

$$x_* - x_*^\mu(b) = \left(A^T A + \mu I\right)^{-1}\left(A^T A + \mu I\right)x_* - \left(A^T A + \mu I\right)^{-1} \underbrace{A^T b}_{A^T A x_*} =$$

$$= \mu\left(A^T A + \mu I\right)^{-1} x_* = \mu \underbrace{\left(A^T A + \mu I\right)^{-1} A^T}_{\Theta(A,\mu)} y_*,$$

$$x_*^{\mu}(b) - x_*^{\mu}(\tilde{b}) = \left(A^T A + \mu I\right)^{-1} A^T \left(b - \tilde{b}\right) = \Theta(A,\mu)\left(b - \tilde{b}\right).$$

Таким образом, необходимо оценить норму оператора $\Theta(A,\mu)$:

$$\left\|\Theta(A,\mu)\right\|_2 = \sqrt{\sup_{\|z\|_2 = 1} \left\|\Theta(A,\mu) z\right\|_2^2}.$$

Пусть $\{e_k\}_k$ – ортонормированный базис пространства из собственных векторов самосопряженного (симметричного) оператора (матрицы) $A^T A$. Тогда по спектральной теореме [58]

$$\left\|\Theta(A,\mu) z\right\|_2^2 = \sum_k \frac{1}{(\lambda_k + \mu)^2} \left\langle e_k, A^T z \right\rangle^2 = \sum_{k:\,\lambda_k > 0} \frac{\lambda_k}{(\lambda_k + \mu)^2} \left\langle \frac{1}{\sqrt{\lambda_k}} A e_k, z \right\rangle^2.$$

Поскольку

$$\left\{ \frac{1}{\sqrt{\lambda_k}} A e_k \right\}_{k:\,\lambda_k > 0}$$

образуют ортонормированную систему (но не обязательно базис), то по неравенству Бесселя имеем

$$\sum_k \left\langle \frac{1}{\sqrt{\lambda_k}} A e_k, z \right\rangle^2 \le \|z\|_2^2.$$

Следовательно,

$$\left\|\Theta(A,\mu) z\right\|_2^2 \le \sup_k \frac{\lambda_k}{(\lambda_k + \mu)^2} \|z\|_2^2.$$

Значит,

$$\left\|\Theta(A,\mu)\right\|_2 \le \sup_k \frac{\sqrt{\lambda_k}}{\lambda_k + \mu} \le \sup_{\lambda \ge 0} \frac{\sqrt{\lambda}}{\lambda + \mu} = \frac{1}{2\sqrt{\mu}}.$$

Таким образом,

$$\left\|x_* - x_*^{\mu}\left(\tilde{b}\right)\right\|_2 \le \left\|x_* - x_*^{\mu}(b)\right\|_2 + \left\|x_*^{\mu}(b) - x_*^{\mu}\left(\tilde{b}\right)\right\|_2 \le \frac{R\sqrt{\mu}}{2} + \frac{\delta}{2\sqrt{\mu}}.$$

Правая часть последнего неравенства минимальная при выборе $\mu(\delta) = \delta/R$. В таком случае

$$\left\|x_* - x_*^{\mu}\left(\tilde{b}\right)\right\|_2 \le \sqrt{R\delta} = \sigma(\delta).$$

К сожалению, на практике $R$ обычно неизвестно. Есть сложности и с предположением об известности $\delta$ (вето Бакушинского [317]). Альтернативным и часто более эффективным способом решения такого рода задач является регуляризация за счет правильного выбора численного метода решения нерегуляризованной задачи

$$\frac{1}{2}\left\|\tilde{A}x - \tilde{b}\right\|_2^2 \to \min_{x\in\mathbb{R}^n},$$

обладающего регуляризирующими свойствами, см. замечание 1.6. ■

# § 5. Универсальный градиентный спуск

Как и в § 2 – § 4, рассмотрим общую задачу выпуклой оптимизации (2.1):

$$f(x) \to \min_{x \in Q}.$$

В данном параграфе, следуя Ю.Е. Нестерову [462] (см. также [7, 18, 30, 95, 285, 463, 585]), будет сделан, пожалуй, самый важный шаг во всем описанном выше подходе – согласованы формулы (2.5), (2.26), (2.27). Как уже ранее отмечалось, для наглядности рассуждения будут проводиться не в максимальной общности (см. § 3).

Прежде всего заметим, что во всех вариантах рассмотренных на данный момент методов градиентного спуска использовался шаг $h = 1/L$, где константа $L$ либо была задана по условию, либо определялась согласно (2.5) с $\delta = \varepsilon/2$ (см., вывод формул (4.6), (4.28)). Рассмотрим следующее (универсальное) обобщение метода (2.28) (описывается $k$-я итерация):

**Универсальный градиентный спуск**

$$L^{k+1} = L^k/2,$$

**While True Do**

$$x^{k+1} = \arg\min_{x \in Q} \left\{ f\left(x^k\right) + \left\langle \nabla f\left(x^k\right), x - x^k \right\rangle + L^{k+1} V\left(x, x^k\right) \right\}.$$

**If** $\left\{ f\left(x^{k+1}\right) \le f\left(x^k\right) + \left\langle \nabla f\left(x^k\right), x^{k+1} - x^k \right\rangle + L^{k+1} V\left(x^{k+1}, x^k\right) + \frac{\varepsilon}{2} \right\}$

Перейти на следующую итерацию: $k \to k+1$

**Else**

$$L^{k+1} := 2L^{k+1}.$$

Для такого метода формула (4.1) перепишется следующим образом:

$$\frac{1}{L^{k+1}} f\left(x^{k+1}\right) \le \frac{1}{L^{k+1}} \left\{ f\left(x^k\right) + \left\langle \nabla f\left(x^k\right), x - x^k \right\rangle \right\} + V\left(x, x^k\right) - V\left(x, x^{k+1}\right) + \frac{\varepsilon}{2L^{k+1}}, \tag{5.1}$$

где константа $L^{k+1}$ подбирается в (5.1) согласно описанной выше процедуре. Причем согласно (2.32), эта константа оценивается сверху константой $L^{k+1}$, подбираемой из соотношения (2.26), в котором $\delta = \varepsilon/2$:

$$f\left(x^{k+1}\right) \le f\left(x^k\right) + \left\langle \nabla f\left(x^k\right), x^{k+1} - x^k \right\rangle + \frac{L^{k+1}}{2}\left\|x^{k+1} - x^k\right\|^2 + \frac{\varepsilon}{2}.$$

То есть автоматически происходит подбор на рассматриваемом отрезке $\left[x^k, x^{k+1}\right]$ двух параметров $\nu$ и $L_\nu$ в (2.27):

$$\left\|\nabla f(y) - \nabla f(x)\right\|_* \le L_\nu \left\|y - x\right\|^\nu, \nu \in [0,1], L_0 < \infty$$

так, чтобы (2.26) выполнилось. Согласно (2.5) и описанной выше процедуре подбора $L^{k+1}$ выполняется, что

$$L^{k+1} \le L_\nu \cdot \left[\frac{L_\nu}{\varepsilon}\frac{1-\nu}{1+\nu}\right]^{\frac{1-\nu}{1+\nu}}.$$

Подчеркнем, что не мы сами решаем задачу подбора $\nu$ и $L_\nu$ – это делает метод за счет описанной процедуры. Свойства гладкости функции $f(x)$ на отрезке $\left[x^k, x^{k+1}\right]$ характеризуются континуальным набором чисел $\{L_\nu\}_{\nu\in[0,1]}$, часть из которых может равняться бесконечности. Мы можем ничего не знать о $\{L_\nu\}_{\nu\in[0,1]}$ – этого и не требуется при универсальном подходе. Тем не менее описанная выше процедура гарантирует, что метод подберет такое $\nu \in [0,1]$, что соответствующая этому $\nu$ константа Гёльдера $L_\nu$ порождает (по формуле (2.26) с $\delta = \varepsilon/2$) на отрезке $\left[x^k, x^{k+1}\right]$ минимально возможную (с точностью до множителя не больше 2) константу $L^k$, которая явно используется в методе. Подобно оценкам (4.6), (4.28) для универсального градиентного спуска можно получить следующий результат.

**Теорема 5.1.** *Пусть нужно решить задачу* (2.1) *в условиях* (2.27). *Для универсального градиентного спуска после*[59]

---

[59] Причем в (5.2) константы $L_\nu$ выбираются не худшие по всему пространству, а «средние» на траектории метода. При получении оценки (5.2) предполагалось, что $L^0$ выбиралось так, что при подстановке $k = 0$ в условие «**If**» рассматриваемого универсального метода неравенство в этом условии будет неверным при $L^{k+1} = L^1 < L^0$. На самом деле, если предположение о выборе $L^0$ не имеет место, то это приведет лишь к небольшому ухудшению числового множителя в оценке (5.2).

$$N = \inf_{\nu \in [0,1]} \left( \frac{2L_\nu R^{1+\nu}}{\varepsilon} \right)^{\frac{2}{1+\nu}} \qquad (5.2)$$

*итераций имеет место следующая оценка*:

$$f\left(\bar{x}^N\right) - f\left(x_*\right) \le f\left(\bar{x}^N\right) -$$

$$- \frac{1}{\sum_{k=0}^{N-1} 1/L^{k+1}} \min_{x \in B_{R,Q}\left(x^0\right)} \left\{ \sum_{k=0}^{N-1} \frac{1}{L^{k+1}} \left[ f\left(x^k\right) + \left\langle \nabla f\left(x^k\right), x - x^k \right\rangle \right] \right\} \le \varepsilon,$$

*где* (*в данном случае*)

$$\bar{x}^N = \frac{1}{\sum_{k=1}^{N} 1/L^k} \sum_{k=1}^{N} \frac{x^k}{L^k},$$

$R^2 = V\left(x_*, x^0\right)$. *Если решение* $x_*$ *не единственно, то оценка* (5.2) *справедлива и для того решения* $x_*$, *которое доставляет минимум* $R^2$.

Еще раз подчеркнем, что задачу минимизации, возникающую в оценке (5.2), нам не нужно решать, равно как и знать хоть что-то о гладкости $f\left(x\right)$, кроме того, что $L_0 < \infty$, чтобы универсальный метод сходился. В негладком случае (когда только $L_0 < \infty$) имеем $N = 4L_0^2 R^2 / \varepsilon^2$, что с точностью до множителя соответствует нижней оценке, см. (2.37).

Рассуждая аналогично [462], несложно показать, что описанный выше универсальный метод на каждой итерации запрашивает один раз $\nabla f\left(x^k\right)$ и в среднем (по итерациям) около трех раз значение функции $f\left(x\right)$. Действительно, среднее число вычислений значения функции $f\left(x\right)$ есть

$$\frac{1}{N} \sum_{k=0}^{N-1} \left( 2 + \log_2 \left( L^{k+1} / \left( L^k / 2 \right) \right) \right) = 3 + \frac{1}{N} \log_2 \left( L^N / L^0 \right).$$

**Замечание 5.1 (универсальный метод для вариационных неравенств и седловых задач [24, 263]).**[60] Аналогом градиентного метода для

[60] В концепции модели функции из примера 1 § 3 на базе этого замечания можно получить универсальный вариант популярного сейчас *метода Chambolle–Pock'a* [197, 421, 541]. Отметим, что в работах [419, 421] предложен альтернативный способ (forward–backward–forward) построения адаптивных методов решения вариационных неравенств. Этот способ используют только одну проекцию на

*вариационных неравенств* (ВН) и седловых задач является *экстраградиентный метод* Г. М. Корпелевич [13, § 15, гл. 5], [60, 142]. Интересные результаты по решению ВН можно найти в работах А. С. Антипина [131], см. также обзор в [144]. Далее рассмотрим один современный вариант экстраградиентного метода, а именно *проксимальный зеркальный метод* А. С. Немировского [441]. Пусть задано векторное поле $g(x)$, в частности $g(x)=\nabla f(x)$. Предположим, что существуют такие $L$ и $\delta$, что для всех $x$, $y$, $z$ из выпуклого компактного множества $Q$ имеет место неравенство

$$\langle g(y)-g(x), y-z\rangle \le LV(y,x)+LV(y,z)+\delta .$$

Тогда для проксимального зеркального метода

$$y^{k+1}=\arg\min_{x\in Q}\left\{\left\langle g\left(x^k\right), x-x^k\right\rangle + LV\left(x,x^k\right)\right\},$$

$$x^{k+1}=\arg\min_{x\in Q}\left\{\left\langle g\left(y^{k+1}\right), x-x^k\right\rangle + LV\left(x,x^k\right)\right\}$$

имеет место следующая оценка:

$$\frac{1}{N}\sum_{k=1}^{N}\left\langle g\left(y^k\right), y^k-x\right\rangle \le \frac{LV\left(x,x^0\right)-LV\left(x,x^N\right)}{N}+\delta. \qquad (5.3)$$

С помощью текста, написанного в конце п. 4.6 [180], несложно построить универсальный вариант такого метода (описывается $k$-я итерация):

**Универсальный проксимальный зеркальный метод**

$$L^{k+1}=L^k/2,$$

**While True Do**

$$y^{k+1}=\arg\min_{x\in Q}\left\{\left\langle g\left(x^k\right), x-x^k\right\rangle + L^{k+1}V\left(x,x^k\right)\right\},$$

$$x^{k+1}=\arg\min_{x\in Q}\left\{\left\langle g\left(y^{k+1}\right), x-x^k\right\rangle + L^{k+1}V\left(x,x^k\right)\right\}.$$

**If** $\left\{\left\langle g\left(y^{k+1}\right)-g\left(x^k\right), y^{k+1}-x^{k+1}\right\rangle \le L^{k+1}V\left(y^{k+1},x^k\right)+L^{k+1}V\left(y^{k+1},x^{k+1}\right)+\frac{\varepsilon}{2}\right\}$

Перейти на следующую итерацию: $k\to k+1$

**Else**

$$L^{k+1}:=2L^{k+1}.$$

---

каждой итерации при этом имеет аналогичные оценки скорости сходимости. Похожий на описанный здесь алгоритм для более общего класса задач (монотонных включений) описан в работе [229].

Если, подобно (2.27), векторное поле $g(x)$ удовлетворяет условию

$$\|g(y)-g(x)\|_* \le L_\nu \|y-x\|^\nu, \nu\in[0,1],\ x,y\in Q,\ L_0<\infty,$$

то, используя неравенство (верное для любых $a,b,L_\nu,\delta>0$, $\nu\in[0,1]$)

$$L_\nu a^\nu b \le L_\nu \cdot \left(\frac{L_\nu}{\delta}\right)^{\frac{1-\nu}{1+\nu}} \left(\frac{a^2}{2}+\frac{b^2}{2}\right)+\delta$$

с $\delta=\varepsilon/2$, можно получить (см. [24]), что для достижения

$$\frac{1}{\sum_{k=1}^{N} 1/L^k} \max_{x\in Q}\left\{\sum_{k=1}^{N}\frac{1}{L^k}\left\langle g\left(y^k\right), y^k-x\right\rangle\right\}\le\varepsilon \tag{5.4}$$

достаточно (следует сравнить с (5.2))

$$N = \inf_{\nu\in[0,1]}\left(\frac{2L_\nu R^{1+\nu}}{\varepsilon}\right)^{\frac{2}{1+\nu}} \tag{5.5}$$

итераций универсального проксимального зеркального метода. Здесь $R^2=\max_{x\in Q} V\left(x,x^0\right)$. При этом среднее число вычислений значений векторного поля $g(x)$ на одной итерации приближенно равно трем. Оценка (5.5) с точностью до числового множителя оптимальна для ВН и для седловых задач при[61] $N\le \dim x$. К сожалению, точной ссылки на обоснование оптимальности не удалось найти, однако различные частные случаи могут быть сведены к разобранному в работах [70, 311, 438, 439, 480].

◊ Следует сравнить (5.4) при $g(x)=\nabla f(x)$ с (2.25), (3.4). Если $g(x)=\left(\nabla_u f(u,w),-\nabla_w f(u,w)\right)$, $x=(u,w)$, $Q=Q_u\otimes Q_w$, где функции $f(u,w)$ выпуклая по $u$ и вогнутая по $w$, то из (5.4) следует, что

$$0\le \max_{w\in Q_w} f\left(\bar{u}^N,w\right)-\min_{u\in Q_u} f\left(u,\bar{w}^N\right)\le\varepsilon.$$

Отметим при этом, что для седловой точки $(u_*,w_*)$ имеем

---

[61] Если это условие не выполняется, то решение монотонного ВН можно свести к задаче негладкой выпуклой оптимизации, которую можно решать методами типа центров тяжести [46, 70], см. указание к упражнению 1.4. Отметим также, что метод эллипсоидов можно применять для решения ВН и непосредственно [443].

$$\max_{w\in Q_w} f(u_*,w)=\min_{u\in Q_u} f(u,w_*). \ \Diamond$$

Заметим, что для *монотонных вариационных неравенств*[62]

$$\langle g(y)-g(x),y-x\rangle\ge 0,\ x,y\in Q,$$

имеем

$$\langle g(x),y^k-x\rangle=\langle g(y^k),y^k-x\rangle+\underbrace{\langle g(x)-g(y^k),y^k-x\rangle}_{\le 0}\le\langle g(y^k),y^k-x\rangle.$$

В этой связи формулу (5.4) можно переписать как[63]

$$\max_{x\in Q}\langle g(x),\bar{y}^N-x\rangle\le\varepsilon, \tag{5.6}$$

где

$$\bar{y}^N=\frac{1}{\sum_{k=1}^{N}1/L^k}\sum_{k=1}^{N}\frac{y^k}{L^k}.$$

В литературе обычно используют именно этот критерий качества решения монотонных ВН, см., например, [71, гл. 3], [441].

Подобно слабой квазивыпуклости из замечания 2.1 можно ослабить условия монотонности ВН, во многом сохранив результаты [215]. ■

**Замечание 5.2 (негладкие задачи и рандомизированные методы).** Как уже отмечалось, основным достоинством универсального подхода является автоматическая и адаптивная настройка на гладкость задачи. И даже если задача заведомо негладкая, универсальный подход может давать существенные преимущества по сравнению с оптимальными мето-дами, настроенными на негладкие задачи, см., например, [7]. Однако у универсального подхода есть несколько минусов. Во-первых, это не адаптивный подход в том смысле, что в метод явным образом зашита желаемая точность $\varepsilon$ (к чему это приводит, см., например, в (2.34)).[64] Отказавшись от универсальности, можно избавиться и от этого ограничения, используя *метод двойственных усреднений* Ю. Е. Нестерова [455].[65] Во-вторых, обоснование метода требует возможности проведения выкладок хотя бы в общности § 2, однако, в действительности, для негладких задач достаточно общности (1.32), что заметно упрощает и вывод основных оценок, и обоснование возможности последующего обобщения на стохастические постановки [92]. Отметим тем не менее, что здесь речь идет только о простоте выводов, но не о потенциальных возможностях вывода. В-третьих, привязка к (1.32) позволяет переносить результаты, полученные непосредственно для негладких выпуклых задач, т. е. без универсализации, на онлайн-постановки, в том числе стохастические и сильно выпуклые[66], см., например, [319, 478]. В-четвертых, для негладких задач концепция неточного оракула (см., например, (2.3), (3.1)) может быть заменена на более простую и менее ограничительную концепцию $\delta$-субградиента (см., например, [82, п. 5 § 1 и п. 3 § 3, гл. 5]), в которой отсутствуют правые неравенства в (2.3), (3.1).[67] В-пятых, при перенесении универсальных методов на стохастические постановки задач [30], в которых случайность искусственно ввели мы сами (это, как правило, называется *рандомизацией* метода) при вычислении градиента или проектировании, чтобы сократить вычислительную сложность этих операций, взамен на увеличение их числа, из-за погони за универсальностью могут теряться некоторые свойства дешевизны этих операций [18]. Связано это, прежде всего, с тем, что при универсальном подходе необходимо рассчитывать значение функции, что может быть намного дороже расчета значения ее стохастического градиента. Вот простой пример [351]:

$$f(x) = \frac{1}{2}\langle x, Ax \rangle,$$

где $A \succ 0$ – плотно заполненная неотрицательно определенная матрица $n \times n$, $x \in S_n(1)$. Несмещенный стохастический градиент этой функции

$$\nabla_x f(x, j) = A^j, \text{ где } P(j = i) = x_i,\ i = 1,\dots,n.$$

---

[62] В случае $g(x)=\nabla f(x)$ это условие соответствует выпуклости $f(x)$.

[63] Собственно под слабым решением вариационного неравенства обычно понимают следующую задачу [144]: найти такую точку $x_*\in Q$, что для всех $x\in Q$ $\langle g(x),x_*-x\rangle\le 0$, т. е. $\max_{x\in Q}\langle g(x),x_*-x\rangle\le 0$. Отсюда становится ясным смысл условия (5.6). Подробнее о ВН см., например, [8, 34], [46, Ч. 3], [71, гл. 3]. В частности, полезно заметить, что ВН можно решать и в сильном смысле [144]: найти такую точку $x_*\in Q$, что для всех $x\in Q$ $\langle g(x_*),x_*-x\rangle\le 0$, т.е. $\max_{x\in Q}\langle g(x_*),x_*-x\rangle\le 0$. Если векторное поле $g(x)$ – непрерывное и монотонное (ВН монотонное), то слабые и сильные решения ВН совпадают [46, часть 3, Лемма 2.2.1] и можно говорить просто о решении ВН. В общем же случае сильное решение всегда будет слабым, но не наоборот.

---

[64] Впрочем, адаптивность по желаемой точности $\varepsilon$ можно получить, рестартуя универсальный метод с новым значением точности $\varepsilon := \varepsilon/2$ в момент, когда достигнута точность $\varepsilon$ и т.д.

[65] Этот метод близок другому популярному методу решения негладких задач оптимизации – *методу зеркального спуска* А. С. Немировского, см., например, [122], [158, гл. 5], [350] и упражнение 2.6.

[66] Отметим, что конструкция рестартов (упражнение 2.3) в онлайн-постановках уже не работает. Более того, для сильно выпуклых задач онлайн-оптимизации оценка (2.37) уже не достижима. Необходима ее корректировка в части $L_0^2/(\mu N) \to L_0^2 \ln N/(\mu N)$. Такая нижняя оценка уже будет достижима [319, 321].

[67] Работая с $\delta$-субградиентами где $\delta = O(\varepsilon/N(\varepsilon))$ (по-видимому, эту оценку можно улучшить до $\delta \sim \varepsilon$), вместо настоящих субградиентов можно получать, например, оценки из указания к упражнению 1.4, изначально установленные для работы с настоящими субградиентами [70, 130], см. также упражнение 5.5.

Ясно, что

$$E_j\left[\nabla_x f(x, j)\right] = Ax = \nabla f(x).$$

Также ясно, что на подсчет $f(x)$ уходит время $O\left(n^2\right)$, в то время как на подсчет $\nabla_x f(x, j)$ – время $O(n)$. Впрочем, если для решения задачи используется *минибатчинг,* что более характерно для решения задач стохастической оптимизации, чем при рандомизации детерминированных процедур (см. приложение), то такой проблемы не возникает. ■

Вернемся к оценке (5.1). Попробуем с помощью этой оценки и *техники рестартов* (упражнение 2.3) распространить описанный выше универсальный градиентный спуск на сильно выпуклые задачи. Ввиду теоремы 1.1 отметим, что существуют и другие способы того, как это можно сделать. Однако выбранный здесь способ представляется наиболее удачным в методическом плане своей общеприменимостью.

Итак, подобно выводу (4.5) из (4.1), из (5.1) можно получить

$$f\left(\bar{x}^N\right) - f(x_*) \le \frac{\bar{L}_N V\left(x_*, x^0\right)}{2N} + \frac{\varepsilon}{2}, \tag{5.7}$$

где

$$\bar{L}_N = \frac{N}{\sum_{k=1}^{N} 1/L^k},\quad \bar{x}^N = \frac{\bar{L}_N}{N}\sum_{k=1}^{N}\frac{x^k}{L^k}.$$

Пусть $f(x)$ – $\mu$-сильно выпуклая функция в норме $\|\ \|$, согласованной с дивергенцией Брэгмана $V(y, x)$ (см. § 2, в частности, формулу (2.32)). Пусть $d\left(x - x^0\right) \le C_n \left\|x - x^0\right\|^2$ (можно считать $C_n = O(\ln n)$, см. п. 2) упражнения 2.3). Тогда из (1.14) и (5.7) следует, что

$$\frac{\mu}{2}\left\|\bar{x}^{N_1}-x_*\right\|^2 \le f\left(\bar{x}^{N_1}\right)-f\left(x_*\right) \le \frac{\bar{L}_{N_1}V\left(x_*,x^0\right)}{2N_1}+\frac{\varepsilon}{2} \le \frac{C_n\bar{L}_{N_1}\left\|x^0-x_*\right\|^2}{2N_1}+\frac{\varepsilon}{2}.$$

Далее используется схема рассуждений, аналогичная той, что была изложена в указании к упражнению 2.3. А именно, из соотношения

$$\frac{\mu}{2}\left\|\bar{x}^{N_1}-x_*\right\|^2 \le \frac{C_n\bar{L}_{N_1}\left\|x^0-x_*\right\|^2}{2N_1}+\frac{\varepsilon}{2} \tag{5.8}$$

выбираем наименьшее такое $N_1$, при котором

$$\left\|\bar{x}^{N_1}-x_*\right\|^2 \le \frac{1}{2}\left\|x^0-x_*\right\|^2. \tag{5.9}$$

Для этого не надо знать $x_*$, можно просто воспользоваться соотношением (5.8). Однако при этом необходимо знать $\mu$. В итоге можно показать, что для такого метода[68] из оценки (5.2) получится оценка вида [22, 27, 265]:

$$N = \mathrm{O}\left(\underbrace{C_n \inf_{\nu\in[0,1]}\left(\frac{L_\nu^2}{\mu^{1+\nu}\varepsilon^{1-\nu}}\right)^{\frac{1}{1+\nu}}}_{\text{число итераций на одном рестарте}}\underbrace{\left\lceil \ln\left(\frac{\mu\left\|x_*-x^0\right\|^2}{\varepsilon}\right)\right\rceil}_{\text{число рестартов}}\right), \tag{5.10}$$

где здесь и далее $\lceil a \rceil = \max\{1,a\}$. Заметим, что оценки (1.24), (2.37) соответствуют (5.10) с точностью до логарифмического множителя при $\nu=1$, $\nu=0$ соответственно.

Формулу (5.10) можно проверить с помощью регуляризации (см. замечание 4.1). А именно, исходную выпуклую задачу всегда можно сделать сильно выпуклой с константой сильной выпуклости $\mu \simeq \varepsilon/R^2$. Подставляя $\mu \simeq \varepsilon/R^2$ в (5.10) с точностью до $C_n$, получим оценку (5.2).

---

[68] На $(k+1)$-м рестарте следует выбирать точку старта как $x^0=\bar{x}^{N_k}$ – среднее арифметическое точек, полученных по ходу работы метода на $k$-м рестарте, а $d(x):=d\left(x-x^0\right)=d\left(x-\bar{x}^{N_k}\right)$ [30, 351, 354]. Натуральный логарифм в (5.10) выбран потому, что именно так наиболее часто выбирают в данном контексте в литературе. Отметим также, что подписи сомножителей в (5.10) (и далее в аналогичных формулах) приводятся с точностью до числовых множителей.

◊ В работе [354] (см. также [30, 71]) отмечается, что число обусловленности (отношение константы Липшица градиента $L^1$ к константе сильной выпуклости $\mu^1$), например, для квадратичных функций

$$f(x)=\frac{1}{2}\langle x,Ax\rangle-\langle b,x\rangle,\ x\in\mathbb{R}^n,$$

посчитанное, в 1-норме не может быть меньше $n$. В то время как число обусловленности, посчитанное в евклидовой норме (2-норме), может при этом равняться 1. Действительно, пусть $\xi=(\xi_1,....,\xi_n)$, где $\xi_k$ – независимые одинаково распределенные случайные величины

$$P(\xi_k=1/n)=P(\xi_k=-1/n)=1/2.$$

Тогда, учитывая, что $\|\xi\|_1\equiv 1$, получим

$$\mu^1(f)\le E_\xi\left[\xi^T A\xi\right]=\frac{1}{n^2}\operatorname{tr}(A)\le\frac{1}{n}\max_{i,j=1,...,n}|A_{ij}|=\frac{1}{n}L^1(f).$$

Отсюда напрашивается вывод, что при решении сильно выпуклых задач естественно выбирать евклидову норму. Как правило, это действительно так. Однако, во-первых, отмеченное выше наблюдение совсем не означает, что обусловленность задачи в 1-норме всегда хуже, чем в 2-норме. Например, числа обусловленности энтропии $f(x)=\sum_{i=1}^n x_i\ln x_i$ на множестве $\left\{x\in\mathbb{R}^n:\sum_{i=1}^n x_i=1;\ x_i\ge\delta,\ i=1,...,n\right\}$ в 1-норме и 2-норме совпадают и равны $1/\delta\ge n$. Во-вторых, для задач композитной оптимизации [446] (см. также пример 3.1) $L^1$ можно брать только от гладкого слагаемого, а $\mu^1$ от всей функции. Таким образом можно решать, например, задачу [27]

$$f(x)=\frac{1}{2}\|Ax-b\|_2^2+\mu\sum_{i=1}^n x_i\ln x_i\to\min_{x\in S_n(1)}.\ ◊$$

Основная проблема с реализацией описанного выше подхода – это явное использование в нем, как правило, неизвестного параметра $\mu$. Довольно естественный способ борьбы с неизвестностью параметра $\mu$ состоит в рестартах по невязке $f(x^N)-f(x_*)$. Если $f(x_*)$ известно, то рестарты можно делать, контролируя эту невязку по функции [265, 281, 507]. Отметим, что если известно $f(x_*)$, то можно также предложить адаптивный по $\mu$ (ускоренный) градиентный метод, не требующий рестартов [146].

**Замечание 5.3 (контроль нормы градиента).** Для неускоренного (универсального) градиентного спуска также можно использовать рестарты по норме градиента (градиентного отображения – см. упражнение 3.4).[69] Ограничимся далее для наглядности задачей безусловной выпуклой оптимизации (1.1) в условиях (1.4). Тогда для метода (1.22)

$$x^{k+1} = x^k - \frac{1}{L}\nabla f\left(x^k\right)$$

имеет место следующая не улучшаемая для данного метода оценка [552] (следует сравнить с оценкой (1.23) – оптимальной для класса гладких невыпуклых задач):

$$\left\|\nabla f\left(x^N\right)\right\|_2 \le \frac{LR}{N+1},$$

где $R = \left\|x^0 - x_*\right\|_2$. Несложно показать, что рестарты (на базе формул (1.14), (1.15)) с данной оценкой приводят к правильным порядкам скорости сходимости (в том числе по функции и по аргументу) для неускоренных методов в сильно выпуклом случае (1.24) [552]. Ситуация меняется для ускоренных методов, см. указание к упражнению 1.3. Для ускоренных (быстрых, моментных) градиентных спусков на данный момент удается получить (см. (1.8)) лишь оценки вида [371, 447, 552]:[70]

$$\left\|\nabla f\left(x^N\right)\right\|_2 \le \sqrt{2L\cdot\left(f\left(x^N\right) - f\left(x_*\right)\right)} \le \frac{2LR}{N+1},\ \min_{k=0,\dots,N}\left\|\nabla f\left(x^k\right)\right\|_2 \le \frac{8LR}{N^{3/2}},$$

которые не улучшаемы больше, чем на числовой множитель для рассмотренных на данный момент классов (ускоренных) методов первого порядка [552]. До недавнего времени открытым оставался вопрос, можно ли добиться для какого-нибудь из методов первого порядка сходимости вида $\left\|\nabla f\left(x^N\right)\right\|_2 \sim LR/N^2$ [447, 552]. Для такого метода можно было бы делать рестарты по норме градиента без риска потерять оптимальность. Тем не менее если регуляризовать исходную задачу (см. замечание 4.1) и решать регуляризованную задачу любым вариантом быстрого градиентного метода, уже настроенного на сильно выпуклую постановку [30, 71, 72, 180, 228], то с точностью до $\ln N$ получается желаемая оценка.[71] К сожалению, для всех известных сейчас вариантов таких методов в размер шага явно входит константа сильной выпуклости, неизвестная в данном контексте, см. указание к упражнению 1.3. К тому же не совсем понятно, какой смысл бороться здесь за оптимальную оценку, чтобы предложить на ее основе с помощью рестартов оптимальный метод для гладких сильно выпуклых задач, если получается, что вся эта конструкция, в свою очередь, сама базируется на таком методе.

---

[69] Следует сопоставить с тем, что ранее отмечалось в комментариях к упражнению 2.3 относительно невозможности использовать такой критерий для рестартов в общем случае (например, для ускоренных методов).

[70] Последняя оценка получается для сочетания сначала $N/2$ шагов быстрого градиентного метода, затем $N/2$ шагов обычного градиентного метода [447].

В 2018 г. в работе [372] был предложен такой вариант быстрого градиентного метода, для которого

$$\left\|\nabla f\left(x^N\right)\right\|_2^2 = \mathrm{O}\left(\frac{L\cdot\left(f\left(x^0\right)-f\left(x_*\right)\right)}{N^2}\right).$$

Излагаемая далее конструкция (два абзаца текста) была нам сообщена Ю.Е. Нестеровым. Если сначала запустить обычный быстрый градиентный метод на $N$ итераций, то получим

$$f\left(x^N\right)-f\left(x_*\right)=\mathrm{O}\left(\frac{LR^2}{N^2}\right).$$

Запустив затем метод из работы [372], получим

$$\left\|\nabla f\left(x^{2N}\right)\right\|_2^2 = \mathrm{O}\left(\frac{L\cdot\left(f\left(x^N\right)-f\left(x_*\right)\right)}{N^2}\right)=\mathrm{O}\left(\frac{L^2R^2}{N^4}\right).$$

Приведенные оценки хорошо согласуются с результатами [439] (см. также замечание 1.6), и в общем случае не улучшаемы.

Отметим, однако, одну сложность, не позволяющую для метода из работы [372] использовать норму градиента в качестве критерия останова. Методу из работы [372] требуется давать на вход в качестве одного из параметров число итераций $N$, которые метод должен сделать. То есть приведенная выше оценка скорости сходимости с параметром $N$ не предполагает, что метод не зависит от этого параметра. Для каждого значения $N$, получается свой метод. Тем не менее, такой метод можно использовать с конструкцией рестартов при известном значении $L$, но не известном значении $\mu$, более эффективно, чем обычный быстрый гради-ентный метод. Для этого нужно брать в качестве критерия окончания рестарта условие уполовинивания нормы градиента на этом рестарте, а в качестве числа итераций на $k$ -м рестарте: $N^k = \sqrt{L/\mu^k}$ , где $\mu^k = 2\mu^{k-1}$ , и далее подбираем параметр $\mu^k$ с возможным понижением: $\mu^k := \mu^{k-1}/2$ , до тех пор, пока не выполнится критерий окончания рестарта.

---

[71] Упражнения 4.2 и данное предложение проясняют необходимость использования регуляризации при решении двойственной задачи ускоренными методами (см. упражнение 4.4) в качестве альтернативы к использованию прямодвойственных методов. Двойственная задача регуляризуется, чтобы оптимально (с точностью до логарифмических множителей) восстанавливать по найденному приближенному решению двойственной задачи решения прямой задачи.

Отметим, что с помощью техники, аналогичной технике работы [372] (см. также замечание 1.5 и последующий поясняющий текст) для выпукло-вогнутых седловых задач (см. замечание 5.2) с $L$ -липшицевым градиентом

$$\min_u \max_w f(u,w)$$

пока можно получить лишь оценку на число вычислений $\nabla_u f(u,w)$ и $\nabla_w f(u,w)$ [369]

$$\left\|\nabla f\left(u^N, w^N\right)\right\|_2^2 = \tilde{O}\left(\frac{L^2R^2}{N^2}\right).$$

Можно ли здесь убрать логарифмический множитель (заменить $\tilde{O}(\ )$ на $O(\ )$)? Ответ на этот вопрос, насколько мы знаем, по-прежнему является открытым. Отметим также, что конструкция упомянутой работы [369], по-видимому, позволяет в некотором смысле перенести идею каталиста на седловые задачи. ■

Для задач безусловной выпуклой оптимизации замечание 5.3 также дает возможность контролировать только норму градиента в качестве критерия останова неускоренного градиентного спуска и его универсального варианта.

Отметим также, что недавно была обнаружена возможность (не связанная с использованием прямодвойственности) избавления от знания значения параметра $\mu$ и в рестартах для ускоренных вариантов градиентных методов [270].

В заключение обратим внимание, что для задач оптимизации на неограниченных множествах не могут одновременно выполняться оба используемых выше неравенства: $\mu > 0$ и $L_\nu < \infty$, где $\nu \in [0,1)$. Однако, в виду (3.16) можно быть уверенным в том, что последовательность точек, сгенерированных методом универсального градиентного спуска с рестартами, будет лежать в компактном множестве, задаваемом (3.16), на котором и стоит определять константы $\mu > 0$ и $L_\nu < \infty$, входящие в итоговую оценку числа итераций (5.10).

**Упражнение 5.1.** Предложите универсальный градиентный спуск в общности § 3, т. е. используя общую концепцию модели функции в точке. Покажите, что оценки (5.2), (5.5) сохранят свой вид, если допускать неточности $\delta = O(\varepsilon)$ и $\tilde{\delta} = O(\varepsilon)$ (см. обозначения § 3). Поясните, как следует понимать эти неточности для подхода из замечания 5.1, приводящего к оценке (5.5).

**Упражнение 5.2.** Попробуйте сформулировать и доказать утверждение аналогичное утверждению из упражнения 2.2 для универсального градиентного спуска.

**Упражнение 5.3.** Попробуйте распространить упражнение 5.1 на сильно выпуклые постановки задач.

**Указание.** См. [541]. Оптимизируемая (целевая) функция $f(x)$, как и раньше, предполагается выпуклой на выпуклом множестве $Q$. Условие гладкости и сильной выпуклости $f(x)$ при наличии шума следует понимать следующим образом (3.1) [228, 257]:

$$\frac{\mu}{2}\|y-x\|^2 \le f(y)-\left(f_\delta(x)+\psi_\delta(y,x)\right) \le \frac{L}{2}\|y-x\|^2+\delta.$$

Впрочем, используя технику рестартов, можно исходить и из обычной концепции модели функции (3.1), предполагая сильную выпуклость у модели функции (3.1) $\psi_\delta(y,x)$ как функции $y$ или (еще более общий случай) предполагая сильную выпуклость только у $f(x)$ [30, 265]. ■

**Упражнение 5.4.** С помощью техники рестартов (см. упражнение 2.3) и формулы (5.3) (или как-то по-другому [87]) попробуйте перенести результаты замечания 5.1 на *сильно монотонные вариационные неравенства* и сильно выпукло/вогнутые седловые задачи [71, гл. 3]. Отметим одно затрудняющее такой перенос обстоятельство: в формуле (5.5) используется $R^2 = \max_{x\in Q} V(x,x^0)$, что не дает возможности формально применять технику рестартов. Можно ли получить для универсального проксимального зеркального метода, использующегося для решения ВН и седловых задач оценки, подобные (2.19)?

**Указание.** Для наглядности будем считать $V(x,y)=\|x-y\|_2^2/2$. Обозначим через $x_*$ решение ВН. Заметим, что для всех $y^k \in Q$

$$\left\langle g(x_*), y^k - x_* \right\rangle \ge 0 .$$

Из формулы (5.3) следует, что

$$\frac{1}{N}\sum_{k=1}^{N}\left\langle g\left(y^k\right), y^k - x_*\right\rangle \le \frac{L\left\|x_* - x^0\right\|_2^2}{2N}.$$

Объединяя эти два неравенства, получим

$$\frac{1}{N}\sum_{k=1}^{N}\left\langle g\left(y^k\right) - g\left(x_*\right), y^k - x_*\right\rangle \le \frac{L\left\|x_* - x^0\right\|_2^2}{2N}.$$

По предположению ВН – сильно монотонное, т.е. существует такое $\mu > 0$, что

$$\left\langle g(y) - g(x), y - x\right\rangle \ge \mu\left\|y - x\right\|_2^2,\ x, y \in Q.$$

Учитывая это и выпуклость функции $\left\|x\right\|_2^2$, получим

$$\mu\left\|\bar{y}^N - x_*\right\|_2^2 \le \mu\sum_{k=1}^{N}\frac{1}{N}\left\|y^k - x_*\right\|_2^2 \le \frac{1}{N}\sum_{k=1}^{N}\left\langle g\left(y^k\right) - g\left(x_*\right), y^k - x_*\right\rangle \le \frac{L\left\|x_* - x^0\right\|_2^2}{2N},$$

где

$$\bar{y}^N = \frac{1}{N}\sum_{k=1}^{N} y^k.$$

На основе последнего неравенства уже можно организовать процедуру рестартов. Получив, следующую оценку на общее (суммарное) число итераций

$$\mathrm{O}\left(\frac{L}{\mu}\ln\left(\frac{\mu R^2}{\varepsilon}\right)\right).$$

Данная оценка уже не может быть улучшена для рассматриваемого класса ВН никакими другими методами [480, 586]. Заметим также, что приведенную оценку можно получить и в модельной общности [541, 542]. ■

**Упражнение 5.5.** Используя принцип множителей Лагранжа [176, гл. 5], распространите проксимальный зеркальный метод из замечания 5.1 на общие задачи выпуклого программирования [62, п. 2], предварительно компактифицировав двойственные переменные (см., например, упражнение 4.1 и замечание 4.2). Рассмотрите альтернативный подход к решению возникшей седловой задачи в случае, когда общее число аффинных ограничений вида равенств и выпуклых ограничений вида неравенств мало. В качестве альтернативного подхода предлагается решать двойственную задачу прямодвойственным методом эллипсоидов [211, 443]. Для этого потребуется число итераций (см. указание к упражнению 1.4) $N_{ellips}(\varepsilon) \sim \ln\left(\varepsilon^{-1}\right)$. Причем на каждой итерации вместо настоящего субградиента можно использовать $\delta$-субградиент (см. замечание 5.2), где $\delta \sim \varepsilon / N_{ellips}(\varepsilon)$, который вычисляется по формуле Демьянова–Данскина из решения с относительной точностью (по функции) $\delta$ вспомогательной задачи выпуклой оптимизации, получающейся из рассматриваемой седловой задачи при фиксации двойственных переменных [82, п. 5 § 1, гл. 5]. Проведите описанные выше рассуждения более аккуратно, прорабатывая детали.

**Указание.** Следует сопоставить данное упражнение с упражнениями 4.3, 4.6, 5.6 и примером 3.2. ■

**Упражнение 5.6.** Предложите способ решения задачи минимизации достаточно гладкой выпуклой функции при наличии, вообще говоря, негладкого скалярного сильно выпуклого ограничения вида неравенства. Причем решение задачи доставляет в этом неравенстве равенство.

**Указание.** Следует воспользоваться упражнением 5.5. При этом на каждой итерации метода в двойственном пространстве на вспомогательную задачу оптимизации следует смотреть как на задачу композитной сильно выпуклой оптимизации (см. пример 3.1), чтобы негладкость ограничения не учитывалась, а его сильная выпуклость, напротив, позволяла решать вспомогательную задачу за линейное время, используя, например, технику рестартов, см. концовку § 5 и [27].

Отметим, что если в условии задачи имеется несколько сильно выпуклых ограничений вида неравенств $h_1(x) \le 0$, …, $h_m(x) \le 0$, то их можно заменить скалярным негладким сильно выпуклым ограничением: $h(x) = \max\{h_1(x), \dots, h_m(x)\} \le 0$ – см. [3], [82, п. 3 § 3, гл. 10], [460, 467] и замечание 3.1.

Заметим, что если функция, задающая сильно выпуклое ограничение, гладкая (имеет липшицев градиент), то нет необходимости ее считать композитно-дружественной. При этом на вспомогательную задачу можно смотреть, как на задачу обычной (не композитной) оптимизации. ■

**Упражнение 5.7.** Предложите способ распространения проксимального зеркального метода из замечания 5.1 на бесконечномерные задачи, например, дифференциальные игры [267].

**Упражнение 5.8.** Определите, какие из результатов, описанных выше (во всем пособии), могут быть перенесены с обычного (неускоренного) градиентного метода на ускоренные (быстрые, моментные) градиентные методы?

**Указание.** На метод подобных треугольников из упражнения 3.7 переносятся все результаты, кроме результатов, собранных в замечаниях 1.1, 1.3, 5.1, 5.3 и кроме результатов, под которые лучше подходит ускоренный проксимальный метод Монтейро–Свайтера или его аналоги, см. замечания 1.6, 3.3. В случае замечания 5.1 частичное ускорение при некоторых дополнительных предположениях оказывается возможным [201, 202, 247].

Отметим, что результаты, связанные с относительной гладкостью из § 3 переносятся на ускоренные методы лишь при дополнительных обременительных предположениях на дивергенцию Брэгмана [238, 313, 315, 327, 415]. По-видимому, аналогичные сложности не позволяют предложить полноценный проксимальный ускоренный градиентный метод с неевклидовым проксом, см. замечание 3.3 и текст сразу после него.

Результаты, касающиеся $\alpha$-слабой квазивыпуклости переносятся на ускоренные методы [308, 329, 463]. ■

**Упражнение 5.9 (Ю. Е. Нестеров, 2014).** Задача поиска такого $x_*$, что $Ax_* = b$ сводится к задаче выпуклой гладкой оптимизации:

$$f(x) = \|Ax - b\|_2^2 \to \min_{x\in\mathbb{R}^n}.$$

Нижняя оценка при $N \le n$ на скорость решения такой задачи (см. также упражнение 1.3 и замечание 1.6) имеет вид

$$\|Ax^N - b\|_2^2 \ge \frac{L_x R_x^2}{2(2N+1)^2},$$

где $L_x = \sigma_{\max}(A)$, $R_x = \|x_*\|_2$. Если решение $x_*$ не единственное, то в определении $R_x$ можно считать, что используется решение с наименьшей 2-нормой. При этом на каждой итерации разрешено не более двух раз умножать матрицу $A$ на вектор (справа и слева). С другой стороны, рассмотрим задачу

$$\frac{1}{2}\|x\|_2^2 \to \min_{Ax=b}.$$

Построим к ней двойственную задачу [176, гл. 5]:

$$\min_{Ax=b}\frac{1}{2}\|x\|_2^2 = \min_x \max_\lambda \left\{\frac{1}{2}\|x\|_2^2 + \langle b - Ax, \lambda\rangle\right\} = \max_\lambda \min_x \left\{\frac{1}{2}\|x\|_2^2 + \langle b - Ax, \lambda\rangle\right\} =$$

$$= \max_\lambda \left\{\langle b, \lambda\rangle - \frac{1}{2}\|A^T\lambda\|_2^2\right\}.$$

С помощью теоремы 4.1 и упражнений 1.3, 5.8 (альтернативный способ базируется на замечании 5.3) покажите, что после $N$ итераций ускоренного градиентного метода, примененного к двойственной задаче, можно восстановить решение исходной задачи $\breve{x}^N$ со следующей точностью:

$$\|A\breve{x}^N - b\|_2 \le \frac{8L_\lambda R_\lambda}{N^2},$$

где $L_\lambda = \sigma_{\max}\left(A^T\right) = \sigma_{\max}\left(A\right)$, $R_\lambda = \left\|\lambda_*\right\|_2$. Если решение $\lambda_*$ не единственное, то в определении $R_\lambda$ можно считать, что используется решение с наименьшей 2-нормой. При этом общее число умножений матрицы $A$ на вектор не будет превышать $2N$. Поясните, почему последняя оценка не противоречит выписанной ранее нижней оценке.

**Указание.** См. [4, 23, 137, 281, 500] и замечание 1.6. При этом важно заметить, что связь оптимальных значений в прямой и двойственной задаче $x_* = A^T\lambda_*$ (следует из связи прямых и двойственных переменных $x(\lambda) = A^T\lambda$) и условие $R_\lambda = \left\|\lambda_*\right\|_2$ можно понимать в совокупности, как *условие истокопредставимости*, см. также упражнение 4.9.

Заметим также, что поскольку система $Ax = b$ совместна, то по *теореме Фредгольма* [46, п. 2.6. Ч. 1] не существует такого $\lambda$, что $A^T\lambda = 0$ и $\langle b, \lambda\rangle > 0$, следовательно, двойственная задача имеет конечное решение, т. е. существует ограниченное решение двойственной задачи $\lambda_*$. Действительно, по предположению существует такой $x$, что $Ax = b$, поэтому для всех $\lambda$ имеет место: $\langle Ax, \lambda\rangle = \langle b, \lambda\rangle$. Следовательно, $\left\langle x, A^T\lambda\right\rangle = \langle b, \lambda\rangle$. Предположив, что существует такой $\lambda$, что $A^T\lambda = 0$ и $\langle b, \lambda\rangle > 0$, придем к противоречию: $0 = \left\langle x, A^T\lambda\right\rangle = \langle b, \lambda\rangle > 0$. ■

◊ В замечании 1.5 отмечалось, что решение системы линейных уравнений $Ax = b$ является краеугольным камнем не только (вычислительной) линейной алгебры [291, 477, 511, 556], вычислительной математики [91], но и численных методов оптимизации [177]. Напомним также, что класс сложности задач выпуклой безусловной оптимизации в категориях $\mathrm{O}(\ )$ характеризуется классом сложности задач квадратичной оптимизации (см. упражнение 1.3 и замечание 1.6), и что необходимость в решении системы линейных уравнений (обращении матрицы) возникает на каждом шаге метода Ньютона (см. приложение).

В свою очередь задачи квадратичной оптимизации получаются из $Ax = b$, обычно, либо как указано в упражнении 5.9, либо (в случае симметричности матрицы $A$) согласно (1.30), см. также замечание 1.6. Упражнения 1.6, 5.9 (см. также [18, гл. 4]) показывают, что в случае дополнительных предположений можно пытаться решать (разреженные) линейные системы быстрее, чем предписывают нижние оценки, полученные из нижних оценок для задач квадратичной оптимизации. Предполагая вы-полненными некоторые свойства матрицы $A$,[72] также можно решать (используя рандомизированные методы) системы линейных уравнений за время, пропорциональное $m = nnz(A)$ (с точностью до больших степеней логарифмических множителей, зависящих от желаемой точности) – числу ненулевых элементов в матрице $A$ [18, гл. 4], [205, 365, 501, 535]. Отметим, что в основе части отмеченных работ лежит полезное наблюдение [535, 564]: для матрицы Лапласа можно построить за $\tilde{O}(m)$ другую матрицу Лапласа (на основе остовного дерева неориентированного графа исходной матрицы с добавлением еще нескольких ребер – полученный граф называют *ультра-спарсификатором*), которая:

1) хорошо предобуславливает (см., например, [477]) исходную матрицу $A$; обобщенное число обусловленности[73] [160] становится $\tilde{O}(1)$;

2) содержит $\tilde{O}(m)$ ненулевых элементов;

3) эффективно обратима.

Далее итеративно сочетая отмеченное предобуславливание задачи с итерациями метода сопряженных градиентов (см. замечание 1.6) можно получить отмеченный выше результат. Заметим, что метод сопряженных градиентов, гарантированно сходящийся к точному решению в общем случае (без предобуславливания) за $n$ итераций, на каждой своей итерации требует умножения матрицы $A$ на вектор. Таким образом, сложность только одной такой итерации пропорциональна числу ненулевых элементов матрицы $A$. Интересно сравнить отмеченную (как правило, завышенную) оценку сложности метода сопряженных градиентов $O(mn)$ с оцен-кой, основанной на быстром матричном умножении $\tilde{O}\left(n^{2.37}\right)$ (см. также конец приложения) и оценкой $\tilde{O}(m)$ [535]. Тем не менее, в смысле универсальности, легкой имплементируемости и реальной скорости работы на практике, обычный метод сопряженных градиентов является, как правило, более предпочтительным. На данный момент отмеченная технология [535] представляет интерес в большей степени с теоретической точки зрения, в виду возникновения больших степеней у возникающих логарифмических множителей в оценках.

---

[72] Например, симметричность и слабое диагональное доминирование. Это свойство выполняется, в частности, для матрицы Лапласа неориентированного графа, см. пример 4.1.

[73] Чтобы приблизительно понять, чем обобщенное число обусловленности отличается от обычного числа обусловленности, вернемся к комментарию к упражнению 4.7, в котором рассматривалась матрица Лапласа графа со звездной топологией. Было показано, что если эту матрицу умножить на специальную диагональную матрицу, то полученная в результате (уже не симметричная) матрица будет иметь неотрицательный действительный спектр с отношением максимального собственного значения к минимальному ненулевому равным приблизительно 2, в то время как отношение соответствующих (максимального и минимального ненулевого) сингулярных чисел этой матрицы будет иметь порядок $n$. Обобщенное число обусловленности здесь 2, а не $n$.

В связи с изучением задач децентрализованной распределенной оптимизации (см. пример 4.1) описанная выше матрица Лапласа ультраспарсификатора исходного графа (сети) позволяет разрабатывать эффективные конструкции разреженных и хорошо обусловленных коммуникационных сетей. ◊

# Приложение. Обзор современного состояния развития численных методов выпуклой оптимизации

Описанные в § 2 – § 5 конструкции переносятся на ускоренные (быстрые, моментные) градиентные методы [95], например, на *метод подобных треугольников*[74] [122], см. также упражнения 3.7, 5.8. Причем дальнейшее ускорение в общем случае уже невозможно (см. упражнение 1.3). Для ускоренного метода оценки (5.2), (5.10) и условия на допустимый уровень шума $\delta$ преобразуются следующим образом [22, 27, 68, 95, 265, 462]:

$$N(\varepsilon) = \mathrm{O}\left(\inf_{\nu\in[0,1]}\left(\frac{L_\nu R^{1+\nu}}{\varepsilon}\right)^{\frac{2}{1+\nu}}\right) \to \mathrm{O}\left(\inf_{\nu\in[0,1]}\left(\frac{L_\nu R^{1+\nu}}{\varepsilon}\right)^{\frac{2}{1+3\nu}}\right),$$

$$N(\varepsilon) = \mathrm{O}\left(\underbrace{C_n \inf_{\nu\in[0,1]}\left(\frac{L_\nu^2}{\mu^{1+\nu}\varepsilon^{1-\nu}}\right)^{\frac{1}{1+\nu}}}_{\substack{\text{число итераций}\\ \text{на одном рестарте}}}\underbrace{\left[\ln\left(\frac{\mu\|x_* - x^0\|^2}{\varepsilon}\right)\right]}_{\text{число рестартов}}\right) \to$$

$$\to \mathrm{O}\left(\underbrace{C_n \inf_{\nu\in[0,1]}\left(\frac{L_\nu^2}{\mu^{1+\nu}\varepsilon^{1-\nu}}\right)^{\frac{1}{1+3\nu}}}_{\substack{\text{число итераций}\\ \text{на одном рестарте}}}\underbrace{\left[\ln\left(\frac{\mu\|x_* - x^0\|^2}{\varepsilon}\right)\right]}_{\text{число рестартов}}\right),$$

$$\delta = \mathrm{O}(\varepsilon) \to \tilde{\mathrm{O}}\left(\frac{\varepsilon}{N(\varepsilon)}\right).$$

---

[74] Заметим, что у метода линейного каплинга (МЛК) [122] за счёт наличия двух проектирований на каждой итерации имеются некоторые дополнительные свойства (по сравнению с методом подобных треугольников, у которого одно проектирование), обнаруженные недавно [16, 258, 259, 307, 309, 463]. К сожалению, пока не удалось предложить такой вариант МЛК, который мог бы работать с моделью функции из § 3, но при этом обладал бы отмеченными выше дополнительными свойствами. Отметим также, что у МЛК в варианте работы [122] при $Q = \mathbb{R}^n$ Grad-шаг лучше заменить на Mirr-шаг, чтобы в случае неевклидовой прокс-структуры гарантировать ввиду формулы (3.16) равномерную ограниченность последовательности, генерируемой методом [4]. Другими словами, при $Q = \mathbb{R}^n$ (1.35), (1.36) стоит заменять на (2.29) с соответствующим выбором размеров шагов.

Данные оценки являются неулучшаемыми (оптимальными) [75, приложение А.С. Немировского], [225, 311, 438].

◊ В свою очередь, эти оценки можно обобщить на так называемые *промежуточные методы* [225, гл. 6], [257], методы, которые представляют собой выпуклые комбинации неускоренного и ускоренного градиентного метода: в выписанных формулах вместо $\left[\frac{1}{1+\nu};\frac{1}{1+3\nu}\right]$ стоит писать $\frac{1}{1+\nu+2p\nu}$, при этом $\delta = \tilde{O}\left(\varepsilon / N(\varepsilon)^p\right)$, $p\in[0,1]$. Такого рода обобщения могут потребоваться, например, при решении задач оптимизации в гильбертовых пространствах [281]. Детали и дальнейшее обобщение на случай неточного проектирования (см. (3.3) и упражнение 3.7) имеется в работе [265]. ◊

Для большей наглядности далее (если не оговорено противного) рассматривается задача выпуклой безусловной оптимизации:

$$f(x)\to\min_{x\in\mathbb{R}^n}.$$

В качестве нормы выбирается 2-норма. В качестве прокс-функции: $d(x)=\frac{1}{2}\|x\|_2^2$, см. § 2.

◊ Тем не менее стоит отметить, что все написанное далее с точностью до логарифмических множителей (см. константу $C_n$ в упражнении 2.3 и в конце § 5) переносится и на задачи выпуклой оптимизации на множествах простой структуры. Исключением являются неполноградиентные методы для гладких задач выпуклой оптимизации. На данный момент для таких методов не удалось полностью перенести основные известные сейчас результаты для безусловных гладких задач на гладкие задачи оптимизации на множествах простой структуры [16, 23, 25, 258, 259, 468].

Так же, как и в основном тексте пособия, далее можно считать, что все константы, характеризующие оптимизируемую функцию, относятся не ко всему пространству, а только к шару с центром в точке старта и радиусом, равным (с точностью до логарифмического множителя) расстоянию от точки старта до (ближайшего) решения [18, 23, 254, 258, 259, 264, 401]. ◊

Изложенные выше результаты (в том числе и в ускоренном случае) с помощью *минибатчинга* (mini-batching'а) переносятся на задачи стоха-стической оптимизации[75] [30, 51, 70, 82, 92, 93, 180, 225, 257, 250]. Заметим, что конструкция минибатчинга позволяет не только переносить оптимальные методы для детерминированных задач на задачи стохастической оптимизации с сохранением свойства оптимальности, но и получать оптимальные методы (по числу вызовов оракула) для задач стохастической оптимизации из неоптимальных методов для детерминированных задач. Опишем вкратце в простейшем случае суть конструкции. В задачах стохастической оптимизации вместо градиента функционала $\nabla f(x)$ оракул выдает его несмещенную оценку (*стохастический градиент*) $\nabla_x f(x,\xi)$ с конечной дисперсией $D$:[76]

$$E_\xi\left[\nabla_x f(x,\xi)\right] \equiv \nabla f(x),\ E_\xi\left[\left\|\nabla_x f(x,\xi)-\nabla f(x)\right\|_2^2\right] \le D.$$

Конструкция *минибатчинга* заключается в подстановке в метод вместо неизвестного градиента $\nabla f(x)$ его оценки (вычислимой)

$$\overset{r}{\nabla}_x f\left(x,\{\xi^l\}_{l=1}^r\right) = \frac{1}{r}\sum_{l=1}^r \nabla_x f\left(x,\xi^l\right),$$

где $\{\xi^l\}_{l=1}^r$ – независимые одинаково распределенные (так же, как $\xi$) случайные величины, и правильного выбора параметра $r$. Выбрать этот па-

---

[75] Можно обойтись и без минибатчинга, также можно рассмотреть и седловые задачи (и даже вариационные неравенства), см., например, [18, 80, 144, 158, 180, 225, 248, 331, 318, 335, 549]. Отметим, что в работе [331] описывается метод, который может работать и в условиях отсутствия точных знаний о параметре $D$ (см. также [144, 223, 335, 364, 405]).

[76] Заметим, что приводимые далее результаты можно распространить на более общую концепцию стохастического шума [295, 367, 481, 549, 565] (обобщенное *условие слабого и сильного роста*):

$$E_\xi\left[\left\|\nabla_x f(x,\xi)\right\|_2^2\right] \le D + c_1\left\|\nabla f(x)\right\|_2^2 + c_2 L\cdot\left(f(x)-f(x_*)\right).$$

Отметим, что в машинном обучении в последнее время достаточно популярны *перепараметризованные модели*, в которых выполняется *условие интерполяции*: $\nabla_x f(x_*,\xi) \cong 0$ тождественно по $\xi$ [565]. В этом случае имеет место условие слабого роста ($D=0$, $c_1=0$) с $c_2=2$. Некоторые частные случаи этой концепции будут далее разобраны (безградиентные методы). Отметим также, что в статистической теории обучения встречаются и другие показатели степени, в частности, $\left\|\nabla f(x)\right\|_2^2 \to \left\|\nabla f(x)\right\|_2^{1/2}$ – отвечает условию малого шума Цыбакова–Массара, Бернштейна [148].

раметр помогают следующие два неравенства (здесь $L = L_1$ в обозначениях (2.4)):

$$E_{\{\xi^l\}_{l=1}^r}\left[\left\|\overset{r}{\nabla}_x f\left(x,\{\xi^l\}_{l=1}^r\right)-\nabla f(x)\right\|_2^2\right]\le\frac{D}{r},$$

$$\left\langle\nabla f(x)-\overset{r}{\nabla}_x f\left(x,\{\xi^l\}_{l=1}^r\right),\mathrm{v}\right\rangle\le\underbrace{\frac{1}{2L}\left\|\overset{r}{\nabla}_x f\left(x,\{\xi^l\}_{l=1}^r\right)-\nabla f(x)\right\|_2^2}_{\delta}+\frac{L}{2}\|\mathrm{v}\|_2^2,$$

и результаты о сходимости исследуемого метода при наличии неточного оракула (подобно § 2). Последнее неравенство можно переписать более удобным образом (см. неравенство (2.3)):

$$f\left(x^{k+1}\right)\le f\left(x^k\right)+\left\langle\overset{r}{\nabla}_x f\left(x^k,\{\xi^{k+1,l}\}_{l=1}^r\right),x^{k+1}-x^k\right\rangle+\frac{2L}{2}\left\|x^{k+1}-x^k\right\|_2^2+\delta^{k+1}.$$

Используя это неравенство в цепочке рассуждений (2.10) – (2.12), придем к следующему аналогу неравенства (2.12):

$$h\left\langle\overset{r}{\nabla}_x f\left(x^k,\{\xi^{k+1,l}\}_{l=1}^r\right),x^k-x\right\rangle\le h\cdot\left(f\left(x^k\right)-f\left(x^{k+1}\right)+\delta^{k+1}\right)+$$
$$+\frac{1}{2}\left\|x-x^k\right\|_2^2-\frac{1}{2}\left\|x-x^{k+1}\right\|_2^2,$$

где $h = 1/(2L)$. Беря от обеих частей этого неравенства условное математическое ожидание $E_{x^{k+1}}\left[\,\cdot\,\middle|x^1,\dots,x^k\right]$, получим

$$f\left(x^k\right)-f(x)\le\left\langle\nabla f\left(x^k\right),x^k-x\right\rangle\le f\left(x^k\right)-E_{x^{k+1}}\left[f\left(x^{k+1}\right)\middle|x^1,\dots,x^k\right]+$$
$$+E_{x^{k+1}}\left[\delta^{k+1}\middle|x^1,\dots,x^k\right]+L\left\|x-x^k\right\|_2^2-E_{x^{k+1}}\left[L\left\|x-x^{k+1}\right\|_2^2\middle|x^1,\dots,x^k\right].$$

Суммируя выписанные неравенства и беря полное математическое ожидание, можно получить при $x = x_*$ аналог неравенства (2.22) с $L := 2L$. Исходя из (2.22), будем выбирать $r$ следующим образом:

$$\frac{D}{2Lr}\simeq E[\delta]=\frac{\varepsilon}{2}\ \Rightarrow\ r\simeq\max\left\{\frac{D}{L\varepsilon},1\right\}.$$

Поскольку всего итераций (см. § 2 и теорему 3.1)

$$\mathrm{O}\left(\frac{LR^2}{\varepsilon}\right),$$

то общее число обращений к оракулу за стохастическим градиентом $\nabla_x f(x,\xi)$ при не малых значениях $D$ будет

$$N(\varepsilon) = \mathrm{O}\left(\frac{DR^2}{\varepsilon^2}\right).$$

Эта же оценка получается и для ускоренных методов. Данная оценка является неулучшаемой оценкой для класса задач выпуклой стохастической оптимизации [70, 105, 527].

◊ В случае если не известны значения $L$ и(или) $D$, можно использовать следующий прием: с помощью подбора $L^{k+1}$ (см. § 5) добиваемся выполнения неравенства

$$f\left(x^{k+1}\right) \le f\left(x^k\right) + \left\langle \overset{r^{k+1}}{\nabla}_x f\left(x^k, \left\{\xi^{k+1,l}\right\}_{l=1}^{r^{k+1}}\right), x^{k+1} - x^k \right\rangle + \frac{2L^{k+1}}{2}\left\|x^{k+1} - x^k\right\|_2^2 + \frac{\varepsilon}{2},$$

с $r^{k+1} \simeq D_0 \big/ \left(L^{k+1}\varepsilon\right)$, где $D_0$ – оценка сверху $D$ ($D_0 \ge D$). Здесь также считаем, что $L \le D_0/\varepsilon$. Описанный способ подбора $L^{k+1}$ позволяет получить, аналогично изложенному выше, следующее неравенство

$$\frac{1}{2L^{k+1}}\left\langle \overset{r^{k+1}}{\nabla}_x f\left(x^k, \left\{\xi^{k+1,l}\right\}_{l=1}^{r^{k+1}}\right), x^k - x \right\rangle \le \frac{1}{2L^{k+1}}\cdot\left(f\left(x^k\right) - f\left(x^{k+1}\right) + \frac{\varepsilon}{2}\right) +$$
$$+\frac{1}{2}\left\|x - x^k\right\|_2^2 - \frac{1}{2}\left\|x - x^{k+1}\right\|_2^2.$$

Далее нужно специальным образом выбрать критерий останова метода. Метод останавливается, когда суммарное (по всем итерациям) число вычислений стохастического градиента достигнет некоторого заданного уровня $\tilde{N}(\varepsilon) = \mathrm{const}\cdot D_0 R^2 / \varepsilon^2$. На практике так делать не обязательно. На последней итерации в общем случае уже не получится подобрать $r^{k+1}$ по $L^{k+1}$, поскольку $r^{k+1}$ будет определяться критерием останова, поэтому на последней итерации, наоборот, по $r^{k+1}$ будет выбрано $L^{k+1} \simeq D_0 \big/ \left(r^{k+1}\varepsilon\right)$. По выпуклости $f(x)$

$$\sum_k \frac{1}{2L^{k+1}} \left\langle \overset{r^{k+1}}{\nabla}_x f\left(x^k, \left\{\xi^{k+1,l}\right\}_{l=1}^{r^{k+1}}\right), x^k - x \right\rangle \geq \sum_k \frac{1}{2L^{k+1}} \cdot \left(f\left(x^k\right) - f(x)\right) +$$

$$+\sum_k \frac{1}{2L^{k+1}} \left\langle \overset{r^{k+1}}{\nabla}_x f\left(x^k, \left\{\xi^{k+1,l}\right\}_{l=1}^{r^{k+1}}\right) - \nabla f\left(x^k\right), x^k - x \right\rangle \simeq$$

$$\simeq \sum_k \frac{1}{2L^{k+1}} \cdot \left(f\left(x^k\right) - f(x)\right) +$$

$$+\frac{\varepsilon}{2D_0} \sum_{i=0}^{\tilde{N}(\varepsilon)} \left\langle \nabla_x f\left(x^{k(i)}, \xi^i\right) - \nabla f\left(x^{k(i)}\right), x^{k(i)} - x \right\rangle.$$

Если бы последняя сумма была суммой мартингал-разностей, что, в частности, означало бы, что полное математическое ожидание последней суммы равно нулю, то описанный выше способ выбора шага $h^{k+1} = 1/\left(2L^{k+1}\right)$ и размера батча $r^{k+1} \simeq D_0/\left(L^{k+1}\varepsilon\right)$ гарантировали, что метод сойдется за $O\left(LR^2/\varepsilon\right)$ итераций, см., например, [225, chapter 7], [252, 254, 257]. Хотя в численных экспериментах и наблюдалось хорошее соответствие результатов данным оценкам (см. работы [475, 476], в которых также рассматривается ускоренный вариант описанной здесь процедуры), строго доказать все это не удалось.

Заметим, что при практической реализации метода, при увеличении $L^{k+1} := 2L^{k+1}$ (см. § 5), не следует генерировать новые стохастические градиенты в соответствующем (меньшем) количестве $r^{k+1} \simeq D_0/\left(L^{k+1}\varepsilon\right)$. В ходе процедуры $L^{k+1} := 2L^{k+1}$ можно использовать один и тот же изначально посчитанный на данной итерации вектор $\overset{r^{k+1}}{\nabla}_x f\left(x^k, \left\{\xi^{k+1,l}\right\}_{l=1}^{r^{k+1}}\right)$.

Отметим также, что в подавляющем большинстве приложений точные значения $f\left(x^k\right)$, $f\left(x^{k+1}\right)$, как правило, недоступны, если недоступны соответствующие градиенты, см. ниже про автоматическое дифференцирование. Доступны обычно случайные несмещенные реализации[77] значений целевой функции $f\left(x^k, \xi\right)$, $f\left(x^{k+1}, \xi\right)$. В этом случае (при допол-нительном предположении о выпуклости $f(x,\xi)$ по $x$ для всех $\xi$) вместо точных значений $f(x^k)$, $f(x^{k+1})$ в описанную выше адаптивную процедуру следует подставлять их оценки

$$f\left(x,\{\xi^l\}_{l=1}^r\right)=\frac{1}{r}\sum_{l=1}^r f\left(x,\xi^l\right),\ x=x^k,x^{k+1},$$

построенные аналогично оценке стохастического градиента. При этом теперь под $L$ следует понимать другую константу. А именно, не константу Липшица градиента $f(x):=E_\xi\left[f(x,\xi)\right]$, а худшую (по $\xi$) из констант Липшица градиентов (по $x$) функций $f(x,\xi)$. В этом случае можно гарантировать, что неравенство

$$f\left(x^{k+1},\xi\right)\le f\left(x^k,\xi\right)+\left\langle\nabla_x f\left(x^k,\xi\right),x^{k+1}-x^k\right\rangle+\frac{2L^{k+1}}{2}\left\|x^{k+1}-x^k\right\|_2^2+\frac{\varepsilon}{2},$$

а следовательно и

$$f\left(x^{k+1},\left\{\xi^{k+1,l}\right\}_{l=1}^{r^{k+1}}\right)\le f\left(x^k,\left\{\xi^{k+1,l}\right\}_{l=1}^{r^{k+1}}\right)+\left\langle\overset{r^{k+1}}{\nabla}_x f\left(x^k,\left\{\xi^{k+1,l}\right\}_{l=1}^{r^{k+1}}\right),x^{k+1}-x^k\right\rangle+$$

$$+\frac{2L^{k+1}}{2}\left\|x^{k+1}-x^k\right\|_2^2+\frac{\varepsilon}{2},$$

заведомо выполняются, если $L^{k+1}$ в процессе подбора дойдет до значения этой новой константы $L$. Так же как и раньше знать настоящее значение этой константы для реализации метода не обязательно.

---

[77] Часто говорят просто о реализациях. То, что они случайные и несмещенные, подразумевается.

К сожалению, построить законченную теорию здесь пока не получилось. Тем не менее, особо отметим работу [144], в которой предлагается универсальный метод для решения монотонных стохастических вариационных неравенств на базе проксимального зеркального метода (см. замечание 5.1). По сути, используется стандартный проксимальный зеркальный метод, в котором $L$ предлагается выбирать не так как в замечании 5.1, а специальным образом, схожим со способом, использующимся в AdaGrad [221, 248]. Естественно такой метод можно использовать и для решения седловых задач и задач выпуклой оптимизации [364]. С точностью до логарифмических множителей, сходится такой метод (и по числу итераций и по числу параллельных обращений к оракулу) для задач выпуклой оптимизации по оценкам, которые были выписаны выше для стохастического варианта градиентного спуска. Однако этот метод не является полностью адаптивным, поскольку также как и в AdaGrad в стратегии выбора шага существенно используется информация о размере решения. Полностью адаптивный метод решения гладких стохастических монотонных вариационных неравенств был построен (с небольшими оговор-ками) в работе [335]. Данная работа представляется достаточно интересной в плане возможности перенесения полученных в ней результатов на ускоренные методы решения задач гладкой стохастической выпуклой оптимизации.

Отметим также работы [221,571], в которой небольшая модификация AdaGrad исследуется с точки зрения глобальной сходимости к экстремуму для невыпуклых задач. ◊

Для $\mu$-сильно выпуклой в 2-норме функции $f(x)$ приведенную оценку

$$N(\varepsilon) = \mathrm{O}\left(\frac{DR^2}{\varepsilon^2}\right)$$

можно улучшить с помощью рестартов (см. указание к упражнению 2.3 и конец § 5):

$$N(\varepsilon) = \mathrm{O}\left(\min\left\{\frac{DR^2}{\varepsilon^2}, \frac{D}{\mu\varepsilon}\right\}\right).$$

Данная оценка является неулучшаемой оценкой для класса задач сильно выпуклой стохастической оптимизации [70, 105]. Заметим, что не сильно выпуклую часть оценки можно получить из сильно выпуклой с помощью регуляризации $\mu \simeq \varepsilon / R^2$ (см. замечание 4.1). На более специальных классах задач приведенные оценки допускают уточнения. Например, сходимость может характеризоваться не дисперсией стохастического градиента (одно число), а его корреляционной матрицей [430, 491, 509].

◊ Заметим, что для $\mu$-сильно выпуклой в 2-норме функции $f(x)$ имеет место достаточно неожиданный результат о сходимости обычного градиентного спуска, в котором вместо градиента использует стохастический градиент (т.е. размер батча $r = 1$), с той же скоростью, что и обычный градиентный спуск (см. § 1) в $\sqrt{D/(L\mu)}$ -окрестность решения [253, 379, 380, 430, 481]. На базе этого результата и стандартных приемов (минибатчинга, каталиста, регуляризации) также можно получать приведенные выше оценки. ◊

Негладкий случай (см. (2.4) с $\nu = 0$) с помощью искусственного введения неточности в оракул (см. § 2) можно свести к гладкому случаю с $L \sim L_0^2/\varepsilon$. Поэтому (также оптимальную) оценку на число обращений к оракулу за стохастическим (суб-)градиентом в негладком случае можно записать в виде (следует сравнить с оценками из упражнений 2.1, 2.3)

$$N(\varepsilon) = \mathrm{O}\left( \min\left\{ \frac{\left(L_0^2 + D\right)R^2}{\varepsilon^2}, \frac{L_0^2 + D}{\mu\varepsilon} \right\} \right).$$

Заметим, что *метод стохастического зеркального спуска* (2.19) (см. также упражнение 2.6) с $h = \varepsilon / M^2$, где $M^2 = L_0^2 + D$, и заменой $\nabla f\left(x^k\right) \to \nabla_x f\left(x^k, \xi^k\right)$, сходится в чезаровском смысле согласно первому аргументу минимума в приведенной оценке [351]. В адаптивном варианте шаг требуется выбирать более изощренным образом, чем просто заменой градиента на стохастический градиент в формулах упражнения 2.6. Адаптивный вариант получил название *AdaGrad* [221, 248]. Различные варианты этого метода [223, 405] являются сейчас одними из основных алгоритмов обучения глубоких нейронных сетей [39]. С помощью рестартов, подобно упражнению 2.3, можно получить вариант метода стохастического зеркального спуска для сильно выпуклых задач [351, 354]. Впрочем, для евклидовой прокс-структуры существуют и прямые (без рестартов) варианты метода зеркального спуска для сильно выпуклых задач, см., например, [318] и цитированную там литературу.

Резюмируем приведенные выше результаты в виде табл. 2. В таблице приведены оценки на число итераций $N(\varepsilon)$ (вызовов (стохастического) градиентного оракула), необходимых для решения задачи (в среднем) с точностью $\varepsilon$ по функции:

$$E\left[ f\left(x^N\right) \right] - f\left(x_*\right) \le \varepsilon .$$

Во всех этих оценках под $R$ можно понимать евклидово расстояние от точки старта, до ближайшего к точке старта решения. Эти оценки являются нижними оценками (не могут быть улучшены) с одними числовыми множителями и достигаются на описанных выше методах (их ускоренных вариантах) с другими числовыми множителями.

Таблица 2

| $N(\varepsilon)$ | $E\left[\left\|\nabla_x f(x,\xi)\right\|_2^2\right] \le M^2$ | $\left\|\nabla f(y) - \nabla f(x)\right\|_2 \le$ $\le L\left\|y - x\right\|_2$ | $\left\|\nabla f(y) - \nabla f(x)\right\|_2 \le L\left\|y - x\right\|_2$ $E\left[\left\|\nabla_x f(x,\xi) - \nabla f(x)\right\|_2^2\right] \le D$ |
|---|---|---|---|
| $f(x)$ выпуклая | $\frac{M^2 R^2}{\varepsilon^2}$ | $\sqrt{\frac{LR^2}{\varepsilon}}$ | $\max\left\{\sqrt{\frac{LR^2}{\varepsilon}}, \frac{DR^2}{\varepsilon^2}\right\}$ |
| $f(x)$ $\mu$-сильно выпуклая | $\frac{M^2}{\mu\varepsilon}$ | $\sqrt{\frac{L}{\mu}}\left\lceil \ln\left(\frac{\mu R^2}{\varepsilon}\right)\right\rceil$ | $\max\left\{\sqrt{\frac{L}{\mu}}\left\lceil \ln\left(\frac{\mu R^2}{\varepsilon}\right)\right\rceil, \frac{D}{\mu\varepsilon}\right\}$ |

В табл. 2 числовые множители опущены. Сильно выпуклая строчка таблицы получается из выпуклой с помощью рестартов (см. упражнение 2.3 и конце параграфа 5), в обратную сторону переход осуществляется с помощью регуляризации, см. замечание 4.1. Негладкий столбец таблицы будет иметь такой же вид и в случае детерминированного оракула, выдающего (суб-)градиент. Под $M$ следует в этом случае понимать константу Липшица функционала. Негладкий столбец получается из гладких столбцов с помощью конструкции универсального метода, см. § 5. Стохастический гладкий столбец получается из детерминированного гладкого с помощью минибатчинга, см. приложение выше. Ускоренные оценки последних двух столбцов (с точностью до логарифмического множителя, зависящего от желаемой точности) могут быть получены из неускоренных с помощью конструкции каталист, см. замечание 3.3. Результаты, приведенные в табл. 2 можно обобщить на случай других норм, при этом $R^2 = 2V\left(x_*, x^0\right)$ (если решение не единственно в этой формуле выбирается то из решений $x_*$, которое минимизирует правую часть), а в последней строчке и последнем столбце в оценках скорости сходимости используемых методов появятся дополнительные логарифмические множители $\ln n$, см. конец § 2, указание к упражнению 2.3 и неравенство больших уклонений в неевклидовом случае ниже. Насколько нам известно, вопрос о том, возможно ли убрать эти множители остается открытым. Скорее всего, ответ тут будет отрицательным, т.е. эти множители в общем случае убрать нельзя.

Заметим, что для задач выпуклой (стохастической) оптимизации таблицу, аналогичную табл. 2, можно построить для неполноградиентных методов [9, 16, 249, 258, 259, 264, 395, 524, 525]; в гладком случае для методов с целевым функционалом: норма градиента целевого функционала [192, 273, 372, 447] (или норма градиентного отображения для целевого функционала[78] [114]); можно построить для методов параллельной опти-мизации [183, 230, 575]; можно построить для методов распределенной оптимизации [136, 165, 252, 294] (на число коммуникационных шагов). Насколько нам известно, на данный момент открытым остается вопрос об оптимальных оценках на число вызовов оракула (при оптимальных оценках на число коммуникаций) на один узел для задач распределенной децентрализованной оптимизации. Имеется гипотеза, дающая ответ на этот вопрос [252]. Открытым остается вопрос и об оптимальных оценках для параллельных алгоритмов в централизованной оптимизации. В работе [575] было показано, что если $r$ узлов (машин) параллельно могут решать исходную задачу стохастической оптимизации, синхронизируясь (через центральный узел) не чаще, чем через каждые $l$ итераций (вычислений стохастического градиента), то после $T$ таких $l$-циклов для широкого класса задач гладкой выпуклой стохастической оптимизации и широкого класса допустимых алгоритмов имеет место следующая нижняя оценка на качество выдаваемого решения $\tilde{x}^T$ (числовые множители опущены для наглядности):

$$E\left[f\left(\tilde{x}^T\right)\right]-f\left(x_*\right)\geq\frac{LR^2}{\left(lT\right)^2}+\sqrt{\frac{DR^2}{rlT}}.$$

---

[78] В данной работе [114] предложена оригинальная техника итеративной регуляризации (см. замечание 4.1, а также работу [217], в которой отмеченная техника проинтерпретирована с помощью неускоренного проксимального градиентного спуска) исходной задачи в сочетании с горячим стартом (warm start). Огрубляя детали для большей наглядности, опишем вкратце эту технику в несильно выпуклом случае. Регуляризуем задачу с помощью $\sigma_1 \left\| x - x^0 \right\|_2^2$, где $\sigma_1 \sim \tilde{O}\left(\varepsilon^2\right)$. Решаем оптимальным методом (с оракулом, выдающим стохастический градиент) регуляризованную задачу (см. табл. 2) до момента, когда невязка по функции уменьшится в два раза. То что получаем на выходе обозначим через $x^1$. Затем по индукции: регуляризуем (дополнительно) задачу с

Активные исследования ведутся в направлении разработки таких параллельных централизованных алгоритмов, которые как можно лучше приближаются к этой нижней границе [576]. На данный момент не известны методы, работающие в общем случае согласно приведенной нижней оценке (см. также нижнюю оценку в немного более общем случае

---

прошлой итерации с помощью $\sigma_{k+1}\left\|x-x^k\right\|_2^2$, где $\sigma_{k+1}=2\sigma_k$ и решаем ее, стартуя из $x^k$, с удвоенной точностью по функции. После $K=\mathrm{O}\left(\log_2\left(C/\varepsilon^2\right)\right)$ таких итераций (перазапусков) получим точку $x^K$, которая (с точностью до размерного множителя) будет $\tilde{\mathrm{O}}\left(\varepsilon^2\right)$-решением исходной задачи по функции (на самом деле в этом месте мы намеренно существенно огрубили описание конструкции работы [114], что, впрочем, практически не влияет на сам метод), а следовательно (см. формулу (1.8)) и $\tilde{\mathrm{O}}\left(\varepsilon\right)$-решением исходной задачи по критерию малости нормы градиента (градиентного отображения). Поскольку, согласно табл. 2 и способу согласованного (пропорционального) увеличения параметра сильной выпуклости и уменьшения невязки по функции, на каждой итерации требуется $\tilde{\mathrm{O}}\left(\varepsilon^{-2}\right)$ раз вычислять стохастический градиент, то общая трудоемкость описанного подхода также будет $\tilde{\mathrm{O}}\left(\varepsilon^{-2}\right)$. Несложно в данном случае понять, что оценка оптимальная [114] (с точностью до логарифмических множителей).

[577]). Основная мотивация здесь – это приложения к Federated Learning [363]. Сюда можно отнести и набирающие популярность исследования по local SGD [362, 368, 374, 538, 576, 577]. Особо отметим работу [374], в которой построена общая теория сходимости неускоренного стохастического градиентного спуска (с консенсусом) при самых разнообразных предположениях об архитектуре (параллельная, типа local SGD, gossip, децентрализованная распределенная, меняющаяся со временем и т.д.).

**Замечание 1.** Во многих задачах (в частности, в задачах анализа данных [39, 522, 579]) функционал имеет вид суммы большого числа слагаемых:

$$f(x) = \frac{1}{m}\sum_{l=1}^{m} f_l(x) \to \min_{x\in Q}.$$

◊ Такие задачи, в частности, возникают, если использовать метод Монте-Карло для задач стохастической оптимизации [526] (минимизации среднего риска)

$$f(x) = E_\xi\left[f(x,\xi)\right] \to \min_{x\in Q\subseteq\mathbb{R}^n}.$$

Подход метода Монте-Карло заключается в замене исходной задачи стохастической оптимизации следующей задачей (минимизации эмпирического риска)

$$\frac{1}{m}\sum_{l=1}^{m} f\left(x,\xi^l\right) \to \min_{x\in Q},$$

где с.в. $\xi^l$ – независимы и распределены также как и $\xi$. Заметим, что для того, чтобы гарантировать, что абсолютно точное решение этой новой задачи является $\varepsilon$-приближенным по функции решением исходной задачи в общем случае может потребоваться взять $m$ порядка [305, 526] (с точностью до числового множителя и оговорок относительно поправок на вероятности больших уклонений; далее для краткости такие оговорки, как правило, будем опускать):[79]

[79] Отметим, что даже для сильно выпуклых задач ( $f(x)$ – сильно выпуклая функция) в общем случае нельзя гарантировать равномерную аппроксимацию среднего риска эмпирическим риском в локальной окрестности решения исходной

$$\frac{DR^2 n \ln\left(L_0 R/\varepsilon\right)}{\varepsilon^2}.$$

Эта наблюдение хорошо поясняет, что подход, связанный с усреднением случайности за счет самого метода, как правило (но не всегда [251]), более предпочтителен, чем замена задачи стохастической оптимизации ее стохастической аппроксимацией [520, 526].[80] Более предпочтителен не только тем, что допускает адаптивность постановки и легко переносится на онлайн модификации исходной задачи, но, прежде всего, лучшей приспособленностью к большим размерностям [537]. Впрочем, при возможности организации распределенных вычислений и(или) невозможности организации онлайн доступа к слагаемым, минимизация (регуляризованной [520]) стохастической аппроксимации может быть более уместна [294]. Отметим также, что задачу минимизации

---

задачи [520]. При этом, для не сильно выпуклых задач (просто выпуклых), если правильно регуляризовать функционал (см. замечание 4.1), то можно гарантировать, что точное решение задачи минимизации соответствующего эмпирического риска будет $\varepsilon$-приближенным по функции решением исходной задачи минимизации среднего риска, если взять число слагаемых $m$ равным $\left(L_0^2 + D\right)R^2 / \varepsilon^2$, см. [520]. Следует сравнить эту оценку с выписанной оценкой на $m$ в нерегуляризованном случае $\left(L_0^2 + D\right)R^2 n \ln\left(L_0 R/\varepsilon\right) / \varepsilon^2$, которая более чем в $n$ раз больше. Отметим также, что в выпуклом случае (без регуляризации) при $m = O\left(\left(L_0^2 + D\right)R^2 / \varepsilon^2\right)$ оптимальное значение эмпирического риска отличается от оптимального значения среднего риска не более чем на $\varepsilon$ [305].

[80] Особенно ярко это проявляется в случае, бесконечномерных пространств (можно говорить и о пространствах, размерность которых много больше планируемого числа вызовов оракула), возникающих в *статистической теории обучения* (SLT, Statistical Learning Theory) [495, 521, 537]. Попытка обучиться (в классе задач с выпуклой структурой) за счет минимизации эмпирического риска может не дать состоятельной оценки/решающего правила, в то время как соответствующий стохастический зеркальный спуск дает состоятельную оценку (здесь также спасает регуляризация эмпирического риска). Отметим, что в работе [537] приводится достаточно интересный общий результат: в задачах обучения с выпуклой структурой (в частности, в задачах SLT, математической статистики и онлайн обучения) способ получения оптимальных (с точностью до логарифмических множителей) оценок/решающих правил (или, другими словами, способ наискорейшего обучения) базируется на применении соответствующего метода зеркального спуска. Правда, найти «соответствующий метод» / прокс структуру, в свою очередь, представляет собой непростую задачу. Отметим также, что отмеченная выше «выпуклая структура», не часто встречается на практике.

эмпирического риска можно решать со сложностью, независящей от значения[81] $m$, см. далее. ◊

Если (объем выборки) $m$ – очень большое число, то вместо честного и дорогого вычисления градиента вычисляют стохастический градиент, случайно (равновероятно) выбирая $r \ll m$ слагаемых $\left\{\xi^l\right\}_{l=1}^r$ и формируя стохастический градиент (несмещенную оценку градиента) по формуле

$$\nabla f\left(x,\left\{\xi^l\right\}_{l=1}^r\right)=\frac{1}{r}\sum_{l=1}^r \nabla f_{\xi^l}(x).$$

Такой подход называют *методом рандомизации суммы*, см., например, [23]. С другой стороны, описанную конструкцию можно понимать и как минибатчинг, если посмотреть на исходную постановку задачи следующим образом:

$$f(x)=E_\xi\left[f(x,\xi)\right]\to\min_{x\in Q},\ f(x,\xi)=f_\xi(x),\ \nabla_x f(x,\xi)=\nabla f_\xi(x),$$

$$P(\xi=l)=\frac{1}{m},\ l=1,\ \ldots,\ m.$$

Именно в таком ключе обычно смотрят на минибатчинг в глубоком обучении [39, 508]. На практике иногда можно получить дополнительное ускорение, выбирая случайно $\left\{\xi^l\right\}_{l=1}^r$ без повторений [111, 425, 494].

Отметим, что минибатчинг хорошо параллелится в отличие от процедур типа стохастического усреднения (stochastic averaging) [261]. В популярной работе [471] обсуждается альтернативная возможность распараллеливания стохастического градиентного спуска, в которой отсутствует синхронизационная накладка и обсуждается возможность одновременной записи в общую память. Интересные исследования по альтерантивным к минибатчингу архитектурам параллелизации (например, local SGD) описаны в работах [289, 362, 368, 374, 538, 575, 576, 577] и цитированной в них литературе.

Продемонстрируем важную роль стохастической оптимизации (рандомизации) при решении *Big Data* задач следующими двумя примерами.

---

[81] Не считая того, что потребуется многократно генерировать случайную величину, принимающую равновероятно значения $1,\ldots,m$, что каждый раз можно делать, например, с помощью генерации $O(\log_2 m)$ случайных битов [50, гл. 2].

Вернемся к задаче минимизации суммы выше. Ограничимся рассмотрением выпуклого (но не сильно выпуклого) случая. Воспользуемся методом рандомизации суммы с $r = 1$. Пусть все функции $f_l(x)$ в определении функционала $f(x)$ имеют ограниченные константы $L_0$, см. (2.4). Тогда для решения исходной задачи минимизации суммы с точностью по функции (в среднем) $\varepsilon$ требуется $O\left(L_0^2 R^2/\varepsilon^2\right)$ обращений к оракулу за (суб-)градиентами случайно выбранных слагаемых $f_l(x)$. Здесь $R$ – расстояние в 2-норме от точки старта до (ближайшего к точке старта) решения. В то же время, даже если все слагаемые достаточно гладкие – имеют ограниченные константы $L_1$ (см. (2.4)), для достижения той же точности $\varepsilon$ быстрому градиентному методу потребуется $O\left(m\sqrt{L_1 R^2/\varepsilon}\right)$ обращений к оракулу за градиентами $f_l(x)$.[82] Для задач с большим числом слагаемых при невысоких требованиях к точности решения первый (рандомизированный) способ может оказаться предпочтительнее. Отметим, что выписанная оценка $O\left(m\sqrt{L_1 R^2/\varepsilon}\right)$ ввиду наличия специальной структуры у задачи уже не будет оптимальной. Оптимальна следующая оценка, достигаемая на методах редукции дисперсии [106, 115, 118, 135, 386, 393, 409, 410, 578]: $O\left(m + \sqrt{mL_1 R^2/\varepsilon}\right)$.[83] В невыпуклом случае с критерием мало-

---

[82] Заметим, что из условия $L_1 < \infty$ не следует, что $L_0 < \infty$ (можно рассмотреть пример квадратичной функции), как следствие, нет гарантий, что стохастический градиент (градиент случайно выбранного слагаемого) имеет конечный второй момент (дисперсию). Это создает определенные трудности, которые в общем случае являются препятствием для прямого применения (без использования приема редукции дисперсии) ускоренных методов стохастической оптимизации для минимизации суммы функций [140], см. также упражнение 1.8 3).

[83] Строго говоря, стоит различать константу $L_1$, входящую в оценку $O\left(m\sqrt{L_1 R^2/\varepsilon}\right)$, и константу $L_1$, входящую в оценку $O\left(m + \sqrt{mL_1 R^2/\varepsilon}\right)$. В первом случае под $L_1$ следует понимать константу Липшица градиента оптимизируемой функции $f(x)$. Во втором случае под $L_1$ следует понимать максимальную из констант Липшица градиентов слагаемых $f_l(x)$, $l = 1,\ldots,m$. Во втором случае $L_1$ может быть заметно больше. Например, для функции ($m = n$)

$$f(x) = \frac{1}{n}\|x\|_2^2 = \frac{1}{n}\sum_{l=1}^{n} x_l^2 = \frac{1}{m}\sum_{l=1}^{m} f_l(x)$$

сти нормы градиента см. [113, 114, 116, 117, 120, 386, 392]. Оптимальные оценки были получены в [268], см. также [133, 489, 569, 570] и самый конец приложения. В случае если дополнительно известно, что функция $f(x)$ – $\mu$-сильно выпуклая в 2-норме, то оптимальной оценкой будет [388]: $\tilde{O}\left(m+\sqrt{mL_1/\mu}\right)$. Наиболее современный обзор текущих достижений в данной области см. в [386, 413, 588]. Приведенные выше результаты обобщаются на случай, когда оракул вместо градиента функции $f_l(x)$ выдает её стохастический градиент [25, 379, 380, 394]. Приведенные выше результаты также обобщаются на седловые задачи и вариационные неравенства [483, 587]. Перенесение на модельную общность на данный момент осуществлено лишь частично: для композитных и проксимальных случаев [219, 394].

Приведенные оценки достижимы, только если имеется доступ к градиенту каждого слагаемого в отдельности [132], а не только к целому градиенту оптимизируемой функции, т.е. используется *инкрементальный метод* [106]. Далее в примере будет продемонстрирован способ получения такого типа оценок. Тем не менее даже с учетом этого замечания можно привести конкретные примеры, когда первый способ по-прежнему остается предпочтительнее. Отметим также, что для конкретных примеров приведенные оценки (и их обобщения на негладкие выпуклые задачи [118]) могут быть улучшены [206].

Еще раз подчеркнем, что при решении гладких выпуклых задач стохастической оптимизации или при решении гладких выпуклых задач рандомизированными методами с помощью конструкции минибатчинга, удается эффективно распараллеливать вычисления. Скажем, в только что разобранном примере, приведенная ранее оценка $O\left(L_0^2R^2/\varepsilon^2\right)$ на число вычислений градиентов случайно выбранных слагаемых $f_l(x)$ (при дополнительно предположении о гладкости функционала) может быть редуцированы (за счет распараллеливания при минибатчинге) до следующей оценки $O\left(\sqrt{L_1R^2/\varepsilon}\right)$. В данном случае результат вполне ожидаем и достижим без минибатчинга. Просто за счет структуры функционала и возможности параллельно (на $m$ процессорах) вычислять градиент в соот-

---

разница будет в $m=n$ раз [324, 547]. Впрочем, если в методах редукции дисперсии выбирать слагаемые в сумме не равновероятно (случайно), а с специальными вероятностями, зависящими от контсант Липшице градиентов слагаемых, то можно заменить в итоговой оценке сложности худшую константу Липшица градиентов слагаемых на некоторую среднюю величину [314, 429].

ветствующем быстром градиентном методе. Однако конструкция минибатчинга естественным образом распараллеливается во всех случаях, независимо от наличия у задачи дополнительной структуры. Здесь лишь хотелось продемонстрировать эффективность такого распараллеливания на простом примере.

Заметим, что методы редукции дисперсии хотя и дают в теории лучший результат в широком диапазоне параметров, на практике на задачах машиного обучения могут проигрывать обычному методу стохастического градиентного спуска с специально выбранным размером шага и небольшим размером батча [565], особенно для невыпуклых задач [528]. В (сильно) выпуклом случае этому есть простое объяснение, см., например, [374] и цитированную там литературу, в основе которого лежит результат о том, что в оценку скорости сходимости стохастического градиентного спуска (Stochastic Gradient Descent (SGD)) для задачи вида суммы вместо дисперсии входит агрегат, который уменьшается (вплоть до нуля) с увеличением степени схожести слагаемых в сумме друг с другом (например, если выполянется упоминаемое выше условие интерполяции, приводящее к условию слабого роста). Последнее обстоятельство в разной степени присуще многим задачам машинного обучения. Другое объяснение связано с дополнительными возможностями для SGD в выборе больших шагов и с меньшей чувствительностью к неизвестным параметрам задачи, таким, как, например, константа сильной выпуклости [430].

Вернемся к задаче минимизации квадратичной формы на симплексе из замечания 5.2 и упражнения 1.6, считая, что матрица $A$, задающая квадратичную форму, плотно заполненная и все элементы этой матрицы ограничены по модулю числом $M$: $|A_{ij}| \le M$. Если для решения этой задачи использовать быстрый градиентный метод, то оценка общего времени работы метода, необходимого для достижения точности по функции $\varepsilon$ будет

$$\underbrace{O\left(n^2\right)}_{\substack{\text{сложность} \\ \text{итерации}}}\underbrace{O\left(\sqrt{L_1 R^2/\varepsilon}\right)}_{\text{число итераций}} = O\left(n^2\sqrt{(M\ln n)/\varepsilon}\right).$$

В этом примере используется 1-норма и энтропия в качестве проксфункции, см. конец § 2 и упражнение 3.7. Если ту же самую задачу решать с той же точностью $\varepsilon$ (только в среднем) рандомизированным методом с рандомизацией, описанной в замечании 5.2, то оценка общего времени работы будет

$$\underbrace{O(n)}_{\substack{\text{сложность}\\ \text{итерации}}} \underbrace{O\left(\left(L_0^2 R^2\right)/\varepsilon^2\right)}_{\text{число итераций}} = O\left(n \cdot \left(M^2 \ln n\right)/\varepsilon^2\right).$$

Для задач очень больших размеров, при невысоких требованиях к точности решения, второй (рандомизированный) способ может оказаться предпочтительнее. ■

◊ Приведем, следуя [25], [180, item 6.3], [404] для задачи минимизации суммы из замечания 1 с $Q = \mathbb{R}^n$ (следовательно $\nabla f(x_*) = 0$) способ (SVRG + каталист) получения оптимальных (с точностью до логарифмического множителя) оценок на число вычислений градиентов слагаемых. Будем считать, что все функции $f_k(x)$ имеют $L_1$-липшицев градиент и $\mu$-сильно выпуклые (всё относительно 2-нормы). В качестве несмещенной оценки градиента будем брать вектор (*variance reduction* – выделение главной части в методах Монте-Карло [50, п. 3.1.1])

$$\nabla_x f\left(x^{s,k}, \xi^{s,k}\right) = \nabla f_{\xi^{s,k}}\left(x^{s,k}\right) - \nabla f_{\xi^{s,k}}\left(y^s\right) + \nabla f\left(y^s\right),$$

где $y^s = x_{\bar{N}}^{s-1}$, случайная величина $\xi$ принимает равновероятно одно из значений $1,\ldots,m$, и использовать специальным образом рестартованный градиентный спуск с минибатчингом (см. выше)

$$x^{s,k+1} = x^{s,k} - \frac{1}{L_1} \overset{r^s}{\nabla}_x f\left(x^{s,k}, \xi^{s,k}\right),\ k = 0,\ldots,\bar{N}-1,\ r^s \simeq \max\left\{\frac{D^s}{L_1 \varepsilon}, 1\right\},$$

где параметр $\bar{N}$ будет выбран позже как $\bar{N} = O(4L_1/\mu)$,

$$E_\xi\left[\left\|\nabla_x f\left(x^{s,k},\xi\right) - E_\xi\left[\nabla_x f\left(x^{s,k},\xi\right)\right]\right\|_2^2\right] \le E_\xi\left[\left\|\nabla_x f\left(x^{s,k},\xi\right)\right\|_2^2\right] \le$$

$$\le \underbrace{2E_\xi\left[\left\|\nabla f_{\xi^{s,k}}\left(x^{s,k}\right) - \nabla f_{\xi^{s,k}}\left(x_*\right)\right\|_2^2\right]}_{\substack{\overset{(1.8)}{\le} 2E_\xi\left[2L_1\left(f_{\xi^{s,k}}\left(x^{s,k}\right) - f_{\xi^{s,k}}\left(x_*\right) - \left\langle \nabla f_{\xi^{s,k}}\left(x_*\right), x^{s,k} - x_*\right\rangle\right)\right] \overset{\nabla f(x_*)=0}{=} \\ = 4L_1\left(f\left(x^{s,k}\right) - f\left(x_*\right)\right)}} +$$

$$+\underbrace{2E_{\xi}\left\|\nabla f_{\xi^{s,k}}\left(y^{s}\right)-\nabla f_{\xi^{s,k}}\left(x_{*}\right)-\nabla f\left(y^{s}\right)\right\|_{2}^{2}}_{\substack{\overset{\nabla f(x_*)=0}{=} 2E_{\xi}\left\|\nabla f_{\xi^{s,k}}\left(y^{s}\right)-\nabla f_{\xi^{s,k}}\left(x_{*}\right)-E_{\xi}\left[\nabla f_{\xi^{s,k}}\left(y^{s}\right)-\nabla f_{\xi^{s,k}}\left(x_{*}\right)\right]\right\|_{2}^{2}\le \\ \le 2E_{\xi}\left\|\nabla f_{\xi^{s,k}}\left(y^{s}\right)-\nabla f_{\xi^{s,k}}\left(x_{*}\right)\right\|_{2}^{2}\overset{(1.8)}{\le} 4L_1\left(f\left(y^{s}\right)-f\left(x_{*}\right)\right)}} \le$$

$$\le 4L_1\cdot\left(f\left(x^{s,k}\right)-f\left(x_{*}\right)\right)+4L_1\cdot\left(f\left(y^{s}\right)-f\left(x_{*}\right)\right)=\mathrm{O}\left(L_1\Delta f^{s}\right)\overset{def}{=}D^{s},$$ [84]

где $\Delta f^{s}=f\left(y^{s}\right)-f\left(x_{*}\right)$. Здесь по $k$ идет внутренний цикл, а по $s$ внешний.

Возьмем $\bar{N}=\mathrm{O}\left(4L_1/\mu\right)$, тогда

$$\Delta f^{s+1}=\max\left\{\mathrm{O}\left(\Delta f^{s}\exp\left(-\frac{\mu}{L_1}\bar{N}\right)\right),\mathrm{O}\left(\frac{D^{s}}{\mu\bar{N}}\right)\right\},$$

$$\frac{D^{s}}{\mu\bar{N}}=\mathrm{O}\left(\frac{L_1\Delta f^{s}}{\mu\bar{N}}\right)=\mathrm{O}\left(\frac{1}{4}\Delta f^{s}\right),\ \Delta f^{s}\exp\left(-\frac{\mu}{L_1}\bar{N}\right)=\mathrm{O}\left(\frac{1}{4}\Delta f^{s}\right).$$

Получим (все рассуждения здесь и далее можно провести более аккуратно, убрав $\mathrm{O}(\ )$)

$$\Delta f^{s+1}\le\mathrm{O}\left(\frac{1}{2}\Delta f^{s}\right).$$

Таким образом, общее число вычислений $\nabla f_l\left(x\right)$, необходимых для достижения точности $\varepsilon$ по функции, будет

$$\mathrm{O}\left(\left(m+\frac{L_1}{\mu}\right)\cdot\ln\left(\frac{\Delta f^{0}}{\varepsilon}\right)\right).$$

---

[84] С помощью этого соотношения и минибатчинга (см. начало приложения) можно получить обобщение описываемого здесь метода редукции дисперсии на случай, когда вместо градиента слагаемого доступен только его стохастический градиент. Ускоряя такой метод каталистом [379] можно получить оптимальные оценки сложности [380, 393].

Собственно, каталист и был впервые предложен в работе [409] как способ ускорения этой оценки до

$$\tilde{O}\left(\left(m+\sqrt{m\frac{L_1}{\mu}}\right)\cdot\ln\left(\frac{\Delta f^0}{\varepsilon}\right)\right).$$

Напомним, что в каталисте (см. замечание 3.3) идет игра на выборе параметра регуляризации $L \gg \mu$ исходной задачи. Сложность решения (число вычислений $\nabla f_l(x)$) на каждой (внешней) итерации каталиста внутренней задачи с необходимой точностью будет

$$O\left(\left(m+\frac{L_1}{\mu+L}\right)\cdot\ln\left(\frac{\Delta f^0}{\varepsilon}\right)\right)=\tilde{O}\left(m+\frac{L_1}{L}\right),$$

а число внешних итераций каталиста будет (чтобы получить такую оценку метод Монтейро–Свайтера из замечания 3.3 необходимо еще рестартовать, см. конец § 5)

$$O\left(\sqrt{\frac{L}{\mu}}\ln\left(\frac{\Delta f^0}{\varepsilon}\right)\right)=\tilde{O}\left(\sqrt{\frac{L}{\mu}}\right).$$

Выбирая,

$$L=O\left(\frac{L_1}{m}\right),$$

получим желаемый результат.

Идея каталиста также будет продемонстрирована ниже на более простом примере. В качестве хорошего упражнения на использование техники каталист можно попробовать ускорить результаты работ [295, 314]. Причем результаты [295] интересно также было бы попробовать ускорить, исходя из подхода, описанного в замечании 1.5 и тексте после него, см. также [549, 553]. ◊

Все отмеченное выше переносится также на покомпонентные и безградиентные постановки задач[85] [9, 16, 18, 23, 25, 70, 82, 93, 203, 259, 264, 344, 395, 445, 468, 512, 524]. Пусть случайный вектор $e$, например, равномерно распределен на евклидовой сфере в $\mathbb{R}^n$ радиуса 1, т. е. $\|e\|_2 = 1$ или равновероятно среди единичных ортов [38, 51, 534] (*покомпонентная рандомизация*). Тогда

$$\frac{f(x+\tau e)-f(x)}{\tau} \simeq \langle \nabla f(x), e\rangle, \; E_e\underbrace{\left[n\langle \nabla f(x), e\rangle e\right]}_{\text{то, что подставляется в метод вместо градиента}} = \nabla f(x),$$

$$\langle \nabla f(x), \langle \nabla f(x), e\rangle e\rangle = \|\langle \nabla f(x), e\rangle e\|_2^2, \; E_e\left[\|n\langle \nabla f(x), e\rangle e\|_2^2\right] = n\|\nabla f(x)\|_2^2.$$

Выписанные соотношения вкупе с основным неравенством (3.1) позволяют получить, что число итераций (число обращений к оракулу за значением функции или производной по направлению) для таких методов в среднем по порядку будет в $n$ раз больше, чем для полноградиентных аналогов. В общем случае этот результат не может быть улучшен [70, 107]. Впрочем, при дополнительных предположениях улучшения возможны [9, 16, 25, 259, 469], см. также разобранный ниже пример решения задачи квадратичной оптимизации покомпонентным методом. Описанный результат вполне понятен, поскольку, запросив частные производные (или значения функции) по $n$ координатным ортам, можно просто восстановить полный градиент. В этой связи отметим, что если есть доступ к программе (точнее, следует говорить о доступе к тексту программы, см. ниже), вычисляющей значение функции (так бывает далеко не во всех приложениях, см., например, [210]), то, как правило, лучше попробовать использовать *автоматическое дифференцирование* [472, гл. 8], чем просто аппроксимировать градиент (и тем более старшие производные) конечными разностями [42].

◊ Автоматическое дифференцирование (automatic differentiation) – способ по программе, вычисляющей значение функции (дереву вычислений), построить программу, вычисляющую градиент функции и работающую не дольше, чем в 4 раза, по сравнению с исходной. Однако такой способ требует в общем случае большей памяти – необходимо в памяти хранить всю историю (дерево) вычисления функции. Изначально такого рода результаты были получены (Баур–Штрассен) для полиномов, см. [85] и цитированную там литературу. Впоследствии, в начале 80-х годов XX

---

[85] Отметим, что про нижние оценки на скорость накопления малых неточностей неслучайной природы для неполноградиентных методов [502] на данный момент известно меньше, чем для полноградиентных методов [225].

века, две группы в ЦЭМИ РАН [57] и ВЦ РАН, см. [49] и цитированную там литературу для более полного и точного исторического обзора, смогли получить описанный выше результат в наибольшей общности. Подробный обзор имеется также в работе [153]. В работе [448] приведен интересный пример использования автоматического дифференцирования для решения вспомогательных подзадач для методов 3-го порядка (см. ниже) с той же по порядку сложностью, что и для методов 2-го порядка (в частности, метода Ньютона), см. также указание к упражнению 3.10. Заметим, что аналогом автоматического дифференцирования для негладких выпуклых функций является лексикографическое дифференцирование, предложенное в конце 80-х годов XX века Ю. Е. Нестеровым [75, 453]. Интересно также заметить, что в ряде классических работ по нейронным сетям, в которых используется частный случай автоматического дифференцирования – *метод обратного распространения* (back propagation), имеются неточности. Эти неточности связанны как раз с тем, что для негладких функций (негладкость получается за счет использования персептронов / функций активаций ReLu) используется процедура автоматического дифференцирования, обоснованно работающая, только для гладких функций. ◊

Далее с помощью техники рестартов, следуя [25], объясняется более точно, откуда в оценке скорости сходимости появляется множитель $\sim n$. Ключевое наблюдение базируется на формуле (1.8):

$$M^2 \simeq E_e\left[\left\|n\langle\nabla f(x),e\rangle e\right\|_2^2\right] = n\left\|\nabla f(x)\right\|_2^2 \le 2Ln\cdot\left(f(x)-f(x_*)\right)$$

и приведенной выше оценке скорости сходимости градиентного спуска с минибатчингом для задач стохастической оптимизации, см. также табл. 2 (здесь, как и раньше, работаем с точностью до поправки на вероятности больших уклонений или оговорки о том, что сходимость понимается в среднем, см., например, [92, 305]):

$$f(x^N)-f(x_*) = O\left(\sqrt{\frac{M^2R^2}{N}}\right) \le O\left(\sqrt{\frac{2Ln\cdot\left(f(x^0)-f(x_*)\right)R^2}{N}}\right).$$

Рестартуя метод, каждый раз, когда происходит гарантированное (выписанной формулой) уполовинивание невязки по функции, получим следующую формулу для общего числа обращений к оракулу за $\langle\nabla f(x),e\rangle$ или $\{f(x), f(x+\tau e)\}$ (см. также указания к упражнениям 1.3, 2.3):

$$N(\varepsilon) = \mathrm{O}\left(n\frac{8LR^2}{\Delta f}\right) + \mathrm{O}\left(n\frac{8LR^2}{\Delta f/2}\right) + \mathrm{O}\left(n\frac{8LR^2}{\Delta f/4}\right) + \ldots$$

$$\ldots + \mathrm{O}\left(n\frac{8LR^2}{\varepsilon}\right) = n\cdot\mathrm{O}\left(\frac{LR^2}{\varepsilon}\right),$$

где $\Delta f = f\left(x^0\right) - f\left(x_*\right)$. К сожалению, по той же причине – грубость формулы (1.8), по которой не следует использовать для ускоренных методов рестарты по норме градиента для решения гладких сильно выпуклых задач (см. замечание 5.3), здесь не получается похожим образом перенести описанную конструкцию на ускоренные градиентные спуски (см. упражнения 1.3, 5.7). Для ускоренных методов требуются более тонкие рассуждения [25, 264, 445, 468, 469].

Приведенные выше рассуждения можно повторить в случае $\mu$-сильно выпуклой в 2-норме функции $f(x)$. Причем все получится еще проще. Можно не повторять рассуждения, а улучшить приведенную оценку с помощью рестартов (см. указание к упражнению 2.3 и конец § 5):

$$N(\varepsilon) = n\cdot\mathrm{O}\left(\frac{L}{\mu}\left\lceil\ln\left(\frac{\mu R^2}{\varepsilon}\right)\right\rceil\right).$$

Проверить соответствие данной оценки, аналогичной оценке, полученной ранее в несильно выпуклом случае, можно с помощью регуляризации $\mu \simeq \varepsilon/R^2$ (см. замечание 4.1).

До недавнего времени считалось, что так же просто, как это было описано выше в неускоренном случае, не удастся перенести описанные выше конструкции на ускоренные методы. Однако недавно было обнаружено [244, 409, 410] (см. также § 3), как приведенные выше оценки (и многие другие оценки скорости сходимости для неускоренных методов) могут быть единообразно перенесены на ускоренный случай с помощью новой довольно общей и вместе с тем достаточно простой техники *Каталист*, см. замечание 3.3:

$$N(\varepsilon) = n\cdot\mathrm{O}\left(\min\left\{\sqrt{\frac{LR^2}{\varepsilon}}, \sqrt{\frac{L}{\mu}}\left\lceil\ln\left(\frac{\mu R^2}{\varepsilon}\right)\right\rceil\right\}\right).$$

Основным «структурным блоком» по-прежнему является вспомогательная задача (3.21),

$$x^{k+1} = \arg\min_{x\in\mathbb{R}^n}\left\{ f(x) + \frac{\tilde{L}}{2}\left\| x - x^k \right\|_2^2 \right\},$$

которую необходимо решать на каждой итерации. Сложность решения этой задачи (число обращений к оракулу за значением функции / производной по направлению на одной итерации) неускоренным безградиентным / покомпонентным методом (с точностью до логарифмического множителя, см. упражнение 3.1) равна $\tilde{O}\left(n\cdot\left(L+\tilde{L}\right)/\left(\mu+\tilde{L}\right)\right)$. С другой стороны, число внешних итераций ускоренного прокс-метода [95, 409, 410] с точностью до логарифмического множителя равно $\tilde{O}\left(\sqrt{\tilde{L}/\mu}\right)$. Для того чтобы получить такую оценку метод Монтейро–Свайтера из замечания 3.3 необходимо еще рестартовать, см. конец § 5. Оценка общего числа обращений к оракулу будет наилучшей $\tilde{O}\left(n\sqrt{L/\mu}\right)$ при выборе $\tilde{L} \simeq L$, что соответствует в сильно выпуклом случае ранее приведенному результату.

Разобранный пример наглядно демонстрирует общую идею подхода: за счет выбора параметра регуляризации $\tilde{L}$ добиваться близкой к единице обусловленности вспомогательной задачи, тогда неоптимальность используемого для ее решения подхода (напомним, что как раз для внутренней задачи используется неускоренный метод, который хотим ускорить) становится несущественной, и за счет ускоренности внешнего метода происходит общее ускорение рассматриваемой процедуры. Примечательно, что внешний ускоренный прокс-метод остается одним и тем же в данном подходе, в то время как внутренний неускоренный метод (для решения вспомогательной задачи) можно как угодно менять в зависимости от контекста. Отметим также, что выше рассуждения проводились с точностью до логарифмических множителей. К сожалению, если честно их выписывать, то выяснится, что такой подход приводит не к оптимальным оценкам, а к оптимальным с точностью до логарифмических (по желаемой точности) множителей.

Каталист хотя и является универсальным способом ускорения всевозможных неускоренных методов, тем не менее на практике предпочитают использовать прямые ускоренные рандомизированные процедуры, см., например, [115, 393, 394] и замечание 2 ниже. Особенно активно в этом направлении работают З. Аллен-Зу [113] и Дж. Лан [383].

Важно также отметить, что даже в неускоренном случае приведенные выше рассуждения при покомпонентной рандомизации оказываются достаточно грубыми, поскольку не учитывают, что константу $L$ теперь можно считать не по худшему направлению, а «средней» по всем направлениям, что может быть в $\sim \sqrt{n}$ раз меньше [25, 469].

**Замечание 2.** Предположим, что для всех $x \in \mathbb{R}^n$ и $h \in \mathbb{R}$:

$$\left| \partial f\left(x + he_i\right)/\partial x_i - \partial f\left(x\right)/\partial x_i \right| \le L_i h .$$

Пусть $\beta \in [0,1]$. Введем

$$\|x\|^2 = \sum_{i=1}^n L_i^{1-2\beta} x_i^2,\ \breve{R}^2 = \frac{1}{2}\left\|x_* - x^0\right\|^2,\ \left\|\nabla f\left(x\right)\right\|_*^2 = \sum_{i=1}^n L_i^{2\beta-1} \cdot \left(\frac{\partial f\left(x\right)}{\partial x_i}\right)^2 .$$

| Метод линейного каплинга (МЛК) | Покомпонентный вариант МЛК |
|---|---|
| См. указание к упражнению 1.3 и замечание 1.6<br>$x^{k+1} = \tau z^k + (1-\tau) y^k,$<br>$y^{k+1} = x^{k+1} - \frac{1}{L}\nabla f\left(x^{k+1}\right),$<br>$z^{k+1} = z^k - h\nabla f\left(x^{k+1}\right).$ | $x^{k+1} = \tau z^k + (1-\tau) y^k,$<br>Случайно и независимо разыгрываем $i_{k+1} \in [1,...,n]$ по правилу:<br>$P\left(i_{k+1} = i\right) = p_i \overset{def}{=} \frac{L_i^\beta}{\sum_{j=1}^n L_j^\beta},\ i = 1,...,n,$<br>$y_{i_{k+1}}^{k+1} = x_{i_{k+1}}^{k+1} - \frac{1}{L_{i_{k+1}}}\frac{\partial f\left(x^{k+1}\right)}{\partial x_{i_{k+1}}^{k+1}},$<br>$z_{i_{k+1}}^{k+1} = z_{i_{k+1}}^k - \frac{h}{p_{i_{k+1}}}\frac{\partial f\left(x^{k+1}\right)}{\partial x_{i_{k+1}}^{k+1}}.$ |

Для МЛК согласно указанию к упражнению 1.3 имеем

$$\left\langle \nabla f\left(x^{k+1}\right), z^k - x_* \right\rangle \le \frac{1}{2h}\left\|z^k - x_*\right\|_2^2 - \frac{1}{2h}\left\|z^{k+1} - x_*\right\|_2^2 + \frac{h\left\|\nabla f\left(x^{k+1}\right)\right\|_2^2}{2},$$

т. е.

$$\left\langle \nabla f\left(x^{k+1}\right), z^k - x_* \right\rangle \le \frac{1}{2h}\left\|z^k - x_*\right\|_2^2 - \frac{1}{2h}\left\|z^{k+1} - x_*\right\|_2^2 + Lh \cdot \left(f\left(x^{k+1}\right) - f\left(y^{k+1}\right)\right).$$

Для ПМЛК аналогом приведенных неравенств будут

$$\frac{1}{p_{i_{k+1}}}\left\langle\left\langle\nabla f\left(x^{k+1}\right),e_{i_{k+1}}\right\rangle e_{i_{k+1}},z^k-x_*\right\rangle\le$$

$$\le\frac{1}{2h}\left\|z^k-x_*\right\|^2-\frac{1}{2h}\left\|z^{k+1}-x_*\right\|^2+\frac{h\left\|\left\langle\nabla f\left(x^{k+1}\right),e_{i_{k+1}}\right\rangle e_{i_{k+1}}\right\|_*^2}{2p_{i_{k+1}}^2},$$

$$\left\langle\nabla f\left(x^{k+1}\right),z^k-x_*\right\rangle\le\frac{1}{2h}\left\|z^k-x_*\right\|^2-$$

$$-E_{i_{k+1}}\left[\frac{1}{2h}\left\|z^{k+1}-x_*\right\|^2\middle|i_0,...,i_k\right]+\breve{L}h\cdot\left(f\left(x^{k+1}\right)-E_{i_{k+1}}\left[f\left(y^{k+1}\right)\middle|i_0,...,i_k\right]\right),$$

где $\breve{L}=\left(\sum_{j=1}^{n}L_j^{\beta}\right)^2$. Второе неравенство получается из первого взятием условного математического ожидания от обеих частей $E_{i_{k+1}}\left[\,\cdot\,\middle|i_0,...,i_k\right]$. Поскольку согласно указанию к упражнению 1.3 оценка числа итераций (вычислений градиента $f(x)$), необходимых для достижения по функции точности $\varepsilon$, имеет вид: $\mathrm{O}\left(\sqrt{LR^2/\varepsilon}\right)$, то для ПМЛК естественно было бы ожидать, что оценка числа итераций (вычислений частных производных $f(x)$), необходимых для достижения по функции (в среднем) точности $\varepsilon$, будет $\mathrm{O}\left(\sqrt{\breve{L}\breve{R}^2/\varepsilon}\right)$. Так оно в действительности и оказывается [25, 123, 445, 469]. Аналогичные рассуждения можно было провести, взяв за основу метод подобных треугольников вместо МЛК [264].

Метод МЛК и ПМЛК можно осуществлять (с такими же оценками скорости сходимости) и без рестартов [25, 122, 123]. Для этого нужно выбирать (см. также замечание 1.6):

$$\tau_k=2/(k+2) \text{ и } h_k=(k+2)/(2L) \text{ (МЛК)},$$

$$\tau_k=2/(k+2) \text{ и } h_k=(k+2)/\left(2\breve{L}\right) \text{ (ПМЛК)}.$$

Описанный покомпонентный метод имеет естественное блочно-покомпонентное обобщение [445]. При этом допускается, что рассматриваемую задачу оптимизации необходимо решать на множестве простой структуры, имеющей вид прямого произведения множеств, отвечающих различным блокам [25, 264]. ■

Поясним абзац, написанный непосредственно перед замечанием 2, простым примером, в котором сравним время работы быстрого градиентного метода, например, из указания к упражнению 1.3 – МЛК, и его покомпонентного варианта – ПМЛК с $\beta = 1/2$ [25]. То же самое можно было продемонстрировать и для неускоренных методов. Итак, рассматривается задача квадратичной выпуклой оптимизации (1.30) в условиях работы [469]:

$$f(x) = \frac{1}{2}\langle Ax, x\rangle - \langle b, x\rangle \to \min_{x\in\mathbb{R}^n},$$

где $A = \|A_{ij}\|_{i,j=1}^{n} \succ 0$ и $1 \le A_{ij} \le 2$ при $i, j = 1, \dots, n$. Из последнего условия имеем

$$L = \lambda_{\max}(A) \ge \lambda_{\max}\left(1_n 1_n^T\right) = n,$$

поэтому оценка общего времени работы МЛК (оптимального, с точностью до числового множителя, метода для данного класса задач) будет

$$\underbrace{O\left(n^2\right)}_{\substack{\text{стоимость}\\ \text{итерации}}} \cdot \underbrace{O\left(\sqrt{\frac{LR^2}{\varepsilon}}\right)}_{\substack{\text{число}\\ \text{итерации}}} = O\left(\frac{n^{5/2}R}{\varepsilon^{1/2}}\right).$$

При этом если использовать покомпонентную рандомизацию с вероятностью выбрать $i$-орт $p_i \sim \sqrt{L_i}$, где $L_i = A_{ii}$ – константа Липшица $i$-й компоненты градиента вдоль $i$-орта, то в оценке числа итераций вместо $\sqrt{L} \ge \sqrt{n}$ можно ставить

$$\sqrt{\bar{L}} = \frac{1}{n}\sum_{i=1}^{n}\sqrt{L_i} \le \sqrt{2},$$

т. е. число итераций согласно замечанию 1 (см. также [25, 469]) будет

$$O\left(n\sqrt{\frac{\bar{L}R^2}{\varepsilon}}\right) = O\left(\frac{nR}{\varepsilon^{1/2}}\right).$$

При этом стоимость итерации теперь будет $O(n)$. Действительно, если $\tilde{x} = x + he_i$, где $e_i$ – $i$-орт и уже посчитано $Ax$, то $A\alpha\tilde{x} = \alpha Ax + \alpha hA^{\langle i\rangle}$, где $\alpha$ – заданное число, может быть посчитано за $O(n)$. Отсюда с помощью указания к упражнению 1.3 следует, что и итерацию можно осуще-ствить за время $O(n)$. Таким образом, общее время работы ПМЛК с $\beta = 1/2$ можно оценить следующим образом:

$$\underbrace{O(n)}_{\text{стоимость итерации}} \cdot \underbrace{O\left(n\sqrt{\frac{\bar{L}R^2}{\varepsilon}}\right)}_{\text{число итерации}} = O\left(\frac{n^2 R}{\varepsilon^{1/2}}\right) \ll O\left(\frac{n^{5/2} R}{\varepsilon^{1/2}}\right).$$

Отметим, что если «собирать» из компонент градиента полный градиент, и использовать обычный быстрый градиентный метод, то в последней оценке вместо $\bar{L}$ следовало бы писать просто $L$. Таким образом, в общем случае в задачах с неполной информацией о градиенте не стоит использовать полноградиентные методы «напрямую»: путем полного восстановления градиента по частично наблюдаемой информации. Этот же вывод (уже и по другой причине) можно сделать и для спусков по направлению и безградиентных методов [9, 16, 258, 259].

Кажется, что получается какое-то противоречие с оптимальностью МЛК для класса гладких выпуклых задач и противоречие с рядом тезисов, которые ранее приводились в пособии. На самом деле никаких противоречий нет [70]:

1) Во-первых, нижние оценки были получены на классе детерминированных методов. Впрочем, введение рандомизации принципиально не меняет нижние оценки.
2) Во-вторых, МЛК оптимален на классе всевозможных гладких выпуклых задач с ограниченной константой Липшица градиента. Выше же был рассмотрен более узкий класс задач – квадратичной оптимизации. Но даже при таком сужении нижние оценки принципиально не изменятся (см. упражнение 1.3). Более важным и ограничительным было предположение на элементы матрицы $A$, которое привело к тому, что след матрицы $A$ имеет тот же порядок, что и наибольшее собственное значение. Именно эта «асимметричность» постановки задачи и определила успешность использования рандомизации.
3) В-третьих, МЛК оптимален на классе задач по числу вычислений градиента функционала (в нашем случае – по числу произведений матрицы на вектор), а не по времени работы (общей трудоемкости). По времени работы не существует теории нижних оценок и вряд ли в ближайшее время можно ожидать ее появления.
4) Наконец, ранее отмечалось, что если есть возможность считать градиент, например, с помощью автоматического дифференцирования, то стоит именно градиент и использовать в методе, а не восстанавливать компоненты градиента по посчитанным значениям функции. Однако в последнем примере используется не процедура расчета значений функции, а процедура их пересчета. Действительно, вычислив самый первый раз $f(x)$ за $O(n^2)$, дальше можно уже пересчитывать $f(A\alpha\tilde{x})$, где $\tilde{x} = x + he_i$, передавая $Ax$ с прошлого запуска, за время $O(n)$. Для осуществления шага нужно посчитать только одну (случайно выбранную) компоненту градиента, что может быть приближенно сделано за два пересчета значения функции, т. е. за время $O(n)$.

◊ Заметим, что аналогичный результат (с точностью до логарифмического множителя) можно было бы получить на основе неускоренного покомпонентного метода [180, item 6.4.2] с $\gamma = 1$, [445] и техники каталист, см. замечание 3.3. Причем для такого ускоренного покомпонентного метода (в отличие от всех других, известных нам [25]) уже можно обеспечить учет разреженности матрицы $A$ в оценке стоимости итерации при минимизации функции $\mu \ln\left(\sum_{k=1}^{m} \exp\left(\langle A_k, x\rangle / \mu\right)\right)$. Такого типа функции возникают, например, при решении двойственной задачи к энтропийно-регуляризованной транспортной задаче ЛП, см. упражнение 3.9 и замечание 4.1.

Отмеченную (двойственную) функцию также можно минимизировать с помощью тензорных методов [187, 262, 450] (см. замечание 3.3 и текст далее). При этом для восстановления решения прямой задачи, в виду замечания 5.3 и упражнения 4.2, в качестве критерия следует выбирать норму градиента. Регуляризация двойственной функции (см. замечание 4.1 и упражнение 4.4) позволяет добиться оптимальной скорости сходимости [262] (см. также работу [302], содержащую нижние оценки) по критерию нормы градиента соответствующих тензорных методов без потери логарифмического множителя в скорости сходимости, как это имело место для методов первого порядка, см. замечание 5.3. ◊

Может показаться, что выше слишком много внимания было уделено, на первый взгляд, довольно незначительной оговорке о возможности в специальных (ассиметричных) случаях ускорять время работы методов в $\sim \sqrt{n}$ раз за счет покомпонентной рандомизации используемых методов.

Одна из причин такого внимания к данному примеру связана с так называемым «препроцессингом». А именно, в описанном выше подходе при рандомизации $p_i \sim \sqrt{L_i}$ явно использовалась дополнительная информация о структуре задачи. Эту информацию несложно было получить до начала работы метода. Хотя в данном случае основной эффект достигался в первую очередь за счет самого факта рандомизации, а не за счет ее специфики [25], в общем случае следует признать удачными в оптимизационной практике приемы типа диагонального шкалирования (предобуславливание), пришедшие из вычислительной линейной алгебры [36, 91, 291, 477, 556], и их рандомизированные варианты типа описанного выше. Особенно популярны сейчас (в связи с приложениями к анализу данных) адаптивные варианты таких методов, в которых по ходу работы метода предпринимается попытка без больших усилий (т. е. не как в методе Ньютона) улучшить обусловленность задачи. К таким методам можно отнести уже упоминавшиеся ранее методы с процедурой растяжения пространства [86, 100] и предобуславливание, описанное в упражнении 3.11 и указании к упражнению 5.9. Также популярными сейчас являются *квазиньютоновские методы* [472] (см. также замечание 3) и методы типа *AdaGrad* [39, 221, 248, 508].

Однако главная причина такого внимания к покомпонентным методам описана ниже. Исследования последних лет показывают (см., например, [25] и цитированную там литературу), что именно покомпонентная рандомизация (для прямой или двойственной задачи) лежит в основе большей части современных подходов к решению задач оптимизации, приходящих из анализа данных [522]. В частности, различные варианты метода рандомизации суммы (см. замечание 1) можно понимать, как варианты покомпонентных рандомизаций для двойственной задачи.

**Пример (типичная задача анализа данных [25, 118, 175, 522, 523]).** Рассмотрим задачу выпуклой оптимизации:

$$\sum_{k=1}^{m} f_k\left(A_k^T x\right) + g(x) \to \min_{x\in Q},$$

где $g(x) = \sum_{i=1}^{n} g_i(x_i)$ (это условие нужно только для того, чтобы можно было эффективно построить двойственную задачу – см. ниже), $Q$ – множество простой структуры. Предполагаем, что трудоемкость вычисления $f_k'(z_k)$ равна $O(1)$ и что для всех $k = 1, \ldots, m$ и всех допустимых $w$, $v$:

$$\left| f_k'(w) - f_k'(v) \right| \le L|w - v|.$$

Функция $g(x)$ предполагается сильно выпуклой в $p$-норме с константой $\mu_p$. Вводя матрицу $A = [A_1, \ldots, A_m]^T$ и вспомогательный вектор $z = Ax$, можно переписать эту задачу в «раздутом» пространстве $\tilde{y} = (x, z)$, как задачу типа проектирования на аффинное многообразие [4, 130, 247], см. также (4.7). Эта задача проектирования решается путем перехода к двойственной задаче. Опишем далее довольно общий способ построения двойственной задачи:

$$\min_{x\in Q}\left\{\sum_{k=1}^{m} f_k\left(A_k^T x\right)+g\left(x\right)\right\}=\min_{\substack{x\in Q\\ z=Ax}}\left\{\sum_{k=1}^{m} f_k\left(z_k\right)+g\left(x\right)\right\}=$$

$$=\min_{\substack{x\in Q\\ z=Ax,z'}}\max_{y}\left\{\left\langle z-z',y\right\rangle+\sum_{k=1}^{m} f_k\left(z_k'\right)+g\left(x\right)\right\}=$$

$$=\max_{y\in\mathbb{R}^m}\left\{-\max_{\substack{x\in Q\\ z=Ax}}\left\{\left\langle -z,y\right\rangle-g\left(x\right)\right\}-\max_{z'}\left\{\left\langle z',y\right\rangle-\sum_{k=1}^{m} f_k\left(z_k'\right)\right\}\right\}=$$

$$=\max_{y\in\mathbb{R}^m}\left\{-\max_{x\in Q}\left(\left\langle -A^T y,x\right\rangle-g\left(x\right)\right)-\sum_{k=1}^{m}\max_{z_k'}\left(z_k' y_k-f_k\left(z_k'\right)\right)\right\}=$$

$$=\max_{y\in\mathbb{R}^m}\left\{-g^*\left(-A^T y\right)-\sum_{k=1}^{m} f_k^*\left(y_k\right)\right\}=-\min_{y\in\mathbb{R}^m}\left\{g^*\left(-A^T y\right)+\sum_{k=1}^{m} f_k^*\left(y_k\right)\right\}.$$

В описанную схему погружаются, например, следующие задачи [130]:

1) $\frac{L}{2}\left\|Ax-b\right\|_2^2+\underbrace{\frac{\mu}{2}\left\|x-x_g\right\|_2^2}_{g(x)}\to\min_{x\in\mathbb{R}^n}$, (Ridge regression)

2) $\frac{L}{2}\left\|Ax-b\right\|_2^2+\underbrace{\mu\sum_{k=1}^{n}x_k\ln x_k}_{g(x)}\to\min_{x\in S_n(1)}$. (Minimal mutual information model)

Константа Липшица производной $f_k\left(z\right)=\frac{L}{2}\left(z_k-b_k\right)^2$ равна $L$, константы сильной выпуклости $g\left(x\right)$ (считаются в разных нормах: 1) в 2-норме, 2) в 1-норме) также одинаковы в обоих случаях и равны $\mu$. Для приведенных выше задач получим следующие двойственные задачи:

1) $\frac{1}{2\mu}\left(\left\|x_g-A^T y\right\|_2^2-\left\|x_g\right\|_2^2\right)+\frac{1}{2L}\left(\left\|y+b\right\|_2^2-\left\|b\right\|_2^2\right)\to\min_{y\in\mathbb{R}^m}$,

2) $\mu\ln\left(\sum_{i=1}^{n}\exp\left(\frac{\left[-A^T y\right]_i}{\mu}\right)\right)+\frac{1}{2L}\left(\left\|y+b\right\|_2^2-\left\|b\right\|_2^2\right)\to\min_{y\in\mathbb{R}^m}$.

В общем случае можно утверждать, что $\sum_{k=1}^{m} f_k^*(y_k)$ (композитный член в двойственной задаче, см. пример 3.1) является сильно выпуклым в 2-норме с константой сильной выпуклости равной $L^{-1}$. Для $g^*(-A^T y)$ можно оценить константу Липшица градиента в 2-норме (см. конец § 2, начало § 4, а также указание к упражнению 4.8):

$$\frac{1}{\mu}\max_{\|y\|_2\le 1,\,\|x\|_p\le 1}\left\langle A^T y, x\right\rangle^2 = \frac{1}{\mu}\max_{\|x\|_p\le 1}\|Ax\|_2^2 = \frac{1}{\mu}\begin{cases} 1)\ \lambda_{\max}\left(A^T A\right) \\ 2)\ \max\limits_{j=1,\dots,n}\left\|A^j\right\|_2^2 \end{cases}$$

и получить следующую оценку сверху на константы Липшица всех частных производных $g^*(-A^T y)$ [25]:

$$\frac{1}{\mu}\max_{\|y\|_1\le 1,\,\|x\|_p\le 1}\left\langle A^T y, x\right\rangle^2 = \frac{1}{\mu}\max_{\|x\|_p\le 1}\|Ax\|_\infty^2 = \frac{1}{\mu}\begin{cases} 1)\ \max\limits_{i=1,\dots,m}\left\|A_i\right\|_2^2 \\ 2)\ \max\limits_{\substack{i=1,\dots,m\\ j=1,\dots,n}}\left|A_{ij}\right|^2 \end{cases}.$$

Для ПМЛК с $\beta = 0$ из замечания 2, примененного к двойственной задаче (в случае, когда матрица $A$ плотно заполненная), имеем следующие оценки общего времени работы (трудоемкости):

1) $$T_1 = \tilde{O}\left(n\cdot m\sqrt{\frac{L\max\limits_{i=1,\dots,m}\left\|A_i\right\|_2^2}{\mu}}\right),$$

2) $$T_2 = \tilde{O}\left(n\cdot m\sqrt{\frac{L\max\limits_{i,j}\left|A_{ij}\right|^2}{\mu}}\right).$$

Если теперь посмотреть на исходную прямую задачу (с $L := L/m$):

$$\frac{1}{m}\sum_{k=1}^{m} f_k\left(A_k^T x\right) + g(x) \to \min_{x\in Q},$$

и оценить трудоемкость оптимальных методов согласно оценкам, приведенным в конце замечания 1, то получим уже анонсированное соответствие оптимальных оценок с полученными только что оценками трудоемкости ПМЛК с $\beta = 0$ (при $L := L/m$). Действительно, с учетом того, что константы Липшица градиентов $f_k\left(A_k^T x\right)$, посчитанные в соответствую-щих нормах (соответствующих норме, в которой сильно выпукл композит прямой задачи), равномерно (по $k = 1, ..., m$) оцениваются следующим образом:

1) $L \max_{i=1,..,m} \|A_i\|_2^2$,

2) $L \max_{i,j} |A_{ij}|^2$,

а сложность вычисления $\nabla f_k\left(A_k^T x\right)$ равна $O(n)$, то в типичном случае:

1) $$\tilde{T}_1 = \tilde{O}\left(n \cdot \left(m + \sqrt{m \frac{L \max_{i=1,..,m} \|A_i\|_2^2}{\mu}}\right)\right) = \tilde{O}\left(n \cdot m \sqrt{\frac{(L/m) \max_{i=1,..,m} \|A_i\|_2^2}{\mu}}\right),$$

2) $$\tilde{T}_2 = \tilde{O}\left(n \cdot \left(m + \sqrt{m \frac{L \max_{i,j} |A_{ij}|^2}{\mu}}\right)\right) = \tilde{O}\left(n \cdot m \sqrt{\frac{(L/m) \max_{i,j} |A_{ij}|^2}{\mu}}\right).$$

Таким образом, имеет место полное соответствие (с точностью до опущенных в рассуждениях логарифмических множителей). Интересно заметить, что для первой задачи здесь можно сполна использовать разреженность матрицы $A$ [25]. Более того, эту задачу (также с полным учетом разреженности) можно решать и прямым ПМЛК с $\beta = 0$. Соответствующие оценки согласно замечанию 2 имеют вид (снова возвращаемся к исходному пониманию параметра $L$):

$$T_1^{\text{прям}} = \tilde{O}\left(s \cdot n \sqrt{\frac{L \max_{j=1,...,n} \|A^j\|_2^2}{\mu}}\right),$$

$$T_1^{\text{двойств}} = \tilde{O}\left(\tilde{s} \cdot m \sqrt{\frac{L \max_{i=1,..,m} \|A_i\|_2^2}{\mu}}\right) = \tilde{O}\left(sn \sqrt{\frac{L \max_{i=1,..,m} \|A_i\|_2^2}{\mu}}\right),$$

где $s$ – среднее число ненулевых элементов в столбцах матрицы $A$, а $\tilde{s}$ – в строках.

В действительности, обе эти оценки оказываются завышенными.[86] Более аккуратные рассуждения позволяют заменить максимум в этих

[86] Это легко усмотреть из способа рассуждений, в котором мы заменяем константы Липшица частных производных на худшую из них. Отметим, что описанное

оценках на некоторые средние. Скажем, в транспортных приложениях [18, гл. 1], [130] матрица $A$ не просто разреженная, но еще и битовая (состоит из нулей и единиц). В таком случае приведенные оценки переписываются следующим образом:

$$T_1^{\text{прям}} = \tilde{O}\left( sn\sqrt{\frac{Ls}{\mu}} \right),\ T_1^{\text{двойств}} = \tilde{O}\left( sn\sqrt{\frac{L\tilde{s}}{\mu}} \right).$$

Отсюда можно сделать довольно неожиданный вывод [130]: при $m \ll n$ стоит использовать прямой ПМЛК, а в случае $m \gg n$ – двойственный. Первый случай соответствует приложениям к изучению больших сетей (компьютерных, транспортных). Второй случай соответствует задачам, приходящим из анализа данных.

В заключение отметим, что описанный в данном примере подход используется при построении эффективных асинхронных распределенных алгоритмов [323, 325], т.е. алгоритмов, не требующих синхронизации распределенных устройств. ■

Несмотря на большую активность в последние годы в направлении разработки рандомизированных методов выпуклой оптимизации, здесь по-прежнему, остаются открытые вопросы. Например, вопрос о возможности построения адаптивных ускоренных методов. Неускоренный блочно-покомпонентный градиентный спуск с адаптивным подбором констант Липшица соответствующих блоков компонент градиента был предложен в работе [445]. Полностью адаптивный ускоренный аналог до сих пор не известен. Известны лишь частично адаптивные конструкции, см., например, [271, 339].

◊ Отметим, что описанное выше направление в его современном варианте появилось во многом под влиянием двух препринтов Ю. Е. Нестерова, подготовленными в 2010–2011 гг. Результаты, приведенные в этих препринтах, впоследствии вошли в две статьи [445, 468]. П. Рихтарик (постдок Ю. Е. Нестерова) активно взялся за развитие отмеченных идей (статей) Ю. Е. Нестерова, особенно в части покомпонентных методов, находя им (их параллельным и распределенным вариантам) всевозможные приложения, особенно в задачах огромных размеров (*Big Data*), приходящих из анализа данных [500]. ◊

Для рандомизированных методов в случае возможности вычислять значение оптимизируемой функции $f(x)$ из результатов о сходимости в среднем

---

огрубление также упрощает использование оценки скорости сходимости ПМЛК с $\beta = 0$ из замечания 2.

$$E\left[f\left(x^N\right)\right]-f\left(x_*\right)\le\varepsilon$$

по неравенству Маркова следует

$$P\left(f\left(x^N\right)-f\left(x_*\right)\ge 2\varepsilon\right)\le\frac{E\left(f\left(x^N\right)-f\left(x_*\right)\right)}{2\varepsilon}\le\frac{1}{2}.$$

Запуская параллельно $m=\left\lceil\log_2\left(\sigma^{-1}\right)\right\rceil$ независимых траекторий метода и выбирая $\hat{x}^N=\arg\min_{k=1,..,m}f\left(x^{N,k}\right)$, получим [25]

$$P\left(f\left(\hat{x}^N\right)-f\left(x_*\right)\ge 2\varepsilon\right)\le\sigma.$$

Некоторые тонкости возникают для задач условной оптимизации [18, п. 4.4], [152], в которых точка $x^N$ может быть не допустимой (не удовлетворять ограничениям) и, как следствие, использовать неравенство Маркова нельзя в виду отсутствия гарантий, что

$$f\left(x^N\right)-f\left(x_*\right)\ge 0$$

Аналогичные результаты, об увеличении числа вычислений (в данном случае числа итераций) в $\sim\ln\left(\sigma^{-1}\right)$, можно получить, не предполагая возможности вычислять значения $f\left(x\right)$ (что затруднительно, например, для обычных задач стохастической оптимизации [526]), но предполагая, что стохастические (суб)градиенты несмещенные и имеют финитный носитель. Если хвосты стохастических (суб)градиентов субгауссовские, то вместо $\sim\ln\left(\sigma^{-1}\right)$ следует писать $\sim\ln^2\left(\sigma^{-1}\right)$. В общем случае, см., например, [23, 225, 257, 261, 351, 383, 526]. Отметим, что в ряде случаев возможность параллелизации вычислений без доступа к значениям $f\left(x\right)$ можно сохранить [261, 289], используя просто немного другой способ формирования $\hat{x}^N$ (похожие идеи используются и в подходе local SGD [368, 374, 538, 576]):

$$\hat{x}^N = \frac{1}{m}\sum_{k=1}^{m} x^{N,k} .$$

Если рандомизированный метод линейно сходится в среднем, см., например, [115, 393], т.е.

$$E\left[f\left(x^{N(\varepsilon)}\right)\right] - f\left(x_*\right) \le \varepsilon \text{ и } N(\varepsilon) \sim \ln\left(\varepsilon^{-1}\right),$$

то для получения точных оценок вероятностей больших уклонений также можно использовать грубое неравенство Маркова [130]:

$$P\left(f\left(x^{N(\varepsilon\sigma)}\right) - f\left(x_*\right) \ge \varepsilon\right) \le \sigma \frac{E\left[f\left(x^{N(\varepsilon\sigma)}\right)\right] - f\left(x_*\right)}{\varepsilon\sigma} \le \sigma .$$

◊ Пусть $\left(f_\delta(x); \langle \nabla f_\delta(x), y - x \rangle\right)$ будет $(\delta, L)$-модель функции $f(x)$ в точке $x$ (относительно нормы $\| \ \|$), см. § 3. Пусть (*субгауссовские хвосты* у модели стохастического градиента)[87]

---

[87] Условия субгауссовости хвостов всегда можно добиться за счет «bias variance trade off»: введения смещения (bias) в стохастический градиент за счет обнуления стохастического градиента вне некоторого шара. Радиус шара выбирается таким образом, чтобы смещение в итоговой оценке давало вклад не больше, чем желаемая точность решения задачи. Субгауссовская дисперсия (variance), входящая в условие субгауссовости для так «обрезанного» стохастического градиента, будет заметно больше его обычной дисперсии. Хочется учесть это обстоятельство. Поэтому вместо неравенства Азума–Хефдинга (приведенного далее) следует использовать более тонкое неравенство Бернштейна–Фридмана. Если так действовать, то можно даже в случае тяжелых хвостов (в предположении конечной дисперсии исходного стохастического градиента) получить экспоненциальную концентрацию в оценках скорости сходимости (аналогичных первозданно субгауссовской концентрации) соответствующих методов для гладких задач выпуклой оптимизации. Для неускоренных методов это сделано в работе [66], а для ускоренных с помощью другой техники аналогичный результат был получен в работе [217]. По-видимому, техника работы [66] также позволяет распространить результаты этой работы на ускоренные методы (недавно это было сделано [293, 583]) и на гладкие монотонные вариационные неравенства, см. замечание 5.1.

$$E_{\xi}\left[\nabla_x f_{\delta}(x,\xi)\right] \equiv \nabla f_{\delta}(x),\ E_{\xi}\left[\exp\left(\frac{\left\|\nabla_x f_{\delta}(x,\xi) - \nabla f_{\delta}(x)\right\|_*^2}{D}\right)\right] \le \exp(1).$$

Тогда имеют место следующие неравенства [92], [225, Chapter 7], [305, 348], [352, Theorem 2.1], [391, Lemma 2]

$$f_{\delta}(y) \le f(y) \le f_{\delta}(y) + \delta;$$

$$P\left( f_{\delta}(y) \le f_{\delta}(x) + \left\langle \overset{r}{\nabla}_x f_{\delta}\left(x, \left\{\xi^l\right\}_{l=1}^r\right), y - x \right\rangle + \frac{2L}{2}\|y-x\|^2 + \right.$$
$$\left. + \frac{3\left(2\chi + 4\Omega\chi + 2\Omega^2\right)}{2Lr} + \delta \right) \ge 1 - \exp\left(-\frac{\Omega^2}{3}\right),$$

где $\chi = 1$ для $\| \ \| = \| \ \|_2$ и $\chi = \min\{q-1, 2\ln n\}$ для $\| \ \| = \| \ \|_p$, $1/p + 1/q = 1$, $p \in [0,\infty]$;

$$P\left( \sum_{k=0}^{N-1} c_{k+1} \left\langle \overset{r^{k+1}}{\nabla}_x f_{\delta}\left(x^k, \left\{\xi^{k+1,l}\right\}_{l=1}^{r^{k+1}}\right) - \nabla f_{\delta}\left(x^k\right), x - x^k \right\rangle \right.$$
$$\left. \le \sqrt{3DR^2\Omega \sum_{k=0}^{N-1} \frac{c_{k+1}^2}{r_{k+1}}} \right) \ge 1 - \exp(-\Omega),$$

где по предположению $c_{k+1}$, $r^{k+1}$ и $x^k$ могут зависеть только от $\left\{\left\{\xi^{t+1,l}\right\}_{l=1}^{r^{t+1}}\right\}_{t=0}^{k-1}$, а $\left\|x - x^k\right\|_2 \le R$ при $k = 0,\dots,N-1$, $\Omega \ge 0$. Последнее неравенство вероятности больших уклонений обычно называют *неравенством концентрации (меры) Азума–Хефдинга* [172, 225, 526].

Приведенные выше неравенства в сочетании с упражнением 3.7 позволяют получить оптимальные методы решения гладких выпуклых задач стохастической оптимизации в модельной общности [253]. ◊

В описанной выше общности полученная линейка методов уже будет практически полностью покрывать основной арсенал методов первого порядка (и ниже), использующихся в современных приложениях для решения задач выпуклой оптимизации большого размера.

- В частности, описанные выше в пособии подходы (в особенности, прямодвойственные универсальные ускоренные методы решения седловых задач и задач выпуклой оптимизации с оракулом, выдающем модель функции и их (блочно) покомпонентные варианты) позволяют строить методы, работающие по наилучшим известным сейчас теоретическим оценкам общего времени работы (трудоемкости) практически для всех известных нам классов задач (структурной) выпуклой оптимизации больших размеров.

Большое число исследований во всем мире сейчас сосредоточено на изучении (переборе) конкретных способов сочетания описанных выше приемов с целью определения их наилучшего сочетания для изучаемого класса задач, как правило, возникающего в одном из актуальных приложений. Отметим в этой связи упражнения 3.7, 5.8 и замечание 3.3, «ускоряющие» почти все оценки, приведенные в пособии, и упражнение 5.5, распространяющее действие разобранных в пособии конструкций с задач выпуклой оптимизации на множествах простой структуры на общие задачи выпуклой оптимизации (с аффинными ограничениями вида равенств и выпуклыми ограничениями вида неравенств). Отметим также, что многие популярные на практике приемы типа *альтернативных направлений, метода штрафных функций, метода модифицированной функции Лагранжа*, *ADMM* и др. [13, 46, 53, 82, 88, 151, 154, 155, 218, 386, 389, 390, 413, 472, 479, 559, 582] и их ускоренные варианты: ускоренный метод альтернативных направлений, AADMM [274, 290, 307, 356, 369, 413, 479], не гарантируют в общем случае лучших <u>теоретических</u> оценок, чем те, которые могут быть получены при описанных в пособии подходах. Улучшения могут быть только в негладком случае за счет проксимальной природы ряда описанных процедур, например, метода модифицированной функции Лагранжа [555], см. также § 3, упражнение 4.3 и замечание 4.3. Но сполна этим можно воспользоваться, как правило, только в случае организации распределенных вычислений, см., например, ADMM [175, 479, 582].

Однако не следует думать, что этим уже исчерпываются современные численные методы выпуклой оптимизации и способы их исследования. Выше мало обсуждались вопросы о стоимости итерации, завязанные на автоматическое дифференцирование [49, 153, 472] и возможность быстрого пересчета значений функции, (компонент) градиента [18, 460]. В этой связи достаточно привести несколько примеров [18, п. 1.5.2, 4, 5.1], [124, 126], [180, п. 3.3], [211, 330, 351, 401, 460, 467], см. также упражнение 1.6. Все эти примеры ярко демонстрируют, что на практике за счет дешевых итераций быстрее могут работать методы, которые далеко не оптимальны с точки зрения числа итераций. Также совсем не рассматривались реальные приложения, например, к задачам оптимизации в гиль-бертовых пространствах, в которых оракул типично зашумлен [14, 48, 49, 281, 487]. Наконец, почти ничего не было сказано о многих других методах и подходах, которые часто используются при решении практических задач умеренных размеров [272]: в частности, квазиньютоновских методах [472] и методах второго порядка (ньютоновские методы), интерес к которым в последние годы резко возрос [21, 35, 49, 71, 72, 108, 121, 138, 143, 191, 209, 286, 299, 301, 358, 359, 412, 426, 444, 448, 449, 450, 451, 452 464, 568, 579].

**Замечание 3 (квазиньютоновские методы [224], [472, гл. 6]).** В замечании 1.2 приводится геометрическая интерпретация градиентного спуска, в основу которой положена замена исходной функции параболоидом вращения, касающимся её графика в текущей точке. Точка, доставляющая минимум параболоиду, принимается за новое положение метода. В замечании 1.4 рассматривается наискорейший спуск, заключающийся в подборе кривизны параболоида с помощью решения вспомогательной задачи одномерной минимизации. В *методе Ньютона* вместо параболоида вращения строится квадратичная аппроксимация оптимизируемой функции (на основе доступного гессиана), однако это приводит к необходимости решения на каждом шаге более сложной задачи – минимизации квадратичной формы (1.30). Естественно, возникает идея построения какого-то «промежуточного» метода, с одной стороны, не требующего вычисления (и тем более обращения) гессиана, а с другой стороны все-таки пытающегося как-то аппроксимировать гессиан, исходя из накопленной информации первого порядка (градиентов). В основу *квазиньютоновских методов* положен следующий общий принцип построения квадратичной аппроксимации: квадратичная аппроксимация должна касаться графика оптимизируемой функции в текущей точке и иметь с ним одинаковые градиенты в точке с предыдущего шага (secant equation). Существует много различных способов удовлетворить этим условиям. У этих способов есть различные интерпретации, среди которых отметим понимание квазиньютоновских методов, как *методов переменной метрики* [82, п. 2 § 3, гл. 3] (варианта метода сопряженных градиентов). Наиболее интересными в практическом плане являются способы, которые требуют не более чем квадратичной (по размерности пространства) трудоемкости и памяти. Среди таких способов наиболее удачно себя зарекомендовал способ, приводящий в итоге к методу *BFGS*:

$$h_k = \arg\min_{h \in \mathbb{R}} f\left(x^k - hH_k \nabla f\left(x^k\right)\right),$$

$$x^{k+1} = x^k - h_k H_k \nabla f\left(x^k\right),$$

$$H_{k+1} = H_k + \frac{H_k\gamma_k\delta_k^T + \delta_k\gamma_k^T H_k}{\langle H_k\gamma_k, \gamma_k\rangle} - \beta_k \frac{H_k\gamma_k\gamma_k^T H_k}{\langle H_k\gamma_k, \gamma_k\rangle},$$

где

$$\beta_k = 1 + \frac{\langle \gamma_k, \delta_k\rangle}{\langle H_k\gamma_k, \gamma_k\rangle},\ \gamma_k = \nabla f\left(x^{k+1}\right) - \nabla f\left(x^k\right),\ \delta_k = x^{k+1} - x^k,\ H_0 = I.$$

В отличие от сопряженных градиентов (см. замечание 1.6) в BFGS не обязательно точно осуществлять вспомогательную одномерную оптимизацию. В целом BFGS оказался наиболее устойчивым (к вычислительным погрешностям) вариантом квазиньютоновских методов. Геометрия квазиньютоновских методов (см. замечание 1.2) близка (но не идентична!) геометрии субградиентных методов с процедурой растяжения пространства [82, п. 4 § 4, гл. 5] (метод эллипсоидов, методы Шора). Рисунок 4 демонстрирует типичное «пилообразное» поведение градиентных спусков в окрестности минимума для плохо обусловленных задач. Из этого рисунка видно, что направления $\gamma_k$ и $\delta_k$ «подсказывают» направления растяжение/сжатие и позволяют адаптивно улучшать обусловленность задачи, правильно аккумулируя собранную информацию в $H_k$.

В теоретическом плане по квазиньютоновским методам на данный момент (как и 30 – 40 лет назад) известно не так уж и много (см., например, [72, п. 1.3.1], [224], [472, гл. 6]): глобальная скорость сходимости для гладких задач выпуклой оптимизации в общем случае не выше, чем у обычных (неускоренных) градиентных методов (во всяком случае, только это пока удалось установить), а локальная скорость сходимости в случае невырожденного минимума – *сверхлинейная*, т. е. быстрее, чем линейная. В работах [504, 505] был предложен такой вариант квазиньютоновского метода, для которого в случае невырожденного минимума была доказана оценка локальной скорости сходимости вида $\left\|x^k - x_*\right\|_2 \sim c^{k^2}$, где $c \in (0,1)$ – некоторая константа, связанная с числом обусловленности задачи.

Основным ограничением по использованию квазиньютоновских методов является необходимость в хранении и обновлении плотной квадратной матрицы $H_k$, что требует (в отличие от того, что имеет место для методов типа сопряженных градиентов) квадратичной памяти и квадратичного времени независимо от разреженности задачи. Это обстоятель-ство существенно ограничивает возможности по использованию таких методов для задач оптимизации с десятками тысяч переменных и более. Однако на практике используют в основном варианты таких методов *с ограниченной памятью,* см., например, метод *LBFGS* [472, п. 7.2]. В этом случае в памяти хранится не матрица $H_k$, а векторы, её порождающие. Проблема, однако, тут в том, что с ростом $k$ размер этой памяти линейно растет. Поэтому обычно последовательности векторов $\{\gamma_l\}_{l=0}^k$ и $\{\delta_l\}_{l=0}^k$ хранят только с $q$ последних итераций ($q$ – глубина памяти), и при этом полагают $H_{k-q} = I$. На практике $q$ часто выбирают совсем небольшим: $q \simeq 3-5$.

Есть основания полагать, что в ближайшее время именно в данной области могут появиться наиболее интересные результаты, объясняющие высокую эффективность на практике таких методов, как, например, LBFGS [472, гл. 7] и всевозможных рандомизированных вариантов квазиньютоновских методов [171, 297, 508]. ■

В заключение все же приведем сопоставительный анализ методов первого порядка (градиентных методов) и методов более высокого порядка, которые могут использоваться для решения задач выпуклой оптимизации умеренных размеров ($n \le 10^4$), в условиях отсутствия шума на классе достаточно гладких

$$\left\|\nabla^r f(y) - \nabla^r f(x)\right\|_2 \le M_r \left\|y - x\right\|_2, \ x, y \in \mathbb{R}^n, \ M_r \le \infty, \ r = 0, 1, 2, \dots,$$

$\mu$-сильно выпуклых в 2-норме задач, где $\mu \ge 0$. Условие $x, y \in \mathbb{R}^n$ можно заменить условием $x, y \in \left\{z \in \mathbb{R}^n : \ f(z) \le f\left(x^0\right)\right\}$. При $r = 0,1$ (и при $r \ge 2$ для оптимального тензорного метода из замечания 3.3, см. также [21, 184, 451]) последнее условие, в свою очередь, можно улучшить, см., например, (2.19). Условие липшицевости $\nabla^r f(x)$ можно заменить условием гёльдеровости $\nabla^r f(x)$ с параметром $\nu \in [0,1]$ [303]. При этом приведенные далее оптимальные оценки можно распространить (если параметр $\nu$ известен) и на этот случай [532].

◊ Заметим, что $\nabla^r f(y)$ – тензор ранга $r$. Поэтому следует пояснить, что понимается под 2-нормой в левой части данного неравенства. Ограничимся случаем $r = 2$, тогда

$$\nabla^2 f(x) = \left\|\partial^2 f(x) / \partial x_i \partial x_j\right\|_{i,j=1}^n,$$

$$\left\|\nabla^2 f(y)-\nabla^2 f(x)\right\|_2 = \sup_{\|h_1\|_2\le 1}\sup_{\|h_2\|_2\le 1}\left\langle\left(\nabla^2 f(y)-\nabla^2 f(x)\right)[h_1],h_2\right\rangle =$$

$$= \sup_{\|h_1\|_2\le 1}\sup_{\|h_2\|_2\le 1}\left\langle\left(\nabla^2 f(y)-\nabla^2 f(x)\right)h_1,h_2\right\rangle =$$

$$= \sup_{\|h\|_2\le 1}\left\langle\left(\nabla^2 f(y)-\nabla^2 f(x)\right)h,h\right\rangle =$$

$$= \max\left\{\lambda_{\max}\left(\nabla^2 f(y)-\nabla^2 f(x)\right),\left|\lambda_{\min}\left(\nabla^2 f(y)-\nabla^2 f(x)\right)\right|\right\}.$$

В общем случае см. [143, 448]. Отметим также, что при $r=0$: $\nabla^0 f(x)=f(x)$, а $\|\ \|_2=|\ |$. В связи с последним замечанием, стоит упомянуть, что, в действительности, под $M_0$ можно понимать меньшую константу (а именно, $L_0$ – см., например, (2.4)), которая только в худшем случае совпадает с введенной здесь [454]. Аналогичное замечание имеет место и по методам 2-го порядка [464] (и, вероятно, более высокого порядка). «Правильный» метод $p$-го порядка ($p\ge 1$) на первых итерациях «осуществляет» желаемую редукцию (уменьшение) констант гладкости $\{M_r\}_{r=0}^{p-1}$ за счет попадания в нужную область сходимости метода. Причем часто достаточно одной (первой) итерации [452, 454, 464]. ◊

Для класса методов, у которых на каждой итерации разрешается не более чем $O(1)$ раз обращаться к оракулу (подпрограмме) за $\nabla^r f(x)$, $r\le 1$, оценка числа итераций, необходимых для достижения точности $\varepsilon$ (по функции), будет иметь вид

$$O\left(\min\left\{n\ln\left(\frac{\Delta f}{\varepsilon}\right);\ \frac{M_0^2R^2}{\varepsilon^2},\left(\frac{M_1R^2}{\varepsilon}\right)^{1/2};\ \frac{M_0^2}{\mu\varepsilon},\left(\frac{M_1}{\mu}\right)^{1/2}\left[\ln\left(\frac{\mu R^2}{\varepsilon}\right)\right]\right\}\right),$$

где, как и раньше, $R=\left\|x^0-x_*\right\|_2$, $\Delta f=f\left(x^0\right)-f\left(x_*\right)$. Данная оценка в общем случае не может быть улучшена даже если дополнительно известно, что, $M_2<\infty$, $M_3<\infty$, … [70]. При этом данная оценка достигается [70, 72, 180, 452].

◊ Заметим, что если вместо $r\le 1$ имеет место $r=0$ (класс безградиентных методов), то в приведенной оценке все аргументы минимума следует домножить на размерность пространства $n$ [9, 16, 18, 70, 84, 258, 259, 264, 395, 399, 468]. Отметим также, что у известных сейчас методов, отвечающих (с точностью до логарифмического множителя) первому ар-гументу минимума, достаточно дорогой является составляющая итерации, не связанная с вычислением градиента: $\gg n^2$ (см. указание к упражнению 1.4 и [70, 180, 400]). ◊

Для класса методов, у которых на каждой итерации разрешается не более чем $O(1)$ раз обращаться к оракулу (подпрограмме) за значениями $\nabla^r f(x)$, $r \le p$, $p \ge 2$, оценка числа итераций, необходимых для достижения точности $\varepsilon$ (по функции), будет иметь вид [21] (в работе [186] приведено обобщение выписанной оценки на случай монотонных вариационных неравенств)

$$O\left(\min\left\{n\ln\left(\frac{\Delta f}{\varepsilon}\right); \frac{M_0^2R^2}{\varepsilon^2}, \left(\frac{M_1R^2}{\varepsilon}\right)^{1/2}, \left(\frac{M_2R^3}{\varepsilon}\right)^{2/7}, \ldots, \left(\frac{M_pR^{p+1}}{\varepsilon}\right)^{2/(3p+1)};\right.\right.$$

$$\min\left\{\left(\frac{M_1}{\mu}\right)^{1/2}, \left(\frac{M_2R}{\mu}\right)^{2/7}, \ldots, \left(\frac{M_pR^{p-1}}{\mu}\right)^{2/(3p+1)}\right\} +$$

$$\left.\left.+\min_{r=2,\ldots,p}\log\left[\log\left(\frac{\left(\mu^{r+1}/M_r^2\right)^{1/(r-1)}}{\varepsilon}\right)\right]\right\}\right).$$

Данная оценка в общем случае не может быть улучшена даже если дополнительно известно, что, $M_{p+1} < \infty$, $M_{p+2} < \infty$, … [70, упражнение 2 § 1 гл. 8], [138]. Полезно отметить, что здесь, также как и для градиентных методов (см. упражнения 1.3, 2.1), можно строить «универсальные худшие в мире функции» в классе выпуклых полиномов (от абсолютных значений линейных комбинаций переменных) степени $p+1$ [138, 448]:

$$f_m(x) = \eta_{p+1}(A_m x) - x_1,\ \eta_{p+1}(x) = \frac{1}{p+1}\sum_{i=1}^{n}|x_i|^{p+1},$$

$$A_m = \begin{pmatrix} U_m & 0 \\ 0 & I_{n-m} \end{pmatrix},\ U_m = \underbrace{\begin{pmatrix} 1\,-1 & 0\ldots 0 \\ 0 & 1\,-1\ldots 0 \\ \ldots\ldots\ldots\ldots\ldots & \\ 0 & 0\ldots 1\,-1 \\ 0 & 0\ldots 0\;\;1 \end{pmatrix}}_{m},\ I_{n-m} = \underbrace{\begin{pmatrix} 1\,0\,0\ldots 0 \\ 0\,1\,0\ldots 0 \\ \ldots\ldots\ldots\ldots \\ 0\,0\ldots 0\,1 \end{pmatrix}}_{n-m}.$$

Точка старта выбирается следующим образом $x^0 = (0,...,0)^T$, а класс допустимых методов (алгоритмов) определяется условием

$$x^{k+1} \in x^0 + \sum_{l=0}^{k} S_f\left(x^l\right) \subseteq \left\{x \in \mathbb{R}^n :\ x_i = 0,\ i = k+2,...,n\right\},$$

где

$$S_f(x) = \mathrm{Lin}\left\{y_x\left(a_0,...,a_p;\gamma;q\right) : a_0,...,a_p \in \mathbb{R};\ \gamma > 0;\ q = 1,2,3,...\right\},$$

$$y_x\left(a_0,...,a_p;\gamma;q\right) \in \operatorname*{Arg\,min}_{y \in \mathbb{R}^n} \left\{ \sum_{r=0}^{p} a_r \nabla^r f(x) \underbrace{[y-x,...,y-x]}_{r} + \gamma \|y-x\|_2^q \right\}.$$

В [452, Theorem 4] было показано, что для всех методов из описанного выше класса при $N \le (n-1)/2$ имеет место оценка

$$f_{2N+1}\left(x^N\right) - \min_{x \in \mathbb{R}^n} f_{2N+1}(x) \ge \frac{3p}{2^{p+1}(p+1)!} \frac{M_p R^{p+1}}{(2N+2)^{\frac{3p+1}{2}}},$$

где $M_p = p!$. При этом недавно было подмечено [138, 426], [452, item 4.3], что данная оценка с точностью до числового множителя достигается при $p = 2$. В работах [143, 448] было показано, как при $p \ge 2$ получать методы (подобно ускоренному методу Ньютона с кубической регуляризацией [444]) с почти оптимальной оценкой скорости сходимости. Также в [300, 448, 450] было показано, что при $p = 3$ (при $p = 2$ см. [444, 450, 464]) трудоемкость каждой итерации предложенного метода сопоставима (с точностью до логарифмических множителей) с трудоемкостью итерации метода Ньютона.

В работе [109] приводятся отличные от [138] нижние оценки при более общих условиях.

Отметим также, что при $p = 2$ функция $\eta_3(x)$ не имеет Липшицевых прозводных третьего порядка. Поэтому, если дополнительно предположить, что для метода второго порядка целевая функция имеет Липшицевы производные третьего порядка, то можно ожидать, что оценка сложности такого класса задач будет лучше. Так оно и есть на самом деле! А именно, с точностью до логарифмических множителей сложность класса трижды гладких задач выпуклой оптимизации для методов второго порядка $p = 2$ будет такой же как сложность класса трижды гладких задач выпуклой оптимизации для методов третьего порядка $p = 3$ [449, 461].

Основная идея получения данного результата заключается в использовании оптимального метода третьего порядка, в котором для решения возникающей на каждой итерации вспомогательной задачи оптимизации используется информация только второго порядка (о градиенте и гессиане целевой функции). Как уже отмечалось в § 1, ничего подобного для методов первого порядка $p=1$ получить уже нельзя.

До недавнего времени считалось, что ускорение, напрямую унаследованное от градиентных методов, не дает оптимальных оценок [143, 444, 448]. Тем не менее, в 2018 году задача построения оптимальных методов высокого порядка (работающих согласно ранее приведенной общей оценке) с помощью проксимального ускоренного метода Монтейро–Свайтера была решена, см. замечание 3.3 [21, 184, 257, 280, 346, 452].

Поясним, как получается сильно выпуклая часть последней оценки из не сильно выпуклой.[88] Для этого рассмотрим *метод Ньютона*[89]:

$$x^{k+1} = \arg\min_{x\in\mathbb{R}^n}\left\{ f\left(x^k\right) + \left\langle \nabla f\left(x^k\right), x - x^k \right\rangle + \frac{1}{2}\left\langle \nabla^2 f\left(x^k\right)\left(x - x^k\right), x - x^k \right\rangle \right\} =$$

$$= x^k - \left[\nabla^2 f\left(x^k\right)\right]^{-1} \nabla f\left(x^k\right).$$

[88] Не сильно выпуклая часть оценки (кроме первых двух аргументов минимума) получается из сильно выпуклой с помощью регуляризации $\mu \simeq \varepsilon / R^2$ (см. замечание 4.1).

[89] В методе Ньютона на каждой итерации предлагается минимизировать параболоид (уже, вообще говоря, не вращения), касающийся (до второго порядка включительно, то есть в точке касания совпадают не только градиенты, но и гессианы) графика рассматриваемой функции. Однако этот параболоид, как правило, не мажорирует рассматриваемую функцию. Это обстоятельство не позволяет гарантировать глобальную сходимость метода. Одним из возможных решений проблемы глобальной сходимости такого типа методов второго порядка является добавление к параболоиду кубического слагаемого – кубической регуляризации [143, 444, 452, 464]. Полученная в результате вспомогательная задача (минимизации регуляризованного кубическим членом параболоида) по сложности сопоставима с исходной [194, 464], но при этом метод приобретает глобальную сходимость. Отметим, что локальная сходимость остается по-прежнему сверхлинейной, как у метода Ньютона.

Заметим, что на примере метода Ньютона для поиска корней простейших полиномов (вида $z^3 - 1$) можно продемонстрировать насколько порой причудливыми бывают бассейны притяжения (бассейны/фракталы Ньютона) численных методов оптимизации в многоэкстремальных задачах [334].

Считая, что $M_2 < \infty$, $\mu > 0$, получим

$$\left\|\nabla f\left(x^{k+1}\right)\right\|_2 = \left\|\nabla f\left(x^{k+1}\right) - \nabla f\left(x^k\right) - \nabla^2 f\left(x^k\right)\left(x^{k+1} - x^k\right)\right\|_2 \le M_2 \left\|x^{k+1} - x^k\right\|_2^2 =$$

$$= M_2 \left\|\left[\nabla^2 f\left(x^k\right)\right]^{-1} \nabla f\left(x^k\right)\right\|_2^2 \le \frac{M_2}{\mu^2}\left\|\nabla f\left(x^k\right)\right\|_2^2.$$

С помощью последнего неравенства можно оценить окрестность квадратичной скорости сходимости метода Ньютона:

$$\frac{M_2}{\mu^2}\left\|\nabla f\left(x^k\right)\right\|_2 < 1 \;\Rightarrow\; \frac{\mu}{2}\left\|x^k - x_*\right\|_2^2 \le f\left(x^k\right) - f\left(x_*\right) \le \frac{1}{2\mu}\left\|\nabla f\left(x^k\right)\right\|_2^2 < \frac{\mu^3}{2M_2^2}.$$

Значит, окрестность квадратичной скорости сходимости содержит открытый шар

$$\left\|x - x_*\right\|_2^2 < \frac{\mu^2}{M_2^2}.$$

Оказавшись в этой окрестности, можно достичь желаемой точности за $\mathrm{O}\left(\left\lceil \log\left\lceil \log\left(\left(\mu^3/M_2^2\right)/\varepsilon\right)\right\rceil\right\rceil\right)$ итераций метода Ньютона.

◊ Заметим, что если последовательность положительных чисел $\{c_k\}_{k=0,1,2,\ldots}$ удовлетворяет условию (обеспечивающему *сверхлинейную скорость сходимости*):

$$c_{k+1} \le \mathrm{const}\cdot\left(c_k\right)^{\gamma},\ \gamma > 1,$$

и $c_0$ достаточно мало, то после $N = \mathrm{O}\left(\log\left\lceil\log\left(c_0/\varepsilon\right)\right\rceil\right)$ итераций $c_N \le \varepsilon$. Для метода Ньютона $c_k = \left\|\nabla f\left(x^k\right)\right\|_2$, $\gamma = 2$; для класса чебышёвских методов высокого порядка $c_k = \left\|x^k - x_*\right\|_2$, $\gamma = 3, 4, 5, \ldots$ [49, п. 2.9], [56, п. 9.5.10], [236]; для методов Ньютона с кубической регуляризацией $c_k = f\left(x^k\right) - f\left(x_*\right)$, $\gamma = 4/3$, причем в последнем случае сильную выпуклость можно заменить градиентным доминированием [71, гл. 4], [444, 464] (см. также замечание 1.1). ◊

Чтобы оказаться в этой окрестности, можно использовать технику рестартов (см. упражнение 2.3, § 5, а также [21], [71, гл. 4], [138, 444]), примененную к методам, обеспечивающим несильно выпуклую составляющую рассматриваемой оценки.

Пусть

$$f\left(x^N\right)-f\left(x_*\right)\le \frac{CM_p\left\|x^0-x_*\right\|_2^{p+1}}{N^{(3p+1)/2}}.$$

Из $\mu$-сильно выпуклости функции $f\left(x\right)$ отсюда имеем

$$\frac{\mu}{2}\left\|x^N-x_*\right\|_2^2\le f\left(x^N\right)-f\left(x_*\right)\le \frac{CM_p\left\|x^0-x_*\right\|_2^{p+1}}{N^{(3p+1)/2}}.$$

Выберем $N$ таким образом, чтобы

$$\left\|x^N-x_*\right\|_2^2\le \frac{1}{2}\left\|x^0-x_*\right\|_2^2.$$

Для этого достаточно сделать

$$\left(\frac{4CM_pR^{p-1}}{\mu}\right)^{2/(3p+1)}$$

итераций. Процедуру можно повторить (рестартовать). Рестарты заканчиваются в момент попадания в окрестность квадратичной скорости сходимости для метода Ньютона. Таким образом, число рестартов можно оценить сверху следующим образом

$$\left\lceil \log_2\left(\frac{M_2^2R^2}{\mu^2}\right)\right\rceil.$$

При этом общее число итераций до попадания в окрестность квадратичной сходимости будет (учитываем, что после каждого рестарта $R^2$ уменьшается как минимум в два раза)

$$\mathrm{O}\left(\left(\frac{M_pR^{p-1}}{\mu}\right)^{2/(3p+1)}\right),$$

что и объясняет вторую заключительную (сильно выпуклую) часть рассматриваемой оценки.

Вместо метода Ньютона здесь можно было бы использовать методы более высокого порядка сходимости [48, 464]. При попадании в окрестность сверхлинейной (квадратичной, кубической и т.д.) скорости сходимости таких методов желаемая относительная точность по функции $\tilde{\varepsilon}$ будет достигнута после $\sim \log\log\left(\tilde{\varepsilon}^{-1}\right)$ [48, 70, 138] дополнительных итераций (основания логарифмов зависят от порядка выбранного метода).

К сожалению, на данный момент не известно, как можно в общем случае использовать сильную выпуклость исходной задачи оптимизации, чтобы быстрее (по сравнению с выпуклым случаем) решать вспомогательную подзадачу (минимизации многочлена Тейлора) на каждой итерации рассматриваемых (тензорных) методов, см. указание к упражнению 3.10.

Отмети также, что имеются обнадеживающие результаты относительно падения чувствительности тензорных методов к точности задания старших производных, что может быть полезно для приложений (рандомизированных) тензорных методов к задачам минимизации суммы большого числа слагаемых [234, 416, 568].

**Замечание 4.** Если вместо сильной выпуклости и достаточной гладкости предполагать *самосогласованность* оптимизируемой функции, то, используя при анализе скорости сходимости вместо евклидовой нормы специальную локальную норму Дикина [45]

$$\|h\|_x = \left\langle h, \nabla^2 f(x) h \right\rangle^{1/2},$$

с помощью *демпфированного метода Ньютона* [72, п. 4.1.5]

$$x^{k+1} = x^k - \frac{1}{1 + \left\|\nabla f\left(x^k\right)\right\|_{x^k,*}} \left[\nabla^2 f\left(x^k\right)\right]^{-1} \nabla f\left(x^k\right)$$

можно получить такую оценку числа итераций [72, 158]:

$$\mathrm{O}\left(f\left(x^0\right) - f\left(x_*\right) + \log\left\lceil \log\left(1/\varepsilon\right)\right\rceil\right).$$

Более того, с помощью специального варианта метода отслеживания траектории [452 item 5.2.2] можно улучшить последнюю оценку в части

$$\mathrm{O}\left(f\left(x^0\right) - f\left(x_*\right)\right) \to \tilde{\mathrm{O}}\left(\sqrt{f\left(x^0\right) - f\left(x_*\right)}\right).$$

Введем $g(t) = f(w + tv)$. Самосогласованность $f(x)$ означает, что для любых $w$ и $v$ справедливо неравенство

$$\left|g'''(t)\right| \le 2\left(g''(t)\right)^{3/2}$$

для всех $t$, причем множитель 2 здесь выбран для определенности. Последнее неравенство эквивалентно неравенству

$$\left|\nabla^3 f(x)[h,h,h]\right| \le 2\|h\|_x^3 .$$

В частности, если для $f(x)$ выполняется $M_2 \le 2\mu^{3/2}$, то функция $f(x)$ – самосогласованная. В сильно выпуклом случае последнее неравенство можно обеспечить, например, шкалированием $f(x)$, учитывая, что левая и правая часть неравенства будут по-разному реагировать на такое шкалирование. Заметим, что из ранее полученных оценок и простого наблюдения

$$\frac{\mu}{2}R^2 = \frac{\mu}{2}\left\|x^0 - x_*\right\|_2^2 \le f\left(x^0\right) - f\left(x_*\right) = \Delta f ,$$

т.е. $R^{p-1} \le \left(2\Delta f/\mu\right)^{(p-1)/2}$, следует

$$\mathrm{O}\left(\left(\frac{M_p R^{p-1}}{\mu}\right)^{2/(3p+1)}\right) \to \mathrm{O}\left(\left(\frac{M_p \Delta f^{(p-1)/2}}{\mu^{(p+1)/2}}\right)^{2/(3p+1)}\right).$$

Таким образом, если $M_2 \le 2\mu^{3/2}$, то при $p = 2$

$$\mathrm{O}\left(\Delta f^{1/7}\right) \ll \tilde{\mathrm{O}}\left(\Delta f^{1/2}\right) \; (\Delta f \gg 1),$$

что значительно лучше оценки сложности метода отслеживания траектории для задач оптимизации самосогласованных функций (см. также [266], [452, chapter 5]).

Отметим, что в описанном подходе константы имеют «физическую» размерность, поэтому в отличие от формул, встречавшихся до настоящего момента, в этом замечании и оставшейся части пособия не стоит пытаться проверять корректность формул из соображений размерности [47, гл. 1]. Если дополнительно предполагать, что оптимизируемая функция является *$\nu$-самосогласованным барьером*, т. е.

$$\nu\nabla^2 f(x) \succ \nabla f(x)\nabla f(x)^T ,$$

то в классе методов 2-го порядка (см. ниже описние метода Нестерова–Немировского) также можно достичь следующей оценки числа итераций [72, 158, 452]:

$$\mathrm{O}\left(\sqrt{\nu}\ln\left(\nu/\varepsilon\right)\right).$$

Для любого выпуклого замкнутого множества с непустой внутренностью, не содержащего прямых, можно построить (универсальный) самосогласованный барьер с $\nu \simeq n$ [72, теорема 4.3.2], [103, 182]. В подавляющем большинстве интересных для практики случаев удается конструктивно, (явно) построить соответствующий самосогласованный барьер. ■

◊ В продолжение темы, затронутой в замечании 4, опишем, пожалуй, самый известный и хорошо проработанный способ глобализации сходимости методов второго порядка, см. замечание 1.3. Излагаемые далее результаты восходят к работам А. С. Немировского и Ю. Е. Нестерова конца 80-х годов XX века, посвященных развитию *методов внутренней точки*/*методов внутренних штрафов* (*барьеров*) [72, 158, 452, 498]. Заметим, что метод (внешних) штрафов был описан в замечании 4.3.

Пусть нужно решить с точностью по функции $\varepsilon$ следующую задачу:

$$\langle c, x\rangle \to \min_{x\in Q},$$

где множество $Q$ уже не предполагается простой структуры в смысле начала § 2. Однако предполагается, что для этого множества можно построить *ν-самосогласованный барьер*, т. е. построить такую достаточно гладкую строго выпуклую функцию $F(x)$, которая, в частности, стремится к бесконечности при приближении аргумента к точкам границы множества $Q$ изнутри [72, 158, 452]. Далее исходная задача заменяется однопараметрическим семейством задач (*метод продолжения по параметру*, см., например, [53, § 5.3])

$$t\langle c, x\rangle + F(x) \to \min_{x},$$

где в конечном итоге нужно, чтобы $t \to \infty$. В упомянутом цикле работ Ю. Е. Нестерова и А. С. Немировского был предложен следующий вариант метода Ньютона:

$$t^{k+1} = \left(1 + \frac{1}{13\sqrt{\nu}}\right)t^k,\quad x^{k+1} = x^k - \left[\nabla^2 F\left(x^k\right)\right]^{-1}\left(t^{k+1}\cdot c + \nabla F\left(x^k\right)\right),$$

который после $\mathrm{O}\left(\sqrt{\nu}\ln\left(\varepsilon^{-1}\right)\right)$ итераций находит с точностью по функции $\varepsilon$ решение исходной задачи, если начальное приближение было разумно выбрано. Также было показано, что получить «разумное» начальное при-ближение всегда можно, сделав дополнительно не более $O\left(\sqrt{\nu}\ln(\nu)\right)$ итераций такого же типа, см., например, [180, п. 5.3].

Заметим, что исходную постановку задачи на самом деле можно завязать на более привычные по данному пособию постановки задачи. В частности:

$$f(x) \to \min_x$$

можно эквивалентным образом переписать как

$$y \to \min_{f(x)\le y}.$$

Описанный выше подход лежит в основе пакета выпуклой оптимизации CVX [592], о котором упоминалось во введении. Дело в том, что основные функции, возникающие в выпуклой оптимизации, можно единообразно записывать, используя следующее обобщение *представления* (*исключения*) *Фурье–Моцкина* [98, 158, 442], см. также замечание 3.1,

$$f(x) = \min_{y:\, A\begin{bmatrix} x \\ y \end{bmatrix} = b,\, \begin{bmatrix} x \\ y \end{bmatrix} \in K} \langle c, x\rangle + \langle d, y\rangle,$$

где выпуклый конус $K$ – есть прямое произведение конусов трех (канонических) типов: $\mathbb{R}_+^n$, $S_+^n$, $L_2^n$, т. е. неотрицательного ортанта, конуса неотрицательно определенных матриц и лоренцовского конуса (ice cream cone). Причем, как было показано А. С. Немировским и А. Бен-Талем в 1998 г. [158, 159, 442], лоренцовский конус $L_2^n$ можно заменить конусом $\mathbb{R}_+^{O\left(n\ln\left(\varepsilon^{-1}\right)\right)}$. Кроме того, $L_2^n$ можно понимать как прообраз $S_+^n$ при специальном аффинном преобразовании, т.е. конус $L_2^n$ может быть полностью исключен из числа “канонических” (базисных) [158, item 3.2]. Важно отметить, что для большого числа классов задач выпуклой оптимизации удалось найти такое представление. Изначально был заготовлен достаточно большой набор стандартных (библиотечных) выпуклых функций с известными, найденными «вручную», представлениями. Затем были определены правила сочетания, которые позволяют по известным функциям получать новые, и написана программа, которая организует разумный перебор этих правил для поиска представления для новых функций, не присутствующих в библиотеке. Детали см., например, в [298, 442].

Отметим, что самосогласованность, вообще говоря, не подразумевает сильную выпуклость. В частности, для задачи квадратичной оптимизации (возникающей, например, при приближенном представлении функции $f(x) = e^x$ с помощью Лоренцовского конуса [442])

$$x_n \to \min_{x_1 \geq 0,\, x_2 \geq x_1^2,\, x_3 \geq x_2^2,\, \ldots,\, x_n \geq x_{n-1}^2}$$

можно явно построить самосогласованный барьер и наблюдать на практике хорошее соответствие с теорией в части сходимости описанного метода по функции. Однако сходимость по аргументу чрезвычайно замедляется с ростом размерности пространства $n$, в котором происходит оптимизация. Схожие проблемы и дополнительная вычислительная неустойчивость возникает и при работе с конусом неотрицательных полиномов [458], играющим важную роль в теории управлении [312] и возникающим также при разработке оптимальных алгоритмов [515]. Особенностью данного конуса является наличие двойственного представления (конуса), порожденного специальным подклассом неотрицательно определенных ганкелевых матриц. Этот двойственный конус практически не пересекает конус неотрицательно определенных матриц и работа с ним приводит к серьезным вычислительным трудностям. В качестве простого примера укажем, что популярная в задачах анализа данных плохо обусловленная матрица Гильберта [82, п. 5 § 1, гл. 11] является одним из самых безобидных представителей данного класса ганкелевых матриц [560].

Таким образом, удалось автоматически редуцировать подавляющее большинство задач выпуклой оптимизации к единому виду:

$$\langle c, x\rangle + \langle d, y\rangle \to \min_{(x,y):\, A\begin{bmatrix} x \\ y \end{bmatrix} = b,\, \begin{bmatrix} x \\ y \end{bmatrix} \in K}.$$

Чтобы применять вариант метода внутренней точки, описанный выше, осталось только подобрать самосогласованные барьеры для канонических конусов $\mathbb{R}_+^n$, $S_+^n$ и понять, что происходит при их линейных преобразованиях. Но эти задачи уже были успешно решены еще в самых первых работах. В частности, были построены следующие $n$-самосогласованные барьеры:

$$F_{\mathbb{R}_+^n}(x) = -\sum_{i=1}^{n} \ln x_i\,,\quad F_{S_+^n}(X) = -\ln\det(X).$$

Таким образом и появился один из самых популярных сейчас пакетов CVX решения задач выпуклой оптимизации умеренных размеров – до десятков тысяч переменных [592]. На наш взгляд, оболочка CVX имеет наилучшую реализацию для использования вместе с MatLab и не очень удачную реализацию для работы с Python. Отметим также другие оболочки и словеры, которые на протяжении многих лет являются помощниками специалистов в области численных методов оптимизации. Среди оболочек выделим также AMPL [591] и GAMS [594], а среди солверов отметим CPLEX [600], MOSEK [602], GUROBI Optimization [599], Local Solver [601]. Имеется большой сервер солверов [604], на котором на данный момент можно найти более 60 различных оптимизационных солверов, в том числе и отмеченные выше. И, конечно, нельзя не упомянуть открытую программную библиотеку Google, заточенную под задачи машинного обучения: TensorFlow [608], в частности, задачи глубокого обучения (см. также [596, 598, 606]).

Заинтересовавшемуся в замечании 4 и последующем тексте читателю можно рекомендовать ознакомиться с [72, гл. 4], [158, п. 2.3, 2.5, 3.2], [180, п. 5.3] для более глубокого погружения в затронутые здесь темы. Особо отметим то, насколько неожиданными порой бывают выпуклые функции, допускающее явно выписываемое обобщенное представление Фурье–Моцкина [442]. ◊

Также, как и для методов 1-го порядка (см. § 5), для методов 2-го порядка и выше можно рассматривать их универсальные композитные варианты [299, 301], можно рассматривать работу методов в условиях наличия шума и неточностей, возникающих при решении вспомогательных задач на каждой итерации [143, 286], также можно переносить и завязанные на наличие шумов конструкции, например, конструкцию *минибатчинга* [286].

Однако при использовании методов 2-го порядка и выше появляется много новых вопросов относительно сильного проигрыша методам первого порядка (градиентного типа) по стоимости итерации и требуемой памяти. Так, для честного осуществления шага метода Ньютона необходимо обратить матрицу Гессе оптимизируемой функции в текущей точке. Эта задача по сложности эквивалента задаче умножения двух матриц такого же по порядку размера [59, гл. 31], что типично в $n$ раз дороже, чем осуществление шага метода типа градиентного спуска (умножение матрицы на вектор).

◊ На самом деле это не так. Умножение двух матриц $n\times n$ современными алгоритмами может быть осуществлено за время $O\left(n^{2.37}\right)$, см. [85, 278, 365] и цитированную там литературу. Однако такого рода результаты проявляются только при очень больших значениях $n$. ◊

В последнее время было предложено несколько подходов, имеющих своей целью хотя бы частичное устранение такого большого зазора в стоимости итерации между методами 1-го и 2-го порядка. Одна из идей активно используется в машинном обучении, когда функционал имеет вид суммы (среднего арифметического) большого числа однотипных слагаемых. Идея заключается в том, чтобы формировать матрицу Гессе оптимизируемой функции исходя из матриц Гессе относительно небольшого числа случайно выбранных слагаемых [286, 568]. Другая идея заключается в отказе от обращения матрицы Гессе на итерации, и вместо этого предлагается использовать информацию о собственном векторе, отвечающем наименьшему собственному значению [108, 190]. Для приближенного вычисления такого вектора вполне достаточно уметь умножать матрицу Гессе на произвольный вектор [120]

$$\nabla^2 f(x)\mathrm{v} = \nabla\langle \nabla f(x), \mathrm{v}\rangle,$$

$$\nabla^2 f(x)\mathrm{v} \approx \frac{\nabla f(x+\tau \mathrm{v}) - \nabla f(x)}{\tau},$$

что может быть сделано с помощью автоматического дифференцирования за то же по порядку время, что и вычисление градиента [153, 472]. Эта идея сейчас активно развивается в связи с поиском наиболее эффективных методов обучения глубоких нейронных сетей [39, 78, 117, 120, 190].

◊ Одной из основных проблем, возникающих при обучении глубоких нейронных сетей, является «застревание» процедур типа градиентного спуска в «паразитных» экстремумах [39] (седловых точках). В качестве примера, следуя [579], приведем способ «проскальзывания седловых точек» с помощью доступа к $\nabla^2 f(x)\mathrm{v}$. Будем считать, что $M_1, M_2 < \infty$, см. выше. Следуя [464], будем искать такую точку $x^N$ ($(\varepsilon,\delta)$-локальный минимум), что

$$\left\|\nabla f(x^N)\right\|_2 \le \varepsilon,\ \lambda_{\min}\left(\nabla^2 f(x^N)\right) \ge -\delta.$$

Рассмотрим следующую процедуру:

$$x^{k+1} = x^k - \frac{1}{M_1}\nabla f(x^k), \text{ если } \left\|\nabla f(x^k)\right\|_2 > \varepsilon;$$

$$x^{k+1} = x^k + hp^k, \text{ если } \left\|\nabla f(x^k)\right\|_2 \le \varepsilon \text{ и } \lambda_{\min}\left(\nabla^2 f(x^k)\right) < -\delta,$$

где $h = 2\delta/M_2$, $p^k$ – собственный вектор $\nabla^2 f(x^k)$, отвечающий наименьшему собственному значению. Вектор $p^k$ определяется с точностью до произвола в нормировке. Распорядимся этим произволом следующим образом:

$$\langle \nabla f(x^k), p^k\rangle \le 0,\ \left\|p^k\right\|_2 = 1.$$

Из определения $M_2$ следует (см., например, [448]), что если «сработала» вторая альтернатива, то

$$f\left(x^{k+1}\right) \le f\left(x^k\right) + h\underbrace{\left\langle \nabla f\left(x^k\right), p^k\right\rangle}_{\le 0} + \frac{h^2}{2!}\underbrace{\left\langle \nabla^2 f\left(x^k\right)p^k, p^k\right\rangle}_{\le -\frac{h^2\delta}{2}} + \frac{h^3}{3!}\underbrace{M_2\left\|p^k\right\|_2^3}_{=M_2} \le$$

$$\le f\left(x^k\right) - \frac{1}{2}\left(\frac{2\delta}{M_2}\right)^2\delta + \frac{M_2}{6}\left(\frac{2\delta}{M_2}\right)^3 = f\left(x^k\right) - \frac{2}{3}\frac{\delta^3}{M_2^2}.$$

Повторяя далее рассуждения из § 1, которые использовались при выводе формулы (1.10), получим, что метод гарантированно остановится после

$$N = \left(f\left(x^0\right) - f\left(x^{\text{local min}}\right)\right)\cdot \max\left\{\frac{2M_1}{\varepsilon^2}, \frac{3M_2^2}{2\delta^3}\right\}$$

итераций. Заметим, что искать в данной процедуре собственный вектор точно нет необходимости. Достаточно либо определить, что $\lambda_{\min}\left(\nabla^2 f\left(x^k\right)\right) \ge -\delta$, либо найти такой $p^k$ ($\left\|p^k\right\|_2 = 1$), что $\left\langle \nabla^2 f\left(x^k\right)p^k, p^k\right\rangle \le -\delta/2$. Эта задача может быть решена обычным степенным методом (методом простой итерации), тесно связанным с методом градиентного спуска (см. замечание 4.4), за $\mathrm{O}\left(1/\delta\right)$ итераций (на каждой итерации происходит умножение $\nabla^2 f\left(x\right)\mathrm{v}$) и методом Ланцоша, тесно связанным с методом сопряженных градиентов[90] [91, п. 21.4], за $\mathrm{O}\left(1/\sqrt{\delta}\right)$ аналогичных итераций.

---

[90] Совсем грубо это можно пояснить, например, так. Пусть известно (в этом месте и делается грубое допущение, поскольку знать $\lambda_{\min}$, на самом деле, не нужно): $\lambda_{\min} = \lambda_{\min}\left(\nabla^2 f\left(x^k\right)\right) < -\delta$. Рассмотрим вспомогательную задачу

$$\frac{1}{2}\left\langle \nabla^2 f\left(x^k\right)p, p\right\rangle + \frac{\left|\lambda_{\min}\right|}{2}\left\|p\right\|_2^2 \to \min_p.$$

Подобно замечанию 4.4 можно показать, что градиентные методы (в том числе ускоренные) сходятся для данной задачи (вырожденной, но выпуклой!) к собственному вектору, отвечающему $\lambda_{\min}$, который наиболее близок (в 2-норме) к точке старта, которую в этой связи следует выбирать отличной от нуля. Осталось только заметить следующее: для поиска такого $p^k$ ($\left\|p^k\right\|_2 = 1$), что $\left\langle \nabla^2 f\left(x^k\right)p^k, p^k\right\rangle \le -\delta/2\left\|p^k\right\|_2^2$, достаточно решить ускоренным методом вспомогательную задачу с точностью (по функции) $\delta/4\left\|p^k\right\|_2^2$. Отсюда и получается оценка $\mathrm{O}\left(1/\sqrt{\delta}\right)$. Отметим цикл работ [178, 529, 530], в которых при-

Полученную выше оценку на число итераций можно улучшить, в особенности, если целевой функционал имеет вид суммы большого числа слагаемых. В этом случае используются различные инкрементальные (рандомизированные) методы, см. замечание 1 и [114, 117, 120, 268, 589], а также цитированную в этих работах литературу.

Обычно различают две постановки задачи:

1) В предположении, что везде $\lambda_{\min}\left(\nabla^2 f(x)\right) \geq -\sigma$, где $\sigma > 0$ – малая величина, требуется найти такую точку $x^N$, что $\left\|\nabla f\left(x^N\right)\right\|_2 \leq \varepsilon$.[91]

2) В условиях достаточной гладкости целевой функции (слагаемых) требуется найти $(\varepsilon, \delta)$-локальный минимум (см. выше).

По постановке 1) уже построена, более менее, законченная теория в классе офлайн методов:[92] получены нижние оценки и предложены методы, работающие по нижним оценкам с точностью до логарифмических множителей, см. [588] и цитированную там литературу. Отметим, используемую в ряде подходов (см., например, [114] и цитированную там литературу) специальную технику итеративной регуляризации, напоминающую проксимальный подход (см. замечания 3.2, 3.3):

$$x^{k+1} \approx \arg\min_x \left\{ f(x) + \sigma \left\| x - x^k \right\|_2^2 \right\},$$

где вспомогательная задача уже будет $\sigma/2$-сильно выпуклой, следовательно, ее можно решать ускоренными методами [392, 408].

По постановке 2) в последние годы активно ведутся исследования [113, 114, 117, 120, 134]. Наилучшие известные сейчас методы приведены в работах [268, 589, 134] и упражнении 1.8 (при должной доработке [268]),

---

водятся нижние оценки для класса *методов Крылова* [91, 212], использующих только умножение матрицы на вектор, для решения систем линейных уравнений и задачи вычисления собственного вектора матрицы, отвечающего максимальному (минимальному) собственному значению, эквивалентной ей по сложности с точностью до логарифмического множителя. Нижние оценки получаются в полном соответствии с оценками, полученными отмеченной выше редукцией исходной постановки задачи к задаче квадратичной выпуклой оптимизации.

[91] Задача поиска глобального $\varepsilon$-минимума по функции даже в такой постановке остается NP-трудной [432].

[92] Использована терминология [114, 117, 268]. Скорость сходимости офлайн методов может зависеть от числа слагаемых, см. упражнение 1.8. При этом можно ограничиться классом методов для задач минимизации функционала вида суммы, сложность которых не зависит от числа слагаемых, см., например, замечание 1 и упражнение 1.8.

см. также работы [204, 489, 569, 570], содержащие практически более эффективные варианты.

Отметим, интересное обобщение упомянутого ранее метода Ланцоша на случай, когда вместо произведения $\nabla^2 f(x)v$ доступно произведение $\nabla^2 f_i(x)v$ для случайно (равновероятно) выбранного слагаемого. В работе [120] была предложена разновидность алгоритма Ойя (рандомизированный вариант метода Ланцоша), позволяющая решать ту же задачу, что решалась выше с помощью метода Ланцоша, за $\tilde{O}(1/\delta^2)$ произведений вида $\nabla^2 f_i(x)v$. При этом такие произведения не обязательно специально генерировать, их можно «собирать» по ходу работы основного алгоритма, как бы «заодно». Отметим, что развиваемые в этом направлении наивные варианты методов типа редукции дисперсии в случае существенно невыпуклой целевой функции могут вести себя очень плохо на практике [220]. ◊

Перспективной также представляется довольно старая идея глобализации сходимости [42, 53, 209, 472]: спуск в область квадратичной (в общем случае, сверхлинейно) сходимости с помощью методов типа градиентного спуска (с дешевыми итерациями) и последующая квадратичная сходимость с использованием, например, метода Ньютона. Проблема в таком подходе – детектирование момента попадания в нужную окрестность. В качестве возможного решения проблемы можно, например, действовать таким образом: через каждые $\sim \sqrt{n}$ итераций метода типа градиентного спуска проверять условие $\left\|\nabla f\left(x^k\right)\right\|_2 \ll 1$. Если оно выполняется, то делать «пристрелочный» шаг метода Ньютона. Если в результате такого шага выполняется еще и условие

$$\left\|\nabla f\left(x^{k+1}\right)\right\|_2 \ll \left\|\nabla f\left(x^k\right)\right\|_2^{3/2},$$

то продолжать делать шаги метода Ньютона, каждый раз проверяя это условие. Если хотя бы одно из этих условий не выполняется, то вернуться к методу типа градиентного спуска.

Надеемся, что данное пособие вызовет желание поработать в описанных в нем направлениях!

# Литература

Учебное издание

**Гасников** Александр Владимирович

# СОВРЕМЕННЫЕ ЧИСЛЕННЫЕ МЕТОДЫ ОПТИМИЗАЦИИ.

## МЕТОД УНИВЕРСАЛЬНОГО ГРАДИЕНТНОГО СПУСКА